\documentclass[11pt,a4paper,reqno]{amsart}

\usepackage[utf8]{inputenc}

\usepackage{amsmath,amsfonts,amssymb,amsthm,color,bm}
\usepackage{dsfont} %math blackboard boldface 1
\usepackage{mathtools} %DeclarePairedDelimeter
\usepackage{esint} %cool integral
\usepackage{cancel} %slash the formula like a boss
\usepackage{mathrsfs}

\usepackage{tikz-cd}%Category setting

\usepackage{hyperref}
\hypersetup{
    colorlinks=true,
    linkcolor=blue,
    citecolor=blue,
    filecolor=magenta,      
    urlcolor=cyan,
} %Referencing with hyperlinks

\usepackage{graphicx}
\graphicspath{ {./images/} }
\usepackage[left=2.5cm, right=2.5cm]{geometry}

\makeatletter
\def\l@subsection{\@tocline{2}{0pt}{2.5pc}{2.5pc}{}}
\makeatother
\theoremstyle{plain}
\newtheorem{theorem}{Theorem}[section]

\newtheorem{corollary}[theorem]{Corollary}
\newtheorem{lemma}[theorem]{Lemma}
\newtheorem{proposition}[theorem]{Proposition}

\theoremstyle{definition}
\newtheorem{definition}[theorem]{Definition}
\newtheorem{conjecture}[theorem]{Conjecture}

\newtheorem{remark}[theorem]{Remark}

\newtheorem{question}[theorem]{Question}
\newtheorem*{theorem*}{Theorem}

\numberwithin{equation}{section}

\newcommand{\todo}[1]{
}
\newcommand{\MH}[1]{
}
\newcommand{\ex}[1]{
}
\newcommand{\Dil}{\operatorname{Dil}}
\newcommand{\Tr}{\operatorname{Tr}}
\newcommand{\Mod}{\operatorname{Mod}}

\newcommand{\abs}[1]{\left\vert #1 \right\vert}
\newcommand{\nrm}[1]{\left\Vert #1 \right\Vert}
\newcommand{\ang}[1]{\left< #1 \right>}
\newcommand{\br}[1]{\left( #1 \right)}
\newcommand{\mr}[1]{\left[ #1 \right]}
\newcommand{\bR}[1]{\left( #1 \right]}
\newcommand{\Br}[1]{\left[ #1 \right]}
\newcommand{\BR}[1]{\left\{ #1 \right\}}
\newcommand{\dif}[1]{\left. #1 \right \vert}

\newcommand{\C}{\mathbb{C}}
\newcommand{\E}{\mathbb{E}}
\newcommand{\N}{\mathbb{N}}
\newcommand{\I}{\mathbb{I}}
\newcommand{\R}{\mathbb{R}}
\newcommand{\RP}{\mathbb{RP}}
\newcommand{\T}{\mathcal{T}}
\newcommand{\Q}{\mathbb{Q}}
\newcommand{\Z}{\mathbb{Z}}
\newcommand{\1}{\mathds{1}}

\newcommand{\4}[1]{\mathcal{F}^{\br{#1}}}
\newcommand{\fM}[1]{\mathcal{M}^{\br{#1}}}
\newcommand{\fP}[1]{\mathcal{P}^{\br{#1}}}
\newcommand{\fQ}[1]{\mathcal{Q}^{\br{#1}}}
\newcommand{\fR}[1]{\mathcal{R}^{\br{#1}}}
\newcommand{\fT}[1]{\mathcal{T}^{\br{#1}}}

\newcommand{\cA}{{\mathcal{A}}}
\newcommand{\cB}{{\mathcal{B}}}
\newcommand{\cC}{{\mathcal{C}}}
\newcommand{\cD}{{\mathcal{D}}}

\newcommand{\cF}{{\mathcal{F}}}
\newcommand{\cG}{{\mathcal{G}}}
\newcommand{\cH}{{\mathcal{H}}}

\newcommand{\cM}{{\mathcal{M}}}

\newcommand{\cP}{{\mathcal{P}}}
\newcommand{\cQ}{{\mathcal{Q}}}
\newcommand{\cR}{{\mathcal{R}}}
\newcommand{\cS}{{\mathcal{S}}}
\newcommand{\cT}{{\mathcal{T}}}

\DeclareRobustCommand{\rchi}{{\mathpalette\irchi\relax}}
\newcommand{\irchi}[2]{\raisebox{\depth}{$#1\chi$}} % inner command, used by \rchi

\tikzset{
  mor/.style={
    draw=none,
    every to/.append style={
      edge node={node [sloped, allow upside down, auto=false]{$#1$}}}
  }
}
\DeclareFontFamily{U}{mathx}{}
\DeclareFontShape{U}{mathx}{m}{n}{<-> mathx10}{}
\DeclareSymbolFont{mathx}{U}{mathx}{m}{n}
\DeclareMathAccent{\widehat}{0}{mathx}{"70}
\DeclareMathAccent{\widecheck}{0}{mathx}{"71}

\newcommand{\avg}{\operatorname{Avg}}
\newcommand{\bavg}{\operatorname{BAvg}}
\newcommand{\cavg}{\operatorname{CAvg}}
\newcommand{\supp}{\operatorname{supp}}
\newcommand{\dom}{\operatorname{dom}}
\newcommand{\dia}{\operatorname{dia}}
\newcommand{\dist}{\operatorname{dist}}
\newcommand{\img}{\operatorname{img}}
\newcommand{\main}{{\operatorname{main}}}
\newcommand{\err}{{\operatorname{err}}}
\newcommand{\lac}{{\operatorname{lac}}}
\newcommand{\sgn}{\operatorname{sgn}}
\newcommand{\op}{{\operatorname{op}}}
\newcommand{\an}{{\operatorname{an}}}
\newcommand{\pv}{\operatorname{p.v.}}
\newcommand{\chr}{\operatorname{char}}

\newcommand{\vnull}{{\bm{0}}}
\newcommand{\valpha}{{\bm{\alpha}}}
\newcommand{\vbeta}{{\bm{\beta}}}
\newcommand{\vgamma}{{\bm{\gamma}}}
\newcommand{\vxi}{{\bm{\xi}}}
\newcommand{\veta}{{\bm{\eta}}}
\newcommand{\vzeta}{{\bm{\zeta}}}
\newcommand{\vlambda}{{\bm{\lambda}}}
\newcommand{\vepsilon}{{\bm{\epsilon}}}
\newcommand{\vLambda}{{\bm{\Lambda}}}
\newcommand{\vPhi}{{\bm{\Phi}}}
\newcommand{\vA}{{\bm{A}}}

\newcommand{\ve}{{\bm{e}}}

\newcommand{\vf}{{\bm{f}}}

\newcommand{\vt}{{\bm{t}}}

\newcommand{\vx}{{\bm{x}}}
\newcommand{\vy}{{\bm{y}}}
\newcommand{\vz}{{\bm{z}}}
\newcommand{\vchi}{{\bm{\rchi}}}
\newcommand{\vM}{{\bm{M}}}
\newcommand{\vP}{{\bm{P}}}
\newcommand{\vQ}{{\bm{Q}}}

\newcommand{\txi}{{\widetilde{\xi}}}
\newcommand{\teta}{{\widetilde{\eta}}}
\newcommand{\trho}{{\widetilde{\rho}}}

\newcommand{\ts}{{\widetilde{s}}}
\newcommand{\tu}{{\widetilde{u}}}

\newcommand{\tv}{{\widetilde{v}}}
\newcommand{\tlambda}{{\widetilde{\lambda}}}
\newcommand{\tPhi}{{\widetilde{\Phi}}}
\newcommand{\tPsi}{{\widetilde{\Psi}}}

\newcommand{\tvlambda}{{\widetilde{\vlambda}}}

\newcommand{\simto}{\stackrel{\sim}{\longrightarrow}}

\title[corner-type Smoothing inequalities]{Smoothing inequalities for corner-type bilinear averages: geometric characterization and applications }
\author{Martin Hsu, Fred Yu-Hsiang Lin}
\date{}

\begin{document}

\subjclass[2020]{42B15, 42B20, 42B25}
\keywords{Smoothing inequality, Triangular Hilbert transform, Spherical maximal operator, O-minimal theory, corner-type pattern}

\begin{abstract}
We study Sobolev smoothing inequalities for bilinear averages associated with corner-type configurations
\[
\cA_{\gamma,\rho}(f_1,f_2)(x_1,x_2)
=\int f_1(x_1+\gamma_1(t),x_2)f_2(x_1,x_2+\gamma_2(t))\rho(t)\,dt.
\]
In the real-analytic setting, we obtain a complete geometric characterization of the curves $\vgamma=(\gamma_{1},\gamma_{2})$ for which $\cA_{\gamma,\rho}(f_1,f_2)$ satisfies the smoothing inequality
\[
\|\cA_{\gamma,\rho}(f_1,f_2)\|_{L^1}
\lesssim
\|f_1\|_{H^{(-\varepsilon,0)}}
\|f_2\|_{H^{(0,-\varepsilon)}}\, 
\]
for some $\varepsilon >0$.

For general \(C^4\) embedded curves, we establish analogous quantitative statement involving purely geometric conditions that encode certain uniform complexity bounds and allows degeneracies on the various curvature conditions. For definable families of curves in an arbitrary o-minimal expansion of the real field, the relevant complexity parameters are uniformly finite, leading to smoothing inequalities that hold uniformly across the family.

As applications, we obtain bounds for triangular Hilbert transforms along a large family of curves, their associated maximal operators, and corner-type lacunary spherical maximal operators. We further prove the existence of configurations of the form
 $(x,y)$, $(x+\gamma_{1}(t),y)$, $(x,y+\gamma_{2}(t))$ inside sets of positive measure, together with a quantitative lower bound on the gap $t$.

\end{abstract}
\maketitle

\newpage

\tableofcontents

\section{Introduction}
Let $I\subseteq  \mathbb{R}$ be an interval, let \(\vgamma:=\br{\gamma_1,\gamma_2}:I\to\R^2\), and let \(\rho\) be an integrable function with \(\supp \rho\Subset I\). We consider the corner-type averaging operator:
\begin{equation}\label{eq_A_op}
    \cA_{\vgamma,\rho}\br{f_1,f_2}\br{\vx}:=\int f_1\br{x_1+\gamma_1\br{t},x_2}f_2\br{x_1,x_2+\gamma_2\br{t}}\rho\br{t}dt\, .
\end{equation}
For $s_{1},s_{2}\in \mathbb{R}$, define
\begin{equation*}
    \nrm{f}_{H^{\br{s_1,s_2}}}:=
    \nrm{
        \ang{\xi_1}^{s_1}\ang{\xi_2}^{s_2}
        % \br{1+\abs{\xi_1}^2}^{\frac{s_1}{2}}
        % \br{1+\abs{\xi_2}^2}^{\frac{s_2}{2}}
        \widehat{f}\br{\vxi}
    }_{L^2\br{d\vxi}},\quad
    \ang{z}:=\br{1+\abs{z}^2}^{\frac{1}{2}}.
\end{equation*}
We aim to determine when the operator \eqref{eq_A_op} satisfies the Sobolev smoothing inequality:
\begin{equation}\label{eq_soblev_smoothing_std}
    \nrm{\cA_{\vgamma,\rho}\br{f_1,f_2}}_{L^1}\leq C_{\vgamma,\rho,\varepsilon}
    \nrm{f_1}_{H^{\br{-\varepsilon,0}}}
    \nrm{f_2}_{H^{\br{0,-\varepsilon}}},\quad f_1,f_{2}\in L^2\br{\R^2}\, ,
\end{equation}
for some $\varepsilon >0$. The central question of this paper is therefore: \textit{which curves $\vgamma$ admit the smoothing inequality \eqref{eq_soblev_smoothing_std}?} Our first main theorem gives a complete characterization in the real-analytic setting.
\begin{theorem}\label{thm_real_ana_smoothing}
    Given real analytic \(\vgamma:I\to\R^2\) and non-negative \(\rho\in C^2_c\br{\R}\) with \(\int \rho\br{t}dt>0\) and \(\supp\rho\Subset I\), the Sobolev smoothing inequality \eqref{eq_soblev_smoothing_std} holds with some \(0<\varepsilon,C_{\vgamma,\rho,\varepsilon}<\infty\) if and only if \(\img\vgamma\) is not contained in any of the following sets
    \begin{equation}\label{eq_line}
        \BR{
            \br{x,y}\in\R^2
        \::\:
            ax+ by = c
        },\tag{Line}
    \end{equation}
    \begin{equation}\label{eq_exp}
        \BR{
            \br{x,y}\in\R^2
        \::\:
            ae^{Ax}+ by = c
        },\tag{Exp}
    \end{equation}
    \begin{equation}\label{eq_log}
        \BR{
            \br{x,y}\in\R^2
        \::\:
            ax+ be^{By} = c
        },\tag{Log}
    \end{equation}
    \begin{equation}\label{eq_psi_line}
        \BR{
            \br{x,y}\in\R^2
        \::\:
            ae^{Ax}+ be^{By} = c
        },\tag{\(\Psi\)-line}
    \end{equation}
    for \(\br{a,b,c}\in\R^3\setminus\BR{\br{0,0,0}}\) and \(A,B\in\R\setminus\BR{0}\).
\end{theorem}
\begin{remark}
The four exceptional classes are affine lines after applying, independently
in each coordinate, either the identity map or an exponential change of variable. They are genuine
obstructions: each class carries a form of (generalized) modulation symmetry that fixes the left-hand side of \eqref{eq_soblev_smoothing_std} while the negative Sobolev norms of suitable modulated inputs tend to zero. In particular, obstructions similar to \eqref{eq_psi_line}, albeit for smoothing inequalities associated to bilinear averages on \(\R\), also appear in the work \cite{MR4776384} by Christ and Zhou. See \textit{Hypothesis 2} and \textit{Lemma 2.5} in \cite{MR4776384}.
\end{remark}
Beyond the real-analytic characterization, we establish more general results for parametrizations $\vgamma$ of curves \(\Gamma\subseteq \R^2\). These results show that the validity of the Sobolev smoothing inequality \eqref{eq_soblev_smoothing_std} is governed by geometric properties of \(\Gamma\) rather than by the particular parametrization \(\vgamma\), and they quantify the dependence of the
implicit constant in \eqref{eq_soblev_smoothing_std} on the relevant geometric and analytic parameters. We begin by specifying the class of curves under consideration.

\begin{definition}[Curves in \(\R^2\)]
    Let \(N\geq 0\). We say \(\Gamma\subseteq \R^2\) is a \(C^N\) curve if it is an embedded \(C^N\) sub-manifold of \(\R^2\) diffeomorphic to an open interval. By a parametric curve, we mean a \(C^N\) embedding \(\vgamma:I\simto\Gamma\subseteq \R^2\) from an open interval \(I\) to \(\R^2\).
\end{definition}
\begin{remark}
    If a property \(P\) of the parametric curve \(\vgamma:I\simto\Gamma\) depends on the parametrization, we denote it by \(P_\vgamma\); otherwise, we denote it by \(P_\Gamma\).
\end{remark}
Henceforth, we fix a \(C^N\) curve \(\Gamma\subseteq \R^2\) and a parametrization \(\vgamma:I\simto\Gamma\).
The fiber structure of \(\cA_{\vgamma,\rho}\) naturally leads us to apply Fourier projections and related operators separately in each coordinate. We therefore
introduce the following notation.
\begin{definition}[Fiber-wise application]\label{def_fiber_op}
    Given an operator \(\cT:L^p\br{\R}\to L^p\br{\R}\), we define its fiber-wise application for functions \(f\in L^p\br{\R^2}\) on the plane via the following formulas:
    \begin{equation*}
        \fT{1}f\br{\vx}:=\br{\cT f\br{\cdot ,x_2}}\br{x_1},\quad
        \fT{2}f\br{\vx}:=\br{\cT f\br{x_1,\cdot}}\br{x_2}.
    \end{equation*}
\end{definition}

With the above notion, we can express \(\cA_{\vgamma,\rho}\) in terms of its Fourier representation:
\begin{equation*}
    \cA_{\vgamma,\rho}\br{f_1,f_2}\br{\vx}=
    \int
        \4{1}f_1\br{\xi_1,x_2}
        \4{2}f_2\br{x_1,\xi_2}
        m_{\vgamma,\rho}\br{\vxi}
        e\br{\vx\cdot\vxi}
    d\vxi,
\end{equation*}
where \(e\br{z}:=e^{2\pi i z}\) denotes the standard character on \(\R\), \(\cF\) stands for the Fourier transform, and the multiplier/symbol \(m_{\vgamma,\rho}\) takes the following form:
\begin{equation}\label{eq_m_gamma_rho}
    m_{\vgamma,\rho}\br{\vxi}:=
    \int
        e\br{\vxi\cdot \vgamma\br{t}}
        \rho\br{t}
    dt.
\end{equation}
The following theorem for the fiber-wise bilinear Fourier integral operator plays a key role in proving our main theorem for the smoothing inequality:
% above Fourier representation suggest that we consider the following related object:
\begin{theorem}[A universal bound for the fiber-wise bilinear Fourier integral operator]\label{thm_universal_multiplier_bound}
    For a symbol \(m\in \br{L^1\cap L^\infty}\br{\R^2}\), and functions \(F,G\in \br{L^1\cap L^\infty}\br{\R^2}\), consider the fiber-wise bilinear Fourier integral operator
    \begin{equation}\label{carlyT}
        \cT_m\br{F,G}\br{x,y}:=
        \1_{\mr{0,1}^2}\br{x,y}
        \int
            F\br{\xi,y}
            G\br{x,\eta}
            m\br{\xi,\eta}
            e\br{x\xi+y\eta}
        d\xi d\eta.
    \end{equation}
    The following estimate holds:
    \begin{equation*}
        \nrm{\cT_m\br{F,G}}_{L^1}
        \lesssim
        \nrm{m}^{\frac{1}{5}}_u
        \nrm{m}^{\frac{4}{5}}_U
        \nrm{F}_{L^2}
        \nrm{G}_{L^2},
    \end{equation*}
    where the two quantities are defined by the following formulas:\footnote{We adopt the convention
    \(\|f(x,y)\|_{L^{p}_{x}L^{q}_{y}}:=\left\| \|f(x,y)\|_{L^{p}_{x}} \right\|_{L^{q}_{y}}\) for the iterated Lebesgue norms.} 
    \begin{equation}\label{mlittleunorm}
        \|m\|_{u}:=\left\|\4{1}\br{\mathcal{D}_{\br{0,s}}m}\br{x,\eta}  \right\|^{\frac{1}{2}}_{L^{\infty}_{s\eta}L^{2}_{x}([-1,1])} ,
    \end{equation}
    \begin{equation}\label{mcapitalunorm}
    \|m\|_{U}:=\left\|\int_{\mathbb{R}}\mathcal{D}_{\br{0,s}}\mathcal{D}_{\br{u,v}}m(\xi, \eta)d\xi   \right\|^{\frac{1}{4}}_{L^{\infty}_{u\eta}L^{1}_{s}L^{2}_{v}},
    \end{equation}
    and the multiplicative derivative is defined by
    \begin{equation}\label{multideri_def}
(\mathcal{D}_{(y_{1},\cdots,y_{d})}f)(x_{1},\cdots ,x_{d}):=f(x_{1}+y_{1},\cdots ,x_{d}+y_{d})\overline{f(x_{1},\cdots,x_{d})}\, .
    \end{equation}
\end{theorem}

To illustrate the connection between the smoothing inequality and the above theorem, we use orthogonality and introduce the localized multiplier near the frequency \(\vlambda:=\br{\lambda_1,\lambda_2}\in \br{\R\setminus\BR{0}}^2\):
\begin{equation*}
    m_{\vgamma,\rho,\vlambda}\br{\vxi}:=
    \int
        e\br{\vxi\cdot \vgamma\br{t}}
        \rho\br{t}
    dt
    \prod_{j=1,2}
    \psi\br{\xi_j/\lambda_j},
\end{equation*}
where $\rho\in C^\infty_c\br{\R}$ will be chosen to be supported in a short interval around \(0\) and \(\psi\in C^\infty_c\br{\R}\) will be chosen to support in a short interval around \(1\) (both $\rho,\psi$ will be rigorously defined in \eqref{eq_m_vlambda_freq_trunc}). It is enough to prove a power-decay estimate of the form:
\begin{equation*}
    \nrm{
        \cT_{m_{\vgamma,\rho,\vlambda}}\br{F,G}
    }_{L^1}
    \underset{\vgamma,\rho}{\lesssim} 
    \lambda^{-\varepsilon}
    \nrm{F}_{L^2}
    \nrm{G}_{L^2},\quad
    \abs{\lambda_j}\eqsim\lambda \gg 1.
\end{equation*}
By \textbf{Theorem \ref{thm_universal_multiplier_bound}}, this reduces to establishing the two bounds
\begin{equation*}
    \nrm{m_{\vgamma,\rho,\vlambda}}_u
    \underset{\vgamma,\rho}{\lesssim}  
    \lambda^{-c},\quad
    \nrm{m_{\vgamma,\rho,\vlambda}}_U
    \underset{\vgamma,\rho}{\lesssim} 
    \log^C\lambda.
\end{equation*}
Estimating these two quantities is the most technical part of the argument. The main difficulty is the oscillatory
structure of the multiplier. Stationary phase suggests the factorization
\begin{equation*}
    m_{\vgamma,\rho,\vlambda}\br{\vxi}=
    e\br{\vxi\cdot\vgamma\br{t_\ast\br{\vxi}}}
    \Psi_{\vlambda}\br{\vxi}
\end{equation*}
where \(t=t_\ast\br{\vxi}\) solves the critical point equation \(\vxi\cdot \vgamma'\br{t}=0\). Thus, \(\vxi\) is normal to \(\Gamma\) at \(\vgamma\br{t_\ast\br{\vxi}}\). This geometric interpretation allows us to relate the phase function \(\vxi\cdot\vgamma\br{t_\ast\br{\vxi}}\) to an expression related to a modified notion of Gauss map. Such a modification maps a curve \(\Gamma\) to the projective space instead of the usual setting mapping to the \(1\)-sphere.

\begin{definition}[The projective space]
The projective space \(\RP^n\) as a set is given by \(\RP^n=^{\R^{n+1}\setminus\BR{\vnull}}\hspace{-1ex}/_\sim\) where
    \begin{equation*}
        \vxi \sim \veta\iff \exists \lambda\in\R\setminus \BR{0},\, \veta=\lambda\vxi,\quad
        \forall \vxi,\veta\in\R^{n+1}\setminus\BR{0}.
    \end{equation*}
For \(\vxi:=\br{\xi_1,\dots,\xi_{n+1}}\in \R^{n+1}\setminus\BR{\vnull}\), we introduce the following notation:
\begin{equation*}
    \mr{\vxi}:=\mr{\xi_1:\cdots:\xi_{n+1}}:=
    \BR{
        \veta\in\R^{n+1}\setminus\BR{\vnull}
    \::\:
        \vxi\sim\veta
    }.
\end{equation*}
We equip \(\RP^n\) with the usual smooth manifold structure and denote \(\RP:=\RP^1\) for our setting.
\end{definition}

\begin{definition}[The Gauss map]
    For a \(C^1\) curve \(\Gamma\subseteq\R^2\), we define the associated Gauss map:
    \begin{equation*}
        \cG_\Gamma:\Gamma\to \RP:p\mapsto N_p\Gamma\setminus \BR{\vnull}.
    \end{equation*}
\end{definition}

\begin{remark}[Inverse of the Gauss map]\label{rmk_inv_Gauss}
    Note that a \(C^1\) curve \(\Gamma\) is homeomorphic to an open interval \(I\subseteq \R\). Via the order topology on \(I\), we see that the inverse \(\cG^{-1}_\Gamma\) exists on \(\img \cG_\Gamma\subsetneq \RP\) and is continuous whenever \(\cG_\Gamma\) is injective. In other words, an injective Gauss map induces a topological embedding \(\cG_\Gamma:\Gamma\hookrightarrow\RP\).
\end{remark}

\begin{definition}[Natural extension of the inverse Gauss map]\label{def_ext_of_inv_Gauss}
    Assume the inverse \(\cG^{-1}_\Gamma:\img\cG_\Gamma\to\Gamma\) exists. We define the extension as below:
    \begin{equation*}
        \vgamma_\Gamma:\:
        \cC_\Gamma:=
        \bigsqcup_{p\in\Gamma}\cG_\Gamma\br{p} \to\Gamma:\:
        \vxi\mapsto \cG^{-1}_\Gamma\br{\mr{\vxi}}.
    \end{equation*}
\end{definition}

With the above notion, we can now state the key identity \(\vgamma_\Gamma\br{\vxi}=\vgamma\br{t_\ast\br{\vxi}}\) and thus, rewrite:
\begin{equation*}
    m_{\vgamma,\rho,\vlambda}\br{\vxi}=
    e\br{\vxi\cdot\vgamma_\Gamma\br{\vxi}}
    \Psi_{\vlambda}\br{\vxi}.
\end{equation*}
The quantities \(\nrm{m_{\vgamma,\rho,\vlambda}}_u\) and \(\nrm{m_{\vgamma,\rho,\vlambda}}_U\) can be computed through oscillatory integrals.
We estimate them using a variant of van der Corput’s lemma stated in \textsc{Appendix \ref{sec_van_der_corput}}. Although neither the
first nor the second derivative of the phase is uniformly non-vanishing on the entire integration domain, at every
point at least one of the two derivatives is nonzero. The relevant version of van der Corput’s lemma therefore
requires a uniform bound on how often the dominant derivative can change. This complexity is measured by
the following quantity.
\begin{definition}[A uniform finiteness characterization]\label{def_uni_fin_ch}
    Given a \(C^1\) curve \(\Gamma\) with injective Gauss map,
    let \(\vgamma_\Gamma=\br{\gamma_{\Gamma,1},\gamma_{\Gamma,2}}\) be the associated extension of the inverse Gauss map. We define:
    \begin{equation}\label{eq_M_Gamma_def}
        \Xi_\Gamma:=
        \sup_{
            % \br{\ast}
            \substack{
            s,u,v,
            w,\eta\\
            :p\br{s,v,\eta}
            }
        }
        \max
        \begin{pmatrix}
            C
            \br{
                \BR{\xi\in\R
                    \::\:
                    \br{\Delta_{\br{0,s}}\partial_1\gamma_{\Gamma,1}}\br{\xi,\eta}=w
                }
            },\\
            C
            \br{
                \BR{\xi\in\R
                    \::\:
                    \br{\Delta_{\br{s,0}}\partial_2\gamma_{\Gamma,2}}\br{\eta,\xi}=w
                }
            },
            \\
            C
            \br{
                \BR{\xi\in\R
                    \::\:
                    \br{\Delta_{\br{0,s}}\Delta_{\br{u,v}}\partial_1\gamma_{\Gamma,1}}\br{\xi,\eta}=w
                }
            },\\
            C
            \br{
                \BR{\xi\in\R
                    \::\:
                    \br{\Delta_{\br{s,0}}\Delta_{\br{v,u}}\partial_2\gamma_{\Gamma,2}}\br{\eta,\xi}=w
                }
            }\phantom{,}
        \end{pmatrix}
        ,
    \end{equation}
    where \(C\br{B}\) denotes the number of connected components of a set \(B\subseteq \R\) with the convention that \(C\br{\varnothing}=0\),
    the symbol \(\Delta_{\br{u,v}}\) stands for the finite difference given by the formula \(\Delta_{u,v}f\br{x,y}:=f\br{x+u,y+v}-f\br{x,y}\), and the notation \(p\br{s,v,\eta}\) stands for the sign condition expressed by the following formula:\footnote{\eqref{eq_sign_cond} is equivalent to the sign condition \(\sgn\br{\eta}=\sgn\br{\eta+s}=\sgn\br{\eta+v}=\sgn\br{\eta+s+v}\).}
    \begin{equation}\label{eq_sign_cond}
        \br{\eta, \eta+s, \eta+v, \eta+s+v > 0}\lor\br{\eta, \eta+s, \eta+v, \eta+s+v < 0}.
    \end{equation}
\end{definition}

\begin{remark}[Interpretation of the equality]\label{rmk_equal}
    For a function \(f:A\to B\), we interpret the formula \(f\br{a}=b\) as a shorthand for the following formula:
    \begin{equation*}
        \br{a\in A}\land\br{b\in B}\land\br{f\br{a}=b}.
    \end{equation*}
\end{remark}

\begin{remark}[Domains of the functions]\label{rmk_dom_func}
    Given a function \(\vf:\R^m\supset \dom\vf\to\R^n\), a vector \(\vz\in\R^m\), and an index \(j\in\BR{1,\dots,m}\), we interpret the domain of the functions \(\Delta_\vz \vf\) and \(\partial_j \vf\) as:
    \begin{equation}\label{eq_dom_Delta_f}
        \dom\Delta_\vz \vf:=
        \BR{
            \vx\in\R^m
        \::\:
            \vx,\vx+\vz\in \dom \vf
        };
    \end{equation}
    \begin{equation}\label{eq_dom_partial_f}
        \dom\partial_j \vf:=
        \BR{
            \vx\in\R^m
        \::\:
            \exists \vy\in\R^n,\,\forall\epsilon>0,\,\exists \delta>0,\,\forall h\in \R,\,
            D_{\vf,j}\br{\vy,\epsilon,\delta,h,\vx}
        }
    \end{equation}
    with the formula \(D_{\vf,j}\br{\vy,\epsilon,\delta,h,\vx}\) expressing the differentiability as below:
    \begin{equation}\label{eq_diff_ability}
        % \exists L\in\R,\,\forall\epsilon>0,\,\exists \delta>0,\,\forall h\in \R,\,
            % \mr{
            0<\abs{h}<\delta
            \implies
            \br{
                \vx,\vx+h\ve_j\in \dom \vf\land
                \abs{
                    \frac{\vf\br{\vx+h\ve_j}-\vf\br{\vx}}{h}
                    -\vy
                }<\epsilon
            % }
            }.
    \end{equation}
    Under the above interpretation, we shall allow functions with empty domain.
\end{remark}

The oscillatory estimates also require quantitative nondegeneracy conditions on the phase \(\vxi\cdot \vgamma_\Gamma\br{\vxi}\). We
encode these conditions in the following notion.
\begin{definition}[\(A\)-admissibility]\label{def_A_admi_ana}
    We say \(\vgamma:I\simto \Gamma\) is \(A\)-admissible if \(\vgamma\in C^4\br{I}\) and satisfies:
    \begin{equation}\label{eq_gamma_admi_ana}
        \abs{\vgamma^{\br{k}}\br{t}}\leq A;\quad
        A^{-1}\leq 
        \abs{\gamma'_j\br{t}},\,
        \abs{\vgamma'\wedge\vgamma''}\br{t},\,
        \abs{\vchi'_\vgamma\wedge\vchi''_\vgamma}\br{t},\quad
        \forall t\in I,\: j=1,2,\:k=1,2,3,4,
    \end{equation}
    where the function \(\vchi_\vgamma\) is defined by the following formula:
    \begin{equation*}
        \vchi_\vgamma:\:
        \BR{t\in I\::\: \vgamma'\br{t}\wedge\vgamma''\br{t}\neq 0}
        \to
        \R^2:\:
        t
        \mapsto
        \br{\frac{\gamma'_1\gamma'_2}{\vgamma'\wedge\vgamma''}}\br{t}\cdot \vgamma'\br{t}.
    \end{equation*}
    % When the exact number \(1\leq A<\infty\) is irrelevant, we simply call such \(\vgamma\) admissible.
\end{definition}

\begin{remark}\label{rmk_small_img}
    As a direct consequence of the mean value theorem,
    for \(P\in\C\mr{z_1,\dots,z_{11}}\) and \(A\)-admissible \(\vgamma:I\simto \Gamma\) with \(A\gg 1\), the following expression:
    \begin{equation*}
        F:=P\br{
            \gamma_1,\dots\gamma'''_1,
            \gamma_2,\dots\gamma'''_2,
            \frac{1}{\gamma'_1},
            \frac{1}{\gamma'_2},
            \frac{1}{\vgamma'\wedge\vgamma''}
        }
    \end{equation*}
    as a function \(F:I\to \C\) satisfies the following control on the image:
    \begin{equation*}
        \dia\br{F\br{U}}\leq A^{O_P\br{1}}\dia\br{U},\quad
        U\subseteq I.
    \end{equation*}
    Equivalently, it says that \(F\) is a Lipschitz map with the Lipschitz constant \(L_F\leq A^{O_P\br{1}}\).
\end{remark}

With the above notions, we are ready to make the previously mentioned argument precise. If \(\vgamma:I\simto\Gamma\) is \(A\)-admissible, the compactness of \(\supp\rho\) permits a
finite decomposition of \(I\) and \(\rho\), reducing the analysis to the following short-interval regime.
\begin{definition}[Short-interval condition]\label{ass_short_I}
     The tuple \(\br{\vgamma,\rho}\) satisfies the \(\br{A,C}\)-short-interval condition, if \(\vgamma:I:=\br{-A^{-C},A^{-C}}\simto\Gamma\) is \(A\)-admissible, and \(\rho\in C^2_c\br{\R}\) satisfies \(\supp\rho\subseteq A^{-2C}I\) and \(\nrm{\rho}_{C^2}=1\).
\end{definition}
We now formulate the key statement:

\subsection{Decay estimates and smoothing inequalities}
\begin{theorem}\label{thm_4ier_op_est}
    % There is a universal constant \(\sigma\in\br{0,1}\) such that 
    Fix a bump function \(\phi\in C^\infty_c\br{\R}\) satisfying the two properties:
    \begin{equation}\label{eq_phi_uni_even}
        \1_{\mr{-\frac{1}{2},\frac{1}{2}}} \leq \phi \leq \1_{\br{-1,1}},\quad
        \phi\br{-z}=\phi\br{z}.
    \end{equation}
    Define for \(\delta>0\) the auxiliary bump function:
    \begin{equation}\label{eq_psi_delta}
        \psi_\delta\br{\zeta}:=\phi\br{\frac{\zeta-1}{\delta}}.
    \end{equation}
    For \(A,C\gg 1\) and \(\br{\vgamma,\rho}\) satisfying the \(\br{A,C}\)-short-interval condition \ref{ass_short_I},
    consider the symbol given by:
    \begin{equation}\label{eq_m_vlambda_freq_trunc}
        m_{\vgamma,\rho,\vlambda}\br{\vxi}:=
        \int
            e\br{\vxi\cdot \vgamma\br{t}}
            \rho\br{t}
        dt
        \prod_{j=1,2}
        \psi_\delta\br{\xi_j/\lambda_j},\quad \delta:=A^{-4C}
    \end{equation}
    with the parameter \(\vlambda:=\br{\lambda_1,\lambda_2}\) satisfying the following conditions:
    \begin{equation}\label{eq_vlambda_cond}
        \abs{\vlambda}\geq A^C,\quad \mr{\vlambda}\in
        \cG_\Gamma\circ \vgamma\br{A^{-C}I}.
    \end{equation}
    The corresponding fiber-wise bilinear Fourier integral operator satisfies the following estimate:
    \begin{equation}\label{eq_4ier_op_est}
        \nrm{
            \cT_{m_{\vgamma,\rho,\vlambda}}\br{F,G}
        }_{L^1}
        \lesssim 
        A^{O\br{1}} \Xi^{\frac{3}{10}}_\Gamma
        \abs{\vlambda}^{-\frac{1}{20}}
        \log^{\frac{3}{20}}\br{\abs{\vlambda}}
        \nrm{F}_{L^2}
        \nrm{G}_{L^2}.
    \end{equation}
\end{theorem}
\begin{remark}
    Notice that \eqref{eq_vlambda_cond} implies \(1\underset{\vgamma}{\ll} \abs{\lambda_1}\underset{\vgamma}{\eqsim}\abs{\lambda_2}\).
    In other words, the above theorem is addressing the diagonal parts of the Littlewood-Payley decomposition. See \textbf{Lemma \ref{lem_osc_st_or_nst}} for the proof of this fact.
\end{remark}

Once we establish the above decay estimate, it becomes standard procedure to produce our main theorem regarding the Sobolev smoothing inequalities for \(\cA_{\vgamma,\rho}\):

\begin{theorem}\label{mainthm_global}
There are universal constants \(c,\varepsilon>0\) such that for all \(A\)-admissible \(\vgamma:I\simto\Gamma\subseteq\R^2\) and \(\rho\in C^2_c\br{\R}\) with \(\supp \rho\Subset I\), the Sobolev smoothing inequality \eqref{eq_soblev_smoothing_std} holds with:
\begin{equation*}
    C_{\vgamma,\rho,\varepsilon}\lesssim \max\br{A,\frac{1}{\dist\br{\supp \rho,I^c}}}^c
    \abs{I}
    \Xi^{\frac{3}{10}}_\Gamma
    \nrm{\rho}_{C^2}.
\end{equation*}
In other words, we have:
\begin{equation}\label{eq_mainthm_global}
    \nrm{\cA_{\vgamma,\rho}\br{f_1,f_2}}_{L^1}\lesssim \max\br{A,\frac{1}{\dist\br{\supp \rho,I^c}}}^c
    \abs{I}
    \Xi^{\frac{3}{10}}_\Gamma
    \nrm{\rho}_{C^2}
    \nrm{f_1}_{H^{\br{-\varepsilon,0}}}
    \nrm{f_2}_{H^{\br{0,-\varepsilon}}}.
\end{equation}
\end{theorem}

We also provide a geometric interpretation for \textbf{Theorem \ref{mainthm_global}}. Consider the following auxiliary map:
\begin{definition}\label{def_vchi_Gamma}
    Let \(N\geq 2\). Given an oriented \(C^N\) curve \(\Gamma\), we define the following \(C^{N-2}\) map:
    \begin{equation*}
        \vchi_\Gamma:\:
        \BR{
            p\in \Gamma
        \::\:
            \kappa\br{p}\neq 0
        }
        \to
        \R^2:\quad
        p
        \mapsto
        \br{\widehat{\vt}\cdot 
        \frac{\widehat{t}_1 \widehat{t}_2}{\kappa}}\br{p}=
        \widehat{\vt}\br{p}\cdot 
        \frac{\widehat{t}_1\br{p} \widehat{t}_2\br{p}}{\kappa\br{p}}
        ,
    \end{equation*}
    where \(\widehat{\vt}\br{p}:=\br{\widehat{t}_1\br{p},\widehat{t}_2\br{p}}\) denotes the directed unit tangent at \(p\in \Gamma\) and \(\kappa\br{p}\) the signed curvature at \(p\in\Gamma\).
\end{definition}
The following observation formally relates the map \(\vchi_\Gamma\) to the function \(\vchi_\vgamma\):
\begin{equation}\label{eq_char_gamma_2_Gamma}
    \vchi_\vgamma
    =
    \frac{\vgamma'}{\abs{\vgamma'}}\cdot\frac{\gamma'_1}{\abs{\vgamma'}}\cdot\frac{\gamma'_2}{\abs{\vgamma'}}/
    \br{\frac{\vgamma'\wedge\vgamma''}{\abs{\vgamma'}^3}}
    =
    \br{
        \widehat{\vt}\cdot 
        \frac{
            \widehat{t}_1 \widehat{t}_2
        }{
            \kappa
        }
    }\circ \vgamma
    =
    \vchi_\Gamma \circ \vgamma.
\end{equation}
As a consequence, we have \(\chr \Gamma:=\img\vchi_\Gamma=\img \vchi_\vgamma\). 
Furthermore, if \(\Gamma\) has non-vanishing curvature and \(\img\cG_\Gamma\subseteq \RP\setminus\BR{\Br{0:1},\Br{1:0}}\), the map \(\vchi_\Gamma:\Gamma\simto \chr \Gamma \subseteq \R^2\) is an embedding of \(\Gamma\) into \(\R^2\). We thus call \(\chr \Gamma\) the characteristic curve of \(\Gamma\) whenever such an embedding exists. With these notions, we may reformulate \textbf{Theorem \ref{mainthm_global}} and provide a qualitative but purely geometric statement:
\begin{theorem}\label{mainthm_geo}
    There is a universal constant \(\varepsilon>0\) such that given any \(C^4\) curve \(\Gamma\subseteq \R^2\) satisfying the following conditions:
    \begin{enumerate}
        \item \label{item_tang_dir}\(\img\cG_\Gamma\) contains neither \(\mr{1:0}\) nor \(\mr{0:1}\);
        \item \label{item_curvature}both \(\Gamma\) and \(\chr \Gamma\) have non-vanishing curvature;
        \item the uniform finiteness condition \(\Xi_\Gamma<\infty\) holds,
    \end{enumerate}
    the Sobolev smoothing inequality \eqref{eq_soblev_smoothing_std} holds with some finite \(0<C_{\vgamma,\rho,\varepsilon}<\infty\)
    % .
    % \begin{equation}\label{eq_mainthm_geo}
    %     \nrm{\cA_{\vgamma,\rho}\br{f_1,f_2}}_{L^1}\underset{\vgamma,\rho}{\lesssim}
    %     \nrm{f_1}_{H^{\br{-\varepsilon,0}}}
    %     \nrm{f_2}_{H^{\br{0,-\varepsilon}}}
    % \end{equation}
    % holds 
    for any \(C^4\) parametrization \(\vgamma:I\simto\Gamma\) and any bump function \(\rho\in C^2_c\br{I}\).
\end{theorem}
Indeed, take \(\Gamma\) satisfying the above three conditions and any \(C^4\) parametrization \(\vgamma:I\simto\Gamma\). We may relate the three conditions to the notion of \(A\)-admissibility for \(\vgamma\). On the one hand, we note that
\begin{equation*}
    \text{condition \ref{item_tang_dir}}\iff
    \widehat{t}_1\br{p}\neq 0,\,\widehat{t}_2\br{p} \neq 0,\quad \forall p\in\Gamma
    \iff
    \abs{\gamma'_1\br{t}},\,\abs{\gamma'_2\br{t}}>0
    ,\quad \forall t\in I.
\end{equation*}
Furthermore, since \(\vgamma:I\simto\Gamma\) and \(\vchi_\vgamma:I\simto \chr \Gamma\), we have
\begin{equation*}
    \text{condition \ref{item_curvature}}\iff
    \frac{\abs{\vgamma'\wedge\vgamma''}}{\abs{\vgamma'}^3}\br{t},\,
    \frac{\abs{\vchi'_\vgamma\wedge\vchi''_\vgamma}}{\abs{\vchi'_\vgamma}^3}\br{t}>0,\quad \forall t\in I
\end{equation*}
\begin{equation*}
    \iff
    \abs{\vgamma'\wedge\vgamma''}\br{t},\,
    \abs{\vchi'_\vgamma\wedge\vchi''_\vgamma}\br{t}>0,\quad \forall t\in I.
\end{equation*}
On the other hand, for any \(\rho\in C^2_c\br{\R}\) with \(\supp\rho\Subset I\), we may use the compactness of the support to find a finite open interval \(J\subseteq \R\) such that \(\supp \rho\Subset J \Subset I\), and restrict the domain of \(\vgamma\) to \(J\). 
Given the two equivalent statements of condition (1) and (2), the compactness of \(\overline{J}\) implies that there is a large enough \(A\gg 1\) such that \(\dif{\vgamma}_{J}\) is \(A\)-admissible. \textbf{Theorem \ref{mainthm_geo}} thus follows from \textbf{Theorem \ref{mainthm_global}} once we invoke the following key lemma:
\begin{lemma}[Monotonicity of \(\Xi_\Gamma\)]\label{lem_M_Gamma_mono}
    Given a \(C^1\) curve \(\Gamma\subseteq \R^2\) with injective Gauss map, if a curve \(\Gamma'\) is a subset of \(\Gamma\), the following inequality holds:
    \begin{equation*}
        \Xi_{\Gamma'}\leq \Xi_\Gamma+O\br{1}
        \lesssim \Xi_\Gamma.
    \end{equation*}
    Note that the second relation follows from the trivial lower bound \(1\leq \Xi_\Gamma\).
\end{lemma}

Return now to the discussion of the quantitative statement \textbf{Theorem \ref{mainthm_global}}. Due to the polynomial dependence on the \(A\)-admissibility condition, with little effort, we may generalize \textbf{Theorem \ref{mainthm_global}} to cover cases where the \(A\)-admissibility fails but in a controllable manner. One such example is the parametrization \(\vgamma\br{t}=\br{\cos t,\sin t}\) for the unit circular arc near \(t=0\).
To formulate the result, we first weaken the notion of \(A\)-admissibility by introducing the following function:
\begin{definition}[Admissibility function]\label{def_admi_func}
    Given \(C^4\) parametric curve \(\vgamma:I\simto\Gamma\), we define the function:
    \begin{equation}\label{eq_admi_func}
        A_\vgamma\br{t}:=
        \max\br{
            \max_{k=1,2,3,4}
                \abs{\vgamma^{\br{k}}}\br{t},
            \max_{j=1,2}
                \abs{\gamma'_j}^{-1}\br{t},
            \abs{\vgamma'\wedge\vgamma''}^{-1}\br{t},
            \abs{\vchi'_\vgamma\wedge\vchi''_\vgamma}^{-1}\br{t}
        },
    \end{equation}
    where we assign \(\infty\) as the value whenever one of the above expressions is not defined.
\end{definition}

Notice that the notion of \(A\)-admissibility coincides with the condition \(\nrm{A_\vgamma}_{L^\infty\br{I}}\leq A\). This fact suggests that we consider the following natural generalization:

\begin{definition}[Weak type control on \(A_\vgamma\)]\label{def_weak_failure}
    Given \(0<p\leq \infty\), we set:
    \begin{equation}\label{eq_weak_failure}
        \nrm{A_\vgamma}_{L^{p,\infty}}:=
        \nrm{A_\vgamma}_{L^{p,\infty}\br{I}}:=
        \sup_{\lambda>0}\lambda 
        \abs{
            \BR{t\in I
            \::\:
                A_\vgamma\br{t}>\lambda
            }
        }^{\frac{1}{p}}
        .
    \end{equation}
Note that for all \(0<p\leq q\leq \infty\), we always have \(1 \underset{I,p}{\lesssim} \nrm{A_\vgamma}_{L^{p,\infty}}\underset{I,p,q}{\lesssim} \nrm{A_\vgamma}_{L^{q,\infty}}\).
\end{definition}

Since the above quantity provides no structural information about the set of poor admissibility, we remedy this issue via a notion controlling the number of connected components:

\begin{definition}[Fluctuation control on \(A_\vgamma\)]\label{def_A_gamma_fluc}
Given \(0\leq \alpha<\infty\), we define:
\begin{equation}\label{eq_A_gamma_fluc}
    \abs{A_\vgamma}_\alpha:=
    \sup_{\lambda>1}
    \lambda^{-\alpha}
    C\br{
        \BR{
            t\in I\::\:
            A_\vgamma\br{t}>\lambda
        }
    }
    % ,\quad
    % \abs{A_\vgamma}_\infty:=1
    .
\end{equation}
Note that for all \(0\leq \alpha\leq \beta <\infty\), we always have \(1\leq \abs{A_\vgamma}_\beta\leq \abs{A_\vgamma}_\alpha\).
\end{definition}

With the above notions, we provide Sobolev smoothing inequalities adapted to this general setting.

\begin{theorem}\label{mainthm_global_w_vanishing_n_fluc}
For all \(0<p \leq \infty\) and \(0\leq \alpha<\infty\) there are \(c_\alpha>0\) and \(\varepsilon_{p,\alpha}>0\) such that for all \(C^4\) parametric curves \(\vgamma:I\simto\Gamma\subseteq\R^2\) with injective Gauss map, and \(\rho\in C^2_c\br{\R}\) with \(\supp \rho\Subset I\), the Sobolev smoothing inequality \eqref{eq_soblev_smoothing_std} holds with:
    \begin{equation*}
        C_{\vgamma,\rho,\varepsilon_{p,\alpha}}
        \underset{p,\alpha}{\lesssim}
        \max\br{
            \nrm{A_\vgamma}_{L^{p,\infty}},
            \abs{A_\vgamma}_\alpha,
            \frac{1}{\dist\br{\supp \rho,I^c}}
        }^{c_\alpha}
        \ang{\abs{I}}
        \Xi^{\frac{3}{10}}_\Gamma
        \nrm{\rho}_{C^2}
    \end{equation*}
\end{theorem}
\begin{remark}\label{rmk_p_infty_case_comments}
    When \(A:=\nrm{A_\vgamma}_{L^{\infty,\infty}}:=\nrm{A_\vgamma}_{L^\infty\br{I}}<\infty\), the parametric curve \(\vgamma:I\simto\Gamma\) is \(A\)-admissible. Therefore the endpoint \(p=\infty\) case follows directly from \textbf{Theorem \ref{mainthm_global}}. 
\end{remark}

\subsection{Finiteness and uniformity via o-minimal theory}

Up to this point, the finiteness of \(\Xi_\Gamma\) and \(\abs{A_\vgamma}_\alpha\) are assumed either explicitly or implicitly. Although one may check for specific examples: 
\begin{equation*}
    \vgamma\br{t}=\br{t,t^2},\, \br{\cos t,\sin t},\,\br{P_1\br{t},P_2\br{t}}
\end{equation*}
that \(\Xi_\Gamma,\abs{A_\vgamma}_0<\infty\), it is more desirable to have an easily checkable criteria for the finiteness. Below, we shall provide such criteria via notions from o-minimal theory\footnote{For readers' convenience, we include the required background knowledge in \textsc{Appendix \ref{sec_o_min}}.}

Fix an o-minimal expansion \(\cR\) of \(\ang{\R,+,\cdot,<}\). Let \(\vgamma:I\simto\Gamma\) be definable (in \(\cR\)). If \(\cG_\Gamma\) is injective, both \(\Xi_\Gamma\) and \(\abs{A_\vgamma}_0\) are finite.
In fact, we will prove the following stronger statements:

\begin{lemma}[A uniform finiteness property for \(\Xi_\Gamma\)]\label{lem_uni_fini_prop_for_Gamma}
    Consider a definable family \(\BR{\Gamma_\sigma}_{\sigma\in\Sigma}\) of \(C^1\) curves in \(\R^2\).
    If for all \(\sigma\in\Sigma\) the Gauss map \(\cG_{\Gamma_\sigma}\) is injective,
    % induces a topological embedding \(\cG_{\Gamma_\sigma}:\cG_\sigma\hookrightarrow \RP\),
    the following uniform bound holds:
    \begin{equation*}
        \sup_{\sigma\in\Sigma}\Xi_{\Gamma_\sigma}<\infty.
    \end{equation*}
\end{lemma}

\begin{lemma}[A uniform finiteness property for \(\abs{A_\vgamma}_0\)]\label{lem_uni_fini_prop_for_A_vgamma}
    Consider a definable family \(
        \BR{\vgamma_\sigma:\R\supset I_\sigma\simto\Gamma_\sigma\subseteq \R^2}_{\sigma\in\Sigma}
    \) of \(C^4\) parametric curves.
    The following uniform bound holds:
    \begin{equation*}
        \sup_{\sigma\in\Sigma}\abs{A_{\vgamma_\sigma}}_0<\infty.
    \end{equation*}
\end{lemma}

These lemmas demonstrate certain uniformity regarding the boundedness of the two quantities. As a direct consequence, we have the following variants of the Sobolev smoothing inequalities:

\begin{theorem}[Smoothing inequalities for families of curves]\label{thm_smoothing_o_mini_uni}
    Consider a definable family 
    \begin{equation*}
        \mathfrak{F}:=\BR{\vgamma_\sigma:\R\supset I_\sigma\simto\Gamma_\sigma\subseteq \R^2}_{\sigma\in\Sigma}
    \end{equation*}
    of \(A_\sigma\)-admissible parametric curves.
    For \(\rho_\sigma\in C^\infty_c\br{\R}\) satisfying \(\supp \rho_\sigma\Subset I_\sigma\), the Sobolev smoothing inequalities \eqref{eq_soblev_smoothing_std} hold uniformly for all \(\sigma\in \Sigma\) with
    \begin{equation*}
        C_{\vgamma_\sigma,\rho_\sigma,\varepsilon}
        \underset{\mathfrak{F}}{\lesssim}
        \max\br{A_\sigma,\frac{1}{\dist\br{\supp \rho_\sigma,I^c_\sigma}}}^c
        \abs{I_\sigma}
        \cdot
        \nrm{\rho_\sigma}_{C^2}
    \end{equation*}
    for the universal constant \(c,\varepsilon>0\) as in \textbf{Theorem \ref{mainthm_global}}.
\end{theorem}

\begin{theorem}[Strong smoothing inequalities for families of curves]\label{thm_smoothing_o_mini_uni_w_vanishing_n_fluc}
    Consider a definable family \(
        \mathfrak{F}:=\BR{\vgamma_\sigma:\R\supset I_\sigma\simto\Gamma_\sigma\subseteq \R^2}_{\sigma\in\Sigma}
    \) of \(C^4\) parametric curves whose Gauss maps are injective.
    For all \(0<p \leq \infty\) and \(\rho_\sigma\in C^2_c\br{\R}\) with \(\supp \rho_\sigma\Subset I_\sigma\), the Sobolev smoothing inequalities \eqref{eq_soblev_smoothing_std} hold uniformly for all \(\sigma\in\Sigma\) with
    \begin{equation*}
        C_{\vgamma_\sigma,\rho_\sigma,\varepsilon_{p,0}}
        \underset{p,\mathfrak{F}}{\lesssim}
        \max\br{
            \nrm{A_{\vgamma_\sigma}}_{L^{p,\infty}},
            \frac{1}{\dist\br{\supp \rho_\sigma,I^c_\sigma}}
        }^{c_0}
        \ang{\abs{I_\sigma}}
        \nrm{\rho_\sigma}_{C^2}
    \end{equation*}
    for \(c_0\) and \(\varepsilon_{p,0}\) defined as in \textbf{Theorem \ref{mainthm_global_w_vanishing_n_fluc}}.
\end{theorem}

Below, we provide a list of o-minimal expansions of \(\R_{\mathrm{alg}}:=\ang{\R,+,\cdot,<}\) and some examples of curves definable in the individual structures:
\begin{itemize}
    \item \(\R_{\mathrm{alg}}\) consists of semi-algebraic sets. Polynomial curves \(\vgamma\br{t}=\br{P_1\br{t},P_2\br{t}}\) are definable.
    % The o-minimality follows from the Tarski-Seidenberg theorem.
    \item \(\R^\R\) expands with power functions: \(t^\alpha,\:\alpha\in\R\). See \cite{MR1278550} for details. Specifically, \(\vgamma\br{t}:=\br{\abs{t}^{\alpha_1},\abs{t}^{\alpha_2}}\) is definable.
    
    \item \(\R_{\mathrm{an}}\) expands with all truncated real analytic functions. See
    \cite{MR245831,MR425152,MR972342,MR1374342,MR1389958} for details. For example, \(\vgamma\br{t}:=\1_{\mr{-\pi/2,\pi/2}}\br{t}\br{\cos t,\sin t}\) is definable.
    
    \item \(\R_{\mathrm{exp}}\) expands with \(e^t\). See \cite{MR1108621,MR996349,MR1398816} for details. For example, \(\vgamma\br{t}:=\br{\log\br{t},e^t}\) is definable.
    % % or perform the Pfaffian closure
    
    \item \(\R_{\mathrm{an,exp}}\) combines the previous two bullet points. See \cite{MR1289495} for details.
\end{itemize}
%   check refs from slides
Notably, the sufficiency of the condition in \textbf{Theorem \ref{thm_real_ana_smoothing}} is one direct consequence of the o-minimality of \(\R_{\operatorname{an}}\) and \textbf{Theorem \ref{thm_smoothing_o_mini_uni_w_vanishing_n_fluc}}. Details are given in \textsc{Section \ref{subsec_pf_thm_real_ana_smoothing}}

\subsection{Applications}

We next describe four applications of the smoothing estimates: singular integrals, maximal
operators, spherical maximal operators, and the detection of corner-type patterns. To place these results
in context, recall that the systematic study of singular and maximal integrals along submanifolds was initiated
by foundational work of Stein, Wainger, Nagel, Christ, Seeger, and Wright \cite{MR508453,MR714828,MR819558,MR816390,MR1046743,MR1726701}. In parallel, Coifman–Meyer multiplier theory \cite{MR511821,MR518170} and Lacey–Thiele’s resolution of Calderon’s conjecture \cite{MR1491450,MR1689336} established central tools and model problems for modern multilinear harmonic analysis.

\textit{Singular Integral.} A broad class combining multilinearity, singular kernels, and Radon-type behavior is
given by
\begin{equation}
    H(f)(x)=\int_{\mathbb{R}^{d}}\prod_{j=1}^{n}f_{j}(\Pi_{j}(x,t))K(t)dt\, ,
\end{equation}
\begin{question}
   For which exponents $(p,p_{1},\cdots,p_{n})\in \mathbb{R}^{n+1}$ does the estimate hold
    \begin{equation}
        \|H(f)\|_{L^{p}}\lesssim \prod_{j=1}^{n}\|f_{j}\|_{L^{p_{j}}} \quad ?
    \end{equation}
\end{question}
The available methods divide naturally according to whether the projection maps $\Pi_{j}$ are flat(linear) or curved. 

In the flat setting, the natural multilinear extension of Calderon–Zygmund singular integrals is the class of
singular Brascamp–Lieb forms \cite{MR4390229}. For instance, the bilinear Hilbert transform is a special case of a one-dimensional trilinear singular Brascamp-Lieb form. Demeter and Thiele \cite{MR2597511} systematically analyzed the two-dimensional trilinear
forms and established bounds in most cases. A notable special case is the twisted paraproduct, whose boundedness was established by Kovač \cite{MR2990138} and later generalized by Durcik \cite{MR3488377,MR3683098}. Among two-dimensional
trilinear singular Brascamp–Lieb forms, the principal remaining case with no known boundedness at any
exponent is the form associated with the triangular Hilbert transform:
\begin{equation}
   \int_{\mathbb{R}}f_{1}(x+t,y)f_{2}(x,y+t)\frac{dt}{t}\, .
\end{equation}
A broader classification of trilinear singular Brascamp–Lieb forms was obtained by Becker, Durcik, and the
second author \cite{becker2024trilinearsingularbrascampliebintegrals} using quiver representations. At the next level of multilinearity, even the trilinear Hilbert
transform, the basic quadrilinear singular Brascamp–Lieb form, remains open.

In the curved setting, oscillatory integral methods become available. The boundedness of the bilinear Hilbert transform along a parabola was established by Li \cite{MR3068544}, and generalizations to other non-flat curves were subsequently pursued by Lie \cite{MR3337797,MR3763348}. Trilinear Hilbert transform along moment curve was studied by Hu and Lie \cite{hu2023curvedtrilinearhilberttransform,hu2025boundednesscurvedtrilinearhilbert}.

Moving to higher dimension, 
let \(\vgamma:=\br{\gamma_1,\gamma_2}:\R\to \R^2\) be a parametric curve in the plane. We define the triangular Hilbert transform associated to \(\vgamma\) as\footnote{We henceforth omit the \(\operatorname{p.v.}\) notation whenever the proper interpretation is clear from the context.}
\begin{equation*}
    T_\vgamma\br{f_1,f_2}\br{\vx}:=\pv
    % \int 
    %     \prod_{j=1,2}
    %         f_j\br{\vx+\ve_j\gamma_j\br{t}}
    % \frac{dt}{t}.
    \int
        f_1\br{x_1+\gamma_1\br{t},x_2}
        f_2\br{x_1,x_2+\gamma_2 \br{t}}
    \frac{dt}{t}.
\end{equation*}
Under this notation, the triangular Hilbert transform is $T_{\br{t,t}}\br{f_1,f_2}\br{\vx}$.  Christ, Durcik, and Roos \cite{MR4295087} obtained the first
nontrivial boundedness result in this corner-type family by treating the parabolic curve \(\vgamma\br{t}:=\br{t,t^2}\). A crucial component in obtaining these bounds is the smoothing inequality. In \cite{MR4295087}, the smoothing inequality is derived via sublevel set estimates. Chen and Guo \cite{MR4800914} later treated polynomial pairs \(\vgamma:=\br{P_1,P_2}\) with distinct leading and lowest degrees by related methods. Our approach is different: inspired by the sparse–uniform analysis of Gaitan and Lie \cite{MR4916715}, it applies to broad geometric classes of curves.

We will utilize the uniformity statement for our Sobolev smoothing inequality provided by \textbf{Theorem \ref{thm_smoothing_o_mini_uni}} to generalize the above results.

To formulate our statement, we first introduce the following preliminary notions:

\begin{definition}[Symmetric punctured neighborhood of \(0\) and \(\infty\)]
    For \(r>0\), we define:
    \begin{equation*}
        U_0\br{r}:=\BR{t\in\R\::\: 0<\abs{t}<r},\quad
        U_\infty\br{r}:=\BR{t\in\R\::\: r<\abs{t}<\infty}.
    \end{equation*}
    Sets of the form \(U_0\br{r}\) (resp. \(U_\infty\br{r}\)) will be referred to as symmetric punctured neighborhoods of \(0\) (resp. \(\infty\)). When the specific value of \(r>0\) is irrelevant, we simply write \(U_0\) and \(U_\infty\).
\end{definition}

\begin{definition}[Positive homogeneity]\label{def_pos_homo}
    A function \(h:\R\setminus\BR{0}\to\R\) is positively homogeneous of degree \(\alpha\in\R\) if \(h\br{st}=s^\alpha h\br{t}\) for all \(t\neq 0\) and \(s>0\).
\end{definition}

\begin{definition}[Asymptotic homogeneity]\label{def_asym_homogeneous}
    Let \(\ast\in\BR{0,\infty}\). A function \(f:\R\setminus\BR{0}\to\R\) is asymptotically homogeneous near \(\ast\)
    if it can be approximated 
    by a positively homogeneous function \(h:\R\setminus\BR{0} \to \R\setminus\BR{0}\) of non-zero degree via the relation \(f\br{t}/h\br{t}=1+O(\abs{t}^{\delta_\ast})\)
    for some \(\delta_\ast\in\R\) satisfying \(\lim_{\abs{t}\to\ast}\abs{t}^{\delta_\ast}=0\).
\end{definition}

\begin{definition}\label{def_asym_def}
    Let \(\ast\in\BR{0,\infty}\). Given an o-minimal expansion \(\cR\) of \(\ang{\R,+,\cdot,<}\), we say a function \(f:\R\setminus\BR{0}\to\R\) is definable near \(\ast\) if \(\dif{f}_{U_\ast}:U_\ast\to\R\) is definable for some \(U_\ast=U_\ast\br{r}\). In other words, \(f\) is germ equivalent to a definable function at \(\ast\).
\end{definition}

Given \(\vgamma=\br{\gamma_1,\gamma_2}:\R\setminus\BR{0}\to\R^2\), we introduce the following two easily checkable conditions:
\begin{equation}\label{eq_THT_ass_curvature_neq_0}
    \forall \ast\in\BR{0,\infty},\quad
    \liminf_{\abs{t}\to\ast}
    \abs{t}/
    \abs{
        \frac{
            \gamma'_1\gamma'_2
        }{ 
            \vgamma'\wedge\vgamma''
        }
    }
    \br{t}
    >0
\end{equation}

\begin{equation}\label{eq_THT_ass_growth_match}
    \forall \ast\in\BR{0,\infty},\quad
    \liminf_{\abs{t}\to \ast}\abs{\gamma_1\br{t}}=\infty \iff \liminf_{\abs{t}\to \ast}\abs{\gamma_2\br{t}}=\infty.
\end{equation}
With the above notions, our result regarding the operator \(T_\vgamma\) can be formulated as below:

\begin{theorem}[Triangular Hilbert transform along asymptotically homogeneous definable curves]\label{thm_THT_asym_homo_curves}
    Let \(\cR\) be an o-minimal expansion of \(\ang{\R,+,\cdot,<}\). For \(j=1,2\), also let \(\gamma_j:\R\setminus\BR{0}\to\R\) be asymptotically homogeneous and definable near both \(0\) and \(\infty\). If \eqref{eq_THT_ass_curvature_neq_0} and \eqref{eq_THT_ass_growth_match} hold, the operator \(T_\vgamma\) satisfies the following estimates:
    \begin{equation}\label{eq_THT_Lp_bdds}
        \nrm{T_\vgamma\br{f_1,f_2}}_{L^{p_3}}\underset{p_1,p_2,\vgamma}{\lesssim}
        \nrm{f_1}_{L^{p_1}}
        \nrm{f_2}_{L^{p_2}},\quad
        \forall
        p_1,p_2\in\br{1,\infty},\quad
        \frac{1}{2}<\frac{1}{p_3}=\frac{1}{p_1}+\frac{1}{p_2}\leq 1.
    \end{equation}
    % extends to a bounded operator \(L^{p_1}\otimes L^{p_2}\to L^p\) for \(p_1,p_2\in\br{1,\infty}\) satisfying \(\frac{1}{2}<\frac{1}{p}=\frac{1}{p_1}+\frac{1}{p_2}\leq 1\).
\end{theorem}

The asymptotic notions involved in \textbf{Theorem \ref{thm_THT_asym_homo_curves}} are motivated by the following observation: the boundedness of \(T_\vgamma\) is dictated by the asymptotic behavior of the parametric curve \(\vgamma\) near \(0\) and \(\infty\). To be precise, the truncated operator satisfies the following trivial bound:
\begin{equation}\label{eq_T_Gamma_log_triv}
    \nrm{
        \int_{r_0\leq\abs{t}\leq r_\infty}
            f_1\br{x_1+\gamma_1\br{t},x_2}
            f_2\br{x_1,x_2+\gamma_2\br{t}}
        \frac{dt}{t}
    }_{L^{p_3}\br{d\vx}}
    \leq
    2
    \log\br{r_\infty/r_0}
    \nrm{f_1}_{L^{p_1}}\nrm{f_2}_{L^{p_2}}
\end{equation}
for \(p_j\in\mr{1,\infty}\) with \(\frac{1}{p_3}=\frac{1}{p_1}+\frac{1}{p_2}\).
As a result, studying the boundedness of \(T_\vgamma\) reduces to understanding truncated operators:
\begin{equation}\label{eq_T_Gamma_trunc_approx}
    T_{\vgamma,\ast}\br{f_1,f_2}\br{\vx}\approx
    \int_{U_\ast\br{r_\ast}} 
        f_1\br{x_1+\gamma_1\br{t},x_2}
        f_2\br{x_1,x_2+\gamma_2\br{t}}
    \frac{dt}{t},\quad \ast\in \BR{0,\infty}
\end{equation}
of which the exact definition will be introduced in \textsc{Section \ref{subsec_tht_along_curve}}.

In the cases \(\vgamma:=\br{P_1,P_2}\) considered in \cite{MR4800914} with \(P_j\) being polynomials of the form:
\begin{equation}\label{eq_poly_form_Pj}
    P_j\br{t}= 
    a_{d_j}t^{d_j}+
    a_{d_j-1}t^{d_j-1}+
    \cdots +
    a_{e_j+1}t^{e_j+1}+
    a_{e_j}t^{e_j},\quad
    a_{d_j},a_{e_j}\neq 0,
\end{equation}
we see that \(P_j\br{t}\) is asymptotically homogeneous near \(0\) (resp. \(\infty\)) of degree \(e_j\) (resp. \(d_j\)).
With \textbf{Theorem \ref{thm_THT_asym_homo_curves}}, we allow real-valued exponents in \eqref{eq_poly_form_Pj} given suitable interpretation of the power function \(t^\alpha\) for \(\alpha\in\R\) and some other well-behaved perturbations.

More generally, the asymptotic homogeneity assumption on \(\gamma_j\) is natural in the case where the o-minimal expansion \(\cR\) of \(\ang{\R,+,\cdot,<}\) is polynomially bounded (see \textbf{Definition \ref{def_omin_poly_bdd}}). The detailed discussion is included in \textsc{Section \ref{subsubsec_poly_bdd_THT}}. As one direct application, we provide the boundedness result analogous to \textbf{Theorem \ref{thm_THT_asym_homo_curves}} for real analytic \(\vgamma\) (see \textbf{Corollary \ref{cor_THT_real_ana}}.)

\textit{Maximal Operator.} The bilinear maximal operator
\begin{equation}
    \underset{r>0}{\operatorname{sup}}\frac{1}{2r}\int_{-r}^{r}|f_{1}(x+t)f_{2}(x-t)|dt
\end{equation}
was initially studied by Lacey \cite{MR1745019}, and subsequently generalized and investigated by Demeter, Tao, and Thiele \cite{MR2403711} in the connection with pointwise ergodic theorems. Curved analogues have been studied by Li, and Xiao \cite{MR3538147} and Gaitan and Lie \cite{MR4131966}. In our corner-type setting, for \(\vgamma:\R\setminus\BR{0}\to\R^2\), define 
\begin{equation}
    M_\vgamma\br{f_1,f_2}\br{\vx}:=
    \sup_{r>0}\abs{
    \fint^{r}_{-r}
        f_1\br{x_1+\gamma_1\br{t},x_2}
        f_2\br{x_1,x_2+\gamma_2\br{t}}
    dt}.
\end{equation}
\begin{theorem}\label{thm_max}
     Let \(\vgamma\) be as given in \textbf{Theorem \ref{thm_THT_asym_homo_curves}}. The operator \(M_\vgamma\) satisfies the following estimates
     \begin{equation*}
         \nrm{M_\vgamma\br{f_1,f_2}}_{L^{p_3}}\underset{\vgamma,p_j}{\lesssim}\nrm{f_1}_{L^{p_1}}\nrm{f_2}_{L^{p_2}}
     \end{equation*}
     for the following Banach range
     \begin{equation}\label{eq_M_vgamma_range}
         1<p_1,p_2\leq\infty,\quad
         \frac{1}{p_3}=\frac{1}{p_1}+\frac{1}{p_2}\leq 1.
     \end{equation}
\end{theorem}
The endpoint $p_{3}=1$ is already nontrivial. Determining the optimal quasi-Banach range $p_{3}<1$ lies beyond
the scope of the present paper.

\textit{Spherical Maximal Operator.} Another maximal operator of interest arises from taking certain supremum over the spherical averages:
\begin{equation}\label{eq_spher_avg}
    \avg_t\br{f}\br{\vx}:=\int_{S^{d-1}} f\br{\vx-t\vy}d\cH^{d-1}\br{\vy},\quad
    \vx\in \R^d.
\end{equation}
Stein \cite{MR420116} introduced the classical formulation of the spherical maximal operator
\begin{equation}\label{eq_spher_max}
    \cS\br{f}\br{\vx}:=\sup_{t>0}\abs{\avg_t\br{f}\br{\vx}},\quad f\in L^p\br{\R^d}
\end{equation}
and proved the \(L^p\) boundedness result for $d\geq 3$ and $p>\frac{d}{d-1}$. The more difficult $d=2$ case was resolved a decade later by Bourgain \cite{MR812567}. These \(L^p\) boundedness results have direct connections to almost everywhere convergence of solutions to certain PDEs for \(L^p\) initial data. Subsequent work has considered multilinear averages and restrictions of the dilation
parameter to lacunary or fractal sets.

By replacing \(\avg_t\) with the following bilinear analog:
\begin{equation}\label{eq_bi_spher_avg}
    \bavg_t\br{f_1,f_2}\br{\vx}:=
    \int_{S^{2d-1}}
        f_1\br{\vx- t\vy_1}
        f_2\br{\vx-t\vy_2}
    d\cH^{2d-1}\br{\vy_1,\vy_2},\quad \vx\in\R^d
\end{equation}
one obtains the bilinear spherical maximal operator:
\begin{equation}\label{eq_bi_spher_max}
    \cB\cS\br{f_1,f_2}\br{\vx}:=
    \sup_{t>0}
        \bavg_t\br{f_1,f_2}\br{\vx},\quad
    f_j\in L^{p_j}\br{\R^d}.
\end{equation}
The operator \eqref{eq_bi_spher_max} was first studied in \cite{MR3188026}. 
Following \cite{MR3917731,MR4299160}, Jeong and Lee \cite{MR4103874}   obtained the sharp range for $d\geq 2$ by a slicing argument. The one-dimensional case was resolved independently by Christ and Zhou \cite{MR4776384} and Dosidis and Ramos \cite{MR4832983}.

Restricting the dilation parameter to dyadic scales gives the lacunary spherical maximal operator
\begin{equation}\label{eq_spher_max_lac}
    \cS_\lac \br{f}\br{\vx}:=
    \sup_{t\in 2^\Z}
        \avg_t\br{f}\br{\vx},\quad
    \cB\cS_\lac\br{f_1,f_2}\br{\vx}:=
    \sup_{t\in 2^\Z}
        \bavg_t\br{f_1,f_2}\br{\vx}\, .
\end{equation}

% A variant of the spherical maximal operator is obtained by taking the supremum over dyadic number which is called lacunary spherical maximal operator. 
Compared with their full counterparts \eqref{eq_spher_max}, \eqref{eq_bi_spher_max}, lacunary spherical maximal operators typically admit a larger boundedness
range. The linear theory was developed in \cite{MR537803,bams/1183540517}. In the bilinear setting, Borges and Foster \cite{MR4800832} established
bounds for $d\geq 2$, while Christ and Zhou \cite{MR4776384} treated $d\geq 1$. 

Bhojak, Choudhary, and Shrivastava \cite{bhojak2024twistedbilinearsphericalmaximal} subsequently introduced twisted/corner-type variants by replacing \(\avg_t\) with
\begin{equation}
    \cavg_t\br{f_1,f_2}\br{\vx_1,\vx_2}
    :=\int_{S^{2d-1}}
        f_1\br{\vx_1-t\vy_1,\vx_2}
        f_2\br{\vx_1,\vx_2-t\vy_2}
    d\cH^{2d-1}\br{\vy_1,\vy_2},\quad
    \vx_j\in\R^d.
\end{equation}
They study both the full and lacunary maximal operators associated with these averages:
\begin{equation}\label{eq_corner_spher_max_ful}
    \cC\cS\br{f_1,f_2}\br{\vx_1,\vx_2}:=
    \sup_{t>0}
        \cavg_t\br{f_1,f_2}\br{\vx_1,\vx_2},
\end{equation}
\begin{equation}\label{eq_corner_spher_max_lac}
    \cC\cS_\lac\br{f_1,f_2}\br{\vx_1,\vx_2}:=
    \sup_{t\in 2^\Z}
        \cavg_t\br{f_1,f_2}\br{\vx_1,\vx_2}.
\end{equation}
% Later, the boundedness of the full and lacunary corner type  bilinear spherical maximal operators were study by Bhojak, Choudhary, and Shrivastava \textcolor{red}{cite}. 
Related work treats more general dilation sets \cite{MR1635404,MR1329463,MR1471145,MR1955209,MR4229592,MR4801507}. In the linear setting, Roos and Seeger \cite{MR4621382}
characterized the $L^{p}$ improving region up to endpoints. Variants associated with surfaces other than
spheres appear in \cite{MR155146,MR4375750,MR4747299,MR4792549,MR4801507,MR4363760,MR4609073}.
For a recent survey of bilinear spherical averages, see Borges \cite{borges2026surveybilinearsphericalaverages}.

Our smoothing inequalities yield the following corner-type lacunary maximal estimate.
Given \(\vgamma:I\to \R^2\), define the following maximal operator
\begin{equation*}
    \cM_\vgamma\br{f_1,f_2}\br{\vx}:=
    \sup_{k\in\Z}
    \abs{
    \int
        f_1\br{x_1+2^k\gamma_1\br{t},x_2}
        f_2\br{x_1,x_2+2^k\gamma_2\br{t}}
        \rho\br{t}
    dt}.
\end{equation*}
\begin{theorem}\label{thm_spherMax}
    Let \(\vgamma:I\to\R^2\) be a real analytic function with \(\img \vgamma\) not contained in the four classes \eqref{eq_line}, \eqref{eq_exp}, \eqref{eq_log}, or \eqref{eq_psi_line}. The operator \(\cM_\vgamma\) satisfies the following estimate
    \begin{equation*}
        \nrm{\cM_\vgamma\br{f_1,f_2}}_{L^{p_3}}
        \underset{\vgamma,p_j}{\lesssim}
        \nrm{f_1}_{L^{p_1}}
        \nrm{f_2}_{L^{p_2}}
    \end{equation*}
    for the same range as in \eqref{eq_M_vgamma_range}.
\end{theorem}

\textit{Corner-Type Pattern.} Finally, the smoothing inequality can be used to detect configurations of the form $(x,y),(x+\gamma_{1}(t),y),(x,y+\gamma_{2}(t))$ in subsets of Euclidean space. This problem belongs to the broader study
of arithmetic progression initiated by Roth’s and Szemeredi’s theorems. Major quantitative
advances include work of Gowers \cite{MR1844079}, Green and Tao \cite{MR3731312}, Bloom and Sisask \cite{bloom2023improvementkelleymekaboundsthreeterm}, Kelly and Meka \cite{MR4720301}, Leng, Sah, and Sawhney \cite{leng2024improvedboundsszemeredistheorem}.

A complementary direction replaces linear progressions by nonlinear patterns $x,x+\gamma_{1}(t),\cdots,x+\gamma_{m}(t)$. Bergelson and Leibman \cite{MR1325795} obtained a qualitative result for polynomial progressions, and a substantial recent literature
has developed quantitative versions in several ambient settings.

Over finite field, Bourgain and Chang \cite{MR3704938} obtained a quantitative bound for $x,x+t,x+t^{2}$. The result was extended to more general three-term polynomial progression by Peluse \cite{MR3874848}, Dong, Li, and Sawin \cite{MR4179774}. Later, a breakthrough was made by Peluse \cite{MR3934588}, giving a quantitative bound for arbitrary length of polynomial progression $x,x+\gamma_{1}(t),\cdots,x+\gamma_{m}(t)$ where $\gamma_{1},\cdots,\gamma_{m}$ are linearly independent polynomials.

Over the integers, Peluse and
Prendiville \cite{MR4824730} established the first quantitative three-term result, and Peluse \cite{MR4199235} extended it to arbitrary-length
polynomial progressions with distinct degrees. 

In the real setting, Bourgain \cite{MR853455} obtained a quantitative bound for $x,x+t,x+t^{d}$. This result was extended by Durcik, Guo and Roos \cite{MR3939567} to $x,x+t,x+\gamma_{2}(t)$ for polynomial $\gamma_{2}$ and by Chen, Guo, and Li \cite{MR4243369} to $x,x+\gamma_{1}(t),x+\gamma_{2}(t)$ where $\gamma_{1},\gamma_{2}$ have distinct degrees. Krause, Mirek, Peluse, and Wright \cite{MR4831142} recently obtained
arbitrary-length results with distinct polynomial degrees over general local fields. 

The PET induction and degree-lowering mechanisms in \cite{MR4831142} are powerful enough to handle arbitrary multilinearity,
but they exploit the algebraic structure of polynomial curves. Their extension to general analytic
curves is therefore not immediate. Our method has a complementary strength and limitation: in the bilinear
corner setting, it replaces polynomial structure by geometric conditions on the curve, but it does not presently
extend to arbitrary multilinearity. Understanding how these method relates to each other and whether these approaches can be combined is beneficial and a natural
direction for future work. 

Recently, there are also several results considering finding nonlinear progression in fractal sets \cite{zhu2024quadraticroththeoremsets,MR4609784,MR4966567,MR4418720,MR3672917,krause2025polynomialprogressionsinsidesets,hong2025polynomialszemeredisetslarge}. 

In higher dimensions, the corresponding corner-type configuration is $x,x+e_{1}\gamma_{1}(t),\cdots,x+e_{m}\gamma_{m}(t)$ for $x\in \mathbb{R}^{m}$ where $e_{i}$ is the $i$-th standard unit vector.  In the two-dimensional case, our smoothing inequality
gives the following quantitative existence result.
\begin{theorem}\label{thm_pattern}
    Let \(\vgamma:=\br{\gamma_1,\gamma_2}\) be real analytic in a neighborhood \(I\) of \(0\) with \(\vgamma\br{0}=\vnull\). If none of the four classes of forbidden sets \eqref{eq_line}, \eqref{eq_exp}, \eqref{eq_log}, and \eqref{eq_psi_line} contains \(\vgamma\br{I}\), there exists \(C_\vgamma,\,c_\vgamma>0\) such that for any \(E\subseteq \Br{0,1}^2\) with positive measure, there exists \(\br{x_1,x_2}\in\R^2\) such that
    \begin{equation}\label{eq_corner_pattern_in_E}
        \br{x_1,x_2},\quad\br{x_1+\gamma_1\br{t},x_2},\quad
        \br{x_1,x_2+\gamma_2\br{t}}\in E
    \end{equation}
    for some \(t\in I\) satisfying the gap estimate
    \begin{equation}\label{eq_corner_pattern_gap_est}
        t>\exp\br{-\exp\br{C_\vgamma \abs{E}^{-c_\vgamma}}}.
    \end{equation}
    % For $\varepsilon \in (0,1)$  and $S\subseteq  [0,1]^{2}$ with $|S|>\varepsilon$, there exist $(x,y),(x+\gamma_{1}(t),y),(x,y+\gamma_{2}(t))\in S$ with $t>\operatorname{exp}(-\varepsilon^{-c})$ for some $c>0$.
\end{theorem}

\begin{remark}
    We note that whenever
    \begin{equation*}
        \lim_{t\to 0}\frac{
            \log\gamma'_1\br{t}-\log\gamma'_2\br{t}
        }{
            \log t
        }\neq 0,
    \end{equation*}
    or equivalently \eqref{eq_THT_ass_curvature_neq_0} with \(\ast=0\):
    \begin{equation*}
        \lim_{t\to 0}\frac{
            t\vgamma'\br{t}\wedge\vgamma''\br{t}
        }{
            \gamma'_1\br{t}\gamma'_2\br{t}
        }\neq 0,
    \end{equation*}
     the gap estimate \eqref{eq_corner_pattern_gap_est} can be further improved with one less exponentiation
     \begin{equation*}
         t>\exp\br{-C_\vgamma\abs{E}^{-c_\vgamma}}.
     \end{equation*}
\end{remark}

\section{Decay estimate for Fourier integral operators}
\begin{proposition}[Control by U-norm of the multiplier]\label{Prop_gowersnormcontrol}
For $\cT_{m}(F,G)$ defined as in \eqref{carlyT}, we have the estimates:
\begin{equation}\label{bddgower2}
     \nrm{\cT_m\br{F,G}}_{L^1}
     \leq \nrm{m}_u\nrm{F}_{L^2}\nrm{G\br{x,\eta}}_{L^1_\eta L^2_x},
\end{equation}
\begin{equation}\label{bddgower3}
     \nrm{\cT_m\br{F,G}}_{L^2}\leq 
     \nrm{m}_U
     \nrm{F\br{\xi,y}}_{L^2_\xi L^4_y}
    \br{
    \nrm{G\br{x,\eta}}_{L^2_\eta L^4_x}
    \nrm{G}_{L^4}
    }^{\frac{1}{2}}.
\end{equation}
\end{proposition}

\begin{proposition}
[Mixed Norm Interpolation]\label{Prop_interpolation}
Let $c_0,c_1,c_2\geq 0$. If a bilinear operator
\begin{equation*}
    T: 
    \br{L^1\cap L^\infty}\br{\R\times\mr{0,1}}
    \times
    \br{L^1\cap L^\infty}\br{\mr{0,1}\times \R}
    \to 
    L^\infty\br{\mr{0,1}^2}
\end{equation*}
satisfies the following support size \(L^0\) control:
\begin{equation}\label{eq_interpolation_L0_ctrl}
    \supp F\subseteq \R\times P,\quad
    \supp G\subseteq Q\times \R
    \implies
    \abs{\supp T\br{F,G}}\leq c_0
    \abs{P}\cdot\abs{Q}
\end{equation}
and the following two mixed-norm estimates:
\begin{equation}\label{eq_interpolation_L2_L12_ctrl}
    \nrm{
        T\br{F,G}
    }_{L^1}\leq c_1
    \nrm{F}_{L^2}\nrm{G\br{x,\eta}}_{L^1_\eta L^2_x},
\end{equation}
\begin{equation}\label{eq_interpolation_L24_L24_L4_ctrl}
    \nrm{
        T\br{F,G}
    }_{L^2}
    \leq c_2
    \nrm{F\br{\xi,y}}_{L^2_\xi L^4_y}
    \br{
    \nrm{G\br{x,\eta}}_{L^2_\eta L^4_x}
    \nrm{G}_{L^4}
    }^{\frac{1}{2}},
\end{equation}
then the following estimate holds:
\begin{equation}\label{eq_interpolation_L0_mix_nrm_result}
    \nrm{
        T\br{F,G}
    }_{L^1}
    \lesssim
    c^{\frac{2}{5}}_0
    c^{\frac{1}{5}}_1
    c^{\frac{4}{5}}_2
    \nrm{F}_{L^2}
    \nrm{G}_{L^2}.
\end{equation}
\end{proposition}

\begin{proposition}[Factorization of $m_{\vgamma,\rho}$]\label{prop_m_fac}
Assume \(\br{\vgamma,\rho}\) satisfies the \(\br{A,C}\)-short-interval assumption  \ref{ass_short_I}.
For all \(\vxi\in \cC_\Gamma\) (i.e. \(\mr{\vxi}\in\img\cG_\Gamma\)), the oscillatory expression \(m_{\vgamma,\rho}\br{\vxi}\) \eqref{eq_m_gamma_rho} admits the following factorization
\begin{equation}\label{eq_m_fac}
    m_{\vgamma,\rho}\br{\vxi}=
    e\br{\vxi\cdot\vgamma_\Gamma\br{\vxi}}
    \cdot
    \Psi\br{\vxi}
\end{equation}
into a purely oscillatory term \(e\br{\vxi\cdot\vgamma_\Gamma\br{\vxi}}\) and an amplitude term \(\Psi\in C^1\br{\cC_\Gamma}\) satisfying the estimates:
\begin{equation}\label{eq_m_fac_amp}
    \abs{\Psi\br{\vxi}}\leq A^{O\br{1}} \abs{\vxi}^{-\frac{1}{2}}
    ,\quad
    \abs{\nabla\Psi\br{\vxi}}
    \leq A^{O\br{1}}
    \abs{\vxi}^{-\frac{3}{2}}.
\end{equation}
\end{proposition}
\begin{definition}[Chirality]
    We define formally the following quantity:
    \begin{equation}\label{eq_def_chirality}
        \Theta\br{\vgamma}:=\Theta\br{\gamma_1,\gamma_2}:=\frac{3}{2}\cdot\frac{\gamma''_1\gamma'_2+\gamma''_2\gamma'_1}{\vgamma'\wedge\vgamma''}
        -\frac{\gamma'_1\gamma'_2}{\vgamma'\wedge\vgamma''}\cdot\frac{\vgamma'\wedge\vgamma'''}{\vgamma'\wedge\vgamma''}.
    \end{equation}
\end{definition}

The quantity detects a certain chirality of the curve since \(\Theta\br{\gamma_2,\gamma_1}=-\Theta\br{\gamma_1,\gamma_2}\). Moreover, when \(\vgamma:I=\br{-A^{-C},A^{-C}}\simto \Gamma\) is \(A\)-admissible for some \(A,C\gg 1\), we deduce from the estimate \(\abs{\Theta\br{\gamma_1,\gamma_2}'\br{t}}\lesssim A^{O\br{1}}\) and the mean value theorem a control on the diameter of the image\footnote{See \textbf{Remark \ref{rmk_small_img}}.}
\begin{equation*}
    \dia \Theta\br{\gamma_1,\gamma_2}\br{I}
    \lesssim
    A^{O\br{1}}\dia I= A^{O\br{1}-C}\ll 1.
\end{equation*}
As a result, when \(A,C\gg 1\), at least one of the following must hold:
\begin{equation}\label{eq_vgamma_chiral_ass}
    \Theta\br{\gamma_1,\gamma_2}\br{t}\leq \frac{1}{4},\quad
    \forall t\in I;
\end{equation}
\begin{equation*}
    \Theta\br{\gamma_2,\gamma_1}\br{t}\leq \frac{1}{4},\quad
    \forall t\in I.
\end{equation*}
Henceforth, after potentially interchanging the order of the pairing \(\br{\gamma_1,\gamma_2}\), we will always assume \eqref{eq_vgamma_chiral_ass}.

\begin{proposition}[Gowers norm calculation for \(m_{\vgamma,\rho,\vlambda}\)]\label{prop_m_uni_norm}
    Let \(\br{\vgamma,\rho}\) satisfy the \(\br{A,C}\)-short-interval condition.\footnote{Recall \textbf{Definition \ref{ass_short_I}}.} Given the two conditions \eqref{eq_vlambda_cond} and \eqref{eq_vgamma_chiral_ass},
    % if further the following inequality is satisfied:
    % \begin{equation}\label{eq_vgamma_chiral_ass}
    %     \br{\frac{3}{2}\cdot\frac{\gamma''_1\gamma'_2+\gamma''_2\gamma'_1}{\vgamma'\wedge\vgamma''}
    %     -\frac{\gamma'_1\gamma'_2}{\vgamma'\wedge\vgamma''}\cdot\frac{\vgamma'\wedge\vgamma'''}{\vgamma'\wedge\vgamma''}}\br{t}\leq \frac{1}{4},\quad \forall t\in I,
    % \end{equation}
    the symbol \(m_{\vgamma,\rho,\vlambda}\) given by \eqref{eq_m_vlambda_freq_trunc} satisfies the following estimates:
    \begin{equation}\label{eq_m_vlambda_small_uni}
        \nrm{m_{\vgamma,\rho,\vlambda}}_u
        \leq A^{O\br{1}} \Xi^{\frac{1}{2}}_\Gamma \abs{\vlambda}^{-\frac{1}{4}}
        \log^{\frac{1}{4}}\br{\abs{\vlambda}};
    \end{equation}
    \begin{equation}\label{eq_m_vlambda_big_uni}
        \nrm{m_{\vgamma,\rho,\vlambda}}_U
        \leq A^{O\br{1}} \Xi^{\frac{1}{4}}_\Gamma
        \log^{\frac{1}{8}}\br{\abs{\vlambda}}
    \end{equation}
\end{proposition}
% \begin{remark}
%     \eqref{eq_vgamma_chiral_ass} is a harmless chirality assumption. 
%     Indeed, by \textbf{Remark \ref{rmk_small_img}}, either the curve itself satisfies \eqref{eq_vgamma_chiral_ass}, or an interchange of indices \(\br{\gamma_1,\gamma_2}\mapsto\br{\gamma_2,\gamma_1}\) guarantees \eqref{eq_vgamma_chiral_ass}. As a result, for the purposes of proving \textbf{Theorem \ref{thm_4ier_op_est}}, we may always assume \eqref{eq_vgamma_chiral_ass}.
% \end{remark}

With \textbf{Proposition \ref{Prop_gowersnormcontrol}, \ref{Prop_interpolation}, \ref{prop_m_fac}, \ref{prop_m_uni_norm}} above, we will prove \textbf{Theorem \ref{thm_universal_multiplier_bound}, \ref{thm_4ier_op_est}} in the next subsection.

\subsection{Proof of Theorem \ref{thm_universal_multiplier_bound} and Theorem \ref{thm_4ier_op_est}}

We begin with the proof of \textbf{Theorem \ref{thm_universal_multiplier_bound}}. 
Observe that the fiber structure and the localization in \eqref{carlyT} suggest that it suffices to consider the case with \(\supp F\subseteq \R\times\mr{0,1}\) and \(\supp G\subseteq \mr{0,1}\times \R\). Moreover, since \(\cT_m\) satisfies the relation:
\begin{equation*}
    \supp F\subseteq \R\times P,\quad
    \supp G\subseteq Q\times \R
    \implies
    \supp \cT_m\br{F,G}\subseteq Q\times P,
\end{equation*}
we deduce that \(\cT_m\) satisfies \eqref{eq_interpolation_L0_ctrl} with \(c_0=1\). Applying now \textbf{Proposition \ref{Prop_gowersnormcontrol} and \ref{Prop_interpolation}}, we deduce:
\begin{equation*}
    \nrm{\cT_m\br{F,G}}_{L^1}
    \lesssim
    1^\frac{2}{5}
    \nrm{m}^{\frac{1}{5}}_u
    \nrm{m}^{\frac{4}{5}}_U
    \nrm{F}_{L^2}
    \nrm{G}_{L^2}.
\end{equation*}
This concludes the proof of \textbf{Theorem \ref{thm_universal_multiplier_bound}}.

To prove \textbf{Theorem \ref{thm_4ier_op_est}}, we recall that it is harmless to assume \eqref{eq_vgamma_chiral_ass}. By \textbf{Proposition \ref{prop_m_uni_norm}}, we deduce the two estimates \eqref{eq_m_vlambda_small_uni} and \eqref{eq_m_vlambda_big_uni}. Apply now \textbf{Theorem \ref{thm_universal_multiplier_bound}}, we deduce:
\begin{equation*}
    \nrm{\cT_{m_{\vgamma,\rho,\vlambda}}\br{F,G}}_{L^1}\lesssim 
    \br{
        A^{O\br{1}}\Xi^{\frac{1}{2}}_\Gamma
        \abs{\vlambda}^{-\frac{1}{4}}\log^{\frac{1}{4}}\br{\abs{\vlambda}}
    }^{\frac{1}{5}}
    \br{
        A^{O\br{1}}\Xi^{\frac{1}{4}}_\Gamma
        \log^{\frac{1}{8}}\br{\abs{\vlambda}}
    }^{\frac{4}{5}}
    \nrm{F}_{L^2}\nrm{G}_{L^2}.
\end{equation*}
By direct calculation, we obtain the desired estimate. This concludes the proof of \textbf{Theorem \ref{thm_4ier_op_est}}.

\subsection{Proof of Proposition \ref{Prop_gowersnormcontrol}: Control via U-norm}

We begin with proving \eqref{bddgower2}. Applying Fubini's theorem, we first integrate over $\xi$, followed by $\eta$, and apply a $L^{1}$ bound in $\eta$ to obtain
\begin{equation}
     |\T_{m}(F,G)(x,y)|\leq \left\|\int_{\mathbb{R}}F(\xi ,y)m(\xi ,\eta)e(x \xi)d\xi  \cdot G(x,\eta)    \right\|_{L^{1}_{\eta}}\, .
\end{equation}
Note that taking the $L^{1}(\R^2)$ norm of $\T_{m}(F,G)$ is the same as the $L^1([0,1]^2)$ norm. Applying Fubini theorem once more, we reorganize the iterated integral as
\begin{equation}
     \|\T_{m}(F,G)\|_{L^{1}}\leq \left\| \left\|  \left\|\int_{\mathbb{R}}F(\xi ,y)m(\xi ,\eta)e(x \xi)d\xi  \cdot G(x,\eta)    \right\|_{L^{1}_{y}([0,1])}  \right\|_{L^{1}_{\eta}}  \right\|_{L^{1}_{x}([0,1])} \, .
\end{equation}   
Factoring out $G(x,\eta)$ from the innermost integral and using the nesting properties of the $L^p$ norms on the
space $[0,1]$, then applying H\"older's inequality, we may further estimate
\begin{equation}\label{bddgower2_2}
    \left\|  \left\| \left\| \int_{\mathbb{R}}F(\xi ,y)m(\xi ,\eta)e(x \xi)d\xi \right\|_{L^{2}_{y}([0,1])}  \cdot G(x,\eta)    \right\|_{L^{1}_{\eta}}  \right\|_{L^{1}_{x}([0,1])}    
\end{equation}
\begin{equation}\label{bddgower2_2_1}
      \leq \left\|   \left\| \left\| \int_{\mathbb{R}}F(\xi ,y)m(\xi ,\eta )e(x \xi)d\xi \right\|_{L^{2}_{y}([0,1])} \right\|_{L^{\infty}_{\eta}}  \cdot \left\| G(x,\eta)    \right\|_{L^{1}_{\eta}}  \right\|_{L^{1}_{x}([0,1])}.
\end{equation}

 We introduce a measurable function  $\eta :[0,1] \rightarrow \mathbb{R}$ to linearize the norm $\|\cdot \|_{L^{\infty}_{\eta}}$ and equate \eqref{bddgower2_2_1} with
\begin{equation}\label{bddgower2_3}
     \left\|   \left\| \int_{\mathbb{R}}F(\xi ,y)m(\xi ,\eta (x))e(x \xi)d\xi \right\|_{L^{2}_{y}([0,1])}  \cdot \left\| G(x,\eta)    \right\|_{L^{1}_{\eta}}  \right\|_{L^{1}_{x}([0,1])}    \, .
\end{equation}

Applying H\"older inequality once more, we bound \eqref{bddgower2_3} by
\begin{equation}
    \left\| \int_{\mathbb{R}}F(\xi ,y)m(\xi ,\eta (x))e(x \xi)d\xi \right\|_{L^{2}_{xy}([0,1]^{2})}\cdot \left\| G(x,\eta)  \right\|_{L^{1}_{\eta}L^{2}_{x}}\, .
\end{equation}

Next, we build up several properties related to multiplicative derivative. Beside the definition \eqref{multideri_def}, we also define the partial multiplicative derivative
\begin{equation}\label{multideri_def_partial}
\mathcal{D}^{(x_{k})}_{y}\left( f(x_{1},\cdots,x_{k},\cdots ,x_{d}) \right):=f(x_{1},\cdots, x_{k}+y,\cdots, x_{d})\overline{f(x_{1},\cdots ,x_{k},\cdots ,x_{d})}\, .
    \end{equation}
\begin{lemma}
For $u,v\in \mathbb{R}$, $1\leq p\leq \infty$, the following properties of multiplicative derivative hold:

\begin{equation}\label{square_to_multiplicativederivative}
        \left|\int_{\mathbb{R}}f(x,y)dx \right|^{2}=\int_{\mathbb{R}^{2}}\mathcal{D}_{u}^{(x)}(f(x,y))dudx \, ,
    \end{equation}
 \begin{equation}\label{lptomulti}
        \left\|f(x,y) \right\|_{L^{p}_{x}}^{2}= \left\| \mathcal{D}_{u}^{(x)}(f(x,y))\right\|_{L^{p}_{xu}}\, ,
    \end{equation}
\begin{equation}\label{twodef_multi_eqiv}
    \mathcal{D}_{u}^{(x)}(f(x,y))=(\mathcal{D}_{(u,0)}f)(x,y)\, ,
\end{equation}
\begin{equation}\label{l2_to_multiplicativederivative}
    \left\| \int_{\mathbb{R}}f(x,y)dx \right\|_{L^{2}_{y}(A)}^{2}=\int_{A}\int_{\mathbb{R}^{2}}\mathcal{D}_{u}^{(x)}(f(x,y))dudxdy\, ,
\end{equation}
\begin{equation}\label{mult_distribute_law}
    \mathcal{D}_{u}^{(x)}\left( \prod_{k=1}^{n}f_{k}(x,y) \right)=\prod_{k=1}^{n} \mathcal{D}_{u}^{(x)}(f_{k}(x,y))\, ,
\end{equation}
\begin{equation}\label{mult_commutative_law}
\mathcal{D}^{(x)}_{u}\mathcal{D}^{(y)}_{v}f(x,y)=\mathcal{D}^{(y)}_{v}\mathcal{D}^{(x)}_{u}f(x,y)\, ,
\end{equation}
\begin{equation}\label{multi_cs}
    \left\| \mathcal{D}_{u}^{(x)}(f(x,y))\right\|_{L^{p}_{y}}\leq \mathcal{D}_{u}^{(x)}\left\|f(x,y) \right\|_{L^{2p}_{y}}\, .
\end{equation}
    
\end{lemma}

\begin{proof}
    For identity \eqref{square_to_multiplicativederivative}, expanding the square and making the change of variables $x'=x+u$, we obtain
    \begin{equation}
         \left|\int_{\mathbb{R}}f(x,y)dx \right|^{2}=\int_{\mathbb{R}}f(x',y)dx' \cdot \overline{\int_{\mathbb{R}}f(x,y)dx }
    \end{equation}
    \begin{equation}
=\int_{\mathbb{R}^{2}}f(x+u,y)\overline{f(x,y)}dudx=\int_{\mathbb{R}^{2}}\mathcal{D}_{u}^{(x)}(f(x,y))dudx 
    \end{equation}
For identity \eqref{lptomulti}, expanding the right-hand side of \eqref{lptomulti} and with the fact that integration is invariant under translation , we obtain
    
    \begin{equation}
        \left\| \mathcal{D}_{u}^{(x)}(f(x,y))\right\|_{L^{p}_{xu}}=\left\|f(x+u,y)\overline{f(x,y)}\right\|_{L^{p}_{xu}}
    \end{equation}
     \begin{equation}
=\left\|\left\|f(x+u,y)\overline{f(x,y)}\right\|_{L^{p}_{u}}\right\|_{L^p(x)}
    \end{equation}
    \begin{equation}
       = \left\| \left\|f(x+u,y) \right\|_{L^{p}_{u}}\cdot f(x,y) \right\|_{L^{p}_{x}}=\left\| \left\|f(u,y) \right\|_{L^{p}_{u}}\cdot f(x,y) \right\|_{L^{p}_{x}}
    \end{equation}
    \begin{equation}
        =\left\|f(u,y) \right\|_{L^{p}_{u}}\cdot \left\|f(x,y) \right\|_{L^{p}_{x}}=\left\|f(x,y) \right\|_{L^{p}_{x}}^{2}\, .
    \end{equation}
     For identity \eqref{twodef_multi_eqiv}, by the definitions of multiplicative derivative \eqref{multideri_def} and partial multiplicative derivative \eqref{multideri_def_partial}, we obtain
    \begin{equation}
        \mathcal{D}_{u}^{(x)}(f(x,y))=f(x+u,y)\overline{f(x,y)}=f(x+u,y+0)\overline{f(x,y)}=(\mathcal{D}_{(u,0)}f)(x,y)\, .
    \end{equation}
    For identity \eqref{l2_to_multiplicativederivative}, expanding the $L^{2}$ norm and by \eqref{square_to_multiplicativederivative}, we obtain
\begin{equation}
     \left\| \int_{\mathbb{R}}f(x,y)dx \right\|_{L^{2}_{y}(A)}^{2}=\int_{A}\left|\int_{\mathbb{R}}f(x,y)dx \right|^{2}dy=\int_{A}\int_{\mathbb{R}^{2}}\mathcal{D}_{u}^{(x)}(f(x,y))dudxdy\, .
\end{equation}
For identity \eqref{mult_distribute_law}, expanding the definition of multiplicative derivative, we obtain
\begin{equation}
     \mathcal{D}_{u}^{(x)}\left( \prod_{k=1}^{n}f_{k}(x,y) \right)=\left( \prod_{k=1}^{n}f_{k}(x+u,y) \right)\overline{\left( \prod_{k=1}^{n}f_{k}(x,y) \right)}
\end{equation}
\begin{equation}
    =\prod_{k=1}^{n}f_{k}(x+u,y)\overline{f_{k}(x,y)}=\prod_{k=1}^{n} \mathcal{D}_{u}^{(x)}(f_{k}(x,y))\, .
\end{equation}
For identity \eqref{mult_commutative_law}, expanding the definition of multiplicative derivative, we obtain
\begin{equation}
\mathcal{D}^{(x)}_{u}\mathcal{D}^{(y)}_{v}f(x,y)=\mathcal{D}^{(x)}_{u}\left(f(x,y+v)\overline{f(x,y)}\right)
\end{equation}
\begin{equation}
=\left(f(x+u,y+v)\overline{f(x+u,y)}\right)\overline{\left(f(x,y+v)\overline{f(x,y)}\right)}
\end{equation}
\begin{equation}
=\left(f(x+u,y+v)\overline{f(x,y+v)}\right)\overline{\left(f(x+u,y)\overline{f(x,y)}\right)}
\end{equation}
\begin{equation}
=\mathcal{D}^{(y)}_{v}\left(f(x+u,y)\overline{f(x,y)}\right)=\mathcal{D}^{(y)}_{v}\mathcal{D}^{(x)}_{u}f(x,y)\, .
\end{equation}
For inequality \eqref{multi_cs}, expanding the definition of multiplicative derivative and by H\"older's inequality, we obtain
\begin{equation}
       \left\| \mathcal{D}_{u}^{(x)}(f(x,y))\right\|_{L^{p}_{y}}=\|f(x+u,y)\overline{f(x,y)}\|_{L^{p}_{y}}
    \leq \|f(x+u,y)\|_{L^{2p}}\cdot \|f(x,y)\|_{L^{2p}}=\mathcal{D}_{u}^{(x)}\left\|f(x,y) \right\|_{L^{2p}_{y}}\, .
\end{equation}
\end{proof}

To show \eqref{bddgower2}, it suffices to show the following estimate
\begin{equation}\label{bddgower2_4}
    \left\| \int_{\mathbb{R}}F(\xi ,y)m(\xi ,\eta (x))e(x \xi)d\xi \right\|_{L^{2}_{xy}([0,1]^{2})}\leq \|m\|_{u}\cdot \|F\|_{L^{2}}\, .
\end{equation}
By \eqref{l2_to_multiplicativederivative}, and \eqref{mult_distribute_law}, the square of the left hand side of \eqref{bddgower2_4} is equal to
\begin{equation}
\int_{0}^{1}\int_{0}^{1}\int_{\mathbb{R}^{2}}\mathcal{D}^{(\xi)}_{u} \left( F(\xi ,y)m(\xi ,\eta (x))e(x \xi)  \right)d\xi dudxdy
\end{equation}
\begin{equation}
=\int_{0}^{1}\int_{0}^{1}\int_{\mathbb{R}^{2}}\mathcal{D}_{u}^{(\xi)}(F(\xi ,y))\cdot \mathcal{D}_{u}^{(\xi)}(m(\xi ,\eta (x))e(x\xi))d\xi du dx dy
\end{equation}
By Fubini theorem, we rewrite this as

\begin{equation}\label{bddgower2_5}
    =\int_{\mathbb{R}^{2}} \left(\int_{0}^{1}\mathcal{D}_{u}^{(\xi)}(F(\xi ,y))dy \right) \cdot \left( \int_{0}^{1}\mathcal{D}_{u}^{(\xi)}(m(\xi ,\eta (x))e(x\xi))dx \right) d\xi du \, .
\end{equation}
By $L^{2}_{\xi u}\times L^{2}_{\xi u} \rightarrow L^{1}_{\xi u}$ H\"older inequality, \eqref{bddgower2_5} is bounded by
\begin{equation}\label{bddgower2_6}
    \left\| \int_{0}^{1}\mathcal{D}_{u}^{(\xi)}(F(\xi ,y))dy  \right\|_{L^{2}_{\xi u}}\cdot \left\| \int_{0}^{1}\mathcal{D}_{u}^{(\xi)}(m(\xi ,\eta (x))e(x\xi))dx  \right\|_{L^{2}_{\xi u}}\, .
\end{equation}
By \eqref{multi_cs} and \eqref{l2_to_multiplicativederivative}, the first term in \eqref{bddgower2_6} can be estimated by 
\begin{equation}\label{bddgower2_5_1}
    \leq \left\| \mathcal{D}_{u}^{(\xi)}\left(\|F(\xi ,y)\|_{L^{2}_{y}}\right)\right\|_{L^{2}_{\xi u}}= \left\| \|F(\xi ,y)\|_{L^{2}_{y}}\right\|_{L^{2}_{\xi }}^{2}=\|F\|_{L^{2}}^{2}\, .
\end{equation}
By \eqref{l2_to_multiplicativederivative} and \eqref{mult_commutative_law}, the square of the second term in \eqref{bddgower2_6} is equal to
\begin{equation}
    \int_{\mathbb{R}^{4}} \mathcal{D}^{(x)}_{v}\left(\1_{[0,1]}(x) \mathcal{D}^{(\xi)}_{u}\left( m(\xi ,\eta (x))e(x \xi) \right)\right)dvdxd\xi du
\end{equation}
\begin{equation}\label{bddgower2_6_2}
     = \int_{\mathbb{R}^{4}} \mathcal{D}^{(x)}_{v}\mathcal{D}^{(\xi)}_{u}\left(\1_{[0,1]}(x)  m(\xi ,\eta (x))e(x \xi) \right)dvdxd\xi du
\end{equation}
\begin{equation}
 = \int_{\mathbb{R}^{4}} \mathcal{D}^{(\xi)}_{u}\mathcal{D}^{(x)}_{v}\left(\1_{[0,1]}(x)  m(\xi ,\eta (x))e(x \xi) \right)dvdxd\xi du
\end{equation}
\begin{equation}\label{bddgower2_6_1}
    = \left\|\int_{\mathbb{R}} \mathcal{D}^{(x)}_{v}\left(\1_{[0,1]}(x)  m(\xi ,\eta (x))e(x \xi) \right)d\xi  \right\|_{L^{2}_{vx}}^{2}\, .
\end{equation}
Expanding the multiplicative derivative, \eqref{bddgower2_6_1} equals to
\begin{equation}
    \left\|\int_{\mathbb{R}}\1_{[0,1]}(x+v)\1_{[0,1]}(x)\cdot m(\xi,\eta(x+v))\overline{m(\xi,\eta(x))} e((x+v)\xi-x\xi)d\xi \right\|_{L^{2}_{xv}}^{2}
\end{equation}
\begin{equation}
    =\left\| \int_{\mathbb{R}}(\mathcal{D}_{(0,\eta (x+v)-\eta (x))} m)(\xi ,\eta (x))e(v \xi)  d\xi \right\|^{2}_{L^{2}_{v}([-x,1-x])L^{2}_{x}([0,1])}
\end{equation}
\begin{equation}\label{bddgower2_7}
     =\left\| (\4{1}\mathcal{D}_{(0,\eta (x+v)-\eta (x))} m)(v ,\eta (x))\right\|^{2}_{L^{2}_{v}([-x,1-x])L^{2}_{x}([0,1])}
\end{equation}
Hence \eqref{bddgower2_7} is bounded by
\begin{equation}
     \left\| \left\|(\4{1}\mathcal{D}_{(0,s)} m)(v ,\eta ) \right\|_{L^{\infty}_{s\eta}}\right\|^{2}_{L^{2}_{v}([-x,1-x])L^{2}_{x}([0,1])}
    \leq \left\| \left\|(\4{1}\mathcal{D}_{(0,s)} m)(v ,\eta ) \right\|_{L^{\infty}_{s\eta}}\right\|^{2}_{L^{2}_{v}([-1,1])}=\|m\|_{u}^{4}\, ,
\end{equation}
which completes the proof of \eqref{bddgower2}.

Last, we prove the inequality \eqref{bddgower3}. By \eqref{l2_to_multiplicativederivative}, \eqref{mult_distribute_law}, we obtain
\begin{equation}
    \|\T_{m} (F,G)\|_{L^{2}}^{2}=\int_{\mathbb{R}^{6}}(\mathcal{D}_{(u,0)}F)(\xi ,y)(\mathcal{D}_{(0,v)}G)(x,\eta)(\mathcal{D}_{(u,v)}m)(\xi ,\eta)e(xu+yv)dxdy d\xi d\eta du dv
\end{equation}
\begin{equation}\label{bddgower3_1}
=\int_{\mathbb{R}^{4}}\left(\br{\4{2}}^{-1}\mathcal{D}_{(u,0)}F  \right)(\xi ,v)\left(\br{\4{1}}^{-1}\mathcal{D}_{(0,v)}G  \right)(u ,\eta)(\mathcal{D}_{(u,v)}m)(\xi ,\eta)d\xi d\eta du dv
\end{equation}
By $L^{2}_{\xi uv}\times L^{2}_{\xi uv} \rightarrow L^{1}_{\xi uv}$ H\"older inequality, \eqref{bddgower3_1} is bounded by
\begin{equation}\label{bddgower3_2}
    \left\| \left(\br{\4{2}}^{-1}\mathcal{D}_{(u,0)}F  \right)(\xi ,v)  \right\|_{L^{2}_{\xi uv}}\cdot \left\| \int_{\mathbb{R}} \left(\br{\4{1}}^{-1}\mathcal{D}_{(0,v)}G  \right)(u ,\eta)(\mathcal{D}_{(u,v)}m)(\xi ,\eta) d\eta \right\|_{L^{2}_{\xi uv}}\, .
\end{equation}
By Plancherel identity and \eqref{lptomulti}, the first component in \eqref{bddgower3_2} equals 
\begin{equation}
    \left\| \left(\mathcal{D}_{(u,0)}F  \right)(\xi ,y)  \right\|_{L^{2}_{\xi uy}}=\left\| \left\|\mathcal{D}^{(\xi)}_{u}(F(\xi,y)) \right\|_{L^{2}_{\xi u}} \right\|_{L^{2}_{y}}=\left\| \left\|  F(\xi ,y) \right\|_{L^{2}_{\xi}}^{2}   \right\|_{L^{2}_{y}}=\left\| F(\xi ,y)  \right\|_{L^{2}_{\xi}L^{4}_{y}}^{2}\, .
\end{equation}
Next, by \eqref{l2_to_multiplicativederivative}, \eqref{mult_distribute_law}, the square of the second component in \eqref{bddgower3_2} equals
\begin{equation}\label{bddgower3_3}
\int_{\mathbb{R}^{5}}\left(\mathcal{D}_{(0,s)}\br{\4{1}}^{-1}\mathcal{D}_{(0,v)}G\right) (u,\eta)  \cdot \left(\mathcal{D}_{(0,s)}\mathcal{D}_{(u,v)}m \right) (\xi ,\eta) d\eta ds d\xi du dv \, .
\end{equation}
Then by $L^{1}_{u\eta}\times L^{\infty}_{u \eta} \rightarrow L^{1}_{u \eta}$ H\"older inequality,  \eqref{bddgower3_3} is bounded by
\begin{equation}\label{bddgower3_4}
    \left\| \left\|\left(\mathcal{D}_{(0,s)}\br{\4{1}}^{-1}\mathcal{D}_{(0,v)}G\right) (u,\eta)   \right\|_{L^{1}_{u\eta}} \cdot \left\|\int_{\mathbb{R}}\left(\mathcal{D}_{(0,s)}\mathcal{D}_{(u,v)}m \right) (\xi ,\eta) d\xi  \right\|_{L^{\infty}_{u\eta}}   \right\|_{L^{1}_{vs}}
\end{equation}
Note that by Cauchy-Schwarz inequality and Plancherel identity, the first component in \eqref{bddgower3_4} is bounded by
\begin{equation}
\left\|\left(\mathcal{D}_{(0,s)}\br{\4{1}}^{-1}\mathcal{D}_{(0,v)}G\right) (u,\eta)   \right\|_{L^{1}_{u\eta}}=\left\| \left(\br{\4{1}}^{-1}\mathcal{D}_{(0,v)}G  \right)(u+s ,\eta) \cdot \overline{\left(\br{\4{1}}^{-1}\mathcal{D}_{(0,v)}G  \right)(u ,\eta)}
  \right\|_{L^{1}_{u\eta}}
\end{equation}
\begin{equation}
    \leq \left\| \left(\br{\4{1}}^{-1}\mathcal{D}_{(0,v)}G  \right)(u+s ,\eta) \right\|_{L^{2}_{u\eta}} \cdot \left\| \left(\br{\4{1}}^{-1}\mathcal{D}_{(0,v)}G  \right)(u ,\eta) \right\|_{L^{2}_{u\eta}}
\end{equation}
\begin{equation}
    =\left\| \left(\br{\4{1}}^{-1}\mathcal{D}_{(0,v)}G  \right)(u ,\eta) \right\|_{L^{2}_{u\eta}}^{2}
    =\left\| \left(\mathcal{D}_{(0,v)}G  \right)(x ,\eta) \right\|_{L^{2}_{x\eta}}^{2}\, .
\end{equation}
Hence, \eqref{bddgower3_4} is bounded by 
\begin{equation}\label{bddgower3_5}
    \left\| \left\| \left(\mathcal{D}_{(0,v)}G  \right)(x ,\eta) \right\|_{L^{2}_{x\eta}}^{2}\cdot   \left\|\int_{\mathbb{R}}\left(\mathcal{D}_{(0,s)}\mathcal{D}_{(u,v)}m \right) (\xi ,\eta) d\xi  \right\|_{L^{\infty}_{u\eta}L^{1}_{s}} \right\|_{L^{1}_{v}}\, .
\end{equation}
Then by $L^{2}_{v}\times L^{2}_{v}\rightarrow L^{1}_{v}$ H\"older inequality,  \eqref{bddgower3_5} is bounded by
\begin{equation}\label{bddgower3_6}
     \left\| \left\| \left(\mathcal{D}_{(0,v)}G  \right)(x ,\eta) \right\|_{L^{2}_{x\eta}}^{2}\right\|_{L^{2}_{v}}\cdot   \left\|\left\|\int_{\mathbb{R}}\left(\mathcal{D}_{(0,s)}\mathcal{D}_{(u,v)}m \right) (\xi ,\eta) d\xi  \right\|_{L^{\infty}_{u\eta}L^{1}_{s}} \right\|_{L^{2}_{v}}\, .
\end{equation}
By the definition \eqref{mcapitalunorm}, the second term in \eqref{bddgower3_6} is exactly the quantity $\|m\|_{U}^{4}$. By Minkovski's inequality and $L^{4}_{x}\times L^{4}_{x}\rightarrow L^{2}_{x}$ H\"older's inequality, the first term in \eqref{bddgower3_6} is bounded by
\begin{equation}
     \left\| \left\| \left(\mathcal{D}_{(0,v)}G  \right)(x ,\eta) \right\|_{L^{2}_{x\eta}}^{2}\right\|_{L^{2}_{v}}=\left\| \left\|G(x,\eta +v)\overline{G(x,\eta)}  \right\|_{L^{2}_{x\eta}}  \right\|_{L^{4}_{v}}^{2}
      \leq \left\| \left\|G(x,\eta +v)\right\|_{L^{4}_{v}} \cdot \overline{G(x,\eta)}  \right\|_{L^{2}_{x\eta}}^{2}
\end{equation}
\begin{equation}
   =\left\| \left\|G(x,v)\right\|_{L^{4}_{v}} \cdot \|G(x,\eta)\|_{L^{2}_{\eta}}  \right\|_{L^{2}_{x}}^{2}   
   \leq \|G(x,v)\|_{L^{4}_{xv}}^{2}\cdot \|G(x,\eta)\|_{L^{2}_{\eta}L^{4}_{x}}^{2}\, ,
\end{equation}
which completes the proof of \eqref{bddgower3}.

\subsection{Proof of Proposition \ref{Prop_interpolation}: Interpolation}

% Due to the localization and the fiber structure of the operator \(\cT_m\), we may also assume: 
% \begin{equation*}
%     \supp F\subseteq \R\times\mr{0,1},\quad
%     \supp G\subseteq \mr{0,1}\times\R.
% \end{equation*}
For \(s\in\br{0,\infty}\), we define the following level sets
\begin{equation}\label{eq_Q_s}
    Q^{s}_{<}:=
    \BR{
        \br{x,\eta}\in \mr{0,1}\times\R
    \::\: 
        s\nrm{G\br{x,\eta'}}_{L^2_{\eta'}}
        <\abs{G\br{x,\eta}}  
    },\quad
    Q^{s}_{\geq}:=
    \br{\mr{0,1}\times\R}\setminus
    Q^{s}_{<}.
\end{equation}
% and decompose the function \(G\) accordingly:
% \begin{equation*}
%     G=G^{s<}+G^{s\geq},\quad G^{s<}:=\1_{Q^s}G.
% \end{equation*}
By Chebyshev inequality, we obtain the following estimate
\begin{equation}\label{eq_cheby_for_freq}
    \nrm{\1_{Q^s_{<}}\br{x,\eta}}_{L^2_\eta}
    \leq
    \Big\Vert
        \frac{G\br{x,\eta}}
        {
            s
            \nrm{G\br{x,\eta'}}_{L^2_{\eta'}}
        }
    \Big\Vert_{L^2_\eta}
    \leq s^{-1}.
\end{equation}
As a direct consequence, we obtain the following estimate
\begin{equation*}
    \nrm{\1_{Q^{s}_{<}}G\br{x,\eta}}_{L^1_\eta L^2_x}
    \leq 
    \nrm{
        \nrm{\1_{Q^{s}_{<}}\br{x,\eta}}_{L^2_\eta}
        \nrm{G\br{x,\eta}}_{L^2_\eta}
    }_{L^2_x}
    \leq s^{-1}\nrm{G}_{L^2}.
\end{equation*}
Applying \eqref{eq_interpolation_L2_L12_ctrl}, we deduce
\begin{equation}\label{eq_inter_mini_u_w_spikes}
    \nrm{
    T\br{F,\1_{Q^{s}_{<}}G}
    }_{L^1}
    \leq
    c_1
    \nrm{F}_{L^2}
    \nrm{\1_{Q^{s}_{<}}G\br{x,\eta}}_{L^1_\eta L^2_x}
    \leq c_1 s^{-1}\nrm{F}_{L^2}\nrm{G}_{L^2}.
\end{equation}
It remains to control \(T\br{F,\1_{Q^{s\geq}}G}\) via \eqref{eq_interpolation_L24_L24_L4_ctrl}. To estimate the more complicated mixed-norm expression involved in \eqref{eq_interpolation_L24_L24_L4_ctrl}, we introduce a further dyadic decomposition.
For \(j,k\in\N\), we define
\begin{equation}\label{eq_P_0}
    P_0:=
    \BR{
        y\in\mr{0,1}
    \::\:
        \nrm{F\br{\xi,y}}_{L^2_\xi}
        \leq
        \nrm{F}_{L^2}
    },
\end{equation}
\begin{equation}\label{eq_P_j}
    P_j:=
    \BR{
        y\in\mr{0,1}
    \::\:
        2^{j-1}\nrm{F}_{L^2}
        <\nrm{F\br{\xi,y}}_{L^2_\xi}
        \leq
        2^j\nrm{F}_{L^2}
    },
\end{equation}
\begin{equation}\label{eq_Q_0}
    Q_0:=
    \BR{
        x\in\mr{0,1}
    \::\:
        \nrm{G\br{x,\eta}}_{L^2_\eta}
        \leq
        \nrm{G}_{L^2}
    },
\end{equation}
\begin{equation}\label{eq_Q_k}
    Q_k:=
    \BR{
        x\in\mr{0,1}
    \::\:
        2^{k-1}\nrm{G}_{L^2}
        <\nrm{G\br{x,\eta}}_{L^2_\eta}
        \leq
        2^k\nrm{G}_{L^2}
    }.
\end{equation}
Correspondingly, we define
\begin{equation}\label{eq_FG_decomp_for_uni}
    F=
    \sum^\infty_{j=0}F_j
    ,\quad
    F_j:=\1_{\R\times P_j}F
    ,\quad
    G=
    \sum^\infty_{k=0}G_k
    ,\quad
    G_k:= \1_{Q_k\times\R}G.
\end{equation}
By Chebyshev's inequality, for $j,k\in \mathbb{N}$, we have the estimate
\begin{equation}
    |P_{j}|\leq 
    \left(2^{j-1}\|F\|_{L^{2}}  \right)^{-2}\cdot \int_{P_{j}}\left(\|F(\xi,y)\|_{L^{2}_{\xi}}  \right)^{2}dy\leq 2^{2-2j}\, ,
\end{equation}
\begin{equation}
    |Q_{k}|\leq 
    \left(2^{k-1}\|G\|_{L^{2}}  \right)^{-2}\cdot \int_{Q_{k}}\left(\|G(x,\eta)\|_{L^{2}_{\eta}}  \right)^{2}dx\leq 2^{2-2k}\, .
\end{equation}
Also note that $P_{0},Q_{0}\subseteq  [0,1]$. Hence for $ j,k\in\N\sqcup\BR{0}$,
\begin{equation}\label{eq_cheby_for_phy}
    \abs{P_j}\lesssim 2^{-2j},\quad\abs{Q_k}\lesssim 2^{-2k}\, .
\end{equation}
Below we estimate three relevant mixed-norm expressions in \eqref{eq_interpolation_L24_L24_L4_ctrl}.
By the upper bound estimate in the definition of \eqref{eq_P_0} and \eqref{eq_P_j}, we have
\begin{equation}\label{eq_L2L4_F}
     \nrm{
        F_j\br{\xi,y}
     }_{L^2_\xi L^4_y}
     \leq
     \nrm{
        F_j\br{\xi,y}
     }^{\frac{1}{2}}_{L^2_\xi L^2_y}
     \nrm{
        F_j\br{\xi,y}
     }^{\frac{1}{2}}_{L^2_\xi L^\infty_y}
     \leq
     \nrm{F}^{\frac{1}{2}}_{L^2}
     \br{2^j\nrm{F}_{L^2}}^{\frac{1}{2}}
     =
     2^{j/2}\nrm{F}_{L^2}.
\end{equation}
Similarly, by the upper bound estimate provided by the defining formula of \eqref{eq_Q_0} and \eqref{eq_Q_k}, we have:
\begin{equation}\label{eq_L2L4_G}
     \nrm{
        G_k\br{x,\eta}
     }_{L^2_\eta L^4_x}
     \leq
     \nrm{
        G_k\br{x,\eta}
     }^{\frac{1}{2}}_{L^2_\eta L^2_x}
     \nrm{
        G_k\br{x,\eta}
     }^{\frac{1}{2}}_{L^2_\eta L^\infty_x}
     \leq
     \nrm{G}^{\frac{1}{2}}_{L^2}
     \br{2^k\nrm{G}_{L^2}}^{\frac{1}{2}}
     =
     2^{k/2}\nrm{G}_{L^2}.
\end{equation}
For the \(L^4\) expression, we use the upper bound estimate in \eqref{eq_Q_s} and deduce
\begin{equation*}
    \nrm{\1_{Q^{s}_{\geq}}G_k}_{L^4}
    \leq
    \nrm{
        \nrm{
            \1_{Q^{s}_{\geq}}G_k\br{x,\eta}
        }^{\frac{1}{2}}_{L^2_\eta}
        \nrm{
            \1_{Q^{s}_{\geq}}G_k \br{x,\eta}
        }^{\frac{1}{2}}_{L^\infty_\eta}
    }_{L^4_x\br{Q_k}}
\end{equation*}
\begin{equation}\label{eq_L4_G}
    \leq
    \nrm{
        \nrm{G\br{x,\eta}}^{\frac{1}{2}}_{L^2_\eta}
        \br{
            s
            \nrm{
                G\br{x,\eta}
            }_{L^2_\eta}
        }^{\frac{1}{2}}
    }_{L^4_x\br{Q_k}}
    =s^{1/2}
    \nrm{
        G_k\br{x,\eta}
     }_{L^2_\eta L^4_x}.
\end{equation}
Observe that \eqref{eq_interpolation_L0_ctrl}
implies the following estimate
\begin{equation*}
    \nrm{T\br{F_j,\1_{Q^{s}_{\geq}}G_k}}_{L^1}
    \leq
    \abs{\supp T\br{F_j,\1_{Q^{s}_{\geq}}G_k}}^{\frac{1}{2}}
    \nrm{T\br{F_j,\1_{Q^{s}_{\geq}}G_k}}_{L^2}
    \leq
    c^{\frac{1}{2}}_0
    \abs{P_j}^{\frac{1}{2}}
    \abs{Q_k}^{\frac{1}{2}}
    \nrm{T\br{F_j,\1_{Q^{s}_{\geq}}G_k}}_{L^2}.
\end{equation*}
Applying \eqref{eq_cheby_for_phy} and \eqref{eq_interpolation_L24_L24_L4_ctrl}, the above term is dominated by
\begin{equation*}
    \lesssim 
    c^{\frac{1}{2}}_{0} 2^{-j} 2^{-k}
    c_2
    \nrm{
        F_j\br{\xi,y}
     }_{L^2_\xi L^4_y}
     \br{
         \nrm{
            \1_{Q^{s}_{\geq}}
            G_k\br{x,\eta}
         }_{L^2_\eta L^4_x}
         \nrm{\1_{Q^{s}_{\geq}}G_k}_{L^4}
     }^{\frac{1}{2}}.
\end{equation*}
By the trivial estimate \(\abs{\1_{Q^{s}_{\geq}}
        G_k}\leq
    \abs{G_k}\) and \eqref{eq_L4_G}, we further dominate the above by
\begin{equation*}
    \lesssim 2^{-j-k}
    c^{\frac{1}{2}}_0
    c_2
    \nrm{
        F_j\br{\xi,y}
     }_{L^2_\xi L^4_y}
     \br{
         \nrm{
            G_k\br{x,\eta}
         }_{L^2_\eta L^4_x}
         \cdot
         s^{1/2}
         \nrm{
            G_k\br{x,\eta}
         }_{L^2_\eta L^4_x}
     }^{\frac{1}{2}}.
\end{equation*}
Finally, combining the two estimates \eqref{eq_L2L4_F} and \eqref{eq_L2L4_G}, we conclude
\begin{equation*}
    \nrm{T\br{F_j,\1_{Q^{s}_{\geq}}G_k}}_{L^1}
    \lesssim 2^{-j/2-k/2}
    c^{\frac{1}{2}}_0
    c_2
    s^{1/4}
    \nrm{
        F
     }_{L^2}
     \nrm{
        G
     }_{L^2}.
\end{equation*}
Summing over \(j,k\in\N\sqcup\BR{0}\), we have
\begin{equation}\label{eq_inter_mini_u_wo_spikes}
    \nrm{T\br{F,\1_{Q^{s}_{\geq}}G}}_{L^1}
    \lesssim
    c^{\frac{1}{2}}_0
    c_2
    s^{1/4}
    \nrm{
        F
     }_{L^2}
     \nrm{
        G
     }_{L^2}.
\end{equation}
Combining \eqref{eq_inter_mini_u_w_spikes} and \eqref{eq_inter_mini_u_wo_spikes}, we have
\begin{equation}
    \|T(F,G)\|_{L^{1}}\leq  \nrm{T\br{F,\1_{Q^{s}_{<}}G}}_{L^1}+ \nrm{T\br{F,\1_{Q^{s}_{\geq}}G}}_{L^1}\lesssim (c_{1}s^{-1}+c_{0}^{\frac{1}{2}}c_{2}s^{\frac{1}{4}})\|F\|_{L^{2}}\|G\|_{L^{2}}\, .
\end{equation}
Taking \(s=\br{c^{-1/2}_0 c^{-1}_2c_1}^{4/5}\), we obtain the desired conclusion: \eqref{eq_interpolation_L0_mix_nrm_result}
\begin{equation*}
    \nrm{T\br{F,G}}_{L^1}
    \lesssim 
    \br{
        c_1\br{c^{-1/2}_0 c^{-1}_2c_1}^{-4/5}+
        c^{\frac{1}{2}}_0 c_2
        \br{c^{-1/2}_0 c^{-1}_2c_1}^{1/5}
    }
    \nrm{F}_{L^2}\nrm{G}_{L^2}
    \sim
    c^{\frac{2}{5}}_0
    c^{\frac{1}{5}}_1
    c^{\frac{4}{5}}_2
    \nrm{F}_{L^2}\nrm{G}_{L^2}.\qed
\end{equation*}

\subsection{Proof of Proposition \ref{prop_m_fac}: Decomposition of the multiplier}
We now perform the decomposition \eqref{eq_m_fac}. We begin by analyzing the phase function.
The critical point of the phase is given by the inverse of the function \(\theta:=-\gamma'_1/\gamma'_2\),
\begin{equation}
    t_\ast\br{\vxi}:=
    \theta^{-1}\br{\xi_2/\xi_1},
\end{equation}
which is well-defined for \(\vxi\in\cC_\Gamma\) due to \eqref{eq_Gauss_2_ratio} and \eqref{eq_theta_n_prime_bds}. Note that $t_{\ast}$ solves the equation
\begin{equation}\label{eq_crit_pt}
    \vxi\cdot \vgamma'\br{t_\ast\br{\vxi}}=0.
\end{equation}
With this in hand, we extract the principal oscillation and obtain the decomposition
\begin{equation}
    m_{\vgamma,\rho}\br{\vxi}=
    e\br{\vxi\cdot \vgamma\br{t_\ast\br{\vxi}}}
    \int
        e\br{
            \vxi\cdot 
            \dif{\vgamma}^t_{t_\ast\br{\vxi}}
        }
        \rho\br{t}
    dt\, ,
\end{equation}
where $\dif{\vgamma}^t_{t_\ast\br{\vxi}}:=\vgamma (t)-\vgamma (t_\ast\br{\vxi})$. Due to the relation given by the geometric interpretation of solving \eqref{eq_crit_pt},
\begin{equation}\label{eq_vgamma_to_inv_Gauss}
    \vgamma\br{t_\ast\br{\vxi}}=\vgamma_\Gamma\br{\vxi}.
\end{equation}
It remains to show the amplitude term
\begin{equation}\label{eq_amp_as_osc_int}
    \Psi\br{\vxi}:=
    \int
        e\br{\vxi\cdot \dif{\vgamma}^t_{t_\ast\br{\vxi}}
        }
        \rho\br{t}
    dt
\end{equation}
satisfies the desired estimate whenever \(\vxi\in \cC_\Gamma\). Applying Van der Corput’s lemma \textbf{\ref{prop_oscillatory}}, we obtain the estimate
\begin{equation*}
    \abs{\Psi\br{\vxi}}
    \lesssim
    \br{\inf_{t\in I} \abs{\vxi\cdot \vgamma''\br{t}}}^{-\frac{1}{2}}
    \br{\nrm{\rho}_{L^\infty}+\nrm{\rho'}_{L^1}}
    \lesssim
    \br{\inf_{t\in I} \abs{\vxi\cdot \vgamma''\br{t}}}^{-\frac{1}{2}}.
\end{equation*}
It remains to estimate the term within the infimum.

\begin{lemma}[Non-degeneracy of curvature]\label{lem_xi_gamma_d2_stuff}
    For \(\vxi\in\cC_\Gamma\) we have the following identities:
    \begin{equation}\label{eq_xi_gamma_d2_id}
        \vxi\cdot \vgamma''\br{t_\ast\br{\vxi}}=
        \br{-1}^i\xi_i
        \br{\frac{\vgamma'\wedge\vgamma''}{\gamma'_j}}\br{t_\ast\br{\vxi}},\quad\BR{i,j}=\BR{1,2};
    \end{equation}
    \begin{equation}\label{eq_xi_gamma_d2_id_abs}
        \abs{\vxi\cdot \vgamma''\br{t_\ast\br{\vxi}}}=
        \abs{\vxi}
        \br{\frac{\abs{\vgamma'\wedge\vgamma''}}{\abs{\vgamma'}}}\br{t_\ast\br{\vxi}}.
    \end{equation}
    Moreover, under \textbf{Assumption \ref{ass_short_I}}, the following estimate holds:
    \begin{equation}\label{eq_xi_gamma_d2_non_degen}
        \abs{\vxi\cdot \vgamma''\br{t}}\geq A^{-O\br{1}}\abs{\vxi},\quad \forall t\in I.
    \end{equation}
\end{lemma}

Given \textbf{Lemma \ref{lem_xi_gamma_d2_stuff}}, the first estimate in \eqref{eq_m_fac_amp} follows. We now prove the lemma.

\begin{proof}[Proof of \textbf{Lemma \ref{lem_xi_gamma_d2_stuff}}]
    Recall \eqref{eq_crit_pt}. Simplifying the following expression
    \begin{equation*}
        \vxi\cdot\vgamma''\br{t_\ast\br{\vxi}}
        =
        \vxi\cdot\vgamma''\br{t_\ast\br{\vxi}}
        -\vxi\cdot \vgamma'\br{t_\ast\br{\vxi}}\cdot
        \br{\frac{\gamma''_j}{\gamma'_j}}\br{t_\ast\br{\vxi}}
    \end{equation*}
    gives exactly \eqref{eq_xi_gamma_d2_id}. As for \eqref{eq_xi_gamma_d2_id_abs}, we observe that \eqref{eq_crit_pt} implies also the relation
    \begin{equation*}
        \vxi=\pm \abs{\vxi} \br{\frac{\br{\gamma'_2,-\gamma'_1}}{\abs{\vgamma'}}}
        \br{t_\ast\br{\vxi}},\quad
    \end{equation*}
    for the correct sign.
    Taking now the inner product with \(\vgamma''\br{t_\ast\br{\vxi}}\), we deduce
    \begin{equation*}
        \vxi\cdot\vgamma''\br{t_\ast\br{\vxi}}
        =
        \pm \abs{\vxi} \br{\frac{\br{\gamma'_2,-\gamma'_1}}{\abs{\vgamma'}}}
        \br{t_\ast\br{\vxi}}\cdot\vgamma''\br{t_\ast\br{\vxi}}
        =
        \mp \abs{\vxi}
        \br{
            \frac{\vgamma'\wedge\vgamma''}{\abs{\vgamma'}}
        }\br{t_\ast\br{\vxi}}.
    \end{equation*}
    The identity \eqref{eq_xi_gamma_d2_id_abs} thus follows.
    Finally, to show \eqref{eq_xi_gamma_d2_non_degen}, we start with the triangle inequality
    \begin{equation*}
        \abs{\vxi\cdot\vgamma''\br{t}}\geq
        \abs{\vxi\cdot \vgamma''\br{t_\ast\br{\vxi}}}
        -\abs{\vxi}\cdot\abs{\vgamma''\br{t}-\vgamma''\br{t_\ast\br{\vxi}}}.
    \end{equation*}
    By \eqref{eq_xi_gamma_d2_id_abs} and \textbf{Remark \ref{rmk_small_img}}, we bound the above from below with
    \begin{equation*}
        \geq 
        \abs{\vxi}
        \br{\frac{\abs{\vgamma'\wedge\vgamma''}}{\abs{\vgamma'}}}\br{t_\ast\br{\vxi}}
        - 
        \abs{\vxi}\cdot \nrm{\vgamma'''}_{L^\infty} \abs{I}
        \geq
        A^{-2}
        \abs{\vxi}
        - A^{1-C}\abs{\vxi}\geq A^{-O\br{1}}\abs{\vxi}
    \end{equation*}
    and conclude the proof of \eqref{eq_xi_gamma_d2_non_degen}.
    % Applying now the identity \eqref{eq_xi_gamma_d2_id}, we bound the above from below with
    % \begin{equation*}
    %     \geq A^{-2}\abs{\xi_i}
    %     - A^{O\br{1}-C}\abs{\vxi}.
    % \end{equation*}
    % Using now the two fact \eqref{eq_Gauss_2_ratio} and \eqref{eq_theta_n_prime_bds}, we deduce \(A^{-2}\leq\abs{\xi_2/\xi_1}\leq A^2\)
    % and conclude \eqref{eq_xi_gamma_d2_non_degen}:
    % \begin{equation*}
    %     \abs{\vxi\cdot\vgamma''\br{t}}
    %     \geq A^{-O\br{1}}\abs{\vxi}
    %     - A^{O\br{1}-C}\abs{\vxi}
    %     \geq A^{-O\br{1}}\abs{\vxi}.
    % \end{equation*}
\end{proof}

Next, to prove the second estimate in \eqref{eq_m_fac_amp}, we pass the derivative inside the integral
\begin{equation}
    \nabla \Psi\br{\vxi}=
    \int
        e\br{
            \vxi\cdot 
            \dif{\vgamma}^t_{t_\ast\br{\vxi}}
        }
        \cdot
        2\pi i 
        \mr{
            \nabla_\vxi
            \br{
                \vxi\cdot 
                \dif{\vgamma}^t_{t_\ast\br{\vxi}}
            }
        }
        \rho\br{t}
    dt.
\end{equation}
Note that the characterizing equation of the critical points \eqref{eq_crit_pt} implies that
\begin{equation}
    \nabla_\vxi
    \br{
        \vxi\cdot 
        \dif{\vgamma}^t_{t_\ast\br{\vxi}}
    }
    =
    \dif{\vgamma}^t_{t_\ast\br{\vxi}}
    +\br{\vxi\cdot
    \vgamma'\br{t_\ast\br{\vxi}}}
    \nabla t_\ast\br{\vxi}
    =
    \dif{\vgamma}^t_{t_\ast\br{\vxi}}.
\end{equation}
We may thus simplify the formula for \(\nabla \Psi\) and write
\begin{equation*}
    \nabla \Psi\br{\vxi}
    =
    \int
        e\br{
            \vxi\cdot 
            \dif{\vgamma}^t_{t_\ast\br{\vxi}}
        }
        \cdot
        2\pi i 
        \dif{\vgamma}^t_{t_\ast\br{\vxi}}
        \rho\br{t}
    dt.
\end{equation*}
Apply integration by-parts and the fact \eqref{eq_crit_pt}. We equate the
above with
\begin{equation}\label{eq_nabla_amp}
    -
    \int
        e\br{
            \vxi\cdot 
            \dif{\vgamma}^t_{t_\ast\br{\vxi}}
        }
        \cdot
        \partial_t
        \br{
            \frac{
                \dif{\vgamma}^t_{t_\ast\br{\vxi}}
            }{
                \vxi \cdot \dif{\vgamma'}^t_{t_\ast\br{\vxi}}
            }
            \cdot \rho\br{t}
        }
    dt.
\end{equation}
This suggests that we consider the function
\begin{equation}\label{eq_nabla_amp_inner_amp}
    \vA_\vxi\br{t}
    :=
    \frac{
        \dif{\vgamma}^t_{t_\ast\br{\vxi}}
    }{
        \vxi \cdot \dif{\vgamma'}^t_{t_\ast\br{\vxi}}
    }
    =
    \frac{
        \int^1_0
            \vgamma'\br{
                t_\ast\br{\vxi}+
                \theta\br{t-t_\ast\br{\vxi}}
            }
        d\theta
    }{
        \int^1_0
            \vxi\cdot
            \vgamma''\br{
                t_\ast\br{\vxi}+
                \theta\br{t-t_\ast\br{\vxi}}
            }
        d\theta
    }\, .
\end{equation}
Applying Van der Corput lemma \textbf{\ref{prop_oscillatory}} to \eqref{eq_nabla_amp}, together with \eqref{eq_xi_gamma_d2_non_degen}, we obtain the estimate
\begin{equation}
    \abs{
        \nabla\Psi\br{\vxi}
    }
    \lesssim
    \br{\inf_{t\in I} \abs{\vxi\cdot \vgamma''\br{t}}}^{-\frac{1}{2}}
    \br{
        \nrm{
            \br{\vA_\vxi \rho}'
        }_{L^\infty}
        +
        \nrm{
            \br{\vA_\vxi \rho}''
        }_{L^1}
    }
    \lesssim
        A^{O\br{1}}
        \abs{\vxi}^{-\frac{1}{2}}
        \nrm{\vA_\vxi}_{C^2\br{I}}.
\end{equation}
To complete the proof of the second estimate in \eqref{eq_m_fac_amp}, it remains to show the following estimate
\begin{equation}\label{eq_nabla_amp_inner_amp_est}
    \nrm{\vA_\vxi}_{C^2\br{I}}\leq A^{O\br{1}}\abs{\vxi}^{-1},\quad \forall \vxi\in \cC_\Gamma.
\end{equation}
Recall \eqref{eq_nabla_amp_inner_amp}. By direct computation, we obtain
\begin{equation*}
    \vA'_\vxi\br{t}
    =
    \frac{
        \int^1_0
            \vgamma''\br{
                t_\ast\br{\vxi}+
                \theta\br{t-t_\ast\br{\vxi}}
            }
        \theta
        d\theta
    }{
        \int^1_0
            \vxi\cdot
            \vgamma''\br{
                t_\ast\br{\vxi}+
                \theta\br{t-t_\ast\br{\vxi}}
            }
        d\theta
    }
    -\vA_\vxi\br{t}\cdot
    \frac{
        \int^1_0
            \vxi\cdot
            \vgamma'''\br{
                t_\ast\br{\vxi}+
                \theta\br{t-t_\ast\br{\vxi}}
            }
        \theta
        d\theta
    }{
        \int^1_0
            \vxi\cdot
            \vgamma''\br{
                t_\ast\br{\vxi}+
                \theta\br{t-t_\ast\br{\vxi}}
            }
        d\theta
    };
\end{equation*}
\begin{equation*}
    \vA''_\vxi\br{t}
    =
     \frac{
        \int^1_0
            \vgamma'''\br{
                t_\ast\br{\vxi}+
                \theta\br{t-t_\ast\br{\vxi}}
            }
        \theta^2
        d\theta
    }{
        \int^1_0
            \vxi\cdot
            \vgamma''\br{
                t_\ast\br{\vxi}+
                \theta\br{t-t_\ast\br{\vxi}}
            }
        d\theta
    }
    -\vA_\vxi\br{t}\cdot
    \frac{
        \int^1_0
            \vxi\cdot
            \vgamma''''\br{
                t_\ast\br{\vxi}+
                \theta\br{t-t_\ast\br{\vxi}}
            }
        \theta^2
        d\theta
    }{
        \int^1_0
            \vxi\cdot
            \vgamma''\br{
                t_\ast\br{\vxi}+
                \theta\br{t-t_\ast\br{\vxi}}
            }
        d\theta
    }
\end{equation*}
\begin{equation*}
    -2\vA'_\vxi\br{t}\cdot
    \frac{
        \int^1_0
            \vxi\cdot
            \vgamma'''\br{
                t_\ast\br{\vxi}+
                \theta\br{t-t_\ast\br{\vxi}}
            }
        \theta
        d\theta
    }{
        \int^1_0
            \vxi\cdot
            \vgamma''\br{
                t_\ast\br{\vxi}+
                \theta\br{t-t_\ast\br{\vxi}}
            }
        d\theta
    }.
\end{equation*}
Applying now the \(A\)-admissibility assumption \eqref{eq_gamma_admi_ana} and \eqref{eq_xi_gamma_d2_non_degen}, we deduce
\begin{equation*}
    \abs{\vA_\vxi\br{t}}\leq A^{O\br{1}}\abs{\vxi}^{-1},\quad
    \abs{\vA'_\vxi\br{t}}\leq A^{O\br{1}}\abs{\vxi}^{-1}+ A^{O\br{1}}\abs{\vA_\vxi\br{t}},
\end{equation*}
\begin{equation*}
    \abs{\vA''_\vxi\br{t}}\leq A^{O\br{1}}\abs{\vxi}^{-1}+ A^{O\br{1}}\br{\abs{\vA_\vxi\br{t}}+\abs{\vA'_\vxi\br{t}}}.
\end{equation*}
Combining the above estimates gives \eqref{eq_nabla_amp_inner_amp_est}. This completes the proof of \eqref{eq_m_fac_amp}.\qed

\subsection{Proof of Proposition \ref{prop_m_uni_norm}: Calculation of the \(u\)- and \(U\)-norms}
We will first prove the inequality \eqref{eq_m_vlambda_small_uni}. Our first step is to obtain the following point-wise estimate
\begin{equation}\label{eq_m_mini_u_pw_est}
    \abs{
        \4{1}\cD_{\br{0,s}}
        m_{\vgamma,\rho,\vlambda}\br{x,\eta}
    }
    \leq A^{O\br{1}}\Xi_\Gamma
    \br{1+\abs{\vlambda x}}^{-\frac{1}{2}}.
\end{equation}
Recall again \eqref{eq_m_vlambda_freq_trunc}. To control the above oscillatory integral expression, we shall first extract the main oscillatory mode of \(m_{\vgamma,\rho,\vlambda}\). This can be achieved by \textbf{Proposition \ref{prop_m_fac}}, once we know for \(\vlambda\) satisfying \eqref{eq_vlambda_cond}, the following implication holds:
\begin{equation*}
    \psi_\delta\br{\xi/\lambda_1}
    \psi_\delta\br{\eta/\lambda_2}\neq 0\implies
    \mr{\xi:\eta}\in\img \cG_\Gamma.
\end{equation*}
Below, we provide the following lemma which not only ensures the above property but also describes the image of Gauss map and the dichotomy between the stationary and the non-stationary phase contributions:
\begin{lemma}\label{lem_osc_st_or_nst}
    Let \(\vlambda:=\br{\lambda_1,\lambda_2}\in\br{\R\setminus\BR{0}}^2\). We have the following implication
    \begin{equation}\label{eq_tight_diag}
        \mr{\vlambda}\in \img\cG_\Gamma \implies A^{-2}\leq \abs{\lambda_2/\lambda_1}\leq A^2,
    \end{equation}
    and for \(\vxi\) satisfying 
    \(\abs{\xi_j-\lambda_j}\leq \delta\abs{\lambda_j}\) the following dichotomy:
    \begin{equation}\label{eq_osc_st}
        \mr{\vlambda}
        \in\cG_\Gamma\circ\vgamma\br{A^{-C}I}
        % \implies
        % \mr{\vxi}\in \cG\circ\vgamma\br{I}
        \implies
        \exists t\in I,\;
            \vxi\cdot\vgamma'\br{t}=0
        \br{i.e.\;
        \mr{\vxi}\in\img\cG_\Gamma
        };
    \end{equation}
    \begin{equation}\label{eq_osc_nst}
        \mr{\vlambda}
        \notin\cG_\Gamma\circ\vgamma\br{A^{-C}I}
        \implies
        \forall t\in \supp \rho,\;
        % \subseteq A^{-2C}I,\;
        \abs{
            \vxi\cdot\vgamma'\br{t}
        }\geq A^{-O\br{C}}\abs{\vlambda}.
    \end{equation}
\end{lemma}
\begin{proof}
    Via rescaling, it suffices to assume \(\vlambda:=\br{1,\lambda}\) for some \(\lambda \in\R\setminus\BR{0}\). As an immediate consequence of the given assumption on \(\vxi\). We have:
    \begin{equation}\label{eq_lambda_close_2_ratio}
        \abs{\frac{\xi_2}{\xi_1}-\lambda}
        \leq\abs{\frac{\xi_2-\lambda}{\xi_1}+\frac{\lambda-\lambda\xi_1}{\xi_1}}
        \leq \frac{2\delta}{1-\delta} \cdot \abs{\lambda}
        \leq A^{O\br{1}-4C}\abs{\lambda}.
    \end{equation}
    Next, we describe the image of \(\cG_\Gamma\circ\vgamma\). Let \(U\subseteq I\) and set \(\theta:=-\gamma'_1/\gamma'_2\). Due to the following identity
    \begin{equation*}
        \cG_\Gamma\circ\vgamma\br{U}=
        \BR{
            \mr{\gamma'_2\br{t}:-\gamma'_1\br{t}}
            =
            \mr{
                1:
                \theta\br{t}
            }
        \::\:
            t\in U
        },
    \end{equation*}
    we can identify the image of \(\cG_\Gamma\circ\vgamma\) with the image of \(\theta\) via the relation
    \begin{equation}\label{eq_Gauss_2_ratio}
        \mr{1:\lambda}\in \cG_\Gamma\circ\vgamma\br{U}
        \iff
        \lambda\in \theta\br{U}.
    \end{equation}
    Notice that, by direct computation, the \(A\)-admissibility assumption on \(\vgamma:I\simto\Gamma\) gives
    \begin{equation}\label{eq_theta_n_prime_bds}
        \theta'=
        \frac{\vgamma'\wedge\vgamma''}{\br{\gamma'_2}^2},\quad
        A^{-2}\leq
        \abs{
            \theta
        }
        \leq A^2
        ,\quad
        A^{-3}\leq
        \abs{\theta'}\leq A^4.
    \end{equation}
    This has various implications. On the one hand, combining \eqref{eq_Gauss_2_ratio} and \eqref{eq_theta_n_prime_bds} yields immediately \eqref{eq_tight_diag}
    \begin{equation*}
        \mr{\vlambda}:=\mr{1:\lambda}\in\img\cG_\Gamma=\cG_\Gamma\br{I}\implies
        \lambda\in \theta\br{I}\implies A^{-2}\leq \abs{\lambda}=\abs{\lambda_2/\lambda_1}\leq A^2.
    \end{equation*}
    On the other hand, we deduce that \(\theta\) must be monotone and satisfy a refinement of \textbf{Remark \ref{rmk_small_img}}:
    \begin{equation*}
        B\br{\theta\br{t},A^{-3} r}\subseteq \theta\br{B\br{t,r}}
        \subseteq B\br{\theta\br{t},A^4 r},\quad B\br{t,r}\subseteq I.
    \end{equation*}
    As an immediate consequence, we have the following nested relations:
    \begin{equation}\label{eq_inner_theta}
        \theta\br{\supp\rho}\subseteq
        \theta\br{A^{-2C}I}\subseteq
        B\br{\theta\br{0},A^{4-3C}}
    \end{equation}
    \begin{equation}\label{eq_middle_theta}
            \subseteq
        B\br{\theta\br{0},A^{-3-2C}}
            \subseteq 
        \theta\br{A^{-C}I}
            \subseteq 
        B\br{\theta\br{0},A^{4-2C}}
    \end{equation}
    \begin{equation}\label{eq_outer_theta}
        \subseteq 
        B\br{\theta\br{0},A^{-3-C}}
            \subseteq 
        \theta\br{I}
    \end{equation}
    whenever \(C\gg 1\). We are now ready to prove the two statements \eqref{eq_osc_st} and \eqref{eq_osc_nst}.

    Starting with \eqref{eq_osc_st}, let \(\mr{\vlambda}:=\mr{1:\lambda}\in\cG_\Gamma\circ\vgamma\br{A^{-C}I}\). By \eqref{eq_Gauss_2_ratio} and \eqref{eq_middle_theta}, we deduce that:
    \begin{equation*}
        \lambda\in\theta\br{A^{-C}I}\subseteq B\br{\theta\br{0},A^{4-2C}}.
    \end{equation*}
    By \eqref{eq_lambda_close_2_ratio} and \eqref{eq_theta_n_prime_bds}, triangle inequality gives:
    \begin{equation*}
        \abs{\frac{\xi_2}{\xi_1}-\theta\br{0}}
        \leq\abs{\frac{\xi_2}{\xi_1}-\lambda}
        +\abs{\lambda-\theta\br{0}}
        \leq A^{4-2C}+A^{O\br{1}-4C}=A^{O\br{1}-2C}.
    \end{equation*}
    Together with \eqref{eq_outer_theta}, this shows that:
    \begin{equation*}
        \frac{\xi_2}{\xi_1}
        \in B\br{\theta\br{0},A^{O\br{1}-2C}}
        \subseteq B\br{\theta\br{0},A^{-3-C}}
        \subseteq \theta\br{I}
    \end{equation*}
    for \(C\gg 1\). Finally, by \eqref{eq_Gauss_2_ratio} again, we conclude \eqref{eq_osc_st}:
    \begin{equation*}
        \mr{\vxi}=\mr{1:\frac{\xi_2}{\xi_1}}
        \in \cG_\Gamma\circ \vgamma\br{I}=\img\cG_\Gamma.
    \end{equation*}
    
    To show \eqref{eq_osc_nst}, we let \(\mr{\vlambda}:=\mr{1:\lambda}\notin\cG_\Gamma\circ\vgamma\br{A^{-C}I}\) and \(t\in \supp\rho\).
    % Since
    % \begin{equation*}
    %     \abs{\vxi}\eqsim \abs{\br{1,\lambda}}\eqsim \max\br{1,\abs{\lambda}},
    % \end{equation*}
    % it suffices to show
    We observe that:
    \begin{equation*}
        \abs{\vxi\cdot\vgamma'\br{t}}=
        \abs{\frac{\xi_2}{\xi_1}-\theta\br{t}}\cdot\abs{\xi_1}\cdot\abs{\gamma'_2\br{t}}\geq 
        \abs{\frac{\xi_2}{\xi_1}-\theta\br{t}}\cdot
        \frac{1-\delta}{A}
        \geq
        A^{-2}
        \abs{\frac{\xi_2}{\xi_1}-\theta\br{t}}.
    \end{equation*}
    % Pairing with the trivial fact: \(\max\br{1,\abs{\lambda}}\geq \abs{\vlambda}/\sqrt{2}\),
    It suffices to prove the following estimate:
    \begin{equation}\label{eq_osc_nst_red}
        \abs{\frac{\xi_2}{\xi_1}-\theta\br{t}}\geq A^{-O\br{C}}
        \abs{\vlambda}
        % \max\br{1,\abs{\lambda}}
    \end{equation}
    to conclude \eqref{eq_osc_nst}.
    To show \eqref{eq_osc_nst_red},
    we start with the triangle inequality:
    \begin{equation*}
        \abs{\frac{\xi_2}{\xi_1}-\theta\br{t}}\geq
        \abs{\lambda-\theta\br{0}}-\abs{\theta\br{t}-\theta\br{0}}-\abs{\frac{\xi_2}{\xi_1}-\lambda}.
    \end{equation*}
    By \eqref{eq_lambda_close_2_ratio} and \eqref{eq_inner_theta}, we bound the above from below with:
    \begin{equation}\label{eq_dist_lower_bd}
        \geq \abs{\lambda-\theta\br{0}}-A^{4-3C}-A^{O\br{1}-4C}\abs{\lambda}.
    \end{equation}
    We now investigate two cases:
    \begin{itemize}
        \item the \(\abs{\lambda}\leq A^C\) case: since \(\lambda\notin\theta\br{A^{-C}I}\) by \eqref{eq_Gauss_2_ratio}, the relation \eqref{eq_middle_theta} implies:
        \begin{equation*}
        \eqref{eq_dist_lower_bd}\geq 
        A^{-3-2C}-A^{4-3C}-A^{O\br{1}-3C}
        =A^{-O\br{1}-2C}.
        \end{equation*}
        In other words, we conclude \eqref{eq_osc_nst_red}:
        \begin{equation*}
            \abs{\frac{\xi_2}{\xi_1}-\theta\br{t}}\geq
            A^{-O\br{1}-2C}\geq
            A^{-4C}\abs{\vlambda}
            % \max\br{1,\abs{\lambda}}
        \end{equation*}
        for large enough \(C\gg 1\).
        
        \item the \(\abs{\lambda}\geq A^C\) case:
        we perform another triangle inequality:
        \begin{equation*}
            \eqref{eq_dist_lower_bd}\geq
            \abs{\lambda}-\abs{\theta\br{0}}
            -A^{4-3C}-A^{O\br{1}-4C}\abs{\lambda}
        \end{equation*}
        and utilize the estimate \eqref{eq_theta_n_prime_bds} to conclude \eqref{eq_osc_nst_red}:
        \begin{equation*}
            \abs{\frac{\xi_2}{\xi_1}-\theta\br{t}}
            \geq A^{-O\br{1}}\abs{\lambda}-A^{O\br{1}}\geq A^{-O\br{1}}\abs{\lambda}
            =
            A^{-O\br{1}}\abs{\vlambda}
            % \max\br{1,\abs{\lambda}}
        \end{equation*}
        for large enough \(C\gg 1\).
    \end{itemize}
    This completes the proof of \eqref{eq_osc_nst}.
\end{proof}

With \textbf{Lemma \ref{lem_osc_st_or_nst}} and \textbf{Proposition \ref{prop_m_fac}}, we may write 
    \begin{equation*}
        \4{1}\cD_{\br{0,s}}
        m_{\vgamma,\rho,\vlambda}\br{x,\eta}
        =
        \int
            e\br{
                \Delta_{\br{0,s}}
                \big(
                    \br{\xi,\eta}\cdot
                    \vgamma_\Gamma\br{\xi,\eta}
                \big)
                +x\xi
            }
            \cD^{\br{\eta}}_s
            \br{
                \psi_\delta\br{\xi/\lambda_1}
                \psi_\delta\br{\eta/\lambda_2}
                \Psi\br{\xi,\eta}
            }
        d\xi.
    \end{equation*}
    We now perform the change of variables: 
\begin{equation}\label{eq_mini_u_CV_rule}
 \xi=\abs{\vlambda}\txi,\quad
 \eta=\abs{\vlambda}\teta,\quad
 s=\abs{\vlambda}\ts,\quad
 \vlambda=\abs{\vlambda}\tvlambda=\abs{\vlambda}\br{\tlambda_1,\tlambda_2}.
\end{equation}
This produces the following identity:
\begin{equation}\label{eq_m_mini_u_post_cv}
 \4{1}\cD_{\br{0,s}}m_{\vgamma,\rho,\vlambda}\br{x,\eta}
 =
 \int
     e\br{
        \abs{\vlambda}
        \cdot
        \tPhi\br{x,\ts,\txi,\teta}
     }
    \tPsi\br{\vlambda,\ts,\txi,\teta}
d\txi,
\end{equation}
where 
\begin{equation}\label{eq_tPhi_mini_u}
\tPhi\br{x,\ts,\txi,\teta}
:=
\Delta_{\br{0,\ts}}
\br{
    \br{\txi,\teta}\cdot
    \vgamma_\Gamma\br{\txi,\teta}
}
+x\txi
;
\end{equation}
\begin{equation}\label{eq_tPsi_mini_u}
\tPsi\br{\vlambda,\ts,\txi,\teta}
:=
\abs{\vlambda}
\cD^{\br{\teta}}_\ts
\br{
    \psi_\delta
        \br{\txi/\tlambda_1}
    \psi_\delta
        \br{\teta/\tlambda_2}
    \Psi\br{
        \abs{\vlambda}\txi,
        \abs{\vlambda}\teta
    }
}.
\end{equation}
Our next step is to estimate $\|\widetilde{\Psi}\|_{L^{\infty}}$ and $\|\partial_{\widetilde{\xi}}\widetilde{\Psi}\|_{L^{1}}$. 
For the sake of discussion, we introduce a super set:
\begin{equation}\label{eq_conic_annuls}
    \cA_{\delta,\Gamma}:=
    \BR{
        \vxi\in \br{\R\setminus\BR{0}}^2
    \::\:
        \mr{\vxi}\in\img\cG_\Gamma,\quad
        1-10\delta\leq \abs{\vxi}\leq 1+10\delta
    }
\end{equation}
of the following set:
\begin{equation*}
    \supset
    \bigcup_{\substack{
        \mr{\tvlambda}\in
        \cG_\Gamma\circ\vgamma\br{A^{-C}I}\\
        :\abs{\tvlambda}=1
    }}
    \BR{
        \vxi \in \R^2
    \::\:
        \psi_\delta\br{\xi_1/\tlambda_1}
        \psi_\delta\br{\xi_2/\tlambda_2}
        \neq 0
    }.
\end{equation*}
Recall \eqref{eq_psi_delta}. By the estimates 
\eqref{eq_m_fac_amp} and the fact that \(\big\vert\tvlambda\big\vert=1\), we obtain
\begin{equation}\label{eq_tPsi_mini_u_bd}
    \big\Vert\tPsi\big\Vert_{L^\infty}
    \leq
    \abs{\vlambda}
    \cdot
    \nrm{\phi}_{L^\infty}
    \nrm{
        \Psi\br{
            \abs{\vlambda}\txi,
            \abs{\vlambda}\teta
        }
    }^2_{L^\infty
        \br{
            \cA_{\delta,\Gamma}
        }
    }
    \leq 
    \abs{\vlambda}
    \br{
        A^{O\br{1}}
        \abs{\vlambda}^{-\frac{1}{2}}
    }^2
    =A^{O\br{1}}
\end{equation}
and, by product rule, the following estimate
\begin{equation*}
    \nrm{\partial_{\txi}\tPsi}_{L^1}
    \leq 
    2\abs{\vlambda}
    \cdot
    \nrm{\phi'}_{L^1}
    \nrm{
        \Psi\br{
            \abs{\vlambda}\txi,
            \abs{\vlambda}\teta
        }
    }^2_{L^\infty
        \br{
            \cA_{\delta,\Gamma}
        }
    }
\end{equation*}
\begin{equation*}
    +
    2
    \abs{\vlambda}
    \cdot 
    \delta\tlambda_1
    \cdot
    \nrm{\phi}_{L^1}
    \nrm{
        \Psi\br{
            \abs{\vlambda}\txi,
            \abs{\vlambda}\teta
        }
    }_{L^\infty
        \br{
            \cA_{\delta,\Gamma}
        }
    }
    \abs{\vlambda}\cdot
    \nrm{
        \nabla
        \Psi\br{
            \abs{\vlambda}\txi,
            \abs{\vlambda}\teta
        }
    }_{L^\infty
        \br{
            \cA_{\delta,\Gamma}
        }
    }
\end{equation*}
\begin{equation}\label{eq_tPsi_mini_u_partial_bd}
    \leq
    \abs{\vlambda}
    A^{O\br{1}}\abs{\vlambda}^{-\frac{2}{2}}+
    \abs{\vlambda}^2
    A^{O\br{1}-O\br{C}}
    \abs{\vlambda}^{-\frac{1}{2}}
    A^{O\br{1}}
    \abs{\vlambda}^{-\frac{3}{2}}
    \leq A^{O\br{1}}.
\end{equation}
Hence, we finish the estimates for the amplitude function \(\tPsi\). Next, we establish estimates for the phase function \(\tPhi\). To proceed, we first introduce relevant derivative identities related to \(\vgamma_\Gamma\) and \(t_\ast\).
Recall \eqref{eq_crit_pt} and \eqref{eq_vgamma_to_inv_Gauss}. For the sake of clarity, we will suppress the \(\vxi\) input dependence when it is clear from the context.
\begin{lemma}[Derivative identities]\label{lem_nabla_ids}
Let \(\vxi\in\cC_\Gamma\) (i.e. \(\mr{\vxi}\in \img \cG_\Gamma\)).
We have the following identities:
\begin{equation}\label{eq_nabla_vxi_inv_Gauss_id}
    \nabla\br{\vxi\cdot \vgamma_\Gamma}=\vgamma_\Gamma
    =\vgamma\br{t_\ast}.
\end{equation}
\begin{equation}\label{eq_nabla_t_ast_id}
    \nabla t_\ast= \frac{1}{\xi_1}\cdot
    \br{
        \frac{
            \gamma'_1\gamma'_2
        }{
            \vgamma'\wedge\vgamma''
        },
        \frac{
            \br{\vgamma'_2}^2
        }{
            \vgamma'\wedge\vgamma''
        }
    }
    \br{t_\ast}
    =
    -\frac{1}{\xi_2}
    \br{
        \frac{
            \br{\vgamma'_1}^2
        }{
            \vgamma'\wedge\vgamma''
        }
        ,
        \frac{
            \gamma'_1\gamma'_2
        }{
            \vgamma'\wedge\vgamma''
        }
    }
    \br{t_\ast}
    .
\end{equation}
\end{lemma}
\begin{proof}[Proof of \textbf{Lemma \ref{lem_nabla_ids}}]
    Starting with \eqref{eq_nabla_vxi_inv_Gauss_id}, we utilize the fact \eqref{eq_vgamma_to_inv_Gauss}, apply chain rule, and invoke \eqref{eq_crit_pt}:
    \begin{equation*}
        \nabla\br{\vxi\cdot \vgamma_\Gamma}
        =\vgamma_\Gamma+
        \br{\vxi\cdot \vgamma'\br{t_\ast}}\cdot\nabla t_\ast
        =\vgamma_\Gamma.
    \end{equation*}
    This justifies the identity \eqref{eq_nabla_vxi_inv_Gauss_id}. To prove \eqref{eq_nabla_t_ast_id}, we begin with the derivatives of \eqref{eq_crit_pt}:
    \begin{equation*}
        \vgamma'\br{t_\ast}+ \br{\vxi\cdot \vgamma''\br{t_\ast}}\nabla t_\ast
        =0.
    \end{equation*}
    By the first identity \eqref{eq_xi_gamma_d2_id} in \textbf{Lemma \ref{lem_xi_gamma_d2_stuff}}, the above can be rewritten as:
    \begin{equation*}
        \vgamma'\br{t_\ast}+ 
        \br{-1}^i\xi_i\br{\frac{\vgamma'\wedge\vgamma''}{\gamma'_j}}\br{t_\ast}
        \nabla t_\ast
        =0,\quad
        \BR{i,j}=\BR{1,2}.
    \end{equation*}
    Solving for \(\nabla t_\ast\) gives exactly \eqref{eq_nabla_t_ast_id}. This concludes the proof of \textbf{Lemma \ref{lem_nabla_ids}}.
\end{proof}
Returning to the discussion on estimating \(\tPhi\), we provide the following lower bound estimates
\begin{lemma}\label{lem_tPhi_mini_u_est}
    Let \(\lambda \abs{x} \geq 4\) and \(\br{\xi,\eta},\br{\xi,\eta+s}\in \cA_{\delta,\Gamma}\). Assume also the correct chirality \eqref{eq_vgamma_chiral_ass}. 
    The phase function defined by \eqref{eq_tPhi_mini_u} satisfies the following estimate:
    \begin{equation}\label{eq_tPhi_mini_u_est}
        \max\br{
            \abs{
                \lambda
                \partial_\xi
                \tPhi
                \br{
                    x,s,\xi,\eta
                }
            },
            \abs{
                \lambda
                \partial^2_\xi
                \tPhi
                \br{
                    x,s,\xi,\eta
                }
            }^{\frac{1}{2}}
        }
        \geq A^{-O\br{1}}
        \abs{\lambda x}^{\frac{1}{2}}.
    \end{equation}
\end{lemma}
\begin{proof}[Proof of \textbf{Lemma \ref{lem_tPhi_mini_u_est}}]

To prove \eqref{eq_tPhi_mini_u_est}, 
it suffices to show that the following upper bound:
\begin{equation}\label{eq_tPhi_mini_u_d1_ubd}
    \abs{\partial_\xi\tPhi\br{x,s,\xi,\eta}}
    \leq \abs{x/\lambda}^{\frac{1}{2}}
\end{equation}
implies the following lower bound:
\begin{equation}\label{eq_tPhi_mini_u_d2_lbd}
    \abs{\partial^2_\xi\tPhi\br{x,s,\xi,\eta}}
    \geq A^{-O\br{1}}\abs{x}
\end{equation}
under the chirality assumption \eqref{eq_vgamma_chiral_ass}.
Assume \eqref{eq_tPhi_mini_u_d1_ubd}. We start by rewriting \eqref{eq_tPhi_mini_u_d1_ubd}.
Indeed, we deduce from \eqref{eq_nabla_vxi_inv_Gauss_id} the following:
\begin{equation*}
    \partial_\xi \tPhi\br{x,s,\xi,\eta}
    =
    x+
    \Delta_{\br{0,s}}
    \gamma_1\br{t_\ast\br{\xi,\eta}}
    =
    x+
    \dif{\gamma_1}
    ^{t_\ast\br{\xi,\eta+s}}
    _{t_\ast\br{\xi,\eta}}.
\end{equation*}
As a result, we obtain the following relation:
\begin{equation*}
    \abs{
        x+
        \dif{\gamma_1}
        ^{t_\ast\br{\xi,\eta+s}}
        _{t_\ast\br{\xi,\eta}}
    }
    \leq \abs{x/\lambda}^{\frac{1}{2}}
    \leq \abs{x/\lambda}^{\frac{1}{2}}\cdot \br{\lambda\abs{x}/4}^{\frac{1}{2}}
    = \abs{x}/2.
\end{equation*}
This implies
\begin{equation}\label{eq_tPhi_mini_u_gamma_1_geq_x}
    \abs{
        \dif{\gamma_1}
        ^{t_\ast\br{\xi,\eta+s}}
        _{t_\ast\br{\xi,\eta}}
    }
    \geq\abs{x}/2.
\end{equation}
On the other hand, by chain rule and \eqref{eq_nabla_t_ast_id}, we also obtain the following identities:
\begin{equation*}
    \partial^2_\xi \tPhi\br{x,s,\xi,\eta}
    =
    \frac{1}{\xi}
    \Delta_{\br{0,s}}
    \br{
    \frac{\br{\gamma'_1}^2\gamma'_2}{\vgamma'\wedge\vgamma''}
    }\br{t_\ast\br{\xi,\eta}}
    =
    \frac{1}{\xi}
    \cdot
    \dif{
    \frac{\br{\gamma'_1}^2\gamma'_2}{\vgamma'\wedge\vgamma''}
    }^{t_\ast\br{\xi,\eta+s}}
    _{t_\ast\br{\xi,\eta}}.
\end{equation*}
We now rewrite the above expression:
\begin{equation*}
    \partial^2_\xi \tPhi\br{x,s,\xi,\eta}
    =
    \frac{1}{\xi}
    \cdot
    \br{
        \dif{
        \frac{\br{\gamma'_1}^2\gamma'_2}{\vgamma'\wedge\vgamma''}
        }^{t_\ast\br{\xi,\eta+s}}
        _{t_\ast\br{\xi,\eta}}
    /
        \dif{\gamma_1}
        ^{t_\ast\br{\xi,\eta+s}}
        _{t_\ast\br{\xi,\eta}}
    }
    \cdot \dif{\gamma_1}
    ^{t_\ast\br{\xi,\eta+s}}
    _{t_\ast\br{\xi,\eta}}.
\end{equation*}
By Cauchy's mean value theorem, there is a value \(t_0\) between the two values \(t_\ast\br{\xi,\eta}\) and \(t_\ast\br{\xi,\eta+s}\) such that:
\begin{equation*}
    \partial^2_\xi \tPhi\br{x,s,\xi,\eta}
    =
    \frac{1}{\xi}
    \cdot
    \br{
        \br{\frac{\br{\gamma'_1}^2\gamma'_2}{\vgamma'\wedge\vgamma''}}'
    /
        \gamma'_1
    }\br{t_0}
    \cdot \dif{\gamma_1}
    ^{t_\ast\br{\xi,\eta+s}}
    _{t_\ast\br{\xi,\eta}}.
\end{equation*}
Applying the relation \eqref{eq_tPhi_mini_u_gamma_1_geq_x}, the fact that \(\br{\xi,\eta}\in \cA_{\delta,\Gamma}\), and the property \eqref{eq_tight_diag}, we deduce:
\begin{equation*}
    \abs{
        \partial^2_\xi \tPhi\br{x,s,\xi,\eta}
    }
    \geq
    A^{-O\br{1}}
    \abs{x}\cdot
    \abs{
        \br{\frac{\br{\gamma'_1}^2\gamma'_2}{\vgamma'\wedge\vgamma''}}'
    /
        \gamma'_1
    }\br{t_0}.
\end{equation*}
We perform the explicit computation and invoke the chirality assumption \eqref{eq_vgamma_chiral_ass}:
\begin{equation*}
    \br{\frac{\br{\gamma'_1}^2\gamma'_2}{\vgamma'\wedge\vgamma''}}'
    /
    \gamma'_1
    =
    \frac{
        2\gamma''_1\gamma'_2+\gamma'_1\gamma''_2
    }{
        \vgamma'\wedge\vgamma''
    }
    -
    \frac{\gamma'_1\gamma'_2}{\vgamma'\wedge\vgamma''}
    \cdot
    \frac{\vgamma'\wedge\vgamma'''}{\vgamma'\wedge\vgamma''}
\end{equation*}
\begin{equation*}
    =
    \frac{3}{2}\cdot\frac{\gamma''_1\gamma'_2+\gamma''_2\gamma'_1}{\vgamma'\wedge\vgamma''}
        -\frac{\gamma'_1\gamma'_2}{\vgamma'\wedge\vgamma''}\cdot\frac{\vgamma'\wedge\vgamma'''}{\vgamma'\wedge\vgamma''}
    -\frac{1}{2}=\Theta\br{\gamma_1,\gamma_2}-\frac{1}{2}\leq -\frac{1}{4}.
\end{equation*}
This concludes \eqref{eq_tPhi_mini_u_d2_lbd} and finishes the proof of \textbf{Lemma \ref{lem_tPhi_mini_u_est}}. 
\end{proof}

We are ready to show \eqref{eq_m_mini_u_pw_est}. Consider the two cases: \(\abs{\vlambda x}< 4\) and \(\abs{\vlambda x}\geq 4\). For the former, we recall \eqref{eq_m_mini_u_post_cv}, \eqref{eq_tPsi_mini_u}, and \eqref{eq_tPsi_mini_u_bd} and perform trivial estimate to conclude \eqref{eq_m_mini_u_pw_est}:
\begin{equation*}
    \abs{\eqref{eq_m_mini_u_post_cv}}
    \leq
    \nrm{
        \tPsi\br{\vlambda,\ts,\txi,\teta}
    }_{L^1\br{d\txi}}
    \leq
    \big\Vert\tPsi\big\Vert_{L^\infty}
    \leq A^{O\br{1}}
    \leq A^{O\br{1}}\Xi_\Gamma\br{1+\abs{\vlambda x}}^{-\frac{1}{2}}.
\end{equation*}
For the latter case \(\abs{\vlambda x}\geq 4\), we shall utilize \textbf{Lemma \ref{lem_tPhi_mini_u_est}} and apply \textbf{Proposition \ref{prop_VdC_reform}}, a reformulation of Van der Corput lemma. Indeed, for the purpose of proving \eqref{eq_m_mini_u_pw_est}, it suffices to consider those values of \(\ts,\txi,\teta\) such that \eqref{eq_tPsi_mini_u} does not vanish. We may thus apply \textbf{Lemma \ref{lem_tPhi_mini_u_est}} to deduce:
\begin{equation}\label{eq_tPhi_mini_u_VdC_phase_bd}
    \inf_{\txi\in J}
    \max_{j=1,2}\abs{
        \partial^j_\txi
        \abs{\vlambda}
        \tPhi\br{x,\ts,\txi,\teta}
    }^{\frac{1}{j}}
    \geq A^{-O\br{1}}\abs{\vlambda x}^{\frac{1}{2}},\quad
    J:=\BR{
        \txi\in\R
    \::\:
        \psi_\delta\br{\txi/\tlambda_1}
        \neq 0
    }.
\end{equation}
On the other hand, due to \eqref{eq_nabla_vxi_inv_Gauss_id}, we have the following identity:
\begin{equation*}
% \label{eq_tPhi_mini_u_to_inv_Gauss}
    \partial^2_\txi\tPhi\br{x,\ts,\txi,\teta}
    =\Delta_{\br{0,\ts}}\partial_\txi\gamma_{\Gamma,1}\br{\txi,\teta}.
\end{equation*}
Coupled the above with the sign condition \(\sgn \br{\teta}=\sgn \br{\teta+\ts}\) given by the support of \eqref{eq_tPsi_mini_u}, we deduce:\footnote{Recall \textbf{Definition \ref{def_uni_fin_ch}}.}
\begin{equation}\label{eq_tPhi_mini_u_conect_bd}
    M\br{\partial^2_\txi \abs{\vlambda} \tPhi\br{x,\ts,\cdot,\teta}}
    =
    M\br{\abs{\vlambda}\Delta_{\br{0,\ts}}\partial_\txi\gamma_{\Gamma,1}\br{\cdot,\teta}}
    =
    M\br{\Delta_{\br{0,\ts}}\partial_\txi\gamma_{\Gamma,1}\br{\cdot,\teta}}
    \leq \Xi_\Gamma.
\end{equation}
With \eqref{eq_tPsi_mini_u_bd}, \eqref{eq_tPsi_mini_u_partial_bd}, \eqref{eq_tPhi_mini_u_VdC_phase_bd}, and \eqref{eq_tPhi_mini_u_conect_bd}, \textbf{Proposition \ref{prop_VdC_reform}} yields the desired estimate \eqref{eq_m_mini_u_pw_est}:
\begin{equation*}
    \abs{\eqref{eq_m_mini_u_post_cv}}
    \lesssim 
    \frac{
        M\br{\partial^2_\txi \abs{\vlambda} \tPhi\br{x,\ts,\cdot,\teta}}
    }{
        \displaystyle{\inf_{\txi\in J}
        \max_{j=1,2}}\abs{
            \partial^j_\txi
            \abs{\vlambda}
            \tPhi\br{x,\ts,\txi,\teta}
        }^{\frac{1}{j}}
    }
    \cdot
    \br{
        \nrm{\tPsi}_{L^\infty\br{J}}
        +\nrm{\partial_\txi\tPsi}_{L^1\br{J}}
    }
    \leq
    A^{O\br{1}}
    \Xi_\Gamma
    \br{1+\abs{\vlambda x}}^{-\frac{1}{2}}.
\end{equation*}

Using now \eqref{eq_m_mini_u_pw_est}, we finish the proof of \eqref{eq_m_vlambda_small_uni}:
\begin{equation*}
    \nrm{m_{\vgamma,\rho,\vlambda}}^2_u:=
    \nrm{
        \4{1}\cD_{\br{0,s}}
        m_{\vgamma,\rho,\vlambda}\br{x,\eta}
    }_{L^\infty_{s,\eta}L^2_x\br{\mr{-1,1}}}
\end{equation*}
\begin{equation*}
    \leq A^{O\br{1}}\Xi_\Gamma
    \br{
    \int^1_0
        \frac{1}{1+\abs{\vlambda} x}
    dx
    }^{\frac{1}{2}}
    \leq A^{O\br{1}}
    \Xi_\Gamma
    \abs{\vlambda}^{-\frac{1}{2}}
    \log^{\frac{1}{2}}\br{\abs{\vlambda}}.
\end{equation*}

Finally, we will prove the estimate \eqref{eq_m_vlambda_big_uni}. We start by showing the following point-wise estimate:
\begin{equation}\label{eq_m_big_u_pw_est}
    \abs{
        \int
            \cD_{\br{0,s}}
            \cD_{\br{u,v}}
            m_{\vgamma,\rho,\vlambda}\br{\xi,\eta}
        d\xi
    }
    \leq 
    \frac{A^{O\br{1}}\Xi_\Gamma}{\abs{\vlambda}}
    \cdot
    \br{1+\frac{\abs{s\br{u,v}}}{\abs{\vlambda}}}^{-\frac{1}{2}}.
\end{equation}
Recall again \eqref{eq_m_vlambda_freq_trunc}. By \textbf{Proposition \ref{prop_m_fac}} and \textbf{Lemma \ref{lem_osc_st_or_nst}}, we may write:
\begin{equation*}
    \int
        \cD_{\br{0,s}}
        \cD_{\br{u,v}}
        m_{\vgamma,\rho,\vlambda}\br{\xi,\eta}
    d\xi
    =
    \int
        e\br{
            \Delta_{\br{0,s}}
            \Delta_{\br{u,v}}
            \mr{
                \br{\xi,\eta}\cdot
                \vgamma_\Gamma\br{\xi,\eta}
            }
        }
        \cD^{\br{\eta}}_s
        \cD^{\br{\xi,\eta}}_{\br{u,v}}
        \br{
            \psi_\delta\br{\frac{\xi}{\lambda_1}}
            \psi_\delta\br{\frac{\eta}{\lambda_2}}
            \Psi\br{\xi,\eta}
        }
    d\xi.
\end{equation*}

 We now perform the change of variables \eqref{eq_mini_u_CV_rule} along with
  \begin{equation*}
      u=\abs{\vlambda}\tu,\quad
      v=\abs{\vlambda}\tv.
  \end{equation*}
  This produces the following identity:
\begin{equation}\label{eq_m_big_u_post_cv}
\int
    \cD_{\br{0,s}}
    \cD_{\br{u,v}}
    m_{\vgamma,\rho,\vlambda}\br{\xi,\eta}
d\xi
 =
\int
    e\br{
        \abs{\vlambda}
        \cdot
        \tPhi\br{\ts,\tu,\tv,\txi,\teta}
    }
    \tPsi\br{\vlambda,\ts,\tu,\tv,\txi,\teta}
d\txi,
\end{equation}
where 
\begin{equation}\label{eq_tPhi_big_u}
    \tPhi\br{\ts,\tu,\tv,\txi,\teta}:=
    \Delta_{\br{0,\ts}}
    \Delta_{\br{\tu,\tv}}
    \mr{
        \br{\txi,\teta}\cdot
        \vgamma_\Gamma\br{\txi,\teta}
    };
\end{equation}
\begin{equation}\label{eq_tPsi_big_u}
    \tPsi\br{\vlambda,\ts,\tu,\tv,\txi,\teta}
    :=
    \abs{\vlambda}
    \cD^{\br{\teta}}_\ts
    \cD^{\br{\txi,\teta}}_{\br{\tu,\tv}}
    \br{
        \psi_\delta
            \br{\txi/\tlambda_1}
        \psi_\delta
            \br{\teta/\tlambda_2}
        \Psi\br{
            \abs{\vlambda}\txi,
            \abs{\vlambda}\teta
        }
    }.
\end{equation}
Our next step is to estimate $\|\widetilde{\Psi}\|_{L^{\infty}}$ and $\|\partial_{\widetilde{\xi}}\widetilde{\Psi}\|_{L^{1}}$.
Utilizing \eqref{eq_m_fac_amp}, we perform computation similar to \eqref{eq_tPsi_mini_u_bd}:
\begin{equation}\label{eq_tPsi_big_u_bd}
    \big\Vert\tPsi\big\Vert_{L^\infty}
    \leq
    \abs{\vlambda}
    \cdot
    \nrm{
        \Psi\br{
            \abs{\vlambda}\txi,
            \abs{\vlambda}\teta
        }
    }^4_{L^\infty
        \br{
            \cA_{\delta,\Gamma}
        }
    }
    \leq 
    \abs{\vlambda}
    \br{
        A^{O\br{1}}
        \abs{\vlambda}^{-\frac{1}{2}}
    }^4
    =A^{O\br{1}}/\abs{\vlambda}
\end{equation}
and computation similar to \eqref{eq_tPsi_mini_u_partial_bd}:
\begin{equation*}
    \nrm{\partial_{\txi}\tPsi}_{L^1}
    \leq
    4\abs{\vlambda}
    \cdot
    \nrm{\phi'}_{L^1}
    \nrm{
        \Psi\br{
            \abs{\vlambda}\txi,
            \abs{\vlambda}\teta
        }
    }^4_{L^\infty
        \br{
            \cA_{\delta,\Gamma}
        }
    }
\end{equation*}
\begin{equation*}
    +
    4
    \abs{\vlambda}
    \cdot 
    \delta \tlambda_1
    \cdot
    \nrm{\phi}_{L^1}
    \nrm{
        \Psi\br{
            \abs{\vlambda}\txi,
            \abs{\vlambda}\teta
        }
    }^3_{L^\infty
        \br{
            \cA_{\delta,\Gamma}
        }
    }
    \abs{\vlambda}\cdot
    \nrm{
        \nabla
        \Psi\br{
            \abs{\vlambda}\txi,
            \abs{\vlambda}\teta
        }
    }_{L^\infty
        \br{
            \cA_{\delta,\Gamma}
        }
    }
\end{equation*}
\begin{equation}\label{eq_tPsi_big_u_partial_bd}
    \leq
    \abs{\vlambda}
    A^{O\br{1}}
    \abs{\vlambda}^{-\frac{4}{2}}
    +
    \abs{\vlambda}^2
    A^{O\br{1}-O\br{C}}
    \abs{\vlambda}^{-\frac{3}{2}}
    A^{O\br{1}}
    \abs{\vlambda}^{-\frac{3}{2}}
    \leq A^{O\br{1}}/\abs{\vlambda}.
\end{equation}
This completes the estimates for the amplitude function \(\tPsi\). Next, we establish an estimate for the phase function $\widetilde{\Phi}$.

\begin{lemma}\label{lem_tPhi_big_u_est}
    Let \(\lambda \abs{s\br{u,v}} \geq 1\) and \(U\subseteq  \cA_{\delta,\Gamma}\) convex with \(\dia\br{U}\leq A^{-C}\). Given that
    \begin{equation}\label{eq_tPhi_big_u_sgn_cond}
        \br{\xi,\eta},\,\br{\xi,\eta+s},\,\br{\xi+u,\eta+v},\,\br{\xi+u,\eta+s+v}\in  U,
    \end{equation}
    the phase function defined by \eqref{eq_tPhi_big_u} satisfies the following estimate:
    \begin{equation}\label{eq_tPhi_big_u_est}
        \max\br{
            \abs{
                \lambda
                \partial_\xi
                \tPhi
                \br{
                    s,u,v,\xi,\eta
                }
            },
            \abs{
                \lambda
                \partial^2_\xi
                \tPhi
                \br{
                    s,u,v,\xi,\eta
                }
            }^{\frac{1}{2}}
        }
        \geq A^{-O\br{1}}
        \abs{\lambda s\br{u,v}}^{\frac{1}{2}}.
    \end{equation}
\end{lemma}

\begin{proof}
    We will first prove
    \begin{equation}\label{eq_tPhi_big_u_vec_bd}
        \abs{
            \begin{pmatrix}
                \partial_\xi \tPhi\br{s,u,v,\xi,\eta}\\
                \partial^2_\xi 
                \tPhi\br{s,u,v,\xi,\eta}
            \end{pmatrix}
        }
        \geq A^{-O\br{1}}
        \abs{s\br{u,v}}.
    \end{equation}
    Since derivative commutes with difference, we have
\begin{equation*}
    \begin{pmatrix}
        \partial_\xi \tPhi\br{s,u,v,\xi,\eta}\\
        \partial^2_\xi 
        \tPhi\br{s,u,v,\xi,\eta}
    \end{pmatrix}
    =
    \begin{pmatrix}
        \Delta_{\br{0,s}}\Delta_{\br{u,v}}
        \partial_\xi\mr{\br{\xi,\eta}\cdot \vgamma_\Gamma\br{\xi,\eta}}\\
        \Delta_{\br{0,s}}\Delta_{\br{u,v}}
        \partial^2_\xi\mr{\br{\xi,\eta}\cdot \vgamma_\Gamma\br{\xi,\eta}}
    \end{pmatrix}
    .
\end{equation*}
Moreover, we may utilize the identity \eqref{eq_nabla_vxi_inv_Gauss_id} and equate the above with
\begin{equation}\label{eq_tPhi_big_u_vec_to_inv_Gauss}
    =
    \begin{pmatrix}
        \Delta_{\br{0,s}}\Delta_{\br{u,v}}
        \gamma_{\Gamma,1}\br{\xi,\eta}\\
        \Delta_{\br{0,s}}\Delta_{\br{u,v}}
        \partial_\xi\gamma_{\Gamma,1}\br{\xi,\eta}
    \end{pmatrix}
    .
\end{equation}

We claim that for a \(C^2\) function \(f:A_{\delta,\Gamma}\to \R\), there is \(\br{\xi_1,\eta_1}\in U\) such that the following identity for iterated finite difference holds:
    \begin{equation}
        \Delta_{\br{0,s}}\Delta_{\br{u,v}}f\br{\xi,\eta}=s\cdot 
        \nabla \partial_\eta
        f\br{\xi_1,\eta_1}
        % \begin{pmatrix}
        % \partial_\xi
        % \partial_\eta
        % f\br{\xi_1,\eta_1}&
        % \partial^2_\eta 
        % f\br{\xi_1,\eta_1}
        % \end{pmatrix}
        \cdot 
        \begin{pmatrix}
            u\\v
        \end{pmatrix} .
    \end{equation}
The above claim holds because:
By the mean value theorem, we have
    \begin{equation}
        (\Delta_{(0,s)}\Delta_{(u,v)})f(\xi,\eta)=(\Delta_{(u,v)}f)(\xi,\eta+s)-(\Delta_{(u,v)}f)(\xi,\eta)
    \end{equation}
    \begin{equation}
        =s\cdot \partial_\eta\Delta_{\br{u,v}}
        f\br{\xi,\eta+\theta_1 s}=
        s\cdot \Delta_{\br{u,v}}\partial_\eta f\br{\xi,\eta+\theta_1 s}
    \end{equation}
    for some \(\theta_1\in\br{0,1}\).
    Then by another application of the mean value theorem, we can further equate the above with
    \begin{equation*}
        =s\cdot\br{
        \partial_\eta 
        f\br{\xi+u,\eta+\theta_1 s +v}-
        \partial_\eta
        f\br{\xi,\eta+\theta_1 s }
        }
        % =s\cdot \int_{0}^{1}\nabla\partial_{\eta}f(\xi+\theta u,\widetilde{\eta}+\theta v)\cdot 
        % \begin{pmatrix}
        %     u\\v
        % \end{pmatrix}
        % d\theta
    \end{equation*}
    \begin{equation*}
        =s\cdot 
        \nabla \partial_\eta f
        \br{\xi+\theta_2 u,
        \eta+\theta_1 s+\theta_2 v}
        % \begin{pmatrix}
        %     \partial_\xi\partial_\eta
        %     f\br{\xi+\theta_2 u,\eta+\theta_1 s+\theta_2 v}&
        %     \partial_{\eta}^{2}
        %     f\br{\xi+\theta_2 u,\eta+\theta_1 s+\theta_2 v}
        % \end{pmatrix}
        \cdot 
        \begin{pmatrix}
            u\\v
        \end{pmatrix} 
    \end{equation*}
for some \(\theta_2\in\br{0,1}\). In other words, we may take \(\br{\xi_1,\eta_1}:=\br{\xi+\theta_2 u,\eta+\theta_1 s+\theta_2 v}
\in U\).

With the claim above, We obtain the following identity:
\begin{equation*}
    \begin{pmatrix}
        \partial_\xi \tPhi\br{s,u,v,\xi,\eta}\\
        \partial^2_\xi 
        \tPhi\br{s,u,v,\xi,\eta}
    \end{pmatrix}
    =
    \eqref{eq_tPhi_big_u_vec_to_inv_Gauss}
    =
    s\cdot
    \begin{pmatrix}
        \nabla \partial_\eta \gamma_{\Gamma,1}\br{\xi_1,\eta_1}
        \\
        \partial_\xi\nabla \partial_\eta
        \gamma_{\Gamma,1}\br{\xi_2,\eta_2}
    \end{pmatrix}
    \cdot
    \begin{pmatrix}
        u\\v
    \end{pmatrix}.
\end{equation*}
We may force a main term to have its entries sharing the same input a the cost of a small error term. That is, via telescoping, we identify the above with:
\begin{equation*}
    =
    s\cdot
    \dif{
    \begin{pmatrix}
        \partial_\eta \nabla \gamma_{\Gamma,1}
        \\
        \partial_\xi \partial_\eta \nabla
        \gamma_{\Gamma,1}
    \end{pmatrix}
    }_{\br{\xi_2,\eta_2}}
    \cdot
    \begin{pmatrix}
        u\\v
    \end{pmatrix}
    +
    s\cdot
    \begin{pmatrix}
        \dif{\partial_\eta \nabla \gamma_{\Gamma,1}}^{\br{\xi_1,\eta_1}}_{\br{\xi_2,\eta_2}}
        \\
        \vnull
    \end{pmatrix}
    \cdot
    \begin{pmatrix}
        u\\v
    \end{pmatrix}
    =: 
    s\cdot 
    \vM_\main
    \cdot
    \begin{pmatrix}
        u\\v
    \end{pmatrix}
    +
    s
    \cdot
    \vM_\err
    \cdot
    \begin{pmatrix}
        u\\v
    \end{pmatrix}.
\end{equation*}
We claim that it suffices to show the following upper bounds:
\begin{equation}\label{eq_tPhi_big_u_matrix_op_bds}
    \nrm{\vM_\main}_\op\leq A^{O\br{1}}
    ,\quad
    \nrm{\vM_\err}_\op \leq A^{O\br{1}-C}
\end{equation}
and lower bound
\begin{equation}\label{eq_tPhi_big_u_matrix_det_bd}
    \abs{\det\vM_\main}\geq A^{-O\br{1}}
\end{equation}
to prove \eqref{eq_tPhi_big_u_vec_bd}.
Indeed, assuming both \eqref{eq_tPhi_big_u_matrix_op_bds} and \eqref{eq_tPhi_big_u_matrix_det_bd}, \(\vM_\main\) admits a singular value decomposition:
\begin{equation*}
    \vM_\main= 
    \vP 
    \begin{pmatrix}
        \sigma_1 & 0\\
        0 &\sigma_2
    \end{pmatrix}
    \vQ^\ast,\quad
    \sigma_1 \geq \sigma_2\geq 0.
\end{equation*}
Using the following facts:
\begin{equation*}
    \sigma_1=\nrm{\vM_\main}_\op,\quad
    \sigma_1\sigma_2=\abs{\det \vM_\main},
\end{equation*}
we deduce the following lower bound:
\begin{equation*}
    \abs{
        \vM_\main \cdot
        \begin{pmatrix}
            u\\ v
        \end{pmatrix}
    }
    \geq \sigma_2
    \cdot
    \abs{
    \br{u,v}
    }
    =
    \frac{
        \abs{\det \vM_\main}
    }{
        \nrm{\vM_\main}_\op
    }
    \cdot
    \abs{
    \br{u,v}
    }
    \geq 
    A^{-O\br{1}}
    \cdot
    \abs{
    \br{u,v}
    }
    .
\end{equation*}
Combing this with the trivial estimate:
\begin{equation*}
    \abs{
        \vM_\err
        \cdot
        \begin{pmatrix}
            u\\v
        \end{pmatrix}
    }
    \leq
    \nrm{\vM_\err}_\op
    \cdot
    \abs{
        \br{u,v}
    }
    \leq 
    A^{O\br{1}-C}
    \cdot
    \abs{
        \br{u,v}
    },
\end{equation*}
we conclude the desired estimate \eqref{eq_tPhi_big_u_vec_bd}:
\begin{equation*}
    \abs{
        \begin{pmatrix}
            \partial_\xi \tPhi\br{s,u,v,\xi,\eta}\\
            \partial^2_\xi 
            \tPhi\br{s,u,v,\xi,\eta}
        \end{pmatrix}
    }
    \geq
    \abs{
    s\cdot 
    \vM_\main
    \cdot
    \begin{pmatrix}
        u\\v
    \end{pmatrix}
    }
    -
    \abs{
    s
    \cdot
    \vM_\err
    \cdot
    \begin{pmatrix}
        u\\v
    \end{pmatrix}
    }
\end{equation*}
\begin{equation*}
    \geq  
    \abs{s}\cdot A^{-O\br{1}}\abs{\br{u,v}}
    -\abs{s}\cdot A^{O\br{1}-C}\abs{\br{u,v}}
    \geq A^{-O\br{1}}\abs{s\br{u,v}}.
\end{equation*}
Return to the proof of \eqref{eq_tPhi_big_u_matrix_op_bds} and \eqref{eq_tPhi_big_u_matrix_det_bd}. 
To proceed, we shall perform a few explicit calculation via identities provided by \textbf{Lemma \ref{lem_nabla_ids}}. Let the input variable be \(\vxi:=\br{\xi_1,\xi_2}\in \cA_{\delta,\Gamma}\). We have:
\begin{equation}\label{eq_tPhi_big_u_d1_id}
    \nabla\gamma_{\Gamma,1}
    =
    \gamma'_1\br{t_\ast}\nabla t_\ast
    =
    \frac{\gamma'_1\br{t_\ast}}{\xi_1}\cdot
    \br{
        \frac{
            \gamma'_1\gamma'_2
        }{
            \vgamma'\wedge\vgamma''
        },
        \frac{
            \br{\vgamma'_2}^2
        }{
            \vgamma'\wedge\vgamma''
        }
    }
    \br{t_\ast}
    =
    \frac{1}{\xi_1}\cdot
    \vchi_\vgamma\br{t_\ast};
\end{equation}
\begin{equation}\label{eq_tPhi_big_u_d2_id}
    \partial_2 \nabla \gamma_{\Gamma,1}
    =\frac{1}{\xi_1}\cdot
    \vchi'_\vgamma\br{t_\ast}\partial_\eta t_\ast
    =-\frac{1}{\xi_1\xi_2}
    \cdot
    \vchi'_\vgamma\br{t_\ast}\cdot \br{\frac{\gamma'_1\gamma'_2}{\vgamma'\wedge\vgamma''}}\br{t_\ast};
\end{equation}
\begin{equation*}
    \partial_j \partial_2 \nabla \gamma_{\Gamma,1}
    =
    \frac{1}{\xi_j\xi_1\xi_2}
    \cdot
    \vchi'_\vgamma\br{t_\ast}\cdot \br{\frac{\gamma'_1\gamma'_2}{\vgamma'\wedge\vgamma''}}\br{t_\ast}
    -
    \frac{1}{\xi_1\xi_2}
    \cdot
    \vchi''_\vgamma\br{t_\ast}\cdot \br{\frac{\gamma'_1\gamma'_2}{\vgamma'\wedge\vgamma''}}\br{t_\ast}\partial_j t_\ast
    -
    \frac{1}{\xi_1\xi_2}
    \cdot
    \vchi'_\vgamma\br{t_\ast}\cdot 
    \br{\frac{\gamma'_1\gamma'_2}{\vgamma'\wedge\vgamma''}}'\br{t_\ast}\partial_j t_\ast
\end{equation*}
\begin{equation}\label{eq_tPhi_big_u_d3_id}
    =
    \frac{1}{\xi_j\xi_1\xi_2}
    \cdot
    \mr{
    \vchi'_\vgamma
    +
    \br{-1}^j
    \vchi''_\vgamma\cdot
    \frac{\gamma'_1\gamma'_2}{\vgamma'\wedge\vgamma''}
    +
    \br{-1}^j
    \vchi'_\vgamma\cdot 
    \br{\frac{\gamma'_1\gamma'_2}{\vgamma'\wedge\vgamma''}}'
    }
    \br{t_\ast}
    \cdot
    \br{\frac{\gamma'_1\gamma'_2}{\vgamma'\wedge\vgamma''}}\br{t_\ast}
    .
\end{equation}
Given the three identities \eqref{eq_tPhi_big_u_d1_id}, \eqref{eq_tPhi_big_u_d2_id}, and \eqref{eq_tPhi_big_u_d3_id}, the main term estimate in \eqref{eq_tPhi_big_u_matrix_op_bds} follows from the \(A\)-admissibility assumption on \(\vgamma:I\simto \Gamma\):
\begin{equation*}
    \nrm{
        \vM_\main
    }_\op
    \leq
    \sqrt{2}
    \max\br{
        \abs{\eqref{eq_tPhi_big_u_d2_id}}
        ,\abs{\eqref{eq_tPhi_big_u_d3_id}}
    }
    \leq A^{O\br{1}}.
\end{equation*}
To prove the error term estimate in \eqref{eq_tPhi_big_u_matrix_op_bds}, we apply fundamental theorem of calculus to the entries:
\begin{equation*}
    \dif{\partial_\eta \nabla\gamma_{\Gamma,1}}^{\br{\xi_1,\eta_1}}_{\br{\xi_2,\eta_2}}
    =
    \int^1_0
        \br{\xi_1-\xi_2,\eta_1-\eta_2}
        \cdot
        \dif{
        \begin{pmatrix}
            \partial_\xi
            \partial_\eta \nabla\gamma_{\Gamma,1}\\
            \partial^2_\eta \nabla\gamma_{\Gamma,1}
        \end{pmatrix}
        }_{
            \br{\xi_2,\eta_2}
            +\theta
            \br{\xi_1-\xi_2,\eta_1-\eta_2}
        }
    d\theta.
\end{equation*}
By \eqref{eq_tPhi_big_u_d3_id} and the \(A\)-admissibility assumption on \(\vgamma:I\simto \Gamma\), we deduce:
\begin{equation*}
    \abs{
        \dif{\partial_\eta \nabla\gamma_{\Gamma,1}}^{\br{\xi_1,\eta_1}}_{\br{\xi_2,\eta_2}}
    }
    \leq A^{O\br{1}}
    \abs{\br{\xi_1-\xi_2,\eta_1-\eta_2}}
    \leq A^{O\br{1}}
    \cdot \dia\br{U}
    \leq A^{O\br{1}-C}.
\end{equation*}
This justifies the error term estimate in \eqref{eq_tPhi_big_u_matrix_op_bds}:
\begin{equation*}
    \nrm{\vM_\err}_\op
    \leq
    \abs{
        \dif{\partial_\eta \nabla\gamma_{\Gamma,1}}^{\br{\xi_1,\eta_1}}_{\br{\xi_2,\eta_2}}
    }
    \leq A^{O\br{1}-C}.
\end{equation*}
Finally, to prove \eqref{eq_tPhi_big_u_matrix_det_bd}, we shall utilize the alternating nature of the determinant expressions. Observe that identities \eqref{eq_tPhi_big_u_d2_id} and \eqref{eq_tPhi_big_u_d3_id} give
\begin{equation*}
    \det \vM_\main
    =
    -\frac{1}{\xi\eta}
    \cdot
    \vchi'_\vgamma\br{t_\ast}\cdot \br{\frac{\gamma'_1\gamma'_2}{\vgamma'\wedge\vgamma''}}\br{t_\ast}
    \wedge
    \frac{1}{\xi^2\eta}
    \cdot
    \mr{
    \vchi'_\vgamma
    -
    \vchi''_\vgamma\cdot
    \frac{\gamma'_1\gamma'_2}{\vgamma'\wedge\vgamma''}
    -
    \vchi'_\vgamma\cdot 
    \br{\frac{\gamma'_1\gamma'_2}{\vgamma'\wedge\vgamma''}}'
    }
    \br{t_\ast}
    \cdot
    \br{\frac{\gamma'_1\gamma'_2}{\vgamma'\wedge\vgamma''}}\br{t_\ast}
    .
\end{equation*}
By removing all parallel vector pairs, the above identity can be simplified into:
\begin{equation*}
    \det \vM_\main
    =
    \frac{1}{\xi^3\eta^2}
    \br{\vchi'_\vgamma\wedge \vchi''_\vgamma}\br{t_\ast}
    \cdot
    \br{\frac{\gamma'_1\gamma'_2}{\vgamma'\wedge\vgamma''}}^3\br{t_\ast}.
\end{equation*}
Invoking the \(A\)-admissibility assumption on \(\vgamma:I\simto\Gamma\), we conclude \eqref{eq_tPhi_big_u_matrix_det_bd}
\begin{equation*}
    \abs{
        \det \vM_\main
    }
    \geq A^{-O\br{1}}
    \abs{
        \vchi'_\vgamma\wedge \vchi''_\vgamma
    }
    \br{t_\ast}
    \cdot
    \abs{\frac{\gamma'_1\gamma'_2}{\vgamma'\wedge\vgamma''}}^3\br{t_\ast}
    \geq A^{-O\br{1}}.
\end{equation*}
Now that \eqref{eq_tPhi_big_u_vec_bd} has been proven, we shall now deduce \eqref{eq_tPhi_big_u_est} from \eqref{eq_tPhi_big_u_vec_bd}. Let
\begin{equation*}
    M:=
    \max\br{
        \abs{
            \lambda
            \partial_\xi
            \tPhi
            \br{
                s,u,v,\xi,\eta
            }
        },
        \abs{
            \lambda
            \partial^2_\xi
            \tPhi
            \br{
                s,u,v,\xi,\eta
            }
        }^{\frac{1}{2}}
    }.
\end{equation*}
We have the following trivial estimate:
\begin{equation*}
    \br{M^2+M^4}^{\frac{1}{2}}\geq
    \lambda
    \cdot
    \abs{
        \begin{pmatrix}
            \partial_\xi \tPhi\br{s,u,v,\xi,\eta}\\
            \partial^2_\xi 
            \tPhi\br{s,u,v,\xi,\eta}
        \end{pmatrix}
    }
    \geq A^{-O\br{1}}
    \abs{\lambda s\br{u,v}} \geq A^{-O\br{1}}.
\end{equation*}
This implies the trivial estimate \(M\geq A^{-O\br{1}}\) and thus:
\begin{equation*}
     A^{O\br{1}} M^2\geq 
     \br{M^2+M^4}^{\frac{1}{2}}
     \geq A^{-O\br{1}}
    \abs{\lambda s\br{u,v}}.
\end{equation*}
After rearrangement, we conclude \eqref{eq_tPhi_big_u_est}. This completes the proof of \textbf{Lemma \ref{lem_tPhi_big_u_est}}.
\end{proof}

We return to the proof of \eqref{eq_m_big_u_pw_est}. 
Consider the two cases:
\begin{equation*}
    \frac{\abs{s\br{u,v}}}{\abs{\vlambda}}=\abs{\vlambda}\cdot\abs{\ts\br{\tu,\tv}}\leq 1;\quad
    \frac{\abs{s\br{u,v}}}{\abs{\vlambda}}=\abs{\vlambda}\cdot\abs{\ts\br{\tu,\tv}}\geq 1.
\end{equation*}
For the former, we recall \eqref{eq_m_big_u_post_cv}, \eqref{eq_tPsi_big_u}, and \eqref{eq_tPsi_big_u_bd} and perform trivial estimate to conclude \eqref{eq_m_big_u_pw_est}:
\begin{equation*}
    \abs{\eqref{eq_m_big_u_post_cv}}
    \leq
    \nrm{
        \tPsi\br{\vlambda,\ts,\tu,\tv,\txi,\teta}
    }_{L^1\br{d\txi}}
    \leq
    \big\Vert\tPsi\big\Vert_{L^\infty}
    \leq A^{O\br{1}}/\abs{\vlambda}
    \leq \frac{A^{O\br{1}}\Xi_\Gamma}{\abs{\vlambda}}\cdot\br{1+\frac{\abs{s\br{u,v}}}{\abs{\vlambda}}}^{-\frac{1}{2}}.
\end{equation*}

For the latter case \(\frac{\abs{s\br{u,v}}}{\abs{\vlambda}}=\abs{\vlambda}\cdot\abs{\ts\br{\tu,\tv}}\geq 1\), we shall utilize \textbf{Lemma \ref{lem_tPhi_big_u_est}} and apply \textbf{Proposition \ref{prop_VdC_reform}} again.
Note that for the purpose of proving \eqref{eq_m_big_u_pw_est}, the support of \eqref{eq_tPsi_big_u} naturally imposes condition \eqref{eq_tPhi_big_u_sgn_cond}
\begin{equation}\label{eq_tPhi_big_u_sgn_cond_specific}
    \br{\txi,\teta},\,\br{\txi,\teta+\ts},\,\br{\txi+\tu,\teta+\tv},\,\br{\txi+\tu,\teta+\ts+\tv}\in U
\end{equation}
with the convex set \(U\) given by
\begin{equation}\label{eq_def_convex_set_4_big_U}
    U:=\BR{
        \br{\txi,\teta}\in \R^2
    \::\:
        \psi_\delta\br{\txi/\tlambda_1}
        \psi_\delta\br{\teta/\tlambda_2}\neq 0
    }\subseteq\cA_{\delta,\Gamma}.
\end{equation}
Consider now the interval
\begin{equation*}
    J_\tu:=\BR{
        \txi\in\R
    \::\:
        \psi_\delta\br{\txi/\tlambda_1}
        \psi_\delta\br{\br{\txi+\tu}/\tlambda_1}
        \neq 0
    }.
\end{equation*}
By \textbf{Lemma \ref{lem_tPhi_big_u_est}}, we have:
\begin{equation}\label{eq_tPhi_big_u_VdC_phase_bd}
    \inf_{\txi\in J_{\tu}}
    \max_{j=1,2}\abs{
        \partial^j_\txi
        \abs{\vlambda}
        \tPhi\br{\ts,\tu,\tv,\txi,\teta}
    }^{\frac{1}{j}}\geq A^{-O\br{1}}\abs{\vlambda}^{\frac{1}{2}}
    \cdot\abs{\ts\br{\tu,\tv}}^{\frac{1}{2}}
    =A^{-O\br{1}}
    \cdot
    \frac{\abs{s\br{u,v}}^{\frac{1}{2}}}{\abs{\vlambda}^{\frac{1}{2}}}.
\end{equation}
On the other hand, due to \eqref{eq_nabla_vxi_inv_Gauss_id}, we have the following identity:
\begin{equation*}
    \partial^2_\txi\tPhi\br{\ts,\tu,\tv,\txi,\teta}
    =\Delta_{\br{0,\ts}}\Delta_{\br{\tu,\tv}}\partial_\txi\gamma_{\Gamma,1}\br{\txi,\teta}.
\end{equation*}
Couple the above with the sign condition \(\sgn\br{\teta}=\sgn\br{\teta+\ts}=\sgn\br{\teta+\tv}=\sgn\br{\teta+\ts+\tv}\) given by \eqref{eq_tPhi_big_u_sgn_cond_specific} and \eqref{eq_def_convex_set_4_big_U}, we deduce the following estimate:\footnote{Recall again \textbf{Definition \ref{def_uni_fin_ch}}.}
\begin{equation}\label{eq_tPhi_big_u_conect_bd}
    M\br{\partial^2_\txi \abs{\vlambda} \tPhi\br{\ts,\tu,\tv,\cdot,\teta}}
    =
    M\br{\abs{\vlambda}\Delta_{\br{0,\ts}}\Delta_{\br{\tu,\tv}}\partial_\txi\gamma_{\Gamma,1}\br{\cdot,\teta}}
    =
    M\br{\Delta_{\br{0,\ts}}\Delta_{\br{\tu,\tv}}\partial_\txi\gamma_{\Gamma,1}\br{\cdot,\teta}}
    \leq \Xi_\Gamma.
\end{equation}
With \eqref{eq_tPsi_big_u_bd}, \eqref{eq_tPsi_big_u_partial_bd}, \eqref{eq_tPhi_big_u_VdC_phase_bd}, and \eqref{eq_tPhi_big_u_conect_bd}, \textbf{Proposition \ref{prop_VdC_reform}} yields the desired estimate \eqref{eq_m_big_u_pw_est}:
\begin{equation*}
    \abs{\eqref{eq_m_big_u_post_cv}}
    \lesssim 
    \frac{
        M\br{\partial^2_\txi \abs{\vlambda} \tPhi\br{\ts,\tu,\tv,\cdot,\teta}}
    }{
        \displaystyle{\inf_{\txi\in J_\tu}
        \max_{j=1,2}
        }\abs{
            \partial^j_\txi
            \abs{\vlambda}
            \tPhi\br{\ts,\tu,\tv,\txi,\teta}
        }^{\frac{1}{j}}
    }
    \cdot
    \br{
        \nrm{\tPsi}_{L^\infty\br{J_\tu}}
        +\nrm{\partial_\txi\tPsi}_{L^1\br{J_\tu}}
    }
    \leq
    \frac{A^{O\br{1}}
    \Xi_\Gamma}{\abs{\vlambda}}
    \cdot
    \br{1+\frac{\abs{s\br{u,v}}}{\abs{\vlambda}}}^{-\frac{1}{2}}.
\end{equation*}

Using now \eqref{eq_m_big_u_pw_est}, we finish the proof of \eqref{eq_m_vlambda_big_uni}:
\begin{equation*}
    \nrm{m_{\vgamma,\rho,\vlambda}}^4_U=
    \nrm{
        \int
            \cD_{\br{0,s}}
            \cD_{\br{u,v}}
            m_{\vgamma,\rho,\vlambda}\br{\xi,\eta}
        d\xi
    }_{L^\infty_{u,\eta}L^1_s L^2_v \br{\abs{s},\abs{v}\leq \abs{\lambda}}}
\end{equation*}
\begin{equation*}
    \leq
    \nrm{
        \nrm{
            \frac{A^{O\br{1}}\Xi_\Gamma}{\abs{\vlambda}}
            \cdot
            \br{
            1+
            \frac{\abs{vs}}{\abs{\vlambda}}
            }^{-\frac{1}{2}}
        }_{L^1\br{ds,\abs{s}\leq\abs{\vlambda}}}
    }_{L^2\br{dv, \abs{v}\leq \abs{\vlambda}}}
\end{equation*}
\begin{equation*}
    =
    A^{O\br{1}}\Xi_\Gamma
    \nrm{
        \nrm{
            \br{
            1+
            \abs{v\ts}
            }^{-\frac{1}{2}}
        }_{L^1\br{d\ts,\abs{\ts}\leq 1}}
    }_{L^2\br{dv, \abs{v}\leq \abs{\vlambda}}}
\end{equation*}
\begin{equation*}
    \leq
    A^{O\br{1}}\Xi_\Gamma
    \br{
    \int^{\abs{\vlambda}}_0
        \frac{1}{1+
        v}
    dv
    }^{\frac{1}{2}}
    \leq 
    A^{O\br{1}}\Xi_\Gamma\log^{\frac{1}{2}}\br{\abs{\vlambda}}.
\end{equation*}
This completes the proof of \textbf{Proposition \ref{prop_m_uni_norm}}.

\section{Smoothing inequalities for averaging operators and their variants}
% Before entering the analysis, we introduce the following key lemma that allows us to reduce the argument to the analysis on smaller segments of curves.
% \begin{lemma}[Monotonicity of \(\Xi_\Gamma\)]\label{lem_M_Gamma_mono}
%     Given a \(C^1\) curve \(\Gamma\subseteq \R^2\) with injective Gauss map, if a curve \(\Gamma'\) is a subset of \(\Gamma\), the following inequality holds:
%     \begin{equation*}
%         \Xi_{\Gamma'}\leq \Xi_\Gamma+O\br{1}
%         \lesssim \Xi_\Gamma.
%     \end{equation*}
% \end{lemma}
Before entering the analysis, we shall first the key statement \textbf{Lemma \ref{lem_M_Gamma_mono}} that allows us to reduce the argument to the analysis on smaller segments of curves.
\subsection{Proof of Lemma \ref{lem_M_Gamma_mono}: monotonicity of \(\Xi_\Gamma\)}
\begin{proof}
    Recall \textbf{Definition \ref{def_uni_fin_ch}}.
    Let \(s,u,v,w,\eta\in\R\) be such that satisfy \(p\br{s,v,\eta}\).
    By symmetry it suffices to establish the following two inequalities:
    \begin{equation}
    \label{eq_sub_curve_leq_curve_1}
        C\br{
            \BR{
                \xi\in\R
            \::\:
                \Delta_{\br{0,s}}\partial_1\gamma_{\Gamma',1}\br{\xi,\eta}=w
            }
        }
        \leq
        C\br{
            \BR{
                \xi\in\R
            \::\:
                \Delta_{\br{0,s}}\partial_1\gamma_{\Gamma,1}\br{\xi,\eta}=w
            }
        }
        +O\br{1};
    \end{equation}
    \begin{equation}
    \label{eq_sub_curve_leq_curve_2}
        C\br{
            \BR{
                \xi\in\R
            \::\:
                \Delta_{\br{0,s}}\Delta_{\br{u,v}}\partial_1\gamma_{\Gamma',1}\br{\xi,\eta}=w
            }
        }
        \leq
        C\br{
            \BR{
                \xi\in\R
            \::\:
                \Delta_{\br{0,s}}\Delta_{\br{u,v}}\partial_1\gamma_{\Gamma,1}\br{\xi,\eta}=w
            }
        }
        +O\br{1}.
    \end{equation}
    By invoking \textbf{Remark \ref{rmk_connect_rel}},
    \eqref{eq_sub_curve_leq_curve_1} and \eqref{eq_sub_curve_leq_curve_2} can be deduced by finding sets \(W_1,W_2\subseteq \R\) with each at most \(O\br{1}\) many connected components such that the following two identities hold:
    \begin{equation}
    \label{eq_sub_curve_cap_curve_1}
            \BR{
                \xi\in\R
            \::\:
                \Delta_{\br{0,s}}\partial_1\gamma_{\Gamma',1}\br{\xi,\eta}=w
            }
        =
            \BR{
                \xi\in\R
            \::\:
                \Delta_{\br{0,s}}\partial_1\gamma_{\Gamma,1}\br{\xi,\eta}=w
            }
        \cap W_1;
    \end{equation}
    \begin{equation}
    \label{eq_sub_curve_cap_curve_2}
            \BR{
                \xi\in\R
            \::\:
                \Delta_{\br{0,s}}\Delta_{\br{u,v}}\partial_1\gamma_{\Gamma',1}\br{\xi,\eta}=w
            }
        =
            \BR{
                \xi\in\R
            \::\:
                \Delta_{\br{0,s}}\Delta_{\br{u,v}}\partial_1\gamma_{\Gamma,1}\br{\xi,\eta}=w
            }
        \cap W_2.
    \end{equation}

    To proceed, we first study the various relations between \(\cG_\Gamma\) and \(\cG_{\Gamma'}\).
    Recall \textbf{Remark \ref{rmk_inv_Gauss}}. By assumption, we have a topological embedding \(\cG_\Gamma:\Gamma\hookrightarrow\RP\). Since \(\Gamma'\) is an open connected subset of \(\Gamma\) under the topology equipped by \(\Gamma\), we have not only the identity:
    \begin{equation*}
        \cG_{\Gamma'}=\dif{\cG_\Gamma}_{\Gamma'}:\Gamma'\hookrightarrow\RP
    \end{equation*}
    but also \(\img\cG_{\Gamma'}=\cG_\Gamma\br{\Gamma'}\) as an open connected subset of the open connected subset \(\img\cG_\Gamma\subseteq \RP\) under the topology equipped by \(\RP\). As a direct consequence, their continuous inverses are related in the following way:
    \begin{equation}
    \label{eq_inv_Gauss_rel_to_sub_curve}
        \cG^{-1}_{\Gamma'}=\dif{\cG^{-1}_\Gamma}_{\cG_\Gamma\br{\Gamma'}}: \cG_\Gamma\br{\Gamma'}\simto \Gamma'.
    \end{equation}
    This suggests that we express \(\vgamma_{\Gamma'}\) in terms of \(\vgamma_\Gamma\) via restricting the domain.
    To be precise, we start by introducing the following topological embeddings:
    \begin{equation*}
        \iota :\:\R \simto \RP\setminus\BR{\mr{1:0}}\subseteq \RP
        % \setminus\BR{\mr{1:0}}
        :\:\xi \mapsto \mr{\xi:1}.
    \end{equation*}
    Consider the following open set:
    \begin{equation}\label{eq_curve_restrictors}
        U:=\iota^{-1} \img\cG_\Gamma=
        \iota^{-1} \br{\img\cG_\Gamma\setminus \BR{\mr{1:0}} }\subseteq \R.
    \end{equation}
    Via the embeddings \(\iota\) and the open set \(U\), we rephrase the extension of the inverse Gauss maps \(\cG^{-1}_\Gamma\) as:
    \begin{equation}\label{eq_ext_inv_Gauss_embeded}
        \vgamma_\Gamma\br{\xi,\eta}= \cG^{-1}_\Gamma\circ \dif{\iota}_U\br{\xi/\eta}.
    \end{equation}
    On the other hand, the relation \eqref{eq_inv_Gauss_rel_to_sub_curve} suggests that we also consider the following open subsets:
    \begin{equation}\label{eq_sub_curve_restrictors}
         V:=\iota^{-1} \cG_\Gamma\br{\Gamma'}=
        \iota^{-1}
        \br{\cG_\Gamma\br{\Gamma'}\setminus \BR{\mr{1:0}} }\subseteq U\subseteq \R.
    \end{equation}
    For similar reasoning, we deduce the following expression
    \begin{equation}\label{eq_ext_inv_Gauss_embeded_sub_curve}
        \vgamma_{\Gamma'}\br{\xi,\eta}= \cG^{-1}_{\Gamma'}\circ \dif{\iota}_V\br{\xi/\eta}=
        \cG^{-1}_\Gamma\circ \dif{\iota}_V\br{\xi/\eta}.
    \end{equation}
    As a consequence of the two expressions \eqref{eq_ext_inv_Gauss_embeded} and \eqref{eq_ext_inv_Gauss_embeded_sub_curve}, we have:
    \begin{equation*}
        \vgamma_{\Gamma'}\br{\cdot,\eta}=\dif{\vgamma_\Gamma\br{\cdot,\eta}}_{\eta V}.
    \end{equation*}
    Next, we examine the domain of their partial derivatives. By \eqref{eq_curve_restrictors} and \eqref{eq_sub_curve_restrictors}, we observe that the domains:
    \begin{equation*}
        \dom \vgamma_{\Gamma'}\br{\cdot,\eta}=
        \eta V
        \subseteq \eta U= \dom \vgamma_\Gamma\br{\cdot,\eta}
    \end{equation*}
    are nested open subsets of \(\R\). Under our interpretation \eqref{eq_dom_partial_f} in \textbf{Remark \ref{rmk_dom_func}}, we deduce:
    \begin{equation}\label{eq_ext_inv_Gauss_partial_rel_to_sub_curve}
        \dom \partial_1\gamma_{\Gamma',1}\br{\cdot,\eta}=
        \eta V
        \cap
        \dom \partial_1\gamma_{\Gamma,1}\br{\cdot,\eta},\quad
        \partial_1\gamma_{\Gamma',1}\br{\cdot,\eta}=
        \dif{\partial_1\gamma_{\Gamma,1}\br{\cdot,\eta}}_{\eta V}
    \end{equation}
    by examining the formula  \eqref{eq_diff_ability} for differentiability criterion. By \eqref{eq_ext_inv_Gauss_partial_rel_to_sub_curve}, we deduce the two key identities:
    \begin{equation}\label{eq_ext_inv_Gauss_partial_rel_to_sub_curve_1}
        \Delta_{\br{0,s}}\partial_1\gamma_{\Gamma',1}\br{\cdot,\eta}
        =
        \dif{\Delta_{\br{0,s}}\partial_1\gamma_{\Gamma,1}\br{\cdot,\eta}}_{W_1},\quad
        W_1:=W_{\eta,s}:=\eta V\cap \br{\eta+s}V;
    \end{equation}
    \begin{equation}\label{eq_ext_inv_Gauss_partial_rel_to_sub_curve_2}
        \Delta_{\br{0,s}}\Delta_{\br{u,v}}\partial_1\gamma_{\Gamma',1}\br{\cdot,\eta}
        =
        \dif{\Delta_{\br{0,s}}\Delta_{\br{u,v}}\partial_1\gamma_{\Gamma,1}\br{\cdot,\eta}}_{W_2},\quad
        W_2:=W_{\eta,s}\cap \br{-u+W_{\eta+v,s} },
    \end{equation}
    which produces the two identities \eqref{eq_sub_curve_cap_curve_1} and \eqref{eq_sub_curve_cap_curve_2} as an immediate consequence. It remains to show the two sets \(W_1,W_2\subseteq \R\) defined in \eqref{eq_ext_inv_Gauss_partial_rel_to_sub_curve_1} and \eqref{eq_ext_inv_Gauss_partial_rel_to_sub_curve_2} have at most \(O\br{1}\) many connected components.

    Recall \eqref{eq_sub_curve_restrictors}. Since \(\cG_\Gamma\br{\Gamma'}\subseteq \RP\) is connected, the set \(\cG_\Gamma\br{\Gamma'}\setminus\BR{\mr{1:0}}\) has at most two connected components.
    Via the homeomorphism \(\iota:\R\simto \RP\setminus\BR{\mr{1:0}}\), we deduce the inequality:
    \begin{equation*}
        C\br{V}=C\br{\iota^{-1}\br{\cG_\Gamma\br{\Gamma'}}\setminus\BR{\mr{1:0}}}\leq 2.
    \end{equation*}
    Apply now \textbf{Remark \ref{rmk_connect_rel}}. We conclude:
    \begin{equation*}
        C\br{W_1}=C\br{W_{\eta,s}}=C\br{\eta V\cap\br{\eta+s} V}\leq C\br{\eta V}+C\br{\br{\eta+s}V}= 2C\br{V}\leq 4;
    \end{equation*}
    \begin{equation*}
        C\br{W_2}=C\br{W_{\eta,s}\cap\br{-u+W_{\eta+v,s}}}
        \leq C\br{W_{\eta,s}}+C\br{-u+W_{\eta+v,s}}
        =C\br{W_{\eta,s}}+C\br{W_{\eta+v,s}}
        \leq 8.
    \end{equation*}
    This concludes \(C\br{W_1},C\br{W_2}=O\br{1}\) and thus the proof of \textbf{Lemma \ref{lem_M_Gamma_mono}}.
\end{proof}

\subsection{Proof of Theorem \ref{mainthm_global}}

Suggested by \textbf{Remark \ref{rmk_small_img}}, we shall benefit from restricting the analysis to only smaller segments of the curve. With \textbf{Lemma \ref{lem_M_Gamma_mono}}, we can indeed impose the short-interval condition via the following reduction:

\subsubsection{Reduction to short-interval case}\label{subsubsec_I_leng_red}
Observe that an \(A\)-admissible parametric curve \(\vgamma:I\simto \Gamma\) is also \(A'\)-admissible for every \(A'\geq A\). Therefore, it suffices to prove \eqref{eq_mainthm_global} assuming \(A\geq\frac{1}{\dist\br{\supp \rho ,I^c}}\). For convenience, we also assume \(A\gg 1 \).

We claim that it suffices to prove \textbf{Theorem \ref{mainthm_global}} under the assumption that \(\br{\vgamma,\rho}\) satisfies the \(\br{A,C}\)-short-interval condition \ref{ass_short_I} for some universal \(C\gg 1\).

Indeed, assume that \textbf{Theorem \ref{mainthm_global}} holds for all \(\br{\vgamma,\rho}\) satisfying the \(\br{A,C}\)-short-interval condition \ref{ass_short_I}. Recall the function \(\phi\in C^\infty_c\br{\R}\) that satisfies \eqref{eq_phi_uni_even}. Consider the sub-interval of \(I\):
\begin{equation*}
    J:=
    \BR{
        s\in I
    \::\:
        \dist\br{s,I^c}\geq A^{-1}/2
    }
\end{equation*}
and for \(s\in J\) the interval \(I_s:=\br{s-A^{-C},s+A^{-C}}\), the
function:
\begin{equation*}
    \rho_s\br{t}:=
    \rho\br{t}
    A^{3C}\phi\br{A^{3C}\br{t-s}},
\end{equation*}
and the sub-curve \(\vgamma_s:=\dif{\vgamma}_{I_s}: I_s \simto \Gamma_s\) of \(\vgamma:I\simto \Gamma\). Observe the identities:
\begin{equation*}
    \cA_{\vgamma,\rho}\br{f_1,f_2}\br{\vx}
    =
    \frac{1}{\nrm{\phi}_{L^1}}
    \cdot
    \int_J
        \cA_{\vgamma_s,\rho_s}\br{f_1,f_2}\br{\vx}
    ds
\end{equation*}
\begin{equation}\label{eq_avg_to_avg_of_avg}
    =
    \int_{\BR{s'\in J\::\:\rho_{s'}\not\equiv 0}}
        \frac{\nrm{\rho_s}_{C^2}}{\nrm{\phi}_{L^1}}
        \cdot
        \cA_{\Tr_{-s}\vgamma_s,\Tr_{-s}\rho_s/\nrm{\rho_s}_{C^2}}\br{f_1,f_2}\br{\vx}
    ds.
\end{equation}
Notice that \(\br{{\Tr_{-s}\vgamma_s,\Tr_{-s}\rho_s/\nrm{\rho_s}_{C^2}}}\) satisfies the \(\br{A,C}\)-short-interval condition \ref{ass_short_I} by design. We may thus utilize our assumption and apply \textbf{Theorem \ref{mainthm_global}} to deduce:
\begin{equation*}
    \nrm{
    \cA_{\Tr_{-s}\vgamma_s,\Tr_{-s}\rho_s/\nrm{\rho_s}_{C^2}}\br{f_1,f_2}
    }_{L^1}
    \lesssim A^{O\br{C}} \Xi^{\frac{3}{10}}_{\Gamma_s}
    \nrm{f_1}_{H^{\br{-\varepsilon,0}}}
    \nrm{f_2}_{H^{\br{0,-\varepsilon}}}.
\end{equation*}
By \textbf{Lemma \ref{lem_M_Gamma_mono}}, we further dominate the above with:
\begin{equation*}
    \lesssim A^{O\br{C}} \Xi^{\frac{3}{10}}_\Gamma
    \nrm{f_1}_{H^{\br{-\varepsilon,0}}}
    \nrm{f_2}_{H^{\br{0,-\varepsilon}}}.
\end{equation*}
On the other hand, we deduce via direct computation that:
\begin{equation*}
    \frac{\nrm{\rho_s}_{C^2}}{\nrm{\phi}_{L^1}}
    \lesssim A^{O\br{C}}\nrm{\rho}_{C^2}.
\end{equation*}
Combining \eqref{eq_avg_to_avg_of_avg} with the above two relations, we recover \textbf{Theorem \ref{mainthm_global}} for the general case with the universal constant \(c=O\br{C}\). Henceforth, we shall assume that:
\begin{itemize}
    \item \(\br{\vgamma,\rho}\) satisfies the \(\br{A,C}\)-short-interval condition \ref{ass_short_I}.
\end{itemize}
To prove \textbf{Theorem \ref{mainthm_global}} under the above assumption, we start by analyzing the Fourier representation:
\begin{equation}\label{eq_4ier_rep_of_A}
    \cA_{\vgamma,\rho}\br{f_1,f_2}\br{\vx}=
    \int
        \4{1}f_1\br{\xi_1,x_2}
        \4{2}f_2\br{x_1,\xi_2}
        m_{\vgamma,\rho}\br{\vxi}
        e\br{\vx\cdot\vxi}
    d\vxi,
\end{equation}
where the symbol/multiplier is given by the following oscillatory integral:
\begin{equation}\label{eq_m_vgamma_chi}
    m_{\vgamma,\rho}\br{\vxi}:=
    \int
        e\br{\vxi\cdot \vgamma\br{t}}
        \rho\br{t}
    dt.
\end{equation}
We note that \(m_{\vgamma,\rho}\) may exhibit distinct behaviors depending on the size \(\abs{\vxi}\) and the ratio \(\xi_1:\xi_2\). This motivates the following cone decomposition.

\subsubsection{Frequency decomposition}
Let \(\phi\in C^\infty_c\br{\R}\) be the smooth bump function satisfying \eqref{eq_phi_uni_even}. We consider:
\begin{definition}[Littlewood-Payley decomposition]\label{def_LP_decomp}
Let \(\phi\in C^\infty_c\br{\R}\) be the smooth bump function satisfying \eqref{eq_phi_uni_even}. Given \(\delta>0\), let \(\psi_\delta\) be as defined in \eqref{eq_psi_delta} and:
\begin{equation}\label{eq_varphi_delta}
    % \psi_\delta\br{\zeta}:=\phi\br{\frac{\zeta-1}{\delta}},\quad
    \varphi_\delta\br{\zeta}:=\1_{\BR{0}}\br{\zeta}+
    \frac{1}{c_\delta}
    \int_{\br{-1,1}\setminus\BR{0}} \psi_\delta\br{\zeta/\lambda}\frac{d\lambda}{\abs{\lambda}},\quad
    c_\delta:=\nrm{\psi_\delta\br{\zeta}}_{L^1\br{\frac{d\zeta}{\zeta},\R_+}}
\end{equation}
and Littlewood-Payley projections for \(\lambda\neq 0\) and \(f\in L^p\br{\R}\):
\begin{equation}\label{eq_LP_projs}
    \widehat{\cP_{\delta,\lambda} f}\br{\xi}:=
    \psi_\delta\br{\xi/\lambda}
    \widehat{f}\br{\xi},\quad
    \widehat{\cQ_{\delta,\lambda} f}\br{\xi}:=
    \varphi_\delta\br{\xi/\lambda}
    \widehat{f}\br{\xi}.
\end{equation}
By construction, we have the following identities:
\begin{equation}\label{eq_LP_id}
    f=
    \frac{1}{c_\delta}
    \int_{\R\setminus\BR{0}}
        \cP_{\delta,\lambda} f
    \frac{d\lambda}{\abs{\lambda}},\quad
    \cQ_{\delta,\lambda_0}f=
    \frac{1}{c_\delta}
    \int_{0<\abs{\lambda}<\lambda_0}
    \cP_{\delta,\lambda}f
    \frac{d\lambda}{\abs{\lambda}}
    % =
    % \cQ_{\delta,\lambda_0}f+ 
    % \frac{1}{c_\delta}
    % \int_{\lambda_0\leq \abs{\lambda}}
    %     \cP_{\delta,\lambda} f
    % \frac{d\lambda}{\abs{\lambda}}
    ,\quad \lambda_0> 0.
\end{equation}
\end{definition}

We can now pre-compose the \(\cA_{\vgamma,\rho}\) with the Littlewood-Payley projections and obtain\footnote{Recall \textbf{Definition \ref{def_fiber_op}}}:
\begin{equation}\label{eq_LP_QQ}
    \cA_{\vgamma,\rho}\br{f_1,f_2}= 
    \cA_{\vgamma,\rho}
    \br{
        \fQ{1}_{\delta,\lambda_0}f_1,
        \fQ{2}_{\delta,\lambda_0}f_2
    }
\end{equation}
\begin{equation}\label{eq_LP_QP}
    +\frac{1}{c_\delta}\int_{\lambda_0\leq\abs{\lambda}}
        \cA_{\vgamma,\rho}
        \br{
            \fP{1}_{\delta,\lambda}f_1,
            \fQ{2}_{\delta,\lambda_0}f_2
        }
        +
        \cA_{\vgamma,\rho}
        \br{
            \fQ{1}_{\delta,\lambda_0}f_1,
            \fP{2}_{\delta,\lambda}f_2
        }
    \frac{d\lambda}{\abs{\lambda}}
\end{equation}
\begin{equation}\label{eq_LP_PP}
    +\frac{1}{c^2_\delta}\int_{\lambda_0\leq\abs{\lambda_j}}
        \cA_{\vgamma,\rho}
        \br{
            \fP{1}_{\delta,\lambda_1}f_1,
            \fP{2}_{\delta,\lambda_2}f_2
        }
    \frac{d\lambda_1}{\abs{\lambda_1}}
    \frac{d\lambda_2}{\abs{\lambda_2}}.
\end{equation}
For the rest of the section, we will set:
\begin{equation}\label{eq_ass_for_cone_dec}
    \delta=A^{-4C},\quad
    \lambda_0=A^C,\quad
    % \fP{j}_\lambda:=\fP{j}_{\delta,\lambda},\quad
    % \fQ{j}_\lambda:=\fQ{j}_{\delta,\lambda},\quad
    \fR{j}_{\delta,\lambda}\in
    \BR{\fP{j}_{\delta,\lambda},\fQ{j}_{\delta,\lambda}}
\end{equation}
and suppress the \(\delta\)-dependency whenever it is clear from the context. Suggested by the above decomposition of \(\cA_{\vgamma,\rho}\br{f_1,f_2}\), we introduce the following three propositions:

\begin{proposition}[The trivial case]\label{prop_ARR_trivl_est}
    \begin{equation}\label{eq_ARR_trivl_est}
        \nrm{
            \cA_{\vgamma,\rho}
            \br{
                \fR{1}_{\lambda_1}f_1,
                \fR{2}_{\lambda_2}f_2
            }
        }_{L^1}
        \leq
        A^{-3C}
        \nrm{
            f_1
        }_{L^2}
        \nrm{
            f_2
        }_{L^2}.
    \end{equation}
\end{proposition}

\begin{proposition}[The off-diagonal case]\label{prop_ARP_nst_est}
Suppose \(\abs{\lambda_1},\abs{\lambda_2}\geq 1\) and either of the two pairs of conditions:
    \begin{equation}\label{eq_strip_cond}
        \fR{1}_{\lambda_1}=\fQ{1}_{\lambda_1}
            \quad \text{and}\quad
        \abs{\lambda_2}/\lambda_1\geq A^2,\quad \lambda_1>0;
    \end{equation}
    \begin{equation}\label{eq_off_diag_cond}
        \fR{1}_{\lambda_1}=\fP{1}_{\lambda_1}
            \quad \text{and}\quad
        \mr{\lambda_1:\lambda_2}\notin
        \cG_\Gamma\circ\vgamma\br{A^{-C}I}.
    \end{equation}
    The following estimate holds:
    \begin{equation}\label{eq_ARP_nst_est}
        \nrm{
            \cA_{\vgamma,\rho}
            \br{
                \fR{1}_{\lambda_1}f_1,
                \fP{2}_{\lambda_2}f_2
            }
        }_{L^1}
        \lesssim
        A^{O\br{C}}
        \abs{\br{\lambda_1,\lambda_2}}^{-1}
        \nrm{
            f_1
        }_{L^2}
        \nrm{
            f_2
        }_{L^2}.
    \end{equation}
\end{proposition}

\begin{proposition}[The diagonal case]\label{prop_APP_st_est}
There is a universal constant \(\sigma\in\br{0,1}\) such that whenever \(\mr{\lambda_1:\lambda_2}\in \cG_\Gamma\circ\vgamma\br{A^{-C}I}\), the following estimate holds:
    \begin{equation}\label{eq_APP_st_est}
        \nrm{
            \cA_{\vgamma,\rho}
            \br{
                \fP{1}_{\lambda_1}f_1,
                \fP{2}_{\lambda_2}f_2
            }
        }_{L^1}
        \lesssim
        A^{O\br{C}}
        \Xi^{\frac{3}{10}}_\Gamma
        \abs{\br{\lambda_1,\lambda_2}}^{-\sigma}
        \nrm{
            f_1
        }_{L^2}
        \nrm{
            f_2
        }_{L^2}.
    \end{equation}
\end{proposition}

We shall now assume \textbf{Proposition \ref{prop_ARR_trivl_est}, \ref{prop_ARP_nst_est}, and \ref{prop_APP_st_est}} and finish the proof of \textbf{Theorem \ref{mainthm_global}}.
\begin{proof}[Proof of \textbf{Theorem \ref{mainthm_global}}]
    We aim to dominate the \(L^1\) norms of \eqref{eq_LP_QQ}, \eqref{eq_LP_QP}, and \eqref{eq_LP_PP} with the right-hand side of \eqref{eq_mainthm_global} under the short-interval condition:
    \begin{equation*}
        A^{O\br{C}}\Xi^{\frac{3}{10}}_\Gamma
        \nrm{f_1}_{H^{\br{-\varepsilon,0}}}
        \nrm{f_2}_{H^{\br{0,-\varepsilon}}},
    \end{equation*}
    where \(\varepsilon>0\) only depends on the universal constant \(\sigma>0\) given in \textbf{Proposition \ref{prop_APP_st_est}}. To begin with, observe that \eqref{eq_phi_uni_even}, \eqref{eq_psi_delta}, and \eqref{eq_varphi_delta} give the following projection property:
    \begin{equation}\label{eq_LP_PQ_proj}
        \fP{j}_{\delta,\lambda}
        \fR{j}_{2\delta,\lambda}
        =
        \fP{j}_{\delta,\lambda},\quad
        \fQ{j}_{\delta,\lambda}
        \fQ{j}_{2\delta,\lambda}
        =
        \fQ{j}_{\delta,\lambda}.
    \end{equation}
    As a direct consequence, we may insert the appropriate Littlewood-Payley projections to the right-hand side of \eqref{eq_ARR_trivl_est}, \eqref{eq_ARP_nst_est}, and \eqref{eq_APP_st_est}. This produces the trivial estimate:
    \begin{equation}\label{eq_L2_2_H_neq}
        \nrm{\fR{1}_{2\delta,\lambda} f}_{L^2}
        % \leq
        % \ang{\br{1+2\delta}\lambda}^\varepsilon
        \lesssim\abs{\lambda}^\varepsilon
        \nrm{f}_{H^{\br{-\varepsilon,0}}},\quad
        \nrm{\fR{2}_{2\delta,\lambda} f}_{L^2}
        % \leq
        % \ang{\br{1+2\delta}\lambda}^\varepsilon
        \lesssim\abs{\lambda}^\varepsilon
        \nrm{f}_{H^{\br{0,-\varepsilon}}}
        ,\quad
        0\leq \varepsilon\leq 1,\;
        1\leq \abs{\lambda}.
    \end{equation}
    
    Starting with the treatment of \eqref{eq_LP_QQ}, we apply \eqref{eq_ARR_trivl_est} and \eqref{eq_LP_PQ_proj} to deduce:
    \begin{equation*}
    \nrm{
        \cA_{\vgamma,\rho}\br{
            \fQ{1}_{\delta,\lambda_0}f_1, 
            \fQ{2}_{\delta,\lambda_0}f_2
        }
    }_{L^1}
    \leq
    A^{-3C}
    \nrm{\fQ{1}_{2\delta,\lambda_0}f_1}_{L^2}
    \nrm{\fQ{2}_{2\delta,\lambda_0}f_2}_{L^2}.
    \end{equation*}
    Applying \eqref{eq_L2_2_H_neq}, we further dominate the above with:
    \begin{equation*}
        \lesssim
        A^{-3C}
        \lambda^{2\epsilon}_0
        \nrm{f_1}_{H^{\br{-\varepsilon,0}}}
        \nrm{f_2}_{H^{\br{0,-\varepsilon}}}
        =A^{O\br{C}}
        \nrm{f_1}_{H^{\br{-\varepsilon,0}}}
        \nrm{f_2}_{H^{\br{0,-\varepsilon}}}
    \end{equation*}
    for any \(\varepsilon\leq 1\). This completes the treatment of \eqref{eq_LP_QQ}.

    To treat \eqref{eq_LP_QP}, we start by estimating the expression within the original integral:
    \begin{equation}\label{eq_LP_QP_PQ_inner}
        \cA_{\vgamma,\rho}
        \br{
            \fP{1}_{\delta,\lambda}f_1,
            \fQ{2}_{\delta,\lambda_0}f_2
        }
        +
        \cA_{\vgamma,\rho}
        \br{
            \fQ{1}_{\delta,\lambda_0}f_1,
            \fP{2}_{\delta,\lambda}f_2
        }.
    \end{equation}
    We have two cases. When \(\abs{\lambda}\leq A^2\lambda_0\), we apply \eqref{eq_ARR_trivl_est}, \eqref{eq_LP_PQ_proj}, and \eqref{eq_L2_2_H_neq} to deduce
    \begin{equation}\label{eq_LP_QP_PQ_inner_L_bd}
        \nrm{\eqref{eq_LP_QP_PQ_inner}}_{L^1}\leq A^{O\br{C}}
        \nrm{f_1}_{H^{\br{-\varepsilon,0}}}
        \nrm{f_2}_{H^{\br{0,-\varepsilon}}}.
    \end{equation}
    For the alternative \(\abs{\lambda}\geq A^2\lambda_0\) case, we invoke symmetry of the indices, apply the \eqref{eq_strip_cond} case in \textbf{Proposition \ref{prop_ARP_nst_est}}, and follow it with \eqref{eq_LP_PQ_proj} to obtain:
    \begin{equation*}
        \nrm{\eqref{eq_LP_QP_PQ_inner}}_{L^1}\lesssim
        \abs{\br{\lambda_0,\lambda}}^{-1}A^{O\br{C}}
        \br{
            \nrm{\fP{1}_{2\delta,\lambda} f_1}_{L^2}
            \nrm{\fQ{2}_{2\delta,\lambda_0} f_2}_{L^2}
            +
            \nrm{\fQ{1}_{2\delta,\lambda_0} f_1}_{L^2}
            \nrm{\fP{2}_{2\delta,\lambda} f_2}_{L^2}
        }.
    \end{equation*}
    Utilizing again \eqref{eq_L2_2_H_neq} and the fact that \(\abs{\br{\lambda_0,\lambda}}\geq\abs{\lambda}\), we deduce for \(\varepsilon\leq 1\) the following:
    \begin{equation}\label{eq_LP_QP_PQ_inner_H_bd}
        \nrm{\eqref{eq_LP_QP_PQ_inner}}_{L^1}
        \lesssim
        \abs{\lambda}^{\varepsilon-1}
        A^{O\br{C}}
        \nrm{f_1}_{H^{\br{-\varepsilon,0}}}
        \nrm{f_2}_{H^{\br{0,-\varepsilon}}}.
    \end{equation}
    Suggested by the two estimates \eqref{eq_LP_QP_PQ_inner_L_bd} and \eqref{eq_LP_QP_PQ_inner_H_bd}, we shall split \eqref{eq_LP_QP} into low and high frequency parts:
    \begin{equation*}
        \nrm{\eqref{eq_LP_QP}}_{L^1}
        \leq
        \frac{1}{c_\delta}
        \int_{\lambda_0\leq\abs{\lambda}\leq A^2\lambda_0}
        \nrm{\eqref{eq_LP_QP_PQ_inner}}_{L^1}
        \frac{d\lambda}{\abs{\lambda}}
        +
        \frac{1}{c_\delta}
        \int_{A^2\lambda_0\leq\abs{\lambda}}
        \nrm{\eqref{eq_LP_QP_PQ_inner}}_{L^1}
        \frac{d\lambda}{\abs{\lambda}}
    \end{equation*}
    and dominate the above expression accordingly:
    \begin{equation*}
        \lesssim
        \frac{A^{O\br{C}}}{c_\delta}
        \int^{A^C}_1
        \frac{d\zeta}{\zeta}
        \nrm{f_1}_{H^{\br{-\varepsilon,0}}}
        \nrm{f_2}_{H^{\br{0,-\varepsilon}}}
        +
        \frac{A^{O\br{C}}}{c_\delta}
        \int^\infty_{A^2\lambda_0}
            \zeta^{\varepsilon-1}
        \frac{d\zeta}{\zeta}
        \nrm{f_1}_{H^{\br{-\varepsilon,0}}}
        \nrm{f_2}_{H^{\br{0,-\varepsilon}}}.
    \end{equation*}
    Via direct computation, we dominate the above with
    \begin{equation*}
        \leq \frac{A^{O\br{C}}}{c_\delta\br{1-\varepsilon}}
        \nrm{f_1}_{H^{\br{-\varepsilon,0}}}
        \nrm{f_2}_{H^{\br{0,-\varepsilon}}}
        ,\quad \varepsilon\in\br{0,1}.
    \end{equation*}
    To address the constant \(c_\delta\), we deduce from \eqref{eq_phi_uni_even} the following:
    \begin{equation}\label{eq_c_delta_est}
        A^{-O\br{C}}\leq
        \log\br{\frac{2+\delta}{2-\delta}}=
        \int^{1+\frac{\delta}{2}}_{1-\frac{\delta}{2}}
        \frac{d\zeta}{\zeta}
        \leq c_\delta.
    \end{equation}
    We may thus absorb \(\frac{1}{c_\delta}\) into a \(A^{O\br{C}}\) factor. In combination, we complete the treatment of \eqref{eq_LP_QP}:
    \begin{equation*}
        \nrm{\eqref{eq_LP_QP}}_{L^1}
        \underset{\varepsilon}{\lesssim}
        A^{O\br{C}}
        \nrm{f_1}_{H^{\br{-\varepsilon,0}}}
        \nrm{f_2}_{H^{\br{0,-\varepsilon}}}
    \end{equation*}
    for any choice of \(\varepsilon\in\br{0,1}\).

    To treat \eqref{eq_LP_PP}, we first combine the \eqref{eq_off_diag_cond} case in \textbf{Proposition \ref{prop_ARP_nst_est}} and \textbf{Proposition \ref{prop_APP_st_est}} to deduce:
    \begin{equation}\label{eq_APP_joint_est}
        \nrm{
            \cA_{\vgamma,\rho}
            \br{
                \fP{1}_{\lambda_1}f_1,
                \fP{2}_{\lambda_2}f_2
            }
        }_{L^1}
        \lesssim
        A^{O\br{C}}
        \Xi^{\frac{3}{10}}_\Gamma
        \abs{\br{\lambda_1,\lambda_2}}^{-\sigma}
        \nrm{
            f_1
        }_{L^2}
        \nrm{
            f_2
        }_{L^2}
    \end{equation}
    for \(\abs{\lambda_1},\abs{\lambda_2}\geq\lambda_0\).
    With \eqref{eq_APP_joint_est}, we apply \eqref{eq_LP_PQ_proj} and \eqref{eq_L2_2_H_neq} to deduce:
    \begin{equation*}
        \nrm{\eqref{eq_LP_PP}}_{L^1}
        \lesssim
        \frac{A^{O\br{C}}\Xi^{\frac{3}{10}}_\Gamma}{c^2_\delta}
        \int_{\lambda_0\leq\abs{\lambda_j}}
            \abs{\br{\lambda_1,\lambda_2}}^{-\sigma}
            \abs{\lambda_1}^\varepsilon
            \abs{\lambda_2}^\varepsilon
        \frac{d\lambda_1}{\abs{\lambda_1}}
        \frac{d\lambda_2}{\abs{\lambda_2}}
        \nrm{f_1}_{H^{\br{-\varepsilon,0}}}
        \nrm{f_2}_{H^{\br{0,-\varepsilon}}}.
    \end{equation*}
    Using \eqref{eq_c_delta_est} and the fact that \(\abs{\br{\lambda_1,\lambda_2}}\geq \abs{\lambda_j}\), we dominate the above with
    \begin{equation*}
        \leq
        A^{O\br{C}}\Xi^{\frac{3}{10}}_\Gamma
        \br{
        \int^\infty_{\lambda_0}
            \zeta^{\varepsilon-\frac{\sigma}{2}}
        \frac{d\zeta}{\zeta}
        }^2
        \nrm{f_1}_{H^{\br{-\varepsilon,0}}}
        \nrm{f_2}_{H^{\br{0,-\varepsilon}}}.
    \end{equation*}
    Finally, direct computation completes the treatment of \eqref{eq_LP_PP}:
    \begin{equation*}
        \nrm{\eqref{eq_LP_PP}}_{L^1}
        \underset{\varepsilon,\sigma}{\lesssim}
        A^{O\br{C}}\Xi^{\frac{3}{10}}_\Gamma
        \nrm{f_1}_{H^{\br{-\varepsilon,0}}}
        \nrm{f_2}_{H^{\br{0,-\varepsilon}}}
    \end{equation*}
    for any choice of \(\varepsilon\in\br{0,\frac{\sigma}{2}}\).
    In combination, we complete the proof of \textbf{Theorem \ref{mainthm_global}}:
    \begin{equation*}
        \nrm{\cA_{\vgamma,\rho}\br{f_1,f_2}}_{L^1}
        \leq
        \nrm{\eqref{eq_LP_QQ}}_{L^1}
        +\nrm{\eqref{eq_LP_QP}}_{L^1}
        +\nrm{\eqref{eq_LP_PP}}_{L^1}
        \underset{\varepsilon,\sigma}{\lesssim}
        A^{O\br{C}}\Xi^{\frac{3}{10}}_\Gamma
        \nrm{f_1}_{H^{\br{-\varepsilon,0}}}
        \nrm{f_2}_{H^{\br{0,-\varepsilon}}}.
    \end{equation*}
\end{proof}
% As is discussed in \textsc{Section \ref{subsubsec_I_leng_red}},
% the verification of \textbf{Theorem \ref{mainthm_global}} gives \textbf{Theorem \ref{mainthm_global}} as a consequence.
It remains to treat \textbf{Proposition \ref{prop_ARR_trivl_est}, \ref{prop_ARP_nst_est}, and \ref{prop_APP_st_est}}.
Observe the following identity:
\begin{equation*}
    \4{j}\fR{j}_\lambda f
    \br{\xi\ve_j+x\ve_k}=
    \vartheta_j\br{\xi/\lambda}\4{j}f
    \br{\xi\ve_j+x\ve_k},\quad
    \vartheta_j\in\BR{\psi_\delta,\varphi_\delta},\quad
    \BR{j,k}=\BR{1,2}.
\end{equation*}
for the suitable choice of \(\vartheta_j\). Combining with \eqref{eq_4ier_rep_of_A}, we write:
\begin{equation}\label{eq_4ier_rep_of_A_w_thetas}
    \cA_{\vgamma,\rho}
    \br{
        \fR{1}_{\lambda_1}f_1,
        \fR{2}_{\lambda_2}f_2
    }
    \br{\vx}
    =
    \int
        \4{1}f_1\br{\xi_1,x_2}
        \4{2}f_2\br{x_1,\xi_2}
        \br{
            \prod_{j=1,2}
                \vartheta_j\br{\xi_j/\lambda_j}}
        m_{\vgamma,\rho}\br{\vxi}
        e\br{\vx\cdot\vxi}
    d\vxi.
\end{equation}
The distinction among the formulation of \textbf{Proposition \ref{prop_ARR_trivl_est}, \ref{prop_ARP_nst_est}, and \ref{prop_APP_st_est}} is a consequence of the different behaviors of the localized symbol \eqref{eq_m_vgamma_chi}:
\begin{equation}\label{eq_m_theta_freq_trunc}
    m_{\vlambda,\vartheta_1,\vartheta_2}\br{\vxi}:=m_{\vgamma,\rho,\vlambda,\vartheta_1,\vartheta_2}\br{\vxi}:=m_{\vgamma,\rho}\br{\vxi}
    \prod_{j=1,2}\vartheta_j\br{\xi_j/\lambda_j}
\end{equation}
as \(\vlambda:=\br{\lambda_1,\lambda_2}\in\br{ \R\setminus\BR{0}}^2\) varies. We identify the following three key cases:
\begin{itemize}
    \item the non-oscillatory regime \(\abs{\lambda_j}\leq A^{O\br{C}}\): this corresponds to \textbf{Proposition \ref{prop_ARR_trivl_est}}.
    \item the oscillatory regime \(\abs{\lambda_j}\geq A^{O\br{C}}\) consists of the two scenarios:
    \begin{itemize}
        \item the non-stationary phase case \eqref{eq_strip_cond} and \eqref{eq_off_diag_cond}: this corresponds to \textbf{Proposition \ref{prop_ARP_nst_est}}.
        \item the stationary phase case \(\mr{\lambda_1:\lambda_2}\in\cG_\Gamma\circ\vgamma\br{A^{-C}I}\): this corresponds to \textbf{Proposition \ref{prop_APP_st_est}}.
    \end{itemize}
\end{itemize}

\subsection{Proof of Proposition \ref{prop_ARR_trivl_est}: the non-oscillatory regime}
When \eqref{eq_m_theta_freq_trunc} exhibits no oscillation, we discard the frequency information completely and perform trivial estimate:
\begin{equation*}
    \nrm{
        \cA_{\vgamma,\rho}\br{
            \fR{1}_{\lambda_1}f_1,
            \fR{2}_{\lambda_2}f_2
        }
    }_{L^1}
    \leq 
    \nrm{
        \nrm{
            \fR{1}_{\lambda_1}f_1\br{x_1+\gamma_1\br{t},x_2}
            \fR{2}_{\lambda_2}f_2\br{x_1,x_2+\gamma_2\br{t}}
        }_{L^1\br{d\vx}}
    \rho\br{t}
    }_{L^1\br{dt}}
\end{equation*}
\begin{equation*}
    \leq
    \nrm{\rho}_{L^1}
    \prod_{j=1,2}
        \nrm{
            \fR{j}_{\lambda_j}f_j
        }_{L^2}
    \leq
    A^{-3C}
    \nrm{f_1}_{L^2}
    \nrm{f_2}_{L^2}.
\end{equation*}
This completes the proof of \textbf{Proposition \ref{prop_ARR_trivl_est}}.

\subsection{Proof of Proposition \ref{prop_ARP_nst_est}: the non-stationary phase case}

Recall now \eqref{eq_4ier_rep_of_A_w_thetas}. Via Fourier inversion, we deduce:
\begin{equation*}
    \cA_{\vgamma,\rho}
    \br{
        \fR{1}_{\lambda_1}f_1,
        \fP{2}_{\lambda_2}f_2
    }\br{\vx}=
    \int
        f_1\br{x_1+z_1,x_2}
        f_2\br{x_1,x_2+z_2}
        \widehat{
            m_{\vlambda,\vartheta_1,\psi_\delta}
        }
        \br{\vz}
    d\vz.
\end{equation*}
Due to the trivial estimates:
\begin{equation*}
    \nrm{\cA_{\vgamma,\rho}
    \br{
        \fR{1}_{\lambda_1}f_1,
        \fP{2}_{\lambda_2}f_2
    }
    }_{L^1}
    \leq
    \nrm{
        \widehat{
            m_{\vlambda,\vartheta_1,\psi_\delta}
        }
    }_{L^1}
    \nrm{f_1}_{L^2}
    \nrm{f_2}_{L^2},
\end{equation*}
proving \textbf{Proposition \ref{prop_ARP_nst_est}} reduces to show the following estimate:
\begin{equation}\label{eq_symb_nst_est}
    \nrm{
        \widehat{
            m_{\vlambda,\vartheta_1,\psi_\delta}
        }
    }_{L^1}
    \lesssim A^{O\br{C}}\abs{\vlambda}^{-1}.
\end{equation}
In fact, due to the trivial fact
\begin{equation*}
    m_{\vgamma,\rho}\br{\vxi}=e\br{\vxi\cdot\vgamma\br{0}} m_{\vgamma-\vgamma\br{0},\rho}\br{\vxi},
\end{equation*}
we may translate the \(L^1\)-norm expression in \eqref{eq_symb_nst_est} to impose the normalization \(\vgamma\br{0}=0\).

To proceed, we
assume either \(\vartheta_1=\varphi_\delta\) with \eqref{eq_strip_cond} or \(\vartheta_1=\psi_\delta\) with \eqref{eq_off_diag_cond} and show:
\begin{equation}\label{eq_nst_symb}
    \vartheta_1\br{\xi_1/\lambda_1}
    \psi_\delta\br{\xi_2/\lambda_2}\neq 0
    \implies 
    \forall t\in \supp\rho,\;
    \abs{\vxi\cdot \vgamma'\br{t}}\geq A^{-O\br{C}}\abs{\vlambda}.
\end{equation}
Indeed, by the construction \eqref{eq_psi_delta} and \eqref{eq_varphi_delta}, we have for \(\lambda'>0\) the following:
\begin{equation}\label{eq_varphi_2_psi_supp}
    \varphi_\delta\br{\xi/\lambda'}\neq 0\implies \exists \lambda\in\br{-\lambda',\lambda'}\setminus\BR{0},\;
    \psi_\delta\br{\xi/\lambda}\neq 0
\end{equation}
and for \(\lambda\neq 0\) the following:
\begin{equation}\label{eq_psi_supp_2_lambda}
    \psi_\delta\br{\xi/\lambda}\implies
    \abs{\xi-\lambda}\leq \delta \abs{\lambda}.
\end{equation}
Consider now the case with \(\vartheta_1=\varphi_\delta\) and \eqref{eq_strip_cond}, we combine \eqref{eq_varphi_2_psi_supp} and \eqref{eq_psi_supp_2_lambda} to deduce:
\begin{equation*}
    \varphi_\delta\br{\xi_1/\lambda_1}\psi_\delta\br{\xi_2/\lambda_2}\neq 0\implies 
    \exists \lambda'_1\in\br{-\lambda_1,\lambda_1}\setminus\BR{0},\;
    \abs{\xi_1-\lambda'_1}\leq \delta\abs{\lambda'_1},\quad
    \abs{\xi_2-\lambda_2}\leq \delta\abs{\lambda_2}.
\end{equation*}
By assumption \eqref{eq_strip_cond}, we have:
\begin{equation*}
    \abs{\lambda_2/\lambda'_1}\geq \abs{\lambda_2}/\lambda_1> A^2.
\end{equation*}
By contrapositive of \eqref{eq_tight_diag}, we deduce from the above the following:
\begin{equation*}
    \implies 
    \mr{\lambda'_1:\lambda_2}\notin \img\cG_\Gamma\implies
    \mr{\lambda'_1:\lambda_2}\notin
    \cG_\Gamma\circ\vgamma\br{A^{-C}I}.
\end{equation*}
Applying now \eqref{eq_osc_nst} and the trivial fact \(2\abs{\br{\lambda'_1,\lambda_2}}\geq\abs{\vlambda}\), we deduce from above the conclusion of \eqref{eq_nst_symb}:
\begin{equation*}
    \implies \forall t\in \supp \rho,\;
    \abs{\vxi\cdot\vgamma'\br{t}}\geq A^{-O\br{C}}\abs{\br{\lambda'_1,\lambda_2}}\geq
    A^{-O\br{C}}\abs{\vlambda}.
\end{equation*}
The case with \(\vartheta_1=\psi_\delta\) and \eqref{eq_off_diag_cond} is a direct consequence of \eqref{eq_osc_nst} and \eqref{eq_psi_supp_2_lambda}.
Given now \eqref{eq_nst_symb}, it suffices to study the oscillatory expression \(m_{\vgamma,\rho}\br{\vxi}\) assuming
\begin{equation*}
    \forall t\in\supp\rho,\;\abs{\vxi\cdot\vgamma'\br{t}}\geq A^{-O\br{C}}\abs{\vlambda}.
\end{equation*}

The assumption allows us to perform the standard integration by-parts argument:
\begin{equation*}
    \partial^{\alpha_1}_1
    \partial^{\alpha_2}_2
    m_{\vgamma,\rho}\br{\vxi}
    =
    \int
        e\br{\vxi\cdot\vgamma\br{t}}
        \br{
            \frac{d}{dt}\cdot
            \frac{1}{-2\pi i \vxi\cdot \vgamma'\br{t}}
            \cdot
        }^2
        \br{
            \rho\br{t}
            \prod_{j=1,2}
            \br{2\pi i\gamma_j\br{t}}^{\alpha_j}
        }
    dt
\end{equation*}
and deduce, via \(A\)-admissibility of \(\vgamma:I\simto\Gamma\) and the normalization \(\vgamma\br{0}=0\), the following estimate:
\begin{equation*}
    \abs{
        \partial^{\alpha_1}_1
        \partial^{\alpha_2}_2
        m_{\vgamma,\rho}\br{\vxi}
    }
    \leq A^{O\br{C}}\abs{\vlambda}^{-2}
    ,\quad
    0 \leq \alpha_1,\alpha_2,\alpha_1+\alpha_2\leq 2.
\end{equation*}
Invoking the assumption \(\abs{\lambda_j}\geq 1\), we obtain for the truncated symbol \eqref{eq_m_theta_freq_trunc}
the following:
\begin{equation*}
    \nrm{
        m_{\vlambda,\vartheta_1,\psi_\delta}
    }_{H^2}
    =
    \nrm{
        \br{1-\frac{\Delta}{4\pi^2}}
        m_{\vlambda,\vartheta_1,\psi_\delta}
    }_{L^2}
    \lesssim
    A^{O\br{C}}\abs{\vlambda}^{-2}\abs{\lambda_1\lambda_2}^{\frac{1}{2}}\abs{\supp\vartheta_1}^{\frac{1}{2}}\abs{\supp\psi_\delta}^{\frac{1}{2}}
    \leq A^{O\br{C}}\abs{\vlambda}^{-1}.
\end{equation*}
Using now the Fourier characterization of \(H^2\) norm, we conclude \eqref{eq_symb_nst_est}:
\begin{equation*}
    \nrm{
        \widehat{
            m_{\vlambda,\vartheta_1,\psi_\delta}
        }
    }_{L^1}
    \leq
    \nrm{
        \br{1+\abs{\vz}^2}
        \widehat{
            m_{\vlambda,\vartheta_1,\psi_\delta}
        }\br{\vz}
    }_{L^2\br{d\vz}}
    \nrm{
        \frac{1}{1+\abs{\vz}^2}
    }_{L^2\br{d\vz,\R^2}}
    \lesssim
    \nrm{
        m_{\vlambda,\vartheta_1,\psi_\delta}
    }_{H^2}
    \lesssim
    A^{O\br{C}}\abs{\vlambda}^{-1}.
\end{equation*}
This completes the proof of \textbf{Proposition \ref{prop_ARP_nst_est}}.\qed

\subsection{Proof of Proposition \ref{prop_APP_st_est}: the stationary phase case}

Let \(\vlambda:=\br{\lambda_1,\lambda_2}\) be such that \(\mr{\vlambda}\in\cG_\Gamma\circ\vgamma\br{A^{-C}I}\).
To beat the trivial bound provided by \textbf{Proposition \ref{prop_ARR_trivl_est}}, we may assume \(\abs{\vlambda}\geq A^C\) throughout the rest of the discussion. We will first demonstrate a physically localized version of \textbf{Proposition \ref{prop_APP_st_est}}. We then utilize a standard trick to recover the general statement.

Recall \eqref{carlyT} and \eqref{eq_m_theta_freq_trunc}. We begin with the following observation:
\begin{equation}\label{eq_avg_to_carlyT}
    \1_{\mr{0,1}^2}
    \cA_{\vgamma,\rho}\br{\fP{1}_{\lambda_1}f_1,\fP{2}_{\lambda_2}f_2}
    =\cT_{m_{\vlambda,\psi_\delta,\psi_\delta}}\br{\4{1}f_1,\4{2}f_2}
    =\cT_{m_{\vgamma,\rho,\vlambda}}\br{\4{1}f_1,\4{2}f_2},
\end{equation}
where \(m_{\vgamma,\rho,\vlambda}\) are defined as in \eqref{eq_m_vlambda_freq_trunc}. Under the short-interval condition, we may apply \textbf{Theorem \ref{thm_4ier_op_est}} and deduce from \eqref{eq_avg_to_carlyT} the following:
\begin{equation}\label{eq_APP_st_loc}
    \nrm{
        \cA_{\vgamma,\rho}
        \br{
            \fP{1}_{\lambda_1}f_1,
            \fP{2}_{\lambda_2}f_2
        }
    }_{L^1\br{\mr{0,1}^2}}
    \lesssim 
    A^{O\br{1}}\Xi^{\frac{3}{10}}_\Gamma 
    \abs{\vlambda}^{-\sigma}
    \nrm{f_1}_{L^2}\nrm{f_2}_{L^2}.
\end{equation}
This establishes the physically localized version of \textbf{Proposition \ref{prop_APP_st_est}}. To prove the statement in full generality, we first apply the same translation trick to impose the normalization \(\vgamma\br{0}=0\) and then utilize the small support of \(\rho\) to propagate the physical localization to the two input functions:\footnote{For \(\varphi_2\), we mean \eqref{eq_varphi_delta} with \(\delta=2\).}
\begin{equation}\label{eq_APP_st_loc_ini_id}
    \1_{\mr{0,1}^2}
    \cA_{\vgamma,\rho}
        \br{
            \fP{1}_{\lambda_1}f_1,
            \fP{2}_{\lambda_2}f_2
        }
    =
    \cA_{\vgamma,\rho}
        \br{
            \br{
                \varphi_2\otimes
                \1_{\mr{0,1}}
            }
            \fP{1}_{\lambda_1}
            f_1,
            \br{
                \1_{\mr{0,1}}\otimes
                \varphi_2
            }
            \fP{2}_{\lambda_2}
            f_2
        }
    .
\end{equation}
The following standard trick further passes the outer physical localization through the frequency projection with a well-controlled error term:

\begin{lemma}\label{lem_phy_loc_pass_proj}
Assuming \(\abs{\lambda_1},\abs{\lambda_2}\gg 1\), we have the following point-wise approximation:
    \begin{equation}\label{eq_APP_st_loc_main}
    \1_{\mr{0,1}^2}\br{\vx}
    \cA_{\vgamma,\rho}
        \br{
            \fP{1}_{\lambda_1}f_1,
            \fP{2}_{\lambda_2}f_2
        }\br{\vx}
    =
    \cA_{\vgamma,\rho}
        \br{
            \fP{1}_{\lambda_1}
            \br{
                \br{
                    \varphi_2\otimes
                    \1_{\mr{0,1}}
                }
            f_1},
            \fP{2}_{\lambda_2}
            \br{
                \br{
                    \1_{\mr{0,1}}\otimes
                    \varphi_2
                }
            f_2}
        }\br{\vx}
\end{equation}
\begin{equation}\label{eq_APP_st_loc_err}
    +
    \frac{A^{O\br{C}}}{\min\br{\abs{\lambda_1},\abs{\lambda_2}}}
    O\br{
        \cA_{\vgamma,\abs{\rho}}
        \br{
            \br{
                \varphi_2\otimes
                \1_{\mr{0,1}}
            }
            \fM{1}f_1,
            \br{
                \1_{\mr{0,1}}\otimes
                \varphi_2
            }
            \fM{2}f_2
        }
        \br{\vx}
    },
\end{equation}
where the operator \(\cM\) stands for the standard Hardy-Littlewood maximal operator.
% \begin{equation*}
%     \fM{j}f\br{\vx}
%     :=
%     \sup_{r>0}
%     \fint^r_{-r}
%         \abs{f\br{\vx+z\ve_j}}
%     dz.
% \end{equation*}
\end{lemma}

Applying the estimate \eqref{eq_APP_st_loc} to the first term \eqref{eq_APP_st_loc_main} and trivial bounds to the second term \eqref{eq_APP_st_loc_err} produces:
\begin{equation*}
    \nrm{
        \cA_{\vgamma,\rho}
        \br{
            \fP{1}_{\lambda_1}f_1,
            \fP{2}_{\lambda_2}f_2
        }
    }_{L^1\br{\mr{0,1}^2}}
    \lesssim
    A^{O\br{1}}\Xi^{\frac{3}{10}}_\Gamma\abs{\vlambda}^{-\sigma}
    \nrm{
        \br{
            \varphi_2\otimes
            \1_{\mr{0,1}}
        }
        f_1
    }_{L^2}
    \nrm{
        \br{
            \1_{\mr{0,1}}\otimes
            \varphi_2
        }
        f_2
    }_{L^2}
\end{equation*}
\begin{equation*}
    +\frac{A^{O\br{C}}}{\min\br{\abs{\lambda_1},\abs{\lambda_2}}}
    \nrm{
        \br{
            \varphi_2\otimes
            \1_{\mr{0,1}}
        }
        \fM{1}f_1
    }_{L^2}
    \nrm{
        \br{
            \1_{\mr{0,1}}\otimes
            \varphi_2
        }
        \fM{2}f_2
    }_{L^2}.
\end{equation*}
Using the property \eqref{eq_tight_diag} and the fact that \(0\leq\varphi_2,\1_{\mr{0,1}}\leq \1_{\mr{-2,2}}\), we further deduce:
\begin{equation*}
    \nrm{
        \cA_{\vgamma,\rho}
        \br{
            \fP{1}_{\lambda_1}f_1,
            \fP{2}_{\lambda_2}f_2
        }
    }_{L^1\br{\mr{0,1}^2}}
    \lesssim
    A^{O\br{C}}\Xi^{\frac{3}{10}}_\Gamma\abs{\vlambda}^{-\sigma}
    \prod_{j=1,2}
    \br{
        \nrm{
            f_j
        }_{L^2\br{\mr{-4,4}^2}}
        +
        \nrm{
            \fM{j}f_j
        }_{L^2\br{\mr{-4,4}^2}}
    }.
\end{equation*}
Invoking now the translation symmetry of the formulation, we obtain for \(\vz\in\Z^2\) the following:
\begin{equation*}
    \nrm{
        \cA_{\vgamma,\rho}
        \br{
            \fP{1}_{\lambda_1}f_1,
            \fP{2}_{\lambda_2}f_2
        }
    }_{L^1\br{\vz+\mr{0,1}^2}}
    \lesssim
    A^{O\br{C}}\Xi^{\frac{3}{10}}_\Gamma\abs{\vlambda}^{-\sigma}
    \prod_{j=1,2}
    \br{
        \nrm{
            f_j
        }_{L^2\br{\vz+\mr{-4,4}^2}}
        +
        \nrm{
            \fM{j}f_j
        }_{L^2\br{\vz+\mr{-4,4}^2}}
    }.
\end{equation*}
Finally, by summing over \(\vz\in\Z^2\) and an application of Cauchy-Schwarz inequality, we recover \eqref{eq_APP_st_est}:
\begin{equation*}
    \nrm{
        \cA_{\vgamma,\rho}
        \br{
            \fP{1}_{\lambda_1}f_1,
            \fP{2}_{\lambda_2}f_2
        }
    }_{L^1}
    =\sum_{\vz\in\Z^2}
    \nrm{
        \cA_{\vgamma,\rho}
        \br{
            \fP{1}_{\lambda_1}f_1,
            \fP{2}_{\lambda_2}f_2
        }
    }_{L^1\br{\vz+\mr{0,1}^2}}
\end{equation*}
\begin{equation*}
    \lesssim
    A^{O\br{C}}\Xi^{\frac{3}{10}}_\Gamma\abs{\vlambda}^{-\sigma}
    \sum_{\vz\in\Z^2}
    \prod_{j=1,2}
    \br{
        \nrm{
            f_j
        }_{L^2\br{\vz+\mr{-4,4}^2}}
        +
        \nrm{
            \fM{j}f_j
        }_{L^2\br{\vz+\mr{-4,4}^2}}
    }
\end{equation*}
\begin{equation*}
    \lesssim
    A^{O\br{C}}\Xi^{\frac{3}{10}}_\Gamma\abs{\vlambda}^{-\sigma}
    \prod_{j=1,2}
    \br{
        \nrm{
            \nrm{
                f_j
            }_{L^2\br{\vz+\mr{-4,4}^2}}
        }_{\ell^2\br{\vz\in\Z^2}}
        +
        \nrm{
            \nrm{
                \fM{j}f_j
            }_{L^2\br{\vz+\mr{-4,4}^2}}
        }_{\ell^2\br{\vz\in\Z^2}}
    }
\end{equation*}
\begin{equation*}
    \lesssim
    A^{O\br{C}}\Xi^{\frac{3}{10}}_\Gamma\abs{\vlambda}^{-\sigma}
    \prod_{j=1,2}
    \br{
        \nrm{
            f_j
        }_{L^2}
        +
        \nrm{
            \fM{j}f_j
        }_{L^2}
    }
    \lesssim
    A^{O\br{C}}\Xi^{\frac{3}{10}}_\Gamma\abs{\vlambda}^{-\sigma}
    \nrm{
        f_1
    }_{L^2}
    \nrm{
        f_2
    }_{L^2}.
\end{equation*}
This completes the proof of \textbf{Proposition \ref{prop_APP_st_est}} assuming \textbf{Lemma \ref{lem_phy_loc_pass_proj}}, which we address below.
\begin{proof}[Proof of \textbf{Lemma \ref{lem_phy_loc_pass_proj}}]
    Suggested by the fiber-wise structure of the operator, we shall begin with the analysis of commutator between the multiplication of \(\varphi_2\) denoted as \(M_{\varphi_2}\) and the Littlewood-Payley projection \(\widehat{\cP_\lambda f}\br{\xi}:=\psi_\delta\br{\xi/\lambda}\widehat{f}\br{\xi}\) for functions \(f\in L^2\br{\R}\) and parameter \(\abs{\lambda}\geq 1\).
    % the following commutator:
    % \begin{equation*}
    %     \mr{\varphi_2,\cP_\lambda}f=
    %     \varphi_2\cP_\lambda f- 
    %     \cP_\lambda\br{\varphi_2 f},\quad
    %     \widehat{\cP_\lambda f}\br{\xi}:=\psi_\delta\br{\xi/\lambda}\widehat{f}\br{\xi}
    % \end{equation*}
    % for functions \(f\in L^2\br{\R}\) and parameter \(\abs{\lambda}\geq 1\).
    By direct computation, we have:
    \begin{equation*}
        \mr{M_{\varphi_2}, \cP_\lambda}f\br{x}
        =
        \varphi_2\br{x}\cP_\lambda f\br{x}- 
        \cP_\lambda\br{\varphi_2 f}\br{x}
        =
        \int
            \lambda \widehat{\psi_\delta}\br{\lambda z}
            f\br{x+z}
            \br{\varphi_2\br{x}-\varphi_2\br{x+z}}
        dz.
    \end{equation*}
    Via mean value theorem and the trivial bound \(\abs{\widehat{\psi_\delta}\br{z}}\leq \frac{A^{O\br{C}}}{\ang{z}^3}\), we obtain the following point-wise estimate:
    \begin{equation*}
        \abs{
            \mr{
                M_{\varphi_2},
                \cP_\lambda
            }
            f\br{x}
        }
        \leq
        \nrm{
            \frac{A^{O\br{C}}\lambda}{\ang{\lambda z}^3}
            f\br{x+z}
            \nrm{\varphi'_2}_{L^\infty}
            z
        }_{L^1\br{dz}}
    \end{equation*}
    With some further simplification, we dominate the above with
    \begin{equation*}
        \leq
        \frac{
        A^{O\br{C}}}{\abs{\lambda}}
        \nrm{
            \frac{\lambda}{\ang{\lambda z}^2}
            f\br{x+z}
        }_{L^1\br{dz}}
        \leq
        \frac{
        A^{O\br{C}}}{\abs{\lambda}}
        \cM f\br{x}.
    \end{equation*}
    In other words, we have the following approximation:
    \begin{equation*}
        \varphi_2\br{x}\cP_\lambda f\br{x}=\cP_\lambda \br{\varphi_2f}\br{x}+\frac{A^{O\br{C}}}{\abs{\lambda}} 
        O\br{\cM f\br{x}}.
    \end{equation*}
The same statement holds fiber-wisely for functions \(f_j\in L^2\br{\R^2}\):
\begin{equation*}
   \varphi_2\br{x_j}\fP{j}_{\lambda_j}
    f_j\br{x
    }=
    \br{\cP_\lambda M_{\varphi_2}}^{\br{j}}
    f_j\br{x}+
    \frac{A^{O\br{C}}}{\abs{\lambda}} 
    O\br{\fM{j} f_j\br{x}}.
\end{equation*}
Inserting the above identity into \eqref{eq_APP_st_loc_ini_id} yields:
\begin{equation*}
    \1_{\mr{0,1}^2}\br{\vx}
    \cA_{\vgamma,\rho}
        \br{
            \fP{1}_{\lambda_1}f_1,
            \fP{2}_{\lambda_2}f_2
        }\br{\vx}
    =
    \eqref{eq_APP_st_loc_main}
\end{equation*}
\begin{equation*}
    +\frac{A^{O\br{C}}}{\abs{\lambda_2}}
    O\br{
    \cA_{\vgamma,\abs{\rho}}
        \br{
            M^{\br{1}}_{\varphi_2}
            \fP{1}_{\lambda_1}
            \br{
                \br{
                    \varphi_2\otimes
                    \1_{\mr{0,1}}
                }
            f_1},
            \br{
                \1_{\mr{0,1}}\otimes
                \varphi_2
            }
            \fM{2}f_2
        }\br{\vx}
    }
\end{equation*}
\begin{equation*}
    +\frac{A^{O\br{C}}}{\abs{\lambda_1}}
    O\br{
    \cA_{\vgamma,\abs{\rho}}
        \br{
            \br{
                \varphi_2\otimes
                \1_{\mr{0,1}}
            }
            \fM{1}f_1
            ,
            M^{\br{2}}_{\varphi_2}
            \fP{2}_{\lambda_2}
            \br{
                \br{
                    \1_{\mr{0,1}}\otimes
                    \varphi_2
                }
            f_2}
        }\br{\vx}
    }
\end{equation*}
\begin{equation*}
    +\frac{A^{O\br{C}}}{\abs{\lambda_1\lambda_2}}
    O\br{
    \cA_{\vgamma,\abs{\rho}}
        \br{
            \br{
                \varphi_2\otimes
                \1_{\mr{0,1}}
            }
            \fM{1}f_1
            ,
            \br{
                \1_{\mr{0,1}}\otimes
                \varphi_2
            }
            \fM{2}f_2
        }\br{\vx}
    }.
\end{equation*}
Finally, applying the trivial fact 
\begin{equation*}
    \abs{
        \br{M_{\varphi_2}\cP_\lambda M_{\varphi_2}}^{\br{j}}
        f_j\br{x}
    }
    \leq 
    A^{O\br{C}}
    \br{M_{\varphi_2}\cM M_{\varphi_2} }^{\br{j}}f_j\br{\vx}
    \leq 
    A^{O\br{C}}
    \varphi_2\br{x_j}\fM{j}f_j\br{\vx}
\end{equation*}
and triangle inequalities conclude the proof of \textbf{Lemma \ref{lem_phy_loc_pass_proj}}.
\end{proof}

\subsection{Proof of Theorem \ref{mainthm_global_w_vanishing_n_fluc}}\label{subsec_pf_mainthm_global_w_vanishing_n_fluc}
Assume the finiteness of the following quantities:
\begin{equation*}
    M_{p,\alpha}:=
    \max\br{\nrm{A_\vgamma}_{L^{p,\infty}},\,\abs{A_\vgamma}_\alpha,\frac{1}{\dist\br{\supp\rho,I^c}}},\quad \Xi_\Gamma < \infty
\end{equation*}
and make the normalization \(\nrm{\rho}_{C^2}=1\). By \textbf{Remark \ref{rmk_p_infty_case_comments}}, it suffices to consider the \(0<p<\infty\) case. Recall now the Littlewood-Payley projection \textbf{Definition \ref{def_LP_decomp}}. We take \(\cQ:=\cQ_{1/2,1}\) and \(\cP_\lambda:=\cP_{1/2,\lambda}\) and perform the standard decomposition.
\begin{equation}\label{eq_cA_freq_decomp_L_for_var}
    \cA_{\vgamma,\rho}\br{f_1,f_2}=
    \cA_{\vgamma,\rho}\br{\fQ{1} f_1, \fQ{2} f_2}
\end{equation}
\begin{equation}\label{eq_cA_freq_decomp_HL_for_var}
    +
    \frac{1}{c_{1/2}}
    \int_{1\leq \abs{\lambda}}
    \cA_{\vgamma,\rho}\br{\fQ{1} f_1, \fP{2}_\lambda f_2}
    +
    \cA_{\vgamma,\rho}\br{\fP{1}_\lambda f_1, \fQ{2} f_2}
    \frac{d\lambda}{\abs{\lambda}}
\end{equation}
\begin{equation}\label{eq_cA_freq_decomp_H_for_var}
    +
    \frac{1}{c^2_{1/2}}
    \int_{1\leq \abs{\lambda_j}}
    \cA_{\vgamma,\rho}\br{\fP{1}_{\lambda_1} f_1, \fP{2}_{\lambda_2} f_2}
    \frac{d\lambda_1}{\abs{\lambda_1}}
    \frac{d\lambda_2}{\abs{\lambda_2}}.
\end{equation}
Via the simple estimate \eqref{eq_L2_2_H_neq}, we dominate the first term trivially:
\begin{equation}\label{eq_cA_good_n_defect_trivial_est}
    \nrm{\cA_{\vgamma,\rho}\br{\fQ{1} f_1, \fQ{2} f_2}}_{L^1}
    \leq 
    \nrm{\rho}_{L^1}
    \nrm{\fQ{1} f_1}_{L^2}
    \nrm{\fQ{2} f_2}_{L^2}
    \lesssim
    \abs{I}
    \nrm{f_1}_{H^{\br{-\varepsilon,0}}}
    \nrm{f_2}_{H^{\br{0,-\varepsilon}}}
\end{equation}
for all \(0< \varepsilon \leq 1\). It remains to address \eqref{eq_cA_freq_decomp_HL_for_var} and \eqref{eq_cA_freq_decomp_H_for_var}. To provide a unified treatment, we write:
\begin{equation}\label{eq_cR_4_var}
    \cR_\lambda:=\begin{cases}
        \cQ & \abs{\lambda} \leq 1;\\
        \cP_\lambda & \text{else}.
    \end{cases}
\end{equation}
For each fixed \(\vlambda:=\br{\lambda_1,\lambda_2}\) with \(\abs{\lambda_1},\abs{\lambda_2}\geq 1\), we aim to show:
\begin{equation}\label{eq_cA_freq_decomp_main_goal}
    \nrm{
        \cA_{\vgamma,\rho}\br{\fR{1}_{\lambda_1} f_1, \fR{2}_{\lambda_2} f_2}
    }_{L^1}\lesssim M^{c_\alpha}_{p,\alpha}
    \ang{I}\Xi^{\frac{3}{10}}_\Gamma \abs{\vlambda}^{-c}\nrm{f_1}_{L^2}\nrm{f_2}_{L^2},\quad\text{for some }c=c\br{p,\alpha}>0.
\end{equation}
To prove \eqref{eq_cA_freq_decomp_main_goal},
we decompose the operator into a main term, where \textbf{Theorem \ref{mainthm_global}} applies, and an error term with a \(\abs{\vlambda}\) power decay.
Let \(\lambda:=\abs{\vlambda}\) and \(\delta,d>0\) to be decided later. Consider:
\begin{equation}\label{eq_S_def_defect_set}
    S:=
    \BR{
        t\in I
    \::\:
        A_\vgamma\br{t}>\lambda^\delta\nrm{A_\vgamma}_{L^{p,\infty}}
    }.
\end{equation}
Notice that the condition \(S=\varnothing\) implies that \(\vgamma:I\simto\Gamma\) is \(A_\lambda:=\lambda ^\delta\nrm{A_\vgamma}_{L^p,\infty}\)-admissible. In this situation, \eqref{eq_cA_freq_decomp_main_goal} is a direct consequence of \textbf{Theorem \ref{mainthm_global}} whenever \(0<\delta<\varepsilon\). Henceforth, we assume
\begin{equation}\label{eq_delta_size_S_nonempty}
    0<\delta<\varepsilon,\quad 
    S\neq\varnothing.
\end{equation}
By \textbf{Definition \ref{def_weak_failure} and \ref{def_A_gamma_fluc}}, we have two estimates regarding \(S\):
\begin{equation*}
    \abs{S}\leq \lambda^{-\delta p},\quad
    C\br{S}\leq \lambda^{\delta \alpha}\nrm{A_\vgamma}^\alpha_{L^{p,\infty}}\abs{A_\vgamma}_\alpha.
\end{equation*}
Define now the auxiliary open set:
\begin{equation}\label{eq_S_d_as_union}
    S_d:=\bigcup_{t\in S}\br{t-d,t+d}.
\end{equation}
By construction, we have \(C\br{S_d}\leq C\br{S}\). As a result, the set \(S_d\) can be decomposed into \(N\leq C\br{S}\) many disjoint unions of open intervals \(\BR{\br{a_k,b_k}}^N_{k=1}\) arranged in ascending order \(b_k\leq a_{k+1}\). By choosing:
\begin{equation}\label{eq_d_choice}
    d:=\frac{
        \lambda^{-\delta\br{\alpha+p}}
    }{
        \nrm{A_\vgamma}^\alpha_{L^{p,\infty}}
        \abs{A_\vgamma}_\alpha
    },
\end{equation}
we deduce from \eqref{eq_delta_size_S_nonempty} the following upper bound:
\begin{equation*}
% \label{eq_d_size_in_strong_smoothing}
    d\leq \lambda^{-\delta p}/ C\br{S} \leq  \lambda^{-\delta p} \leq 1
\end{equation*}
and the following size estimate:
\begin{equation}\label{eq_S_d_size}
    \abs{S_d}=\sum^N_{k=1}\abs{I_k}\leq
    \abs{S}+C\br{S}\cdot 2 d
    \lesssim 
    \lambda^{-\delta p}
    +\lambda^{\delta \alpha}
    \nrm{A_\vgamma}^\alpha_{L^{p,\infty}}
    \abs{A_\vgamma}_\alpha d
    \lesssim \lambda^{-\delta p}.
\end{equation}
We now introduce smooth bump functions \(\chi_k\in C^\infty_c\br{\R}\) adapted to the intervals \(\br{a_k,b_k}\) in the following sense:
\begin{equation}\label{eq_smooth_bump_for_defect}
    \1_{\mr{a_k+2d/3,b_k-2d/3}}\leq \chi_k \leq \1_{\br{a_k+d/3,b_k-d/3}},\quad
    \nrm{\chi_k}_{C^2}\lesssim \ang{1/d}^2\eqsim d^{-2}.
\end{equation}
This allows us to extract the parts of the operator where \textbf{Theorem \ref{mainthm_global}} does not apply due to defects of \(\vgamma\br{t}\) for \( t\in S_d\). Explicitly, we consider:
\begin{equation}\label{eq_cA_defect_decomp}
    \cA_{\vgamma,\rho}\br{\fR{1}_{\lambda_1}f_1,\fR{2}_{\lambda_2}f_2}
    =\cA_{\vgamma,\trho}\br{\fR{1}_{\lambda_1}f_1,\fR{2}_{\lambda_2}f_2}
    +\sum^N_{k=1}\cA_{\vgamma,\rho\chi_k}\br{\fR{1}_{\lambda_1}f_1,\fR{2}_{\lambda_2}f_2}
    ,\quad
    \trho:=\rho-\rho\sum^N_{k=1}\chi_k.
\end{equation}
Apply trivial estimates to the summand. We deduce:
\begin{equation*}
    \nrm{\sum^N_{k=1}\cA_{\vgamma,\rho\chi_k}\br{\fR{1}_{\lambda_1}f_1,\fR{2}_{\lambda_2}f_2}}_{L^1}
    \leq 
    \sum^N_{k=1}
    \nrm{\rho\chi_k}_{L^1}
    \nrm{f_1}_{L^2}
    \nrm{f_2}_{L^2}
\end{equation*}
\begin{equation}\label{eq_cA_defect_est}
    \leq
    \abs{S_d}\cdot
    \nrm{f_1}_{L^2}
    \nrm{f_2}_{L^2}
    \lesssim \lambda^{-\delta p}
    \nrm{f_1}_{L^2}
    \nrm{f_2}_{L^2}.
\end{equation}
We now treat the remaining contribution. Due to \eqref{eq_smooth_bump_for_defect}, \(\trho\) can be further decomposed into \(N+1\) many smooth bump functions \(\rho_l\in C^\infty_c\br{\R}\) with disconnected supports. More precisely, we have:
\begin{equation}\label{eq_supp_good_rho_ends}
    \supp \rho_0\subseteq \br{-\infty ,a_1+2d/3}\cap \supp\rho,\quad
    \supp \rho_N\subseteq \br{b_N-2d/3,\infty}\cap \supp\rho;
\end{equation}
\begin{equation}\label{eq_supp_good_rho_mids}
    \supp \rho_l\subseteq \br{b_l-2d/3,a_{l+1}+2d/3}\cap \supp\rho ,\quad
    0<l<N.
\end{equation}
Correspondingly, we decompose the operator:
\begin{equation}\label{eq_cA_good_decomp}
    \cA_{\vgamma,\trho}
        \br{
            \fR{1}_{\lambda_1}f_1,
            \fR{2}_{\lambda_2}f_2
        }
    =
    \sum^N_{l=0}
        \cA_{\vgamma,\rho_l}
            \br{
                \fR{1}_{\lambda_1}f_1,
                \fR{2}_{\lambda_2}f_2
            }
\end{equation}
and consider sub-intervals of \(I\) defined below:
\begin{equation}\label{eq_good_sub_I_ends}
    I_0:=\br{-\infty ,a_1+d}\cap I,\quad
    I_N:=\br{b_N-d,\infty}\cap I;
\end{equation}
\begin{equation}\label{eq_good_sub_I_mids}
    I_l:=\br{b_l-d,a_{l+1}+d}\cap I,\quad 0<l<N.
\end{equation}
As a consequence of the design of \(\rho_l,I_l\) and construction \eqref{eq_S_d_as_union}, we have the following relations:
\begin{equation}\label{eq_good_rho_sub_I_rels}
    \nrm{\rho_l}_{C^2}\lesssim 1+d^{-1}+d^{-2}\eqsim d^{-2},\quad
    \supp \rho_l\subseteq I_l\subseteq I\setminus S,\quad
    \dist \br{\supp \rho_l,I^c_l }\geq \min\br{d/3,\dist\br{\supp\rho,I^c}}
\end{equation}
whenever \(\rho_l\) is not identically zero. Given such an \(l\), we restrict curves to the sub-interval \(I_l\):
\begin{equation}
    \vgamma_l:=\dif{\vgamma}_{I_l}:I_l\simto \Gamma_l\subseteq \Gamma.
\end{equation}
By \eqref{eq_S_def_defect_set} and \eqref{eq_good_rho_sub_I_rels}, we have as a direct consequence that \(\vgamma_l\) is \(A_\lambda:=\lambda^\delta\nrm{A_\vgamma}_{L^{p,\infty}}\)-admissible.
We can thus apply \textbf{Theorem \ref{mainthm_global}} to estimate \eqref{eq_cA_good_decomp}:
\begin{equation*}
    \nrm{
        \eqref{eq_cA_good_decomp}
    }_{L^1}
    \leq
    \sum^N_{l=0}
    \nrm{
        \cA_{\vgamma,\rho_l}
            \br{
                \fR{1}_{\lambda_1}f_1,
                \fR{2}_{\lambda_2}f_2
            }
    }_{L^1}
\end{equation*}
\begin{equation*}
    \lesssim
    \sum^N_{l=0}
    \max\br{
        A_\lambda,
        \frac{1}{\dist\br{\supp \rho_l,I^c_l}}
    }^C
    \abs{I_l}
    \Xi^{\frac{3}{10}}_{\Gamma_l}
    \nrm{\rho_l}_{C^2}
    \nrm{
        \fR{1}_{\lambda_1}f_1
    }_{H^{\br{-\varepsilon,0}}}
    \nrm{
        \fR{2}_{\lambda_2}f_2
    }_{H^{\br{0,-\varepsilon}}}.
\end{equation*}
Recall now the choice \eqref{eq_d_choice} and the relations \eqref{eq_good_rho_sub_I_rels}. We dominate the above with:
\begin{equation*}
    \underset{C}{\lesssim}
    \max
    \br{
        \lambda^\delta\nrm{A_\vgamma}_{L^{p,\infty}},
        \lambda^{\delta\br{\alpha+p}}
        \nrm{A_\vgamma}^\alpha_{L^{p,\infty}}
        \abs{A_\vgamma}_\alpha,
        \frac{1}{\dist\br{\supp \rho,I^c}}
    }^C
    \abs{I}
    \br{\max_{0\leq l\leq N} \Xi_{\Gamma_l}}^{\frac{3}{10}}
\end{equation*}
\begin{equation*}
    \cdot
    \br{
        \lambda^{\delta\br{\alpha+p}}
        \nrm{A_\vgamma}^\alpha_{L^{p,\infty}}
        \abs{A_\vgamma}_\alpha
    }^2
    \lambda^{-\varepsilon}
    \nrm{
        f_1
    }_{L^2}
    \nrm{
        f_2
    }_{L^2}
\end{equation*}
By direct computation and \textbf{Lemma \ref{lem_M_Gamma_mono}}, we further dominate the above with
\begin{equation*}
    \lesssim
    M^{\br{C+3}\max\br{1,\alpha}}_{p,\alpha}
    \lambda^{\delta \br{C+2} \max\br{1,\alpha+p}-\varepsilon}
    \abs{I}
    \Xi^{\frac{3}{10}}_\Gamma
    \nrm{
        f_1
    }_{L^2}
    \nrm{
        f_2
    }_{L^2}.
\end{equation*}
Finally, we take:
\begin{equation*}
    c_\alpha:= \br{C+3}\max\br{1,\alpha}
\end{equation*}
and solve \(\delta>0\) for
\begin{equation*}
    -\sigma_{p,\alpha}:=
    \delta \br{C+2} \max\br{1,\alpha+p}-\varepsilon=-\delta p.
\end{equation*}
we deduce:
\begin{equation}\label{eq_cA_good_est}
    \nrm{
    \cA_{\vgamma,\trho}
        \br{
            \fR{1}_{\lambda_1}f_1,
            \fR{2}_{\lambda_2}f_2
        }
    }_{L^1}
    \lesssim
    M^{c_\alpha}_{p,\alpha}
    \lambda^{-\sigma_{p,\alpha}}
    \abs{I}
    \Xi^{\frac{3}{10}}_\Gamma
    \nrm{
        f_1
    }_{L^2}
    \nrm{
        f_2
    }_{L^2}.
\end{equation}
Combining this with \eqref{eq_cA_defect_est}, we conclude \eqref{eq_cA_freq_decomp_main_goal}.
Finally, direct computation yields:
\begin{equation*}
    \nrm{
        \frac{1}{c_{1/2}}
        \int_{1\leq \abs{\lambda}}
        \cA_{\vgamma,\rho}\br{\fQ{1} f_1, \fP{2}_\lambda f_2}
        +
        \cA_{\vgamma,\rho}\br{\fP{1}_\lambda f_1, \fQ{2} f_2}
        \frac{d\lambda}{\abs{\lambda}}
    }_{L^1}
\end{equation*}
\begin{equation*}
    +
    \nrm{
        \frac{1}{c^2_{1/2}}
        \int_{1\leq \abs{\lambda_j}}
        \cA_{\vgamma,\rho}\br{\fP{1}_{\lambda_1} f_1, \fP{2}_{\lambda_2} f_2}
        \frac{d\lambda_1}{\abs{\lambda_1}}
        \frac{d\lambda_2}{\abs{\lambda_2}}
    }_{L^1}
\end{equation*}
\begin{equation*}
    \lesssim
    M^{c_\alpha}_{p,\alpha}
    \ang{\abs{I}}
    \Xi^{\frac{3}{10}}_\Gamma
    \br{
    \int_{1\leq \abs{\lambda}}
        \lambda^{2\varepsilon_{p,\alpha}-\sigma_{p,\alpha}}
    \frac{d\lambda}{\abs{\lambda}}
    +
    \br{
    \int_{1\leq \abs{\lambda}}
        \lambda^{\varepsilon_{p,\alpha}-\sigma_{p,\alpha}/2}
    \frac{d\lambda}{\abs{\lambda}}
    }^2
    }
    \nrm{f_1}_{H^{\br{-\varepsilon_{p,\alpha},0}}}
    \nrm{f_2}_{H^{\br{0,-\varepsilon_{p,\alpha}}}}
\end{equation*}
\begin{equation*}
    \underset{\sigma_{p,\alpha},\varepsilon_{p,\alpha}}{\lesssim}
    M^{c_\alpha}_{p,\alpha}
    \ang{\abs{I}}
    \Xi^{\frac{3}{10}}_\Gamma
    \nrm{f_1}_{H^{\br{-\varepsilon_{p,\alpha},0}}}
    \nrm{f_2}_{H^{\br{0,-\varepsilon_{p,\alpha}}}},\quad
    \forall \varepsilon_{p,\alpha} \in\br{0,\sigma_{p,\alpha}/2}.
\end{equation*}
Together with the trivial estimate \eqref{eq_cA_good_n_defect_trivial_est}, this concludes the proof of \textbf{Theorem \ref{mainthm_global_w_vanishing_n_fluc}}. \qed

\section{Smoothing inequalities and O-minimal theory}\label{sec_smoothing_o_min}
    Observe that a definable family \(\mathfrak{F}:=\BR{\vgamma_\sigma: \R\supset I_\sigma\simto \Gamma_\sigma\subseteq \R^2}_{\sigma\in\Sigma}\) of \(C^1\) parametric curves induces a definable family \(\BR{\Gamma_\sigma\subseteq \R^2}_{\sigma\in\Sigma}\) of \(C^1\) curves. 
    This is a direct consequence of the fact that projection preserves definability. Indeed, recall \textbf{Definition \ref{def_struct_over_M} and \ref{def_definable_family}}. We associate the family of parametric curves with a definable function \(F:\Sigma\times\R\supset \dom F \to \R^2\) such that the following relations hold:
    \begin{equation}\label{eq_def_fam_parametric_curve}
        \vgamma_\sigma\br{t}=F\br{\sigma,t},\quad
        \forall t\in \dom \vgamma_\sigma=\BR{
            t\in\R
        \::\:
            \br{\sigma,t}\in \dom F
        }=
        \br{\dom F}_\sigma.
    \end{equation}
    Consider now the set:
    \begin{equation*}
        S:=
        \BR{
            \br{\sigma,\vgamma}\in \Sigma\times \R^2
        \::\:
            \exists t\in I_\sigma,\, \vgamma=\vgamma_\sigma\br{t}
        }.
    \end{equation*}
    By construction, we have that:
    \begin{equation*}
        \Gamma_\sigma=\BR{
            \vgamma\in\R^2
        \::\:
            \br{\sigma,\vgamma}\in S
        }=S_\sigma,\quad \forall \sigma\in \Sigma.
    \end{equation*}
    To show that \(\BR{\Gamma_\sigma}_{\sigma\in\Sigma}\) is a definable family, it now suffices to show the definibility of \(S\), which can be seen by the following fact:
    \begin{equation*}
        S=
        \operatorname{Perm}\circ
        \pi
        \BR{
            \br{\vgamma,\sigma,t}\in \R^2\times \dom F 
        \::\:
            F\br{\sigma,t}=\vgamma
        },\quad
        \pi\br{\vgamma,\sigma,t}:=\br{\vgamma,\sigma},\quad
        \operatorname{Perm}\br{\vgamma,\sigma}=\br{\sigma,\vgamma}.
    \end{equation*}
    By \textbf{Definition \ref{def_struct_over_M}},
    the definability of \(S\) thus follows from the definability of \(F\).

Combining the fact that all \(A\)-admissible curves have injective Gauss maps and the above observation with \textbf{Lemma \ref{lem_uni_fini_prop_for_Gamma}} and \textbf{Lemma \ref{lem_uni_fini_prop_for_A_vgamma}},
we deduce the uniformity result \textbf{Theorem \ref{thm_smoothing_o_mini_uni}} from \textbf{Theorem \ref{mainthm_global}} and uniformity result \textbf{Theorem \ref{thm_smoothing_o_mini_uni_w_vanishing_n_fluc}} from \textbf{Theorem \ref{mainthm_global_w_vanishing_n_fluc}}.

\subsection{Proof of Lemma \ref{lem_uni_fini_prop_for_Gamma}}

Recall \textbf{Definition \ref{def_definable_family}}. We may thus associate the family of curves with a definable set \(S\subseteq \Sigma\times\R^2\subseteq \R^{k+2}\) via the following relation
    \begin{equation*}
        \Gamma_\sigma:=\BR{\vgamma\in \R^2\::\: \br{\sigma,\vgamma}\in S },\quad
        \sigma\in \Sigma.
    \end{equation*}
For each \(\sigma\in \Sigma\), the natural extension of the inverse Gauss map \(\vgamma_{\Gamma_\sigma}\) of \(\Gamma_\sigma\) exists. As a result, the function \(F\br{\sigma,\vxi}:=\vgamma_{\Gamma_\sigma}\br{\vxi}\) exists on the natural domain:
\begin{equation*}
    \dom F:=\BR{\br{\sigma,\vxi}\in \Sigma\times\R^2\::\: \vxi\in \dom\vgamma_{\Gamma_\sigma}=\cC_{\Gamma_\sigma}}.
\end{equation*}
We claim that such function is definable. Consider the first-order formula \(\phi\br{\sigma,\vxi,\vgamma}\):
\begin{equation*}
    \forall \epsilon>0,\,\exists \delta>0,\, \forall \vzeta\in\R^2,\quad 
    \Big(
        \big(
            \br{\sigma,\vzeta}\in S\land 
            0<\abs{\vgamma-\vzeta}<\delta
        \big)
    \implies\,
        \abs{\ang{\vxi,\vgamma-\vzeta}}<\epsilon \abs{\vgamma-\vzeta}
    \Big).
\end{equation*}
The following set is definable:
\begin{equation*}
    \BR{
        \br{\sigma,\vxi,\vgamma}\in \Sigma\times\br{\R^2\setminus\BR{\vnull}}\times \R^2\::\:
        \br{\sigma,\vgamma}\in S\land
        \phi\br{\sigma,\vxi,\vgamma}
    }
\end{equation*}
and coincides with the graph of the function \(F\). Consequently, the relation \(F\br{\sigma,\vxi}=\vgamma\) and thus, the function \(F=\br{F_1,F_2}\) are definable.
To control \(\Xi_{\Gamma_\sigma}\), we consider the auxiliary functions \(G_i:=\partial_{k+i}F_i\) for \(i=1,2\). By \textbf{Proposition \ref{prop_func_def_closure} and \ref{prop_diff_def}}, \(G_1,G_2\) are definable. Furthermore, the following four functions:
\begin{equation*}
    h_1\br{\sigma,s,\xi,\eta}:=\br{\Delta_{\br{0,s}}\partial_1\vgamma_{\Gamma_\sigma,1}}\br{\xi,\eta}=G_1\br{\sigma,\xi,\eta+s}-G_1\br{\sigma,\xi,\eta},
\end{equation*}

\begin{equation*}
    h_2\br{\sigma,s,\xi,\eta}:=\br{\Delta_{\br{s,0}}\partial_2\vgamma_{\Gamma_\sigma,2}}\br{\eta,\xi}=G_2\br{\sigma,\eta+s,\xi}-G_2\br{\sigma,\eta,\xi},
\end{equation*}

\begin{equation*}
    H_1\br{\sigma,s,u,v,\xi,\eta}:=\br{\Delta_{\br{0,s}}\Delta_{\br{u,v}}\partial_1\vgamma_{\Gamma_\sigma,1}}\br{\xi,\eta}
    =h_1\br{\sigma,s,\xi+u,\eta+v}-h_1\br{\sigma,s,\xi,\eta},
\end{equation*}

\begin{equation*}
    H_2\br{\sigma,s,u,v,\xi,\eta}:=\br{\Delta_{\br{s,0}}\Delta_{\br{v,u}}\partial_2\vgamma_{\Gamma_\sigma,2}}\br{\eta,\xi}
    =h_2\br{\sigma,s,\xi+u,\eta+v}-h_2\br{\sigma,s,\xi,\eta},
\end{equation*}
are definable by \textbf{Proposition \ref{prop_func_def_closure}} and \textbf{Remarks \ref{rmk_def_func}}. Introduce now the following auxiliary sets:
\begin{equation*}
    A^{\br{i}}:=\BR{
        \br{\sigma,\eta,s,w,\xi}\in\Sigma\times \R^4
    \::\:
        h_i\br{\sigma,s,\xi,\eta}=w
    },
\end{equation*}
\begin{equation*}
    B^{\br{i}}:=\BR{
        \br{\sigma,\eta,s,u,v,w,\xi}\in\Sigma\times \R^6
    \::\:
        H_i\br{\sigma,s,u,v,\xi,\eta}=w
    }.
\end{equation*}
By \textbf{Remark \ref{rmk_def_set_construct}}, \(A^{\br{1}},A^{\br{2}},B^{\br{1}},B^{\br{2}}\) are definable sets. We finish the proof by applying \textbf{Proposition \ref{prop_uni_fin}}:
\begin{equation*}
    \sup_{\sigma\in\Sigma}\Xi_{\Gamma_\sigma}\leq\max_{i=1,2}\sup_{\substack{s,u,v,\\ \sigma,\eta,w}}
    \max\br{
            C\br{A^{\br{i}}_{\br{\sigma,\eta,s,w}}},
            C\br{B^{\br{i}}_{\br{\sigma,\eta,s,u,v,w}}}
    }\leq \max_{i=1,2}\br{M_{{A^{\br{i}}}},M_{B^{\br{i}}}}<\infty.
\end{equation*}

\subsection{Proof of Lemma \ref{lem_uni_fini_prop_for_A_vgamma}}
    By \textbf{Proposition \ref{prop_uni_fin}}, it suffices to demonstrate the following family of sets:
    \begin{equation}
    \label{eq_A_vgamma_def_fam_set_S_sig_lambda}
        \BR{
            S_{\sigma,\lambda}:=\BR{
                t\in I_\sigma
            \::\:
                A_{\vgamma_\sigma}\br{t}>\lambda
            }
            \subseteq \R
        }_{\br{\sigma,\lambda}\in\Sigma\times\br{0,\infty}}
    \end{equation}
    is a definable family. By \textbf{Definition \ref{def_definable_family}}, there is a definable function \(F:\Sigma\times \R\supset S\to\R^2\) such that for all \(\sigma\in\Sigma\), the relations \eqref{eq_def_fam_parametric_curve} hold. Consider now the following definable sets:
    \begin{equation*}
        U:=\bigcup^4_{k=1}
        \BR{
            \br{\sigma,\lambda,t}\in \Sigma\times \br{0,\infty}\times \R
        \::\:
            \abs{
                \partial^k_t F\br{\sigma,t}
            }>\lambda
        },
    \end{equation*}
    \begin{equation*}
        V:=\bigcup_{j=1,2}
        \BR{
            \br{\sigma,\lambda,t}\in \Sigma\times \br{0,\infty}\times \R
        \::\:
            \abs{
                \partial_t F_j\br{\sigma,t}
            }<1/\lambda
        },
    \end{equation*}
    \begin{equation*}
        W:=
        \BR{
            \br{\sigma,\lambda,t}\in \Sigma\times \br{0,\infty}\times \R
        \::\:
            \abs{
                \partial_t F
                \wedge
                \partial^2_t F
            }\br{\sigma,t}
            <1/\lambda
        },
    \end{equation*}
    \begin{equation*}
        X:=
        \BR{
            \br{\sigma,\lambda,t}\in \Sigma\times \br{0,\infty}\times \R\setminus W
        \::\:
            \abs{
                \partial_t \vchi_F
                \wedge
                \partial^2_t \vchi_F
            }\br{\sigma,t}
            <1/\lambda
        }
    \end{equation*}
    with the function \(\vchi_F\) defined on the set \(\BR{\partial_t F
    \wedge
    \partial^2_t F\neq 0}\) with the explicit formula:
    \begin{equation*}
        \vchi_F\br{\sigma,t}:=
        \frac{
            \partial_t F_1
            \partial_t F_2
            \partial_t F
        }{
            \partial_t F
            \wedge
            \partial^2_t F
        }
        \br{\sigma,t}.
    \end{equation*}
    It is easy to check that:
    \begin{equation*}
        S_{\sigma,\lambda}
        =\br{
            U\cup V\cup W\cup X
        }_{\sigma,\lambda}.
    \end{equation*}
    This demonstrates that \(\BR{S_{\sigma,\lambda}}_{\br{\sigma,\lambda}\in\Sigma\times\br{0,\infty}}\) is a definable family and thus, concludes the proof.

\subsection{Proof of Theorem \ref{thm_real_ana_smoothing}}\label{subsec_pf_thm_real_ana_smoothing}

Fix \(\vgamma:I\to\R^2\) real analytic and \(\rho\in C^2_c\br{\R}\) with \(\supp\rho\Subset I\). We begin with the following trivial geometric fact:
\begin{equation}\label{eq_iff_line}
    \vgamma'\wedge\vgamma''\equiv 0 \iff \img \vgamma\subseteq \eqref{eq_line}
\end{equation}
for some \(\br{a,b,c}\in\R^3\setminus\BR{\br{0,0,0}}\). Assuming now that \(\vgamma'\wedge\vgamma''\not\equiv 0\), we shall demonstrate the following slightly more involved geometric statement:
\begin{equation}\label{eq_iff_3else}
    \vchi'_\vgamma\wedge\vchi''_\vgamma\equiv 0 \iff
    \br{
        \img \vgamma \subseteq \eqref{eq_exp} \text{ or }
        \img \vgamma \subseteq \eqref{eq_log} \text{ or }
        \img \vgamma \subseteq \eqref{eq_psi_line}
    }
\end{equation}
for some \(\br{a,b,c}\in\R^3\setminus\BR{\br{0,0,0}}\) and \(A,B\in\R\setminus\BR{0}\). Details are reserved in \textsc{Section \ref{subsub_iff_ode}}.

As a direct consequence of \eqref{eq_iff_line} and \eqref{eq_iff_3else}, whenever \(\img \vgamma\) is not contained in any of the four: \eqref{eq_line}, \eqref{eq_exp}, \eqref{eq_log}, or \eqref{eq_psi_line}, none of the below expressions
\begin{equation}\label{eq_4_nonzero}
    \gamma'_j,\, \vgamma'\wedge \vgamma'',\, \vchi'_\vgamma\wedge\vchi''_\vgamma,\quad j=1,2,
\end{equation}
are zero functions. Since \(\gamma'_1,\gamma'_2,\vgamma'\wedge\vgamma''\) are real analytic, and \(\vchi'_\vgamma\wedge\vchi''_\vgamma\) is real meromorphic, we deduce that when none of \eqref{eq_4_nonzero} are zero functions, they can only vanishes up to some finite order at finitely many points in \(\supp \rho\). After suitable partition of \(\rho\), we may isolate each points where any of the expressions in \eqref{eq_4_nonzero} vanishes. Thus, with the above partition and changes of variables,
it suffices to consider a function \(\vgamma:I:=\br{-1,1}\to\R^2\) with a real analytic extension on a neighborhood of the closed interval \(\Br{-1,1}\);\footnote{This makes sure the function \(\vgamma:I:=\br{-1,1}\to\R^2\) is definable in \(\R_\an\) (see \cite{MR245831,MR425152,MR972342,MR1374342,MR1389958}).} additionally, the following expressions have at most finite vanishing order at the single point \(t=0\):
\begin{equation}\label{eq_ana_bad_at_0}
    \abs{\gamma'_j}\br{t},\,
    \abs{\vgamma'\wedge\vgamma''}\br{t},\,
    \abs{\vchi'_\vgamma\wedge\vchi''_\vgamma}\br{t}
    \gtrsim \abs{t}^n ,\quad
    \forall t\in I:=\br{-1,1},\,j=1,2.
\end{equation}
% when proving the sufficiency of the condition in \textbf{Theorem \ref{thm_real_ana_smoothing}}. 
After the above reduction, a direct computation shows that the admissibility function \(A_\vgamma\) as defined in \textbf{Definition \ref{def_admi_func}} satisfies the estimate
\begin{equation}\label{eq_admin_func_ana_bd}
    \nrm{A_\vgamma}_{L^{1/n,\infty}\br{I}}\lesssim 1.
\end{equation}
Using now the definability of \(\vgamma\) in \(\R_\an\) (see \cite{MR245831,MR425152,MR972342,MR1374342,MR1389958}), we aim to apply \textbf{Theorem \ref{thm_smoothing_o_mini_uni_w_vanishing_n_fluc}} to conclude the sufficiency of the condition in \textbf{Theorem \ref{thm_real_ana_smoothing}}. Unfortunately, \(\vgamma\) may have self-intersection and other undesirable geometric pathologies. To remedy those issues, we perform further surgery on \(\vgamma\). This is the content of \textsc{Section \ref{subsub_suffi_real_ana_smoothing}}.

Finally, we demonstrate the necessity of the condition by explicitly constructing functions \(f_{j,\lambda}\) such that:
\begin{equation}\label{eq_res_obstr}
    \nrm{f_{1,\lambda}}_{H^{\br{-\varepsilon,0}}}
    \nrm{f_{2,\lambda}}_{H^{\br{0,-\varepsilon}}}\lesssim \lambda^{-\varepsilon}\ll
    1\lesssim \nrm{\cA_{\vgamma,\rho}\br{f_{1,\lambda},f_{2,\lambda}}}_{L^1}\quad  \forall \lambda \gg 1,\,\varepsilon>0.
\end{equation}
The existence of these functions is a result of certain generalized resonance phenomenon/modulation symmetry which we elaborate in \textsc{Section \ref{subsub_res_obstr}}.

\subsubsection{Proof of \eqref{eq_iff_3else}: the non-linear ODE}\label{subsub_iff_ode}
Consider the non-linear ODE:
\begin{equation}\label{eq_ode}
    \vchi'_\vgamma\br{t}\wedge\vchi''_\vgamma\br{t}= 0,\quad \forall t\in I.
\end{equation}
Due to \eqref{eq_char_gamma_2_Gamma}, for the purpose of understanding the geometry of the solution to \eqref{eq_ode} under the assumption \(\vgamma'\wedge\vgamma''\not\equiv 0\), it suffices to consider the case \(\vgamma\br{t}=\br{t,\gamma\br{t}}\) to simplify the calculation. By direct computation,
\begin{equation*}
    0\not\equiv\vgamma'\wedge\vgamma''=\gamma'', \quad \vchi_\vgamma=\theta\cdot\br{1,\gamma'},\quad
    \theta:=\gamma'/\gamma''.
\end{equation*}
We thus deduce the following:
\begin{equation*}
    0\equiv
    \vchi'_\vgamma\wedge\vchi''_\vgamma
    =
    2\br{\theta'}^2\gamma''+\theta'\theta\gamma'''-\theta\theta''\gamma''.
\end{equation*}
Since \(\gamma''\not\equiv 0\), the above is equivalent to the following:
\begin{equation}\label{eq_ode_theta_gamma}
    2\br{\theta'}^2+\theta'\theta\frac{\gamma'''}{\gamma''}-\theta\theta''\equiv 0.
\end{equation}
Utilizing now the following identity:
\begin{equation*}
    \theta'= 1- \frac{\gamma'\gamma'''}{\br{\gamma''}^2}=1-\theta\frac{\gamma'''}{\gamma''},
\end{equation*}
we further deduce an equivalent formulation of \eqref{eq_ode_theta_gamma} involving only \(\theta\) and its derivatives:
\begin{equation}\label{eq_ode_theta}
    \theta'\br{\theta'+1}-\theta\theta''\equiv 0.
\end{equation}
As a consequence, two family of solutions to \eqref{eq_ode} arises from solving:
\begin{equation*}
    \theta'\equiv 0 \quad \text{and}\quad \theta' + 1\equiv 0.
\end{equation*}
Recall now \(\theta:=\gamma'/\gamma''\). Standard argument shows that there are \(A,a\neq 0\) and \(c\in\R\) such that:
\begin{equation*}
    \theta'\equiv 0\iff
    \gamma\br{t}=c_1e^{c_2 t}+c_3\text{ for some }c_1,c_2\neq 0,\, c_3\in\R;
\end{equation*}
\begin{equation*}
    \theta'+1\equiv 0 \iff
    \gamma\br{t}=c_1\log\abs{t-c_2}+c_3\text{ for some }c_1\neq 0,\, c_2,c_3\in\R,
\end{equation*}
In other words, we deduce the following geometric conditions:
\begin{equation*}
    \theta'\equiv 0\iff
    \img \vgamma \subseteq \eqref{eq_exp};\quad
    \theta'+1\equiv 0 \iff\img\vgamma \subseteq \eqref{eq_log}.
\end{equation*}
It remains to solve \eqref{eq_ode_theta} assuming \(\theta'\not\equiv 0\) and \(\theta'+1\not\equiv 0\). We may thus divide \eqref{eq_ode_theta} by \(\br{\theta'+1}^2\) and obtain:
\begin{equation*}
    0\equiv \frac{\theta'}{\theta'+1}-\frac{\theta\theta''}{\br{\theta'+1}^2}=\br{\frac{\theta}{\theta'+1}}'.
\end{equation*}
Standard argument solving for \(\theta\) yields:
\begin{equation*}
    \gamma'\br{t}/\gamma''\br{t}=:
    \theta\br{t}= c_1 +c_2 e^{t/c_1}\text{ for some } c_1,c_2\neq 0.
\end{equation*}
Finally, it remains to solve for \(\gamma\):
\begin{equation*}
    \gamma''\br{t} - \frac{\gamma'\br{t}}{
        c_1 +c_2 e^{t/c_1}
    }=0.
\end{equation*}
We multiply the above by the factor \(c_1 e^{-t/c_1}+c_2\) and deduce:
\begin{equation*}
    0= \br{c_1 e^{-t/c_1}+c_2} \gamma''\br{t}
    -e^{-t/c_1}\gamma'\br{t}
    = \br{\br{c_1 e^{-t/c_1}+c_2} \gamma'\br{t}}'.
\end{equation*}
As a result, we solve \(\gamma'\):
\begin{equation*}
    \gamma'\br{t}=\frac{c_3}{c_1 e^{-t/c_1}+c_2}= \frac{\frac{c_3}{c_2} e^{t/c_1}}{\frac{c_1}{c_2}+e^{t/c_1}} \text{ for some }c_1,c_2,c_3\neq 0
\end{equation*}
\begin{equation*}
    = \frac{c_3 e^{c_1t}}{e^{c_1 t}+c_2} \text{ for some other }c_1,c_2,c_3\neq 0
\end{equation*}
and thus
\begin{equation*}
    \gamma\br{t}= c_3\log\abs{e^{c_1t}+c_2}+c_4 
    \text{ for some }c_1,c_2,c_3\neq 0,\,c_4\in\R.
\end{equation*}
This concludes the last case:
\begin{equation*}
    \br{\frac{\theta}{\theta'+1}}'\equiv 0
    \iff \img \vgamma \subseteq \eqref{eq_psi_line}.
\end{equation*}

\subsubsection{Proof of sufficiency: reduction to \textbf{Theorem \ref{thm_smoothing_o_mini_uni_w_vanishing_n_fluc}}}\label{subsub_suffi_real_ana_smoothing}
Let \(\vgamma:I\to\R^2\) be real analytic.
Recall that it suffices to consider the case satisfying \eqref{eq_ana_bad_at_0}. Decompose now the domain into positive and negative parts. By avoiding the pathology at \(t=0\), condition \eqref{eq_ana_bad_at_0} implies the splitting of the domain produces two genuine real analytic parametric curves:
\begin{equation*}
    \vgamma_+:I_+:=\br{0,1}\simto \Gamma_+\subseteq \R^2,\quad
    \vgamma_-:I_-:=\br{-1,0}\simto \Gamma_-\subseteq \R^2
\end{equation*}
with both \(\cG_{\Gamma_+}\) and \(\cG_{\Gamma_-}\) injective. As a direct consequence of \textbf{Theorem \ref{thm_smoothing_o_mini_uni_w_vanishing_n_fluc}}\footnote{We take \(\Sigma:=\BR{0}\) to be singleton here.} and the definability of \(\vgamma_\pm\) in \(\R_\an\), the Sobolev smoothing inequality \eqref{eq_soblev_smoothing_std} holds for some \(0<\varepsilon,C_{\vgamma,\rho,\varepsilon}<\infty\) \textbf{as long as} \(0\notin\supp\rho\).

To address the \(0\in \supp \rho\) case, we recall the frequency decomposition \eqref{eq_cA_freq_decomp_L_for_var}, \eqref{eq_cA_freq_decomp_HL_for_var}, \eqref{eq_cA_freq_decomp_H_for_var}, the trivial bound \eqref{eq_cA_good_n_defect_trivial_est}, and the notation \eqref{eq_cR_4_var}.
Following the discussion in \textsc{Section \ref{subsec_pf_mainthm_global_w_vanishing_n_fluc}}.
It suffices to establish for \(\lambda:=\abs{\vlambda}\gg 1\) the following decay estimate:
\begin{equation}\label{eq_cA_ana_dec}
    \nrm{\cA_{\vgamma,\rho}\br{\fR{1}_{\lambda_1}f_1,\fR{2}_{\lambda_2}f_2}}_{L^1}\underset{\vgamma,\rho}{\lesssim}\lambda^{-c}\nrm{f_1}_{L^2}\nrm{f_2}_{L^2}. 
\end{equation}
Let \(\phi\) be as in \eqref{eq_phi_uni_even} and \(\delta>0\) to be decided later. We decompose \(\rho\) into three \(C^2\) bump functions:
\begin{equation*}
    \rho=\rho_{\lambda,-} + \rho_{\lambda,0}+ \rho_{\lambda,+}
\end{equation*}
with each given by:
\begin{equation*}
    \rho_{\lambda,0}\br{t}:=\phi\br{\lambda^\delta t}\rho\br{t},\quad
    \rho_{\lambda,\pm}:=\dif{\br{\rho-\rho_{\lambda,0}}}_{\R_\pm}.
\end{equation*}
On the one hand, we have the trivial bound:
\begin{equation}\label{eq_cA_ana_0}
    \nrm{\cA_{\vgamma,\rho_{\lambda,0}}\br{\fR{1}_{\lambda_1}f_1,\fR{2}_{\lambda_2}f_2}}_{L^1}
    \lesssim \nrm{\rho_{\lambda,0}}_{L^1}
    \nrm{f_1}_{L^2}
    \nrm{f_2}_{L^2}
    \underset{\rho}{\lesssim} \lambda^{-\delta}\nrm{f_1}_{L^2}
    \nrm{f_2}_{L^2}.
\end{equation}
On the other hand, since \(0\notin \supp \rho_{\lambda,\pm}\Subset I_\pm\), we recall \eqref{eq_admin_func_ana_bd}, calculate:
\begin{equation*}
    \dist\br{\supp_{\rho_{\lambda,\pm}}, I^c_\pm}\underset{\rho}{\gtrsim} \lambda^{-\delta},\quad
    \nrm{\rho_{\lambda,\pm}}_{C^2}\underset{\rho}{\lesssim}\lambda^{2\delta},
\end{equation*}
and apply \textbf{Theorem \ref{thm_smoothing_o_mini_uni_w_vanishing_n_fluc}} to obtain the following quantitative estimate:
\begin{equation*}
    \nrm{\cA_{\vgamma,\rho_{\lambda,\pm}}\br{\fR{1}_{\lambda_1}f_1,\fR{2}_{\lambda_2}f_2}}_{L^1}
    \underset{n,\vgamma}{\lesssim }
    \max\br{
        \nrm{A_{\vgamma_{\pm}}}_{L^{1/n,\infty}},
        \frac{1}{\dist\br{\supp_{\rho_{\lambda,\pm}}, I^c_\pm}}
    }^{c_0}
\end{equation*}
\begin{equation*}
    \cdot\ang{\abs{I_\pm}}
    \nrm{\rho_{\lambda,\pm}}_{C^2}
    \nrm{\fR{1}_{\lambda_1}f_1}_{H^{\br{-\varepsilon_{1/n,0},0}}}
    \nrm{\fR{2}_{\lambda_2}f_2}_{H^{\br{0,-\varepsilon_{1/n,0}}}}
\end{equation*}
\begin{equation}\label{eq_cA_ana_pm}
    \underset{n,\vgamma,\rho}{\lesssim}
    \lambda^{\delta \br{2+c_0} -\varepsilon_{1/n,0}}\nrm{f_1}_{L^2}\nrm{f_2}_{L^2}
\end{equation}
Taking \(\delta=\frac{\varepsilon_{1/n,0}}{3+c_0}>0\), \eqref{eq_cA_ana_0} and \eqref{eq_cA_ana_pm} yield \eqref{eq_cA_ana_dec} and conclude the proof.

\subsubsection{Proof of necessity: resonance obstructions}\label{subsub_res_obstr}

We claim that each of \eqref{eq_line}, \eqref{eq_exp}, \eqref{eq_log}, \eqref{eq_psi_line} corresponds to certain forms of resonance phenomena. The resonance for the \eqref{eq_line} case is quite well-known. See \cite{becker2024trilinearsingularbrascampliebintegrals,MR3482272}. However, for readers convenience, we sketch the standard argument below:

Assume that \(\img \vgamma \subseteq \eqref{eq_line}\). That is,
\begin{equation}\label{eq_line_vgamma}
    a\gamma_1\br{t}+b\gamma_2\br{t}=c,\quad \br{a,b}\neq \br{0,0},\;c\in\R.
\end{equation}
Let \(\phi\in C^\infty_c\br{\R}\) be given as in \eqref{eq_phi_uni_even}. Consider the standard linear modulation:
\begin{equation*}
    \Mod_\vxi f\br{\vx}:=e\br{\vxi\cdot \vx}f\br{\vx}.
\end{equation*}
Due to \eqref{eq_line_vgamma}, the following resonance relation holds
\begin{equation}\label{eq_line_res}
    e\br{c\lambda}\cdot\Mod_{\lambda\br{a,b}}\cA_{\vgamma,\rho}\br{f_1,f_2}=
    \cA_{\vgamma,\rho}\br{\Mod_{\br{a\lambda,0}}f_1,\Mod_{\br{0,b\lambda}}f_2}
\end{equation}
for all \(\lambda\in\R\). To see \eqref{eq_line_res} forbids \eqref{eq_soblev_smoothing_std}, we consider:
\begin{equation*}
    f_{1,\lambda}:=\Mod_{\br{a\lambda,0}}\big\vert \widehat{\phi}\big\vert^2\otimes \big\vert \widehat{\phi}\big\vert^2,\quad
    f_{2,\lambda}:=\Mod_{\br{0,b\lambda}}\big\vert \widehat{\phi}\big\vert^2\otimes \big\vert \widehat{\phi}\big\vert^2.
\end{equation*}
On the one hand, \eqref{eq_line_res} implies
\begin{equation*}
    1\eqsim 
    \nrm{\cA_{\vgamma,\rho}\br{\big\vert \widehat{\phi}\big\vert^2\otimes \big\vert \widehat{\phi}\big\vert^2,\big\vert \widehat{\phi}\big\vert^2\otimes \big\vert \widehat{\phi}\big\vert^2}}_{L^1}
    =\nrm{\cA_{\vgamma,\rho}\br{f_{1,\lambda},f_{2,\lambda}}}_{L^1}.
\end{equation*}
On the other hand, direct computation shows that
\begin{equation*}
    \nrm{f_{1,\lambda}}_{H^{\br{-\varepsilon,0}}}\lesssim \ang{a\lambda}^{-\varepsilon},\quad
    \nrm{f_{2,\lambda}}_{H^{\br{0,-\varepsilon}}}\lesssim \ang{b\lambda}^{-\varepsilon},\quad \forall \lambda\gg 1.
\end{equation*}
This with the condition \(\br{a,b}\neq \br{0,0}\) produces \eqref{eq_res_obstr} and forbids the Sobolev smoothing inequality \eqref{eq_soblev_smoothing_std} for any \(\varepsilon>0\).

As for the three remaining cases, since \eqref{eq_exp} and \eqref{eq_log} are mirror image to each other, it suffices to demonstrate the resonance obstructions for only the two cases: \eqref{eq_exp} and \eqref{eq_psi_line}.

We continue with the more symmetric \eqref{eq_psi_line} case. Assume that \(\img \vgamma\subseteq \eqref{eq_psi_line}\):
\begin{equation}\label{eq_psi_line_vgamma}
    ae^{A\gamma_1\br{t}}+be^{B\gamma_2\br{t}}=c,
    \quad a,b,A,B\neq 0,\,c\in\R.
\end{equation}
What follow are constructions inspired by \cite[\textbf{Lemma 2.3}]{MR4776384}: consider the following generalized modulation:
\begin{equation}\label{eq_exp_mod_joint}
    \Mod_{\lambda\exp_{A,B}}f\br{\vx}:=e\br{\lambda e^{Ax_1+B x_2}}f\br{\vx}.
\end{equation}
Due to \eqref{eq_psi_line_vgamma}, the following resonance relation holds:
\begin{equation}\label{eq_psi_line_res}
    \Mod_{c\lambda\exp_{A,B}}
    \cA_{\vgamma,\rho}\br{f_1,f_2}=\cA_{\vgamma,\rho}\br{\Mod_{a\lambda\exp_{A,B}}f_1,\Mod_{b\lambda\exp_{A,B}}f_2}
\end{equation}
for all \(\lambda\in\R\).
To see \eqref{eq_psi_line_res} causes the failure of \eqref{eq_soblev_smoothing_std}, we take \(\phi\in C^\infty_c\br{\R}\) to be given as in \eqref{eq_phi_uni_even} and consider:
\begin{equation}\label{eq_psi_line_constr}
    f_{1,\lambda}:=\Mod_{ a \lambda\exp_{A,B}}\phi\otimes\phi,\quad
    f_{2,\lambda}:= \Mod_{b\lambda\exp_{A,B}}\phi\otimes\phi.
\end{equation}
Observe that \eqref{eq_psi_line_res} implies:
\begin{equation*}
    1\eqsim \nrm{\cA_{\vgamma,\rho}\br{\phi\otimes\phi,\phi\otimes\phi}}_{L^1} = \nrm{\cA_{\vgamma,\rho}\br{f_{1,\lambda},f_{2,\lambda}}}_{L^1},\quad \forall \lambda\in\R.
\end{equation*}
On the other hand, we claim that
\begin{equation}\label{eq_psi_line_upper_bd}
    \nrm{f_{1,\lambda}}_{H^{\br{-\varepsilon,0}}},\,
    \nrm{f_{2,\lambda}}_{H^{\br{0,-\varepsilon}}}\lesssim \lambda^{-\varepsilon},\quad \forall \lambda\gg 1.
\end{equation}
In combination, the relation \eqref{eq_res_obstr} forbids any Sobolev smoothing inequality \eqref{eq_soblev_smoothing_std} for any \(\varepsilon>0\).

To demonstrate \eqref{eq_psi_line_upper_bd} and for the sake of discussion, we introduce the single variable variant of \eqref{eq_psi_line_constr}:
\begin{equation}\label{eq_phi_exp_mod}
    \phi_{S,\lambda}:=e\br{\lambda e^{Sx}}\phi\br{x}.
\end{equation}
Observe the two identities:
\begin{equation*}
    \nrm{f_{1,\lambda}}_{H^{\br{-\varepsilon,0}}}=\nrm{\nrm{\phi_{A,ae^{Bx_2}\lambda}}_{H^{-\varepsilon}}\phi\br{x_2}}_{L^2\br{dx_2}},\quad
    \nrm{f_{2,\lambda}}_{H^{\br{0,-\varepsilon}}}=\nrm{\phi\br{x_1}\nrm{\phi_{B,be^{Ax_1}\lambda}}_{H^{-\varepsilon}}}_{L^2\br{dx_1}}.
\end{equation*}
Using the fact that \(e^{Ax}, e^{Bx}\eqsim 1\) for \(x\in\br{-1,1}\), we reduces \eqref{eq_psi_line_upper_bd} to the following estimate:
\begin{equation}\label{eq_phi_S_lambda_sobo}
    \nrm{\phi_{S,\lambda}}_{H^{-\varepsilon}}:=\nrm{
    \ang{\xi}^{-\varepsilon} \widehat{\phi_{S,\lambda}}\br{\xi}
    }_{L^2\br{d\xi}}\underset{S}{\lesssim} \abs{\lambda}^{-\varepsilon}, \quad \forall \abs{\lambda}\gg 1.
\end{equation}
To proceed, we shall estimate the Fourier coefficient:
\begin{equation*}
    \widehat{\phi_{S,\lambda}}\br{\xi}=
    \int 
        e\br{\lambda e^{Sx}- \xi x}
        \phi\br{x}
    dx.
\end{equation*}
Suggested by the (non-)stationary phase principle, we consider two cases:
\begin{itemize}
    \item \(\min\br{\abs{\xi}, \abs{\lambda}} \underset{S}{\ll} \max\br{\abs{\xi},\lambda}\) the non-stationary phase case: standard integration by part argument yields:
    \begin{equation}\label{eq_phi_S_lamb_nst}
        \abs{\widehat{\phi_{S,\lambda}}\br{\xi}}\underset{S,N}{\lesssim}
        \br{\max\br{\abs{\xi},\abs{\lambda}}}^{-N}.
    \end{equation}
    \item \(\abs{\xi} \underset{S}{\eqsim} \lambda\) the stationary phase case: since there is no degenerate critical point, we have:
    \begin{equation}\label{eq_phi_S_lamb_st}
        \abs{\widehat{\phi_{S,\lambda}}\br{\xi}}\underset{S}{\lesssim}1/\sqrt{
            \inf_{x\in\br{-1,1}} \lambda \br{e^{Sx}-\xi x}''
        }
        \underset{S}{\eqsim}\abs{\lambda}^{-\frac{1}{2}}.
    \end{equation}
\end{itemize}
Combining \eqref{eq_phi_S_lamb_nst} and \eqref{eq_phi_S_lamb_st}, we conclude \eqref{eq_phi_S_lambda_sobo}:
\begin{equation*}
    \nrm{\phi_{S,\lambda}}_{H^{-\varepsilon}}
    \underset{S,N}{\lesssim}
    \nrm{\ang{\xi}^{-\varepsilon}\abs{\lambda}^{-N}}_{L^2\br{\abs{\xi}\ll \abs{\lambda}}}
    +
    \nrm{\ang{\xi}^{-\varepsilon}\abs{\lambda}^{-\frac{1}{2}}
    }_{L^2\br{\abs{\xi}\eqsim \abs{\lambda}}}
    +\nrm{
        \ang{\xi}^{-\varepsilon}
        \abs{\xi}^{-N}
    }_{L^2\br{\abs{\xi}\gg \abs{\lambda}}}
\end{equation*}
\begin{equation*}
    \underset{N,\varepsilon}{\lesssim}
    \abs{\lambda}^{\frac{1}{2}-N}+\abs{\lambda}^{-\varepsilon-\frac{1}{2}+\frac{1}{2}}+\abs{\lambda}^{-\varepsilon- N+\frac{1}{2}}
    \eqsim\abs{\lambda}^{-\varepsilon}.
\end{equation*}
This concludes the demonstration of the resonance obstruction for the \eqref{eq_psi_line} case.

The \eqref{eq_exp} case exhibits new resonance phenomenon that relies on the corner structure of \(\cA_{\vgamma,\rho}\). Let now \(\img\vgamma\subseteq \eqref{eq_exp}\):
\begin{equation}\label{eq_exp_vgamma}
    ae^{A\gamma_1\br{t}}+b\gamma_2\br{t}=c,\quad a,b,A\neq 0,\,c\in\R.
\end{equation}
Consider the following joint modulation:
\begin{equation*}
    \Mod_{\lambda,A}f\br{\vx}
    =e\br{\lambda x_2 e^{A x_1}}f\br{\vx}.
\end{equation*}
\eqref{eq_exp_vgamma} produces the following more elaborate resonance relation:
\begin{equation}\label{eq_exp_res}
    \Mod_{c\lambda\exp_{A,0}}\Mod_{b\lambda,A}\cA_{\vgamma,\rho}\br{f_1,f_2}
    =
    \cA_{\vgamma,\rho}\br{
        \Mod_{a\lambda\exp_{A,0}}f_1,
        \Mod_{b\lambda,A}f_2
    }
\end{equation}
for all \(\lambda\in\R\). We now consider:
\begin{equation}\label{eq_exp_res_constr}
    f_{1,\lambda}:= \Mod_{a\lambda\exp_{A,0}} \phi\otimes\phi,\quad
    f_{2,\lambda}:=\Mod_{b\lambda, A}\phi\otimes \vert\widehat{\phi}\vert^2.
\end{equation}
Similar to previous discussion on the \eqref{eq_psi_line} case, \eqref{eq_exp_res} implies:
\begin{equation*}
    1\eqsim\nrm{\cA_{\vgamma,\rho}\br{\phi\otimes\phi,\phi\otimes\vert\widehat{\phi}\vert^2}}_{L^1}
    =\nrm{\cA_{\vgamma,\rho}\br{f_{1,\lambda},f_{2,\lambda}}}_{L^1}.
\end{equation*}
Moreover, the analogous statement of \eqref{eq_psi_line_upper_bd} also holds. Indeed, on the one hand, \eqref{eq_phi_S_lambda_sobo} implies the following decay estimate:
\begin{equation*}
    \nrm{f_{1,\lambda}}_{H^{\br{-\varepsilon,0}}}=
    \nrm{\phi_{A,a\lambda}}_{H^{-\varepsilon}}\nrm{\phi}_{L^2}\lesssim\abs{\lambda}^{-\varepsilon},\quad \forall \abs{\lambda}\gg 1.
\end{equation*}
On the other hand, by design, we have:
\begin{equation*}
    \nrm{f_{2,\lambda}}_{H^{\br{0,-\varepsilon}}}
    =\nrm{
        \phi\br{x_1}
        \nrm{
            \Mod_{b\lambda e^{Ax_1}}
            \vert\widehat{\phi}\vert^2
        }_{H^{-\varepsilon}}
    }_{L^2\br{dx_1}}
\end{equation*}
\begin{equation*}
    \leq
    \nrm{
        \phi\br{x_1}
        \nrm{
            \ang{\xi}^{-\varepsilon}
            \1_{\br{-2,2}}\br{\xi-b\lambda e^{Ax_1}}
        }_{L^2\br{d\xi}}
    }_{L^2\br{dx_1}}
    \lesssim \abs{\lambda}^{-\varepsilon},\quad \forall \abs{\lambda}\gg 1.
\end{equation*}
In combination, we deduce \eqref{eq_res_obstr} and thus conclude that no Sobolev smoothing inequality \eqref{eq_soblev_smoothing_std} can hold for any \(\varepsilon>0\) in the \eqref{eq_exp} case.

\section{Applications}\label{sec_application}
Recall the following operators:
\begin{equation}\label{eq_THT_vgamma}
    T_\vgamma\br{f_1,f_2}\br{\vx}:=
    \int 
        \prod_{j=1,2}
            f_j\br{\vx+\gamma_j\br{t}\ve_j}
    \frac{dt}{t},
\end{equation}
% \begin{equation*}
%     A_\vgamma\br{f_1,f_2}\br{\vx}:=
%     \int^1_0 
%         \prod_{j=1,2}
%             f_j\br{\vx+\gamma_j\br{t}\ve_j}
%     dt,
% \end{equation*}
\begin{equation}\label{eq_M_vgamma}
    M_\vgamma\br{f_1,f_2}\br{\vx}:=
    \sup_{R>0}\abs{
    \fint^R_0
        \prod_{j=1,2}
            f_j\br{\vx+\gamma_j\br{t}\ve_j}
    dt},
\end{equation}
\begin{equation}\label{eq_SM_vgamma}
    \cM_\vgamma\br{f_1,f_2}\br{\vx}:=
    \sup_{k\in\Z}
    \abs{
    \int
        \prod_{j=1,2}
            f_j\br{\vx+2^k\gamma_j\br{t}\ve_j}
        \rho\br{t}
    dt}.
\end{equation}
By considering the standard partition of unity:
\begin{equation}\label{eq_part_R+}
    \1_{\R\setminus\BR{0}}=\sum_{k\in\Z}\Dil^\infty_{2^k}\rho=
    \sum_{k\in\Z}
    \sum_{\pm\in\BR{+,-}}
        \Dil^\infty_{2^k}\rho_\pm
    ,\quad \rho:=\Dil^{\infty}_2\phi-\phi
    ,\quad
    \rho_\pm:=\1_{\R_\pm}\rho,
\end{equation}
we obtain for \(f_1,f_2\in \mathscr{S}\br{\R^2}\) the following relations:
\begin{equation}\label{eq_SIO_2_A}
    T_\vgamma\br{f_1,f_2}\br{\vx}=
    \sum_{k\in\Z}
    \cA_{\Dil^\infty_{2^{-k}}\vgamma,\rho/t}\br{f_1,f_2}\br{\vx},
\end{equation}
% \begin{equation}\label{eq_A_rough_2_A}
%     A_\vgamma\br{\1_E,\1_E}\br{\vx}\geq
%     \sum_{k<0}
%     2^{k}
%     \cA_{\Dil^\infty_{2^{-k}}\vgamma,\rho_+}\br{\1_E,\1_E}\br{\vx},
% \end{equation}
\begin{equation}\label{eq_M_2_A}
    M_\vgamma\br{f_1,f_2}\br{\vx}\lesssim
    \sup_{k\in\Z}
    \abs{\cA_{\Dil^\infty_{2^{-k}}\vgamma,\rho_+}\br{f_1,f_2}\br{\vx}},
\end{equation}
\begin{equation}\label{eq_SM_2_A}
    \cM_\vgamma\br{f_1,f_2}\br{\vx}=
    \sup_{k\in\Z}
    \abs{\cA_{2^k\vgamma,\rho}\br{f_1,f_2}\br{\vx}}.
\end{equation}
The occurrence of the corner-type averaging operators suggests a common framework. Under suitable assumptions on \(\vgamma\), we expect the existence of a sequence of diagonal matrices \(\vM:=\BR{\vM_k}_k\):
\begin{equation*}
    \vM_k:=\begin{pmatrix}
    M_{k,1} & 0\\
    0 & M_{k,2}
    \end{pmatrix}
\end{equation*}
whose diagonal entries \(\BR{M_{k,1}}_k\) and \(\BR{M_{k,2}}_k\) are positive lacunary sequences:
\begin{equation}\label{eq_lacu_seqs_condi}
    1<M_\ast\leq \inf_k M_{k+1,j}/M_{k,j},\quad j=1,2
\end{equation}
such that the rescaled parametric curves:
\begin{equation}\label{eq_vgamma_rescaled_via_lac}
    \vgamma_k\br{t}:=\vM^{-1}_k \vgamma\br{2^k t}
\end{equation}
satisfy certain admissibility conditions: \textbf{Definition \ref{def_A_admi_ana}}, \eqref{eq_weak_failure}, \eqref{eq_A_gamma_fluc} in a uniform manner.
Regardless of the above assumptions, we reduce \eqref{eq_SIO_2_A}, \eqref{eq_M_2_A}, and \eqref{eq_SM_2_A} to the following two types of operators:
\begin{equation}\label{eq_SIO_M_models}
    \sum_k \cA_{\vM_k\vgamma_k,\rho}\br{f_1,f_2}\br{\vx},\quad
    \sup_k
    \abs{\cA_{\vM_k\vgamma_k,\rho}\br{f_1,f_2}\br{\vx}}
\end{equation}
for some sequence of \(\BR{\vgamma_k:I\to\R^2}_k\) that may or may not arise from \eqref{eq_vgamma_rescaled_via_lac} and some smooth function \(\rho\in C^\infty_c\br{I}\); depending on the context, \(\rho\) may be required to satisfy the mean zero condition \(\widehat{\rho}\br{0}=0\). From this point on, we shall first study the general form \eqref{eq_SIO_M_models} and later provide specific context to address the three operators \(T_\vgamma,M_\vgamma,\cM_\vgamma\) individually.

% To begin with, we recall \eqref{eq_soblev_smoothing_std} and introduce for \(\BR{\vgamma_k:I\to\R^2}_k\) and \(\rho\in C^\infty_c\br{I}\) the following two quantities
% \begin{equation}\label{eq_cond_4_Hparts}
%     B_0:=\sup_k
%     \nrm{\vgamma_k}_{L^\infty\br{I}}
%     ,\quad
%     C_\varepsilon:=\sup_k C_{\vgamma_k,\rho,\varepsilon}.
% \end{equation}

Fix \(k\in \Z\) for a brief moment. We aim to perform suitable frequency decomposition on \(\cA_{\vM_k\vgamma_k,\rho}\br{f_1,f_2}\) adapted to the scaling given by \(\vM_k\).
We recall \textbf{Definition \ref{def_LP_decomp}}, fix a small \(\delta=1/5\), omit the \(\delta\)-dependency.
Suggested by the dilation relation:
\begin{equation}\label{eq_A_Mk_2_Dil_A}
    \cA_{\vM_k\vgamma_k,\rho}\br{f_1,f_2}
    =\Dil^{p_3}_{\vM_k}\cA_{\vgamma_k,\rho}\br{
        \Dil^{p_1}_{\vM^{-1}_k}f_1,
        \Dil^{p_2}_{\vM^{-1}_k}f_2
    },\quad
    \Dil^p_{\vM}f\br{\vx}:=\abs{\det\br{\vM}}^{-\frac{1}{p}}f\br{\vM^{-1}\vx},
\end{equation} 
we perform the frequency decomposition:
\begin{equation*}
    \cA_{\vM_k\vgamma_k,\rho}\br{f_1,f_2}=
    \frac{1}{c^2_{1/5}}
    \iint
        \Dil^{p_3}_{\vM_k}
        \cA_{\vgamma_k,\rho}
        \br{\fP{1}_{\lambda_1} \Dil^{p_1}_{\vM^{-1}_k}f_1,
        \fP{2}_{\lambda_2} \Dil^{p_2}_{\vM^{-1}_k}f_2}
    \frac{d\lambda_1}{\abs{\lambda_1}}
    \frac{d\lambda_2}{\abs{\lambda_2}}
\end{equation*}
\begin{equation*}
    =
    \frac{1}{c^2_{1/5}}
    \iint
    \cA_{\vM_k\vgamma_k,\rho}\br{\fP{1}_{\lambda_1/M_{k,1}} f_1,\fP{2}_{\lambda_2/M_{k,2}} f_2}
    \frac{d\lambda_1}{\abs{\lambda_1}}
    \frac{d\lambda_2}{\abs{\lambda_2}}
\end{equation*}
We now consider the two frequency regions:
\begin{equation}
    H:=\BR{
        \vlambda\in\br{\R\setminus\BR{0}}^2
    \::\:
        \abs{\lambda_1},\abs{\lambda_2}\geq 1
    },\quad
    L:=\br{\R\setminus\BR{0}}^2\setminus H
\end{equation}
and define the corresponding frequency projected averages:
\begin{equation}
    \cA^U_{\vM_k\vgamma_k,\rho}\br{f_1,f_2}
    =
    \iint_U
    \cA_{\vM_k\vgamma_k,\rho}\br{\fP{1}_{\lambda_1/M_{k,1}} f_1,\fP{2}_{\lambda_2/M_{k,2}} f_2}
    \frac{d\lambda_1}{\abs{\lambda_1}}
    \frac{d\lambda_2}{\abs{\lambda_2}},\quad
    U\in\BR{H,L}.
\end{equation}
We obtain for \eqref{eq_SIO_M_models} the following two preliminary estimates:
\begin{equation*}
    \abs{\sum_k \cA_{\vM_k\vgamma_k,\rho}\br{f_1,f_2}\br{\vx}}
    \lesssim 
    \abs{\sum_k \cA^L_{\vM_k\vgamma_k,\rho}\br{f_1,f_2}\br{\vx}}
    +
    \iint_H
    \nrm{
    \cA_{\vM_k\vgamma_k,\rho}\br{\fP{1}_{\lambda_1/M_{k,1}} f_1,\fP{2}_{\lambda_2/M_{k,2}} f_2}\br{\vx}
    }_{\ell^1\br{k}}
    \frac{d\lambda_1}{\abs{\lambda_1}}
    \frac{d\lambda_2}{\abs{\lambda_2}};
\end{equation*}
\begin{equation*}
    \sup_k\abs{ \cA_{\vM_k\vgamma_k,\rho}\br{f_1,f_2}\br{\vx}}
    \lesssim 
    \sup_k
    \abs{
        \cA^L_{\vM_k\vgamma_k,\rho}\br{f_1,f_2}\br{\vx}
    }
    +
    \iint_H
    \nrm{
    \cA_{\vM_k\vgamma_k,\rho}\br{\fP{1}_{\lambda_1/M_{k,1}} f_1,\fP{2}_{\lambda_2/M_{k,2}} f_2}\br{\vx}
    }_{\ell^1\br{k}}
    \frac{d\lambda_1}{\abs{\lambda_1}}
    \frac{d\lambda_2}{\abs{\lambda_2}}.
\end{equation*}
We claim that the common term shown in the above two estimates admits \(L^{p_1}\times L^{p_2}\to L^{p_3}\) boundedness for all \(p_1,p_2\in\br{1,\infty}\) and \(\frac{1}{p_3}=\frac{1}{p_1}+\frac{1}{p_2}\leq 1\), once we verify the following quantities related to \eqref{eq_soblev_smoothing_std}:

\begin{equation}\label{eq_cond_4_Hparts}
    B_N:=\sup_{\substack{
    0\leq \beta\leq N \\k
    }}
    \nrm{\vgamma^{\br{\beta}}_k}_{L^\infty\br{\supp\rho}}
    ,\quad
    C_\varepsilon:=\sup_k C_{\vgamma_k,\rho,\varepsilon}
\end{equation}
are finite and positive for some \(\varepsilon>0\). This requires two key statements which we state below:
\begin{proposition}\label{prop_HF_tame}
Let \(B_N\) be defined as in \eqref{eq_cond_4_Hparts}. We have the following estimate:
    \begin{equation*}
        \nrm{
            \nrm{
                \cA_{\vM_k\vgamma_k,\rho}
                \br{\fP{1}_{\lambda_1/M_{k,1}}f_1,
                \fP{2}_{\lambda_2/M_{k,2}}f_2}
            }_{\ell^1\br{k}}
        }_{L^{p_3}}
        \underset{p_j,M_\ast}{\lesssim}  
        \log^4
        \br{e+B_0\abs{\vlambda}}
        \nrm{f_1}_{L^{p_1}}\nrm{f_2}_{L^{p_2}}
    \end{equation*}
    for the following range
    \begin{equation}\label{eq_SIO_Lp_range}
        1<p_1,p_2<\infty,\quad\frac{1}{p_3}=\frac{1}{p_1}+\frac{1}{p_2}\leq 1.
    \end{equation}
\end{proposition}
\begin{proposition}\label{prop_HF_dec}
Let \(C_\varepsilon\) be defined as in \eqref{eq_cond_4_Hparts}. We have the following estimate:
\begin{equation*}
    \nrm{
        \nrm{
            \cA_{\vM_k\vgamma_k,\rho}
            \br{\fP{1}_{\lambda_1/M_{k,1}}f_1,
            \fP{2}_{\lambda_2/M_{k,2}}f_2}
        }_{\ell^1\br{k}}
    }_{L^1}
    \underset{M_\ast}{\lesssim} C_\varepsilon
    \ang{\abs{\vlambda}}^{-\varepsilon}
    \nrm{f_1}_{L^2}\nrm{f_2}_{L^2}.
\end{equation*}
\end{proposition}
Indeed, when \(0<\varepsilon,B_0,C_\varepsilon<\infty\), we derive from interpolating \textbf{Proposition \ref{prop_HF_dec} and \ref{prop_HF_tame}} the following:
\begin{equation*}
    \nrm{
        \nrm{
            \cA_{\vM_k\vgamma_k,\rho}
            \br{\fP{1}_{\lambda_1/M_{k,1}}f_1,
            \fP{2}_{\lambda_2/M_{k,2}}f_2}
        }_{\ell^1\br{k}}
    }_{L^{p_3}}
    \underset{C_\varepsilon,B_0,p_j,M_\ast}{\lesssim}
    \ang{\abs{\vlambda}}^{-\varepsilon_{p_1,p_2}}
    \nrm{f_1}_{L^{p_1}}\nrm{f_2}_{L^{p_2}}
\end{equation*}
for all \(p_1,p_2\in\br{1,\infty}\) with \(\frac{1}{p_3}=\frac{1}{p_1}+\frac{1}{p_2}\leq 1\) and some \(\varepsilon_{p_1,p_2}>0\). As a direct consequence, we conclude:
\begin{equation*}
    \nrm{
        \iint_H
        \nrm{
        \cA_{\vM_k\vgamma_k,\rho}\br{\fP{1}_{\lambda_1/M_{k,1}} f_1,\fP{2}_{\lambda_2/M_{k,2}} f_2}
        }_{\ell^1\br{k}}
        \frac{d\lambda_1}{\abs{\lambda_1}}
        \frac{d\lambda_2}{\abs{\lambda_2}}
    }_{L^{p_3}}
    \underset{C_\varepsilon,B_0,p_j,M_\ast}{\lesssim}
    \nrm{f_1}_{L^{p_1}}
    \nrm{f_2}_{L^{p_2}}.
\end{equation*}
The treatment of the remaining two terms:
\begin{equation*}
    \abs{\sum_k \cA^L_{\vM_k\vgamma_k,\rho}\br{f_1,f_2}\br{\vx}},\quad
    \sup_k\abs{ \cA^L_{\vM_k\vgamma_k,\rho}\br{f_1,f_2}\br{\vx}}
\end{equation*}
requires specifications on \(\vM_k,\vgamma_k,\rho\) depending on the associated contexts \eqref{eq_SIO_2_A}, \eqref{eq_M_2_A}, and \eqref{eq_SM_2_A}.
% We summarize the results below:

\begin{proposition}\label{prop_LF_SIO}
    Given the following list of conditions:
    \begin{enumerate}
        \item \(\rho\) has zero mean;
        \item \textbf{Conjecture \ref{conj_twist_vM}} holds for \(\vM:=\BR{\vM_k}_k\);
        \item \(B_{N_\vM+2}\) is finite\footnote{Recall \(N_\vM\in\N\sqcup\BR{0}\) in the statement of \textbf{Conjecture \ref{conj_twist_vM}}.};
        \item there is \(A>0\) such that \(1/A\leq \abs{\gamma'_{k,j}\br{t}}\) for all \(t\in I\), \(k\in\Z\), and \(j=1,2\),
        % \item there is \(\valpha=\br{\alpha_1,\alpha_2}\in\br{0,\infty}^2\) such that \(M_{k,j}=2^{k\alpha_j}\) for all \(k\in\Z\) and \(j=1,2\),
    \end{enumerate}
    the following estimates hold:
    \begin{equation*}
        \nrm{
            \sum_k \cA^L_{\vM_k\vgamma_k,\rho}\br{f_1,f_2}\br{\vx}
        }_{L^{p_3}\br{d\vx}}
        \underset{A,B_{N_\vM+2},M_\ast,p_j}{\lesssim}
        \nrm{f_1}_{L^{p_1}}
        \nrm{f_2}_{L^{p_2}}
    \end{equation*}
    for the range stated as in \eqref{eq_THT_Lp_bdds}.
\end{proposition}
\begin{remark}\label{rmk_Conj_holds_for_good_vM}
    In particular, \textbf{Conjecture \ref{conj_twist_vM}} holds for sequences \(\vM\) given by
\begin{equation*}
    \vM_k:= 
    \begin{pmatrix}
        2^{k \alpha_1 } /\lambda_1 & 0\\
        0 & 2^{k \alpha_2} /\lambda_2
    \end{pmatrix},
    % \quad
    % \text{ for some fixed }
    % \valpha:=\br{\alpha_1,\alpha_2},\vlambda:=\br{\lambda_1,\lambda_2} \in\br{0,\infty}^2
\end{equation*}
for some fixed \(\valpha:=\br{\alpha_1,\alpha_2},\vlambda:=\br{\lambda_1,\lambda_2} \in\br{0,\infty}^2\). Our main result \textbf{Theorem \ref{thm_THT_asym_homo_curves}} for \eqref{eq_THT_vgamma} only uses this specific case. See \eqref{eq_vM_hold_for_twist} and its surrounding discussion.
\end{remark}

\begin{proposition}\label{prop_LF_M}
    Given the following list of conditions:
    \begin{enumerate}
        \item \(B_0\) is finite;
        \item there is \(A>0\) such that \(1/A\leq \abs{\gamma'_{k,j}\br{t}}\) for all \(t\in I\), \(k\in\Z\), and \(j=1,2\),
    \end{enumerate}
    the following estimates hold:
    \begin{equation*}
        \nrm{
            \sup_k \abs{
            \cA^L_{\vM_k\vgamma_k,\rho}\br{f_1,f_2}}\br{\vx}
        }_{L^{p_3}\br{d\vx}}
        \underset{A,B_0,\vM_\ast,p_j}{\lesssim}
        \nrm{f_1}_{L^{p_1}}
        \nrm{f_2}_{L^{p_2}}
    \end{equation*}
    for the following range:
    \begin{equation}\label{eq_maximal_op_range}
        1<p_1,p_2\leq \infty,\quad
        \frac{1}{p_3}=\frac{1}{p_1}+\frac{1}{p_2}.
    \end{equation}
\end{proposition}

\begin{proposition}\label{prop_LF_SM}
    Given that \(\vgamma\) is real analytic, and both \(\gamma_1,\gamma_2\) are non-constant,
    % \begin{enumerate}
    %     \item \(\vgamma\) is real analytic and not constant;
    %     \item \(\vM:=\BR{\vM_k}_k\) is of the form \eqref{eq_vM_valpha_scaling},
    % \end{enumerate}
    the following estimates hold:
    \begin{equation*}
        \nrm{
            \sup_k \abs{
            \cA^L_{\vM_k\vgamma,\rho}\br{f_1,f_2}}\br{\vx}
        }_{L^{p_3}\br{d\vx}}
        \underset{\vgamma,M_\ast,p_j}{\lesssim}
        \nrm{f_1}_{L^{p_1}}
        \nrm{f_2}_{L^{p_2}}
    \end{equation*}
    for the range as in \eqref{eq_maximal_op_range}.
\end{proposition}

\subsection{High-frequency components: proof of Proposition \ref{prop_HF_tame} and \ref{prop_HF_dec}}

We begin with the proof of \textbf{Proposition \ref{prop_HF_tame}}. Recall \eqref{eq_phi_uni_even}, \eqref{eq_psi_delta}, and \textbf{Definition \ref{def_LP_decomp}}.
By the fact that:
\begin{equation*}
    \cP_{\delta,\lambda} f\br{x}=
    \cP_{2\delta,\lambda}
    \cP_{\delta,\lambda}
    f\br{x}=
    \int
    \cP_{\delta,\lambda}
    f\br{x+y}\lambda\widehat{\psi_{2\delta}}\br{\lambda y}
    dy,
\end{equation*}
we deduce the following identity:
\begin{equation}\label{eq_HF_tame_avg_id}
    \cA_{\vM_k\vgamma_k,\rho}\br{\fP{1}_{\lambda_1/M_{k,1}}f_1,\fP{2}_{\lambda_2/M_{k,2}}f_2}\br{\vx}=
    \int
        \prod_{j=1,2}
        \fP{j}_{\lambda_j/M_{k,j}}f_j\br{\vx+y_j\ve_j}
        K_{k,\vlambda}\br{\vy}
    d\vy
\end{equation}
with the kernel defined as below:
\begin{equation*}
    K_{k,\vlambda}\br{\vy}
    :=
    \int
    \br{
        \prod_{j=1,2}
        \frac{\lambda_j}{M_{k,j}}
        \widehat{\psi_{2\delta}}\br{
            \frac{\lambda_j \br{y_j-M_{k,j}\gamma_{k,j}\br{t}}}{M_{k,j}}
        }
    }
    \rho\br{t}
    dt.
\end{equation*}
Using now the trivial estimate \(\abs{\widehat{\psi_{2\delta}}\br{y}}\underset{N}{\lesssim} \ang{y}^{-N}\), we deduce the following point-wise bound
\begin{equation}\label{eq_HF_tame_avg_K_rel}
    \abs{K_{k,\vlambda}\br{\vy}}
    \underset{N}{\lesssim}
    \int
    \br{
        \prod_{j=1,2}
        \Dil^1_{M_{k,j}/\lambda_j}
        \Tr_{\lambda_j \gamma_{k,j}\br{t}}
        \ang{y_j}^{-N}
    }
        \rho\br{t}
    dt.
\end{equation}
With the relation \eqref{eq_HF_tame_avg_K_rel}, we derive the following
\begin{equation}\label{eq_HF_tame_avg_rel_k}
    \abs{\eqref{eq_HF_tame_avg_id}}\underset{N}{\lesssim}
    \int
        \prod_{j=1,2}
        \int
            \abs{
            \fP{j}_{\lambda_j/M_{k,j}}f_j\br{\vx
            +\frac{M_{k,j}}{\lambda_j}\br{y_j+\lambda_j\gamma_{k,j}\br{t}}
            }
            }
            \ang{y_j}^{-N}
        d y_j
    \abs{\rho\br{t}}
    dt.
\end{equation}
This suggests that we introduce a variant of shifted maximal operator \eqref{eq_MN_lac_bd} to further simplify the above expression. More precisely, by setting
\begin{equation}\label{eq_HF_tame_avg_Lambda_Sigma_settings}
    \Lambda_j:=\BR{M_{k,j}/\lambda_j}_k,\quad
    \Sigma_j\br{t}:=\BR{\lambda_j\gamma_{k,j}\br{t}}_k,
\end{equation}
we may reformulate and simplify \eqref{eq_HF_tame_avg_rel_k} into:
\begin{equation*}
    \cA_{\vM_k\vgamma_k,\rho}\br{\fP{1}_{\lambda_1/M_{k,1}}f_1,\fP{2}_{\lambda_2/M_{k,2}}f_2}\br{\vx}
    \lesssim
    \nrm{
        \rho\br{t}
        \prod_{j=1,2}
        \fM{j}_{\Lambda_j,\Sigma_j\br{t}}
            \fP{j}_{\lambda_j/M_{k,j}}f_j
        \br{\vx}
    }_{L^1\br{dt}}
    .
\end{equation*}
Take now the \(\ell^1\br{k}\)-norm, apply Cauchy-Schwarz inequality
\begin{equation}\label{eq_HF_tame_avg_rel_ell1}
    \nrm{
        \cA_{\vM_k\vgamma_k,\rho}\br{\fP{1}_{\lambda_1/M_{k,1}}f_1,\fP{2}_{\lambda_2/M_{k,2}}f_2}\br{\vx}
    }_{\ell^1\br{k}}
    \lesssim
    \nrm{
        \rho\br{t}
        \prod_{j=1,2}
        \nrm{
        \fM{j}_{\Lambda_j,\Sigma_j\br{t}}
            \fP{j}_{\lambda_j/M_{k,j}}f_j
        \br{\vx}
        }_{\ell^2\br{k}}
    }_{L^1\br{dt}},
\end{equation}
and finally apply \(L^{p_3}\br{d\vx}\)-norm and H\"{o}lder's inequality
\begin{equation*}
    \nrm{
        \eqref{eq_HF_tame_avg_rel_ell1}
    }_{L^{p_3}\br{d\vx}}
    \lesssim
    \nrm{
        \rho\br{t}
        \prod_{j=1,2}
        \nrm{
            \nrm{
                \fM{j}_{\Lambda_j,\Sigma_j\br{t}}
                    \fP{j}_{\lambda_j/M_{k,j}}f_j
                \br{\vx}
            }_{\ell^2\br{k}}
        }_{L^{p_j}\br{d\vx}}
    }_{L^1\br{dt}}.
\end{equation*}
Applying \textbf{Proposition \ref{prop_FS_MN_lac_bd}} and Littlewood-Paley theory, we dominate the above with
\begin{equation*}
    \lesssim
    \nrm{
        \rho\br{t}
        \prod_{j=1,2}
        \log^2_{M_\ast}\br{e+\nrm{\lambda_j\gamma_{k,j}\br{t}}_{\ell^\infty\br{k}}}
        \nrm{
            \nrm{\fP{j}_{\lambda_j/M_{k,j}}f_j}_{\ell^2\br{k}}
        }_{L^{p_j}}
    }_{L^1\br{t}}
\end{equation*}
\begin{equation*}
    \lesssim
        \prod_{j=1,2}
        \log^2_{M_\ast}\br{e+\abs{\vlambda}B_0}
        \nrm{
            \nrm{\fP{j}_{\lambda_j/M_{k,j}}f_j}_{\ell^2\br{k}}
        }_{L^{p_j}}
    \underset{M_\ast,p_j}\lesssim
    \log^4\br{e+\abs{\vlambda}B_0}
    \nrm{f_1}_{L^{p_1}}
    \nrm{f_2}_{L^{p_2}}.
\end{equation*}
This concludes the proof of \textbf{Proposition \ref{prop_HF_tame}}.

To prove \textbf{Proposition \ref{prop_HF_dec}}, we utilize the identity \eqref{eq_A_Mk_2_Dil_A} and the following commutation rule:
\begin{equation}\label{eq_HF_dec_LP_Dil_com}
    \Dil^2_{\vM^{-1}_k}\fP{j}_{\lambda_j/M_{k,j}}f_j=
    \fP{j}_{\lambda_j}
    \Dil^2_{\vM^{-1}_k} f_j
\end{equation}
to deduce via \eqref{eq_soblev_smoothing_std} the following chain of relations:
\begin{equation*}
    \nrm{
        \cA_{\vM_k\vgamma_k,\rho}
            \br{\fP{1}_{\lambda_1/M_{k,1}}f_1,
            \fP{2}_{\lambda_2/M_{k,2}}f_2}
    }_{L^1}
    =
    \nrm{
        \cA_{\vgamma_k,\rho}
        \br{
            \bigotimes_{j=1,2}
            \Dil^2_{\vM^{-1}_k}
            \fP{j}_{\lambda_j/M_{k,j}}f_j
        }
    }_{L^1}
\end{equation*}
\begin{equation*}
    =
    \nrm{
        \cA_{\vgamma_k,\rho}
        \br{
            \bigotimes_{j=1,2}
            \fP{j}_{\lambda_j}
            \Dil^2_{\vM^{-1}_k} f_j
        }
    }_{L^1}
    \leq 
    C_{\vgamma_k,\rho,\varepsilon}
    \nrm{\fP{1}_{\lambda_1}
        \Dil^2_{\vM^{-1}_k} f_1}_{H^{\br{-\varepsilon,0}}}
    \nrm{
        \fP{2}_{\lambda_2}
            \Dil^2_{\vM^{-1}_k} f_2
    }_{H^{\br{0,-\varepsilon}}}
\end{equation*}
\begin{equation*}
    \lesssim
    C_\varepsilon
    \ang{\abs{\vlambda}}^{-\varepsilon}
    \prod_{j=1,2}
        \nrm{\fP{j}_{\lambda_j}
            \Dil^2_{\vM^{-1}_k} f_j}_{L^2}
    =
    C_\varepsilon
    \ang{\abs{\vlambda}}^{-\varepsilon}
    \prod_{j=1,2}
        \nrm{
        \fP{j}_{\lambda_j/M_{k,j}}
        f_j}_{L^2}.
\end{equation*}
We finish the proof of \textbf{Proposition \ref{prop_HF_dec}} by taking \(\ell^1\br{k}\)-norm
\begin{equation*}
    \nrm{
        \nrm{
            \cA_{\vM_k\vgamma_k,\rho}
                \br{\fP{1}_{\lambda_1/M_{k,1}}f_1,
                \fP{2}_{\lambda_2/M_{k,2}}f_2}
        }_{\ell^1\br{k}}
    }_{L^1},
\end{equation*}
applying Cauchy-Schwarz, and invoking orthogonality
\begin{equation*}
    \lesssim
    C_\varepsilon
    \ang{\abs{\vlambda}}^{-\varepsilon}
    \prod_{j=1,2}
    \nrm{
        \nrm{
        \fP{j}_{\lambda_j/M_{k,j}}
        f_j}_{L^2}
    }_{\ell^2\br{k}}
    \underset{M_\ast}{\lesssim}
    C_\varepsilon
    \ang{\abs{\vlambda}}^{-\varepsilon}
    \nrm{f_1}_{L^2}
    \nrm{f_2}_{L^2}.
\end{equation*}

\subsection{Low-frequency components: proof of Proposition \ref{prop_LF_SIO}, \ref{prop_LF_M}, and \ref{prop_LF_SM}}
We start with the treatment of \textbf{Proposition \ref{prop_LF_SIO}}. Observe the following identity
\begin{equation}\label{eq_LF_SIO_L}
    \cA^L_{\vM_k\vgamma_k,\rho}\br{f_1,f_2}
    =
    \frac{1}{c_{1/5}}
    \int_{\abs{\lambda}\leq 1}
        \cA_{\vM_k\vgamma_k,\rho}\br{\fQ{1}_{\lambda/M_{k,1}} f_1,\fP{2}_{\lambda /M_{k,2}} f_2}
        +
        \cA_{\vM_
        k\vgamma_k,\rho}\br{\fP{1}_{\lambda/M_{k,1}} f_1,\fQ{2}_{\lambda /M_{k,2}} f_2}
    \frac{d\lambda}{\abs{\lambda}}
\end{equation}
\begin{equation}\label{eq_LF_SIO_H}
    +
    \frac{1}{c_{1/5}}
    \int_{\abs{\lambda}> 1}
        \cA_{\vM_k\vgamma_k,\rho}\br{\fQ{1}_{1/M_{k,1}} f_1,\fP{2}_{\lambda /M_{k,2}} f_2}
        +
        \cA_{\vM_k\vgamma_k,\rho}\br{\fP{1}_{\lambda/M_{k,1}} f_1,\fQ{2}_{1/M_{k,2}} f_2}
    \frac{d\lambda}{\abs{\lambda}}.
\end{equation}
It suffices to prove the following two estimates:
\begin{equation}\label{eq_LF_SIO_L_reduc}
    \nrm{
        \sum_k \cA_{\vM_k\vgamma_k,\rho}
        \br{\fQ{1}_{\lambda_1/M_{k,1}} f_1,\fP{2}_{\lambda_2 /M_{k,2}} f_2}
    }_{L^{p_3}}
    \underset{B_1,M_\ast,p_j}{\lesssim}
    \abs{\vlambda}\cdot
    \nrm{f_1}_{L^{p_1}}
    \nrm{f_2}_{L^{p_2}},\quad \abs{\vlambda}\lesssim A B_1;
\end{equation}
\begin{equation}\label{eq_LF_SIO_H_reduc}
    \nrm{
        \sum_k \cA_{\vM_k\vgamma_k,\rho}
        \br{\fQ{1}_{1/M_{k,1}} f_1,\fP{2}_{\lambda /M_{k,2}} f_2}
    }_{L^{p_3}}
    \underset{\substack{A,B_{N_\vM+2},\\N_\vM,M_\ast,p_j}}{\lesssim}
    \abs{\lambda}^{-1}\cdot
    \nrm{f_1}_{L^{p_1}}
    \nrm{f_2}_{L^{p_2}},\quad \abs{\lambda}\gg A B_1
\end{equation}
for the range stated as in \eqref{eq_THT_Lp_bdds}.
Indeed, by symmetry, \eqref{eq_LF_SIO_L_reduc} and \eqref{eq_LF_SIO_H_reduc} imply the desired estimate\footnote{Due to the integration \(\int \cdots \frac{d\lambda}{\lambda}\) in \eqref{eq_LF_SIO_L} and \eqref{eq_LF_SIO_H}, our proof do not cover the full \(L^p\) range \eqref{eq_T_m_twisted_multiplier_cond}. However, it still suffices for the application we stated.}
\begin{equation*}
    \nrm{
        \sum_k
            \cA^L_{\vM_k\vgamma_k,\rho}\br{f_1,f_2}
    }_{L^{p_3}}
    \underset{\substack{A,B_{N_\vM+2},\\N_\vM,M_\ast,p_j}}{\lesssim}
    \br{\int_{\abs{\lambda}\leq 1} \abs{\lambda}
    \cdot
    \frac{d\lambda}{\abs{\lambda}}+
    \int_{\abs{\lambda}>1}
    \abs{\lambda}^{-1}
    \frac{d\lambda}{\abs{\lambda}}
    }
    \nrm{f_1}_{L^{p_1}}
    \nrm{f_2}_{L^{p_2}}\eqsim 
    \nrm{f_1}_{L^{p_1}}
    \nrm{f_2}_{L^{p_2}}.
\end{equation*}

To show \eqref{eq_LF_SIO_L_reduc} and \eqref{eq_LF_SIO_H_reduc}, we express both sums on the left-hand-side as operators of the form \eqref{eq_T_m_twisted_para} associated with symbol
\begin{equation}\label{eq_symb_seq_struct_lambdas}
    m_\vlambda\br{\vxi}:=
    \sum_k
        \int
            e\br{\vxi^\top \vLambda \vM_k \vgamma_k\br{t}}
            \rho\br{t}
        dt
        \cdot 
        \varphi\otimes\psi
        \br{
            \vM_k
            \vxi
        },\quad
    \vLambda:=
    \begin{pmatrix}
        \lambda_1 & 0\\
        0 & \lambda_2
    \end{pmatrix},
\end{equation}
where, more precisely, we have the following identities via formula \eqref{eq_T_m_twisted_para} and \eqref{eq_T_m_twisted_rescaling_id}
\begin{equation}\label{eq_T_m_vlambda}
    \sum_k \cA_{\vM_k\vgamma_k,\rho}
        \br{\fQ{1}_{\lambda_1 /M_{k,1}} f_1,\fP{2}_{\lambda_2 /M_{k,2}} f_2}
    =T_{\Dil^\infty_\vLambda m_\vlambda}\br{f_1,f_2}=
    \Dil^{p_3}_{\vLambda^{-1}}T_{m_\vlambda}\br{\Dil^{p_1}_\vLambda f_1, \Dil^{p_2}_\vLambda f_2}.
\end{equation}
\eqref{eq_LF_SIO_L_reduc} and \eqref{eq_LF_SIO_H_reduc} can thus be further rephrased respectively as
\begin{equation}\label{eq_LF_SIO_L_T_m_est}
    \nrm{T_{m_\vlambda}\br{f_1,f_2}}_{L^{p_3}}\underset{B_1,M_\ast,p_j}{\lesssim}\abs{\vlambda}\cdot \nrm{f_1}_{L^{p_1}}\nrm{f_2}_{L^{p_2}},\quad \abs{\vlambda}\lesssim AB_1;
\end{equation}
\begin{equation}\label{eq_LF_SIO_H_T_m_est}
    \nrm{T_{m_{\br{1,\lambda}}}\br{f_1,f_2}}_{L^{p_3}}\underset{\substack{A,B_{N_\vM+2},\\N_\vM,M_\ast,p_j}}{\lesssim}\abs{\lambda}^{-1}\cdot \nrm{f_1}_{L^{p_1}}\nrm{f_2}_{L^{p_2}},\quad \abs{\lambda}\gg AB_1.
\end{equation}
Therefore, under the condition that \textbf{Conjecture \ref{conj_twist_vM}} holds for \(\vM:=\BR{\vM_k}_k\), it suffices to check that \(m_\vlambda\) and \(m_{\br{1,\lambda}}\) from \eqref{eq_LF_SIO_L_T_m_est} and \eqref{eq_LF_SIO_H_T_m_est} above can be formulated in terms of \eqref{eq_symb_w_seq_struct} with respective coefficients \(\vepsilon:=\BR{\epsilon_k}_k\) of desired \(\ell^\infty\) sizes. 

Starting with \eqref{eq_LF_SIO_L_T_m_est}, we consider the following function
\begin{equation}\label{eq_def_vPhi_vlambda_k}
    \Phi_{\vlambda,k}\br{\vxi}:=
    \int
        e\br{\vxi^\top \vLambda \vgamma_k\br{t}}
        \rho\br{t}
    dt
    \cdot
    \varphi\otimes \psi\br{\vxi}.
\end{equation}
The mean zero condition on \(\rho\) implies that
\begin{equation}\label{eq_int_mean_0_vlambda_bd}
    \abs{
    \int
        e\br{\vxi^\top \vLambda \vgamma_k\br{t}}
        \rho\br{t}
    dt
    }
    =
    \abs{
    \int
        \br{e\br{\vxi^\top \vLambda \vgamma_k\br{t}}-1}
        \rho\br{t}
    dt
    }
    \lesssim \abs{\vLambda \vxi} B_0,
\end{equation}
and direct computation also shows for \(\vbeta \in \BR{0,1,\dots,N_\vM}^2\setminus\BR{\br{0,0}}\) the trivial bound
\begin{equation}\label{eq_int_trivial_vlambda_power_bd}
    \abs{
    \partial^\vbeta
    \int
        e\br{\vxi^\top \vLambda \vgamma_k\br{t}}
        \rho\br{t}
    dt
    }
    \eqsim
    \abs{
    \int
        e\br{\vxi^\top \vLambda \vgamma_k\br{t}}
        \br{\lambda_1\gamma_{k,1}\br{t}}^{\beta_1}
        \br{\lambda_2\gamma_{k,2}\br{t}}^{\beta_2}
        \rho\br{t}
    dt
    }
    \lesssim \br{\abs{\vlambda}B_0}^{\abs{\vbeta}}.
\end{equation}
In view of \eqref{eq_int_mean_0_vlambda_bd} and \eqref{eq_int_trivial_vlambda_power_bd}, we deduce that 
\begin{equation*}
    \abs{\vlambda}\lesssim AB_1\implies
    \abs{\partial^\vbeta \Phi_{\vlambda,k}\br{\vxi}}\underset{A,B_1,N_\vM}{\lesssim} \abs{\vlambda}\1_{\br{1/2,2}}\br{\abs{\vxi}}
    ,\quad
    \forall k,\,\forall \vbeta \in\BR{0,1,\dots,N_\vM}^2.
\end{equation*}
Therefore, in terms of \eqref{eq_symb_w_seq_struct}, we may write
\begin{equation*}
    m_\vlambda=m_{\vepsilon,\vPhi,\vM},\quad
    \vepsilon=\BR{\abs{\vlambda}}_k,\quad
    \vPhi:=\BR{\Phi_{\vlambda,k}/\abs{\vlambda}}_k
\end{equation*}
and conclude via the validity of \textbf{Conjecture \ref{conj_twist_vM}} for \(\vM:=\BR{\vM_k}\) the desired estimate \eqref{eq_LF_SIO_L_T_m_est}. 

As for \eqref{eq_LF_SIO_H_T_m_est}, we again analyze the derivatives of the oscillatory expression in the defining formula of \(\Phi_{\br{1,\lambda},k}\)
\begin{equation*}
    \abs{
        \partial^\vbeta
        \int
            e\br{\xi_1\gamma_{k,1}\br{t}+
            \lambda\xi_2\gamma_{k,2}\br{t}}
            \rho\br{t}
        dt
    }
    \eqsim
    \abs{
    \int
        e\br{\xi_1\gamma_{k,1}\br{t}+
        \lambda\xi_2\gamma_{k,2}\br{t}}
        \gamma^{\beta_1}_{k,1}\br{t}
        \br{\lambda\gamma_{k,2}\br{t}}^{\beta_2}
        \rho\br{t}
    dt
    }
\end{equation*}
for \(\vbeta\in\BR{0,1,\dots,N_\vM}^2\).
Suggested by non-stationary phase principle, we equate the above with
\begin{equation*}
    =
    \abs{
    \int
        e\br{\xi_1\gamma_{k,1}\br{t}+
        \lambda\xi_2\gamma_{k,2}\br{t}}
        % \gamma^{\beta_1}_{k,1}\br{t}
        \br{\frac{\partial}{2\pi i \partial t}\cdot  \frac{1}{\xi_1 \gamma'_{k,1}\br{t}+\lambda \xi_2\gamma'_{k,2}\br{t}}\cdot }^{\beta_2+1}
        \gamma^{\beta_1}_{k,1}\br{t}
        \br{\lambda\gamma_{k,2}\br{t}}^{\beta_2}
        \rho\br{t}
    dt
    }
\end{equation*}
and dominate the above under the condition that \(\abs{\lambda}\gg AB_1\) and \(\varphi\otimes\psi\br{\vxi}\neq 0\) with 
\begin{equation*}
    \underset{N_\vM}{\lesssim}\max\br{1,A,B_{N_\vM+2}}^{4N_\vM+1}/\abs{\lambda}.
\end{equation*}
Consequently, we deduce 
\begin{equation*}
    \abs{\lambda}\gg AB_1\implies
    \abs{\partial^\vbeta \Phi_{\br{1,\lambda},k}\br{\vxi}}\underset{A,B_{N_\vM+2},N_\vM}{\lesssim} \abs{\lambda}^{-1}\1_{\br{1/2,2}}\br{\abs{\vxi}}
    ,\quad
    \forall k,\,\forall \vbeta \in\BR{0,1,\dots,N_\vM}^2.
\end{equation*}
Therefore, following the previous discussion, we shall set
\begin{equation*}
    m_{\br{1,\lambda}}=m_{\vepsilon,\vPhi,\vM},\quad
    \vepsilon=\BR{\abs{\lambda}^{-1}}_k,\quad
    \vPhi:=\BR{\abs{\lambda}\Phi_{\vlambda,k}}_k
\end{equation*}
and conclude via the validity of \textbf{Conjecture \ref{conj_twist_vM}} for \(\vM:=\BR{\vM_k}\) the desired estimate \eqref{eq_LF_SIO_H_T_m_est}. 
With both \eqref{eq_LF_SIO_L_T_m_est} and \eqref{eq_LF_SIO_H_T_m_est}, we complete the proof of \textbf{Proposition \ref{prop_LF_SIO}}.

The proof of \textbf{Proposition \ref{prop_LF_M} and \ref{prop_LF_SM}} follow similar structure. Consider the general form and the identity
\begin{equation*}
    \cA^L_{\vM_k\vgamma_k,\rho}\br{f_1,f_2}\br{\vx}=
    -\cA_{\vM_k\vgamma_k,\rho}\br{\fQ{1}_{1/M_{k,1}} f_1,\fQ{2}_{1/M_{k,2}} f_2}
\end{equation*}
\begin{equation}\label{eq_LF_off}
    +\cA_{\vM_k\vgamma_k,\rho}
    \br{\fQ{1}_{1/M_{k,1}} f_1,f_2}
    +\cA_{\vM_k\vgamma_k,\rho}
    \br{f_1,\fQ{2}_{1/M_{k,2}} f_2}.
\end{equation}
By triangle inequality and the symmetry between the two terms in \eqref{eq_LF_off}, it suffices to study the two operators
\begin{equation}\label{eq_LF_M_L}
    \sup_k
    \abs{
        \cA_{\vM_k\vgamma_k,\rho}\br{\fQ{1}_{1/M_{k,1}} f_1,\fQ{2}_{1/M_{k,2}} f_2}
        \br{\vx}
    },
\end{equation}
\begin{equation}\label{eq_LF_M_A}
    \sup_k
    \abs{
        \cA_{\vM_k\vgamma_k,\rho}\br{\fQ{1}_{1/M_{k,1}} f_1, f_2}
        \br{\vx}
    }.
\end{equation}
Recall now \eqref{eq_varphi_delta} and \textbf{Definition \ref{def_LP_decomp}}. We utilize the following trivial fact
\begin{equation*}
    \abs{z}\lesssim 1\implies \abs{\cQ_{\delta,1/R}f\br{x+ R z}}\lesssim Mf\br{x},\quad \forall x\in\R,\,R>0
\end{equation*}
to deduce the following point-wise estimate
\begin{equation*}
    \abs{
        \cA_{\vM_k\vgamma_k,\rho}\br{\fQ{1}_{1/M_{k,1}} f_1, f_2}
        \br{\vx}
    }
    \underset{B_0}{\lesssim} 
    M^{\br{1}}f_1\br{\vx} 
    \int
        \abs{f_2\br{\vx+M_{k,2}\gamma_{k,2}\br{t}\ve_2}
        \rho\br{t}}
    dt
\end{equation*}
and consequently the following as well
\begin{equation*}
    \abs{
        \cA_{\vM_k\vgamma_k,\rho}\br{\fQ{1}_{1/M_{k,1}} f_1, \fQ{2}_{1/M_{k,2}} f_2}
        \br{\vx}
    }
    \underset{B_0}{\lesssim} 
    M^{\br{1}}f_1\br{\vx} 
    M^{\br{2}}f_2\br{\vx}
    .
\end{equation*}
Taking supremum over \(k\) yields the following two point-wise estimates
\begin{equation}\label{eq_LF_M_L_pw_bd}
    \abs{\eqref{eq_LF_M_L}}\underset{B_0}{\lesssim}M^{\br{1}}f_1\br{\vx}M^{\br{2}}f_2\br{\vx},
\end{equation}
\begin{equation}\label{eq_LF_M_A_pw_bd}
    \abs{\eqref{eq_LF_M_A}}\underset{B_0}{\lesssim}M^{\br{1}}f_1\br{\vx}
    \sup_k
    \int
        \abs{f_2\br{\vx+M_{k,2}\gamma_{k,2}\br{t}\ve_2}
        \rho\br{t}}
    dt.
\end{equation}
On the one hand \eqref{eq_LF_M_L_pw_bd} implies the full range \eqref{eq_maximal_op_range} boundedness
\begin{equation*}
    \nrm{\sup_k
        \abs{
            \cA_{\vM_k\vgamma_k,\rho}\br{\fQ{1}_{1/M_{k,1}} f_1,\fQ{2}_{1/M_{k,2}} f_2}
        }\br{\vx}
    }_{L^{p_3}\br{d\vx}}\underset{B_0,p_j}{\lesssim}\nrm{f_1}_{L^{p_1}}\nrm{f_2}_{L^{p_2}},
\end{equation*}
which address the contribution from \eqref{eq_LF_M_L} in both \textbf{Proposition \ref{prop_LF_M} and \ref{prop_LF_SM}}.
On the other hand, \eqref{eq_LF_M_A_pw_bd} suggests that we introduce maximal operators of the following form
\begin{equation}\label{eq_Max_Gamma_rho}
    M_{\BR{\Gamma_k}_k,\rho} f\br{x}:=
    \sup_k
    \int
        \abs{
        f\br{x+\Gamma_k\br{t}}
        \rho\br{t}
        }
    dt
    ,\quad
    f\in L^p\br{\R}.
\end{equation}
With \eqref{eq_Max_Gamma_rho} and \eqref{def_fiber_op}, we rephrase \eqref{eq_LF_M_A_pw_bd} as
\begin{equation*}
    \sup_k
    \abs{
        \cA_{\vM_k\vgamma_k,\rho}\br{\fQ{1}_{1/M_{k,1}} f_1, f_2}
        \br{\vx}
    }
    \lesssim
    M^{\br{1}}f_1\br{\vx}
    M^{\br{2}}_{\BR{M_{k,2}\gamma_{k,2}}_k,\rho}
    f_2\br{\vx}.
\end{equation*}
As a result, to address the contribution from \eqref{eq_LF_M_A}, it suffices to provide full \(p\in(1,\infty]\) range \(L^p\) boundedness result for \eqref{eq_Max_Gamma_rho} under the respective settings for \textbf{Proposition \ref{prop_LF_M} and \ref{prop_LF_SM}}.

\subsubsection{\(\BR{\Gamma_k}_k=\BR{M_{k,2}\gamma_{k,2}}_k\) with assumptions given in \textbf{Proposition \ref{prop_LF_M}}}\label{subsubsec_LF_M_lin}
Observe that the condition:
\begin{equation*}
    \abs{\gamma_{k,2}\br{t}}\leq B_0,\quad 1/A\leq \abs{\gamma'_{k,2}\br{t}},\quad \forall t\in\supp\rho
\end{equation*}
allows us to perform the following changes of variables
\begin{equation}\label{eq_CV_of_M_mgamma}
    M_{\BR{M_{k,2}\gamma_{k,2}}_k,\rho} f\br{x}=
    \sup_k
    \int^{B_0}_{-B_0}
    \abs{
        f\br{x+M_{k,2}\tau}
        \rho\br{\gamma^{-1}_{k,2}\br{\tau}}
    }
    \frac{d\tau}{\abs{\gamma'_{k,2}\br{\gamma^{-1}_{k,2}\br{\tau}}}}.
\end{equation}
Via a trivial bound, it follows that
\begin{equation*}
    M_{\BR{M_{k,2}\gamma_{k,2}}_k,\rho} f\br{x}
    \lesssim 
    A B_0
    \nrm{\rho}_{L^\infty}
    \sup_k
    \fint^{B_0}_{-B_0}
    \abs{
        f\br{x+M_{k,2}\tau}
    }
    d\tau
    \leq 
    A B_0
    \nrm{\rho}_{L^\infty} Mf\br{x}.
\end{equation*}
This point-wise estimate provides the desired full range \(L^p\) boundedness.

\subsubsection{\(\BR{\Gamma_k}_k=\BR{M_{k,2}\gamma_2}_k\) with assumptions given in \textbf{Proposition \ref{prop_LF_SM}}}\label{subsubsec_LF_SM_lin} By previous discussion, we are done if \(\gamma'_2\) vanishes nowhere on \(\supp \rho\).  Therefore, it is most interesting when \(\gamma'_2\not\equiv0\) but vanishes somewhere on \(\supp \rho\). By partitioning \(\rho\) and some standard normalization, we may assume that \(\supp\rho=\mr{0,1}\) and that \(\gamma'_2\br{t}\) vanishes only at \(t=0\) with some order \(n\in\N\) such that the following asymptotics hold
\begin{equation*}
    \abs{\gamma'_2\br{t}}\gtrsim \abs{t}^n,\quad
    \abs{\gamma_2\br{t}-\gamma_2\br{0}}\lesssim \abs{t}^{n+1},\quad\forall t\in\mr{0,1}.
\end{equation*}
We thus obtain via changes of variables analogous to \eqref{eq_CV_of_M_mgamma} and triangle inequalities the following
\begin{equation*}
    M_{\BR{M_{k,2}\gamma_2},\rho}f\br{x}
    =
    \sup_k
    \int
    \abs{
        f\br{x+M_{k,2}\br{\gamma_2\br{0}+\tau}}
        \rho\br{\gamma^{-1}_2\br{\gamma_2\br{0}+\tau}}
    }
    \frac{d\tau}{\abs{\gamma'_2\br{\gamma^{-1}_2\br{\gamma_2\br{0}+\tau}}}}
\end{equation*}
\begin{equation*}
    \lesssim
    \nrm{\rho}_{L^\infty}
    \sup_k
    \int_{\abs{\tau}\leq C}
    \abs{
        f\br{x+M_{k,2}\br{\gamma_2\br{0}+\tau}}
    }
    \frac{d\tau}{\abs{\tau}^{\frac{n}{n+1}}}
\end{equation*}
\begin{equation*}
    \underset{\rho}{\lesssim}
    \sum_{\substack{
        j\leq \log_2 C\\
        \pm\in\BR{+,-}
    }}
    2^{-\frac{jn}{n+1}}
    \sup_k
    \int^{2^j}_{2^{j-1}}
    \abs{
        f
        \br{
            x+M_{k,2}\br{\gamma_2\br{0}
            \pm \tau
            }
        }
    }
    d \tau.
\end{equation*}
Via another change of variables, we then dominate the above with shifted maximal operators \eqref{eq_M_lac_bd}:
\begin{equation*}
    \lesssim
    \sum_{\substack{
        j\leq \log_2 C\\
        \pm\in\BR{+,-}
    }}
    2^{\frac{j}{n+1}}
    \sup_k
    \int^1_0
        \abs{
            f\br{
                x+2^{j-1}M_{k,2}
                \br{
                    2^{1-j}\gamma_2\br{0}+ \frac{\pm 3-1}{2}
                    +\tau
                }
            }
        }
    d\tau
\end{equation*}
\begin{equation*}
    =
    \sum_{\substack{
        j\leq \log_2 C\\
        \pm\in\BR{+,-}
    }}
    2^{\frac{j}{n+1}}
    M_{
        \BR{2^{j-1}M_{k,2}}_k,
        \BR{
            2^{1-j}\gamma_2\br{0}+\frac{\pm 3-1}{2}
        }_k
    }
    f\br{x}.
\end{equation*}
Thus, by the non-vector-valued result of \textbf{Proposition \ref{prop_FS_M_lac_bd}}, we conclude the desired full range \(L^p\) boundedness
\begin{equation*}
    \nrm{
        M_{
            \BR{M_{k,2}\gamma_2}_k,\rho
        }
        f
    }_{L^p}
    \underset{p}{\lesssim}
    \sum_{\substack{
        j\leq \log_2 C\\
        \pm\in\BR{+,-}
    }}
    2^{\frac{j}{n+1}}
    \log_{M_\ast}
    \br{
        e+\abs{2^{1-j}\gamma_2\br{0}+\frac{\pm 3-1}{2}}
    }
    \nrm{f}_{L^p}
    \underset{M_\ast,n,\gamma_2\br{0}}{\lesssim}
    \nrm{f}_{L^p}.
\end{equation*}

Finally, with the two cases in \textsc{Section \ref{subsubsec_LF_M_lin} and \ref{subsubsec_LF_SM_lin}} established, we obtain the desired full range \eqref{eq_maximal_op_range} boundedness for \eqref{eq_LF_M_A_pw_bd} in both cases and thus, complete the proof of \textbf{Proposition \ref{prop_LF_M} and \ref{prop_LF_SM}}.

\subsection{Singular integral operator}\label{subsec_tht_along_curve}

Before establishing our main result \textbf{Theorem \ref{thm_THT_asym_homo_curves}} on singular integral operators \eqref{eq_THT_vgamma} associated with curves \(\vgamma\), we first provide some basic facts regarding \(\vgamma:=\br{\gamma_1,\gamma_2}\). Recall \textbf{Definition \ref{def_asym_homogeneous} and \ref{def_asym_def}}. Let \(\ast\in\BR{0,\infty}\), \(j\in\BR{1,2}\). The asymptotic homogeneity assumptions ensure the presence of a homogeneous function \(h_{\ast,j}:\R\setminus\BR{0}\to\R\setminus\BR{0}\) of degree \(\alpha_{\ast,j}\neq 0\) such that 
\begin{equation}\label{eq_gamma_j_asym_homo}
    \gamma_j\br{t}/h_{\ast,j}\br{t}=
    1+O(\abs{t}^{\delta_\ast})
\end{equation}
for some \(\delta_\ast\in\R\) satisfying the vanishing condition \(\lim_{\abs{t}\to\ast}\abs{t}^{\delta_\ast}=0\).
We shall utilize the definability of \(\gamma_j\) to deduce smallness of the error terms and their derivatives. Indeed, we have the trivial estimate:
\begin{equation*}
    \gamma_j\br{t}= h_{\ast,j}\br{t}+ O\br{\abs{t}^{\alpha_{\ast,j}+\delta_\ast}}.
\end{equation*}
% \begin{equation*}
%     \epsilon_{\ast,j}\br{t}:=\frac{\gamma_j\br{t}}{h_{\ast,j}\br{t}}-1=
%     O(\abs{t}^{\delta_\ast}).
% \end{equation*}
To control the derivatives, we begin with the following fact regarding the definability of specific power functions:
\begin{lemma}\label{lem_spec_power_func_def}
    The power function \(\br{0,\infty}\to\br{0,\infty}:s\mapsto s^{\alpha_{\ast,j}}\) is definable.
\end{lemma}
\begin{proof}
    We observe that the definability of \(\gamma_j\) near \(\ast\) and \eqref{eq_gamma_j_asym_homo} guarantees the existence of a \(r_\ast\in\br{0,\infty}\) such that \(\dif{\gamma_j}_{U_\ast\br{r_\ast}}:U_\ast\br{r_\ast}\to\R\setminus\BR{0}\) is definable. Consider now the definable function \(F_j:\br{0,\infty}\times U_\ast\br{r_\ast}\to\R\) given by the formula:
    \begin{equation*}
        F_j\br{s,t}:=
        \1_{U_\ast\br{r_\ast}}\br{st}
        \gamma_j\br{st}/\gamma_j\br{t}.
    \end{equation*}
    By \eqref{eq_gamma_j_asym_homo}, we have:
    \begin{equation*}
        s^{\alpha_{\ast,j}}=
        \lim_{t\to\ast}
            F_j\br{s,t}.
    \end{equation*}
    As a result, the graph \(\Gamma_{pow}\subseteq \br{0,\infty}^2\) of the power function can be expressed as:
    \begin{equation*}
        \Gamma_{pow}:=
        \BR{
            \br{s,u}\in \br{0,\infty}^2
        \::\:
            u=\lim_{t\to\ast}F_j\br{s,t}
        }
    \end{equation*}
    where the limit expression can be expressed in terms of first-order formulas involving only the definable function \(F_j\) and inequalities. This implies that \(\Gamma_{pow}\) is definable and shows that the power function is also definable.
\end{proof}

As a direct consequence, the fact that there is a positively homogeneous function \(C_{\ast,j}:\R\setminus\BR{0}\to\R\setminus\BR{0}\) of degree \(0\) such that the factorization holds
\begin{equation}\label{eq_h_to_g}
    h_{\ast,j}\br{t}=
    C_{\ast,j}\br{t}\abs{t}^{\alpha_{\ast,j}}
    =
    C_{\ast,j}\br{1}\cdot t^{\alpha_{\ast,j}}\1_{\br{0,\infty}}\br{t}
    +
    C_{\ast,j}\br{-1}\cdot \br{-t}^{\alpha_{\ast,j}}\1_{\br{-\infty,0}}\br{t}
\end{equation}
% for some \(C_{\ast,j}:\BR{-1,1}\to\R\setminus\BR{0}\) 
implies the definability of \(h_{\ast,j}\) and thus the definability of the error term \(\gamma_{\ast,j}-h_{\ast,j}\) near \(\ast\). Finally, the following lemma provides control on the derivatives of the error term.

\begin{lemma}\label{lem_poly_bd_to_der_est}
    Let \(f:U_\ast\to\R\) be definable near \(\ast\). If there are \(\alpha\in\R\) and \(\delta_\ast\in\R\) satisfying the vanishing condition \(\lim_{\abs{t}\to\ast}\abs{t}^{\delta_\ast}=0\) such that
    \begin{equation}\label{eq_f_poly_est}
        \abs{f\br{t}}\lesssim
        \abs{t}^{\alpha+\delta_\ast},\quad \forall t\in U_\ast,
    \end{equation}
    for any \(N\gg 1\), there is a positive constant \(r_\ast\in\br{0,\infty}\) such that the restriction \(\dif{f}_{U_\ast\br{r_\ast}}:U_\ast\br{r_\ast}\to\R\) is \(C^N\), and each derivative \(\dif{f^{\br{k}}}_{U_\ast\br{r_\ast}}:U_\ast\br{r_\ast}\to\R\) is definable and satisfies the estimate:
    \begin{equation}\label{eq_f_der_poly_est}
        \abs{f^{\br{k}}\br{t}}\underset{N}{\lesssim} \abs{t}^{\alpha-k+\delta'_\ast}
        ,\quad \forall t\in U_\ast\br{r_\ast},\,k\in\BR{0,1,\dots,N}
    \end{equation}
    for some \(\delta'_\ast\in\R\) strictly between \(\delta_\ast\) and \(0\) that satisfy the vanishing condition \(\lim_{\abs{t}\to\ast}\abs{t}^{\delta'_\ast}=0\).
\end{lemma}

\begin{proof}
    Observe that it suffices to establish and iterate the result for \(N=k=1\). In fact, via pre-composing \(f\) with the reflection map \(t\mapsto-t\), it suffices to establish the one-sided result on \(U_\ast\cap\br{0,\infty}\). Furthermore, via pre-composing \(f\) with the reciprocal map \(t\mapsto 1/t\), we may assume without loss of generality\footnote{The definablility of \(t\mapsto-t\), \(t\mapsto 1/t\), and their composition with \(f\) are guaranteed by \textbf{Remark \ref{rmk_def_func}} and \textbf{Proposition \ref{prop_func_def_closure}}.} that \(\ast=\infty\) and thus, the vanishing condition reads as \(\delta_\ast=\delta_\infty<0\). As a brief summary, it remains to show the estimate
    \begin{equation}\label{eq_f_prime_one_side_reduc_est}
        \abs{f'\br{t}}\lesssim t^{\alpha-1+\delta_\infty'},\quad \forall t\in \br{r_\infty,\infty}
    \end{equation}
    for some \(\delta_\infty'<0\) under the assumption
    \begin{equation}\label{eq_f_one_side_reduc_est}
        \abs{f\br{t}}\lesssim t^{\alpha+\delta_\infty},\quad \forall t\in \br{r,\infty}
    \end{equation}
    for some \(\delta_\infty<0<r\).
    With these reductions, let \(f\) be given with the stated assumptions.
    By \textbf{Proposition \ref{prop_diff_def} and \ref{prop_diff_ae}} and \textbf{Theorem \ref{thm_mono}}, we can take \(r_\infty\) such that \(\dif{f}_{U_\infty\br{r_\infty}}\) is definable, \(C^1\), and with definable derivative. Specifically, both \(\dif{f}_{\br{r_\infty,\infty}}\) and \(\dif{f'}_{\br{r_\infty,\infty}}\) are continuous definable functions. 
    % Furthermore, \textbf{Theorem \ref{thm_mono}} allows us to take \(r_\infty>0\) large enough such that both \(\dif{f}_{\br{r_\infty,\infty}}\) and \(\dif{f'}_{\br{r_\infty,\infty}}\) are monotone.
    Let now \(\delta_\infty'\in\br{\delta_\infty,0}\) be any real number such that \(\beta:=\alpha+\delta_\infty'\in\Q\setminus\BR{0}\). To prove \eqref{eq_f_prime_one_side_reduc_est},
    we shall utilize the L'H\^{o}pital's rule.
    
    First, we show the existence of the following limit
    \begin{equation}\label{eq_lim_quo_exists}
        -\infty \leq \lim_{t\to \infty}\frac{f'\br{t}}{t^{\beta-1}}\leq \infty
    \end{equation}
    in the extended real. Indeed, by the fact \(\beta\in\Q\) and \textbf{Remark \ref{rmk_def_func}}, \(
    t\mapsto t^{\beta-1}
    \) is definable on \(\br{r_\infty,\infty}\), and so is the quotient \(f'\br{t}/t^{\beta-1}\). \textbf{Corollary \ref{cor_lim}} thus guarantees the existence of the limit \eqref{eq_lim_quo_exists}.

    Next, we observe that \eqref{eq_f_one_side_reduc_est} and the fact that \(\delta_\infty'>\delta_\infty\) imply
    \begin{equation*}
        \abs{f\br{t}}\lesssim t^{\alpha+\delta_\infty}\underset{r_\infty}{\lesssim}t^{\alpha+\delta_\infty'}=t^\beta,\quad
        \forall t\in\br{r_\infty,\infty}.
    \end{equation*}
    Using also the fact that \(\beta\neq 0\), we either have \(\lim_{t\to\infty}t^\beta=\infty\) or via the above relation, the following condition
    \begin{equation*}
        \lim_{t\to\infty}f\br{t}= \lim_{t\to\infty}t^\beta=0.
    \end{equation*}
    In both cases, we may apply L'H\^{o}pital's rule to deduce
    \begin{equation*}
        \lim_{t\to \infty}\frac{f'\br{t}}{t^{\beta-1}}=\lim_{t\to\infty} \frac{\beta f\br{t}}{t^\beta}
    \end{equation*}
    Use again the fact \(\delta_\infty'>\delta_\infty\) and the assumption \eqref{eq_f_one_side_reduc_est}. We equate the above with
    \begin{equation*}
         =\beta
        \lim_{t\to\infty}t^{\delta_\infty-\delta_\infty'}\cdot \frac{f\br{t}}{t^{\alpha+\delta_\infty}}=0.
    \end{equation*}
    This implies \eqref{eq_f_prime_one_side_reduc_est} and concludes the proof of the lemma.

\end{proof}

Henceforth, \(\delta_\ast\) will always denote a number satisfying the vanishing condition \(\lim_{\abs{t}\to\ast}\abs{t}^{\delta_\ast}=0\) and be allowed to vary from line to line.
% Recall \eqref{eq_h_to_g}. For simplicity, we set \(\vh_\ast:=\br{h_{\ast,1},h_{\ast,2}}\) and \(\vC_\ast:=\br{C_{\ast,1},C_{\ast,2}}\) and introduce for \(s\in\br{0,\infty}\) and \(\valpha:=\br{\alpha_1,\alpha_2}\in\R^2\) the shorthand notation:
% \begin{equation}\label{eq_power_mat}
%     s^\valpha:=
%     \begin{pmatrix}
%         s^{\alpha_1} & 0\\
%         0 & s^{\alpha_2}
%     \end{pmatrix}.
% \end{equation}
% We thus obtain the following normal form:
% \begin{equation}\label{eq_vgamma_nrm_form}
%     \vgamma\br{t}
%     =\vh_\ast\br{t}+\abs{t}^{\valpha_\ast}\vO\br{\abs{t}^{\delta_\ast}}
%     =\abs{t}^{\valpha_\ast}
%     \br{
%         \vC_\ast\br{\sgn t}+
%         \vO(\abs{t}^{\delta_\ast})
%         % \begin{pmatrix}
%         %     O(\abs{t}^{\delta_\ast})\\
%         %     O(\abs{t}^{\delta_\ast})
%         % \end{pmatrix}
%     },\quad
%     \vO(\abs{t}^{\delta_\ast}):=
%     \br{
%         O(\abs{t}^{\delta_\ast}),
%         O(\abs{t}^{\delta_\ast})
%     }.
% \end{equation}
% Pairing the normal form \eqref{eq_vgamma_nrm_form} with \textbf{Lemma \ref{lem_poly_bd_to_der_est}}, we deduce:
With \textbf{Lemma \ref{lem_poly_bd_to_der_est}}, we deduce:
\begin{lemma}\label{lem_vgamma_def_deriv_approx}
    For all \(N\gg 1\) and \(\ast\in\BR{0,\infty}\), we can find \(r_\ast\in\br{0,\infty}\) close enough to \(\ast\) such that \(\dif{\vgamma}_{U_\ast\br{r_\ast}}\) is definable and \(C^N\) and satisfy the following relations:
    \begin{equation*}
        \vgamma^{\br{k}}_j\br{t}=
        h^{\br{k}}_{\ast,j}\br{t}
        +
        O\br{\abs{t}^{\alpha_{\ast,j}-k+\delta_\ast}},\quad
        \forall t\in U_\ast\br{r_\ast},\, k\in\BR{0,1,\dots,N},\, j\in\BR{1,2},
    \end{equation*}
    where the implicit constants only depends on \(\vgamma\), \(r_\ast\), and \(N\).
\end{lemma}
The above lemma allows us to compute \eqref{eq_THT_ass_curvature_neq_0} explicitly. This gives:
\begin{equation*}
    0<
    \liminf_{\abs{t}\to\ast}
    \abs{t}/
    \abs{
        \frac{
            \gamma'_1\gamma'_2
        }{ 
            \vgamma'\wedge\vgamma''
        }
    }
    \br{t}
    =\abs{\alpha_{\ast,2}-\alpha_{\ast,1}}.
\end{equation*}
Additionally, \eqref{eq_THT_ass_growth_match} yields the two statements:
\begin{equation*}
    \alpha_{0,1}<0\iff \liminf_{\abs{t}\to 0}\abs{\gamma_1\br{t}}=\infty \iff \liminf_{\abs{t}\to 0}\abs{\gamma_2\br{t}}=\infty
    \iff \alpha_{0,2}<0;
\end{equation*}
\begin{equation*}
    \alpha_{\infty,1}>0\iff \liminf_{\abs{t}\to \infty}\abs{\gamma_1\br{t}}=\infty \iff \liminf_{\abs{t}\to \infty}\abs{\gamma_2\br{t}}=\infty
    \iff \alpha_{\infty,2}>0.
\end{equation*}
In other words, the two conditions \eqref{eq_THT_ass_curvature_neq_0} and \eqref{eq_THT_ass_growth_match} can be rephrased as:
\begin{equation}\label{eq_alpha_conds}
    \alpha_{\ast,1}\alpha_{\ast,2}>0,\quad \alpha_{\ast,1}\neq\alpha_{\ast,2},\quad
    \forall\ast\in\BR{0,\infty}.
\end{equation}
% and call \(\vh_{\ast,\pm}\) the characteristic homogeneous curve.

\subsubsection{Proof of \textbf{Theorem \ref{thm_THT_asym_homo_curves}}}\label{subsubsec_tht_asym_homo}
Return now to the analysis of the operator \(T_\vgamma\). Fix \(\delta=1/5\). We recall \eqref{eq_part_R+} and perform the standard decomposition \eqref{eq_SIO_2_A}.
To reiterate, we have
\begin{equation*}
    T_\vgamma\br{f_1,f_2}\br{\vx}=
    \sum_{k\in\Z}
        \cA_{\Dil^\infty_{2^{-k}} \vgamma, \rho/t}       \br{f_1,f_2}\br{\vx}.
\end{equation*}
Since, similar to \eqref{eq_T_Gamma_log_triv}, we have the following trivial estimate
\begin{equation*}
    \nrm{
        \sum_{k\in Z}
        \cA_{\Dil^\infty_{2^{-k}} \vgamma, \rho/t}       \br{f_1,f_2}\br{\vx}.
    }_{L^{p_3}}
    \lesssim
    \# Z\cdot
    \nrm{f_1}_{L^{p_1}}
    \nrm{f_2}_{L^{p_2}},\quad
    Z\subseteq \Z
\end{equation*}
for \(p_j\in\mr{1,\infty}\) with \(\frac{1}{p_3}=\frac{1}{p_1}+\frac{1}{p_2}\), it suffices to control the following pair of operators:
\begin{equation}\label{eq_THT_2_ends_def}
    \sum_{k<N_0}
    \cA_{\Dil^\infty_{2^{-k}} \vgamma, \rho/t}       \br{f_1,f_2}\br{\vx},\quad
    \sum_{k>N_\infty}
    \cA_{\Dil^\infty_{2^{-k}} \vgamma, \rho/t}       \br{f_1,f_2}\br{\vx}
\end{equation}
% \begin{equation}\label{eq_THT_0_end_def}
%     T_{\vgamma,0}\br{f_1,f_2}\br{\vx}
%     :=
%     \sum_{k<-N_0}
%     \cA_{\Dil^\infty_{2^{-k}} \vgamma, \rho/t}       \br{f_1,f_2}\br{\vx}
% \end{equation}
% \begin{equation}\label{eq_THT_infty_end_def}
%     T_{\vgamma,\infty}\br{f_1,f_2}\br{\vx}
%     :=
%     \sum_{k>N_0}
%     \cA_{\Dil^\infty_{2^{-k}} \vgamma, \rho/t}       \br{f_1,f_2}\br{\vx}
% \end{equation}
for some fixed \(N_0<0<N_\infty\). To perform the general strategy laid out in the beginning of \textsc{Section \ref{sec_application}}, we first study the expression \(\Dil^\infty_{2^{-k}}\vgamma\) and in fact, \(\Dil^\infty_{1/s}\vgamma\) for reasons which will be clear in a moment.

Choose \(N_0\) small enough and \(N_\infty\) large enough such that via \textbf{Lemma \ref{lem_vgamma_def_deriv_approx}} and the definability of multiplication in \(\R\) as bivariate function, we deduce the definability of function
\begin{equation*}
    F:\dom F \to\R^2:\br{s,t}\mapsto \vgamma\br{st}=\Dil^\infty_{1/s}\vgamma\br{t}
\end{equation*}
with its domain described as below
\begin{equation*}
    \dom F:=\BR{\br{s,t}\in\br{0,\infty}\times\R\::\: \br{s<2^{N_0}\lor s>2^{N_\infty}} \land \br{ -4<t<-1/4 \lor 1/4< t<4}}
\end{equation*}
and the approximation identities up to \(N\)-th derivatives
\begin{equation}\label{eq_dil_gammas_aprox_0}
    \br{\Dil^\infty_{1/s} \gamma_j}^{\br{l}}\br{\pm t}=s^{\alpha_{0,j}}h^{\br{l}}_{0,j}\br{\pm t}
    +s^{\br{\alpha_{0,j}+\delta_0}}O\br{t^{\alpha_{0,j}-l+\delta_0}},\quad \forall t\in\br{1/4,4},\,\forall
    s\in\br{0,2^{N_0}};
\end{equation}
\begin{equation}\label{eq_dil_gammas_aprox_infty}
    \br{\Dil^\infty_{1/s} \gamma_j}^{\br{l}}\br{\pm t}
    =s^{\alpha_{\infty,j}}h^{\br{l}}_{\infty,j}\br{\pm t}
    +s^{\br{\alpha_{\infty,j}+\delta_\infty}}O\br{t^{\alpha_{\infty,j}-l+\delta_\infty}},\quad \forall t\in\br{1/4,4},\,\forall
    s\in\br{2^{N_\infty},\infty}
\end{equation}
for each coordinate. This suggests that we consider the definable families of functions
\begin{equation}\label{eq_def_fam_4_THT_0_+}
    \mathfrak{F}_{\vgamma,0,+}:=\BR{
        \vgamma_{s,+}:=\br{
            s^{-\alpha_{0,1}}\Dil^\infty_{1/s}\gamma_1, s^{-\alpha_{0,2}}\Dil^\infty_{1/s}\gamma_2
        }
    :\br{1/4,4}\to\R^2
    }_{s\in\br{0,2^{N_0}}},
\end{equation}
\begin{equation}\label{eq_def_fam_4_THT_0_-}
    \mathfrak{F}_{\vgamma,0,-}:=\BR{
        \vgamma_{s,-}:=\br{
            s^{-\alpha_{0,1}}\Dil^\infty_{1/s}\gamma_1, s^{-\alpha_{0,2}}\Dil^\infty_{1/s}\gamma_2
        }
    :\br{-4,-1/4}\to\R^2
    }_{s\in\br{0,2^{N_0}}},
\end{equation}
\begin{equation}\label{eq_def_fam_4_THT_infty_+}
    \mathfrak{F}_{\vgamma,\infty,+}:=\BR{
        \vgamma_{s,+}:=\br{
            s^{-\alpha_{\infty,1}}\Dil^\infty_{1/s}\gamma_1, s^{-\alpha_{\infty,2}}\Dil^\infty_{1/s}\gamma_2
        }:
   \br{1/4,4}\to\R^2
    }_{s\in\br{2^{N_\infty},\infty}},
\end{equation}
\begin{equation}\label{eq_def_fam_4_THT_infty_-}
    \mathfrak{F}_{\vgamma,\infty,-}:=\BR{
        \vgamma_{s,-}:=\br{
            s^{-\alpha_{\infty,1}}\Dil^\infty_{1/s}\gamma_1, s^{-\alpha_{\infty,2}}\Dil^\infty_{1/s}\gamma_2
        }:
   \br{-4,-1/4}\to\R^2
    }_{s\in\br{2^{N_\infty},\infty}}.
\end{equation}
Notice that \eqref{eq_dil_gammas_aprox_0} and \eqref{eq_dil_gammas_aprox_infty} implies the following uniform estimates:
\begin{equation}\label{eq_nrm_gamma_j_est_4_THT_0}
    \abs{\br{s^{-\alpha_{0,j}}\Dil^\infty_{1/s}\gamma_j-h_{0,j}}^{\br{l}}\br{\pm t}}\lesssim 2^{N_0\delta_0},\quad \forall t\in\br{1/4,4},\,\forall
    s\in\br{0,2^{N_0}};
\end{equation}
\begin{equation}\label{eq_nrm_gamma_j_est_4_THT_infty}
    \abs{\br{s^{-\alpha_{\infty,j}}\Dil^\infty_{1/s}\gamma_j-h_{\infty,j}}^{\br{l}}\br{\pm t}}\lesssim 2^{N_\infty\delta_\infty},\quad \forall t\in\br{1/4,4},\,\forall
    s\in\br{2^{N_\infty},\infty}.
\end{equation}
In particular, all of the above estimates hold for all \(s=2^k\) with either \(k<N_0\) or \(k>N_\infty\).
Thus, via setting
\begin{equation}\label{eq_vM_k_4_THT}
    \vM_k:=
    \begin{cases}
        \begin{pmatrix}
            2^{k\alpha_{0,1}} & 0\\
            0 & 2^{k\alpha_{0,2}}
        \end{pmatrix},
        &
        k<N_0,\\
        \begin{pmatrix}
            2^{k\alpha_{\infty,1}} & 0\\
            0 & 2^{k\alpha_{\infty,2}}
        \end{pmatrix},
        &
        k>N_\infty;
    \end{cases}
\end{equation}
\begin{equation}\label{eq_vgamma_s_4_THT}
    \vgamma_s:=
    \begin{cases}
        \br{
            s^{-\alpha_{0,1}}\Dil^\infty_{1/s}\gamma_1, s^{-\alpha_{0,2}}\Dil^\infty_{1/s}\gamma_2
        },
        &
        0<s<2^{N_0},\\
        \br{
            s^{-\alpha_{\infty,1}}\Dil^\infty_{1/s}\gamma_1, s^{-\alpha_{\infty,2}}\Dil^\infty_{1/s}\gamma_2
        },
        &
        2^{N_\infty}<s<\infty,
    \end{cases}
\end{equation}
we may reduce either of \eqref{eq_THT_2_ends_def} to an infinite partial sum in \(k\)
\begin{equation*}
    \sum_{k<N_0} \cA_{\vM_k\vgamma_{2^k},\rho/t}\br{f_1,f_2}\br{\vx},
    \quad
    \sum_{k>N_\infty} \cA_{\vM_k\vgamma_{2^k},\rho/t}\br{f_1,f_2}\br{\vx}.
\end{equation*}
As a consequence of \eqref{eq_nrm_gamma_j_est_4_THT_0} and \eqref{eq_nrm_gamma_j_est_4_THT_infty}, for any prescribed \(N\in\N\) and \(\epsilon>0\), we may choose \(N_0,N_\infty\) such that 
\begin{equation}\label{eq_vgamma_s_to_h_4_THT_0}
    \nrm{
        \vgamma_s-\br{h_{0,1},h_{0,2}}
    }_{C^N\br{
        \br{-4,4}\setminus\Br{-1/4,1/4}
    }}<\epsilon,\quad \forall s\in\br{0,2^{N_0}},
\end{equation}
\begin{equation}\label{eq_vgamma_s_to_h_4_THT_infty}
    \nrm{
        \vgamma_s-\br{h_{\infty,1},h_{\infty,2}}
    }_{C^N\br{
        \br{-4,4}\setminus\Br{-1/4,1/4}
    }}<\epsilon,\quad \forall s\in\br{2^{N_\infty},\infty}.
\end{equation}
From this point on, there is no essential distinction between the two cases \(k<N_0\) and \(k>N_\infty\). We therefore suppress the subscripts \(\br{\cdot}_0\) and \(\br{\cdot}_\infty\) indicating the asymptotic behavior.

Proceed now with the general strategy illustrated in the beginning of \textsc{Section \ref{sec_application}}. To prove \textbf{Theorem \ref{thm_THT_asym_homo_curves}} it remains to verify the required assumptions for applying \textbf{Proposition \ref{prop_HF_tame}, \ref{prop_HF_dec}, and \ref{prop_LF_SIO}}. Recall \eqref{eq_cond_4_Hparts}. We note that by previous discussion, we always have \(B_N<\infty\). To show \(C_\varepsilon<\infty\) for some \(\varepsilon>0\), we observe that
\begin{equation*}
    \cA_{\vgamma_s,\rho/t}\br{f_1,f_2}=
    \cA_{\vgamma_{s,+},\rho_+/t}\br{f_1,f_2}+
    \cA_{\vgamma_{s,-},\rho_-/t}\br{f_1,f_2}.
\end{equation*}
As a direct consequence, we have the trivial estimate
\begin{equation*}
    C_{\vgamma_s,\rho/t,\varepsilon}\leq C_{\vgamma_{s,+},\rho_+/t,\varepsilon}+
    C_{\vgamma_{s,-},\rho_-/t,\varepsilon}.
\end{equation*}
By \eqref{eq_vgamma_s_to_h_4_THT_0} and \eqref{eq_vgamma_s_to_h_4_THT_infty}, it is a routine calculation to check that the functions in the four definable families \eqref{eq_def_fam_4_THT_0_+}, \eqref{eq_def_fam_4_THT_0_-}, \eqref{eq_def_fam_4_THT_infty_+}, and \eqref{eq_def_fam_4_THT_infty_-} are \(A\)-admissible for a fixed \(A>1\) that only depends on the positively homogeneous functions \(h_1,h_2\). Therefore, we can apply \textbf{Theorem \ref{thm_smoothing_o_mini_uni}} to deduce
\begin{equation*}
    \sup_s C_{\vgamma_{s,+},\rho_+/t,\varepsilon}\underset{\mathfrak{F}_{\vgamma,+}}{\lesssim} 1,\quad
    \sup_s C_{\vgamma_{s,-},\rho_-/t,\varepsilon}\underset{\mathfrak{F}_{\vgamma,-}}{\lesssim} 1
\end{equation*}
for some universal \(\varepsilon>0\). In particular, the above holds with \(s\) restricted to dyadic numbers \(\BR{2^k}_{k\in\Z}\). Namely, we have shown that \(C_\varepsilon<\infty\). This addresses the required assumption to apply \textbf{Proposition \ref{prop_HF_tame} and \ref{prop_HF_dec}}. As for the required assumption for \textbf{Proposition \ref{prop_LF_SIO}}, we recall \textbf{Remark \ref{rmk_Conj_holds_for_good_vM}}. In our current setting, the sequence \(\BR{\vM_k}_k\) falls within the scope stated in \textbf{Remark \ref{rmk_Conj_holds_for_good_vM}} except potentially both \(\alpha_1,\alpha_2<0\). This can be addressed simply by changing the index \(k\mapsto -k\) and replacing \(\alpha_j\mapsto-\alpha_j\). Combining these fact, we complete the proof of \textbf{Theorem \ref{thm_THT_asym_homo_curves}}.

\subsubsection{Polynomially bounded settings}\label{subsubsec_poly_bdd_THT}
Recall \textbf{Definition \ref{def_omin_poly_bdd}}. Let now \(\cR\) be polynomially bounded (see \textbf{Definition \ref{def_omin_poly_bdd}}). We claim that in such a setting, the asymptotic homogeneity assumption in our main statement \textbf{Theorem \ref{thm_THT_asym_homo_curves}} is a mild requirement due to the following result:
% o-minimal expansion of \(\ang{\R,+,\cdot,<}\).
% We have:
\begin{proposition}[Chris Miller\cite{MR1195484}]\label{prop_chris_mill_poly_bdd_asym}
    A function \(f:\br{r,\infty}\to \R\) definable in \(\cR\) satisfies the following dichotomy:
    \begin{itemize}
        \item either \(f\) is eventually zero; that is, there is \(t_0\in\br{r,\infty}\) such that \(f\br{t}=0\) whenever \(t>t_0\),
        \item or there are \(C\neq 0\) and \(c\in\R\) such that \(\lim_{t\to \infty} f(t)/t^c=C\), and in particular, the specific power function \(t\mapsto t^c\) on \(\br{0,\infty}\) is definable in \(\mathcal{S}\).
    \end{itemize}
\end{proposition}

As a direct consequence, whenever \(f\) is not eventually zero, there is a definable function \(h:\br{0,\infty}\to\R\setminus\BR{0}\) given by \(h\br{t}:=Ct^c\) such that \(\lim_{t\to\infty}f\br{t}/h\br{t}=1\). Since \(f/h-1\) is also definable, another application of \textbf{Proposition \ref{prop_chris_mill_poly_bdd_asym}} implies that either \(f/h-1\) is eventually zero, or there are \(D\neq 0\) and \(d\in\R\) such that:
\begin{equation*}
    \lim_{t\to\infty} t^d =0
    ,\quad
    \lim_{t\to\infty} \frac{f\br{t}/h\br{t}-1}{t^d}=D.
\end{equation*}
In both cases, we deduce that \(f\br{t}/h\br{t}=1+O(\abs{t}^{\delta_\infty})\) for some \(\delta_\infty\in\R\) satisfying \(\lim_{t\to\infty}\abs{t}^{\delta_\infty}=0\).

Now applying the changes of variables \(t\mapsto 1/t\). We deduce analogous statements as \(t\) approaches \(0\) from above. Finally, via reflection \(t\mapsto-t\), we deduce analogous statements as \(t\) approaches \(-\infty\) and \(0\) from below. In other words, excluding the eventually zero and constant trivial cases,
the asymptotic homogeneity assumption serves to ensure the degrees of homogeneity of the function \(\gamma_j\) on the positive and the negative side of \(\ast\in\BR{0,\infty}\) match and are non-zero.
We further provide the following characterization for the degree of asymptotic homogeneity:
\begin{proposition}\label{prop_homo_deg_ext}
    Let \(\cR\) be a polynomially bounded o-minimal expansion of \(\ang{\R,+,\cdot,<}\). Let also \(\ast\in\BR{0,\infty}\) and \(\pm\in\BR{+,-}\). Given a definable function \(f:U_\ast\to\R\setminus\BR{0}\), the number \(c\in\R\) shown in the following limit:
    \begin{equation}\label{eq_one_side_asym_homo}
        \lim_{t\to \pm\ast} \frac{f\br{t}}{\abs{t}^c}=C,\quad
        C\neq 0
    \end{equation}
    can be computed via the following limit:
    \begin{equation}\label{eq_one_side_asym_homo_deg_ext}
        \lim_{t\to\pm\ast} \frac{tf'\br{t}}{f\br{t}}=c.
    \end{equation}
\end{proposition}

\begin{proof}
    For the same reason as in the proof of \textbf{Lemma \ref{lem_poly_bd_to_der_est}}, it suffices to consider the case \(\ast=\infty\) and \(\pm = +\).
    By \textbf{Proposition \ref{prop_diff_def} and \ref{prop_diff_ae}}, we may assume \(f'\) exists and is definable on \(U_\ast\).
    We now have two cases: \(c\neq 0\) and \(c=0\).
    We start by addressing the \(c\neq 0\) case.
    Recall that \textbf{Proposition \ref{prop_chris_mill_poly_bdd_asym}} guarantees the definability of \(t\mapsto t^c\). By \textbf{Corollary \ref{cor_lim}}, the following limit:
    \begin{equation*}
        -\infty\leq \lim_{t\to\infty} \frac{f'\br{t}}{ct^{c-1}} \leq \infty
    \end{equation*}
    exists in the extended real line.
    By assumption \eqref{eq_one_side_asym_homo}, the condition \(c\neq 0\) implies either \(\lim_{t\to\infty}t^c=\infty\) or
    \begin{equation*}
        \lim_{t\to\infty}f\br{t}=\lim_{t\to\infty}t^c=0.
    \end{equation*}
    We may thus apply the L'H\^{o}pital's rule to deduce that:
    \begin{equation*}
        C=
        \lim_{t\to\infty}
        \frac{f\br{t}}{t^c}
        =
        \lim_{t\to\infty} \frac{f'\br{t}}{ct^{c-1}}.
    \end{equation*}
    As a direct consequence, we conclude:
    \begin{equation*}
        \lim_{t\to\infty}
            \frac{tf'\br{t}}{f\br{t}}
        =c
        \lim_{t\to\infty}
            \frac{f'\br{t}}{ct^{c-1}}/
        \lim_{t\to\infty}
            \frac{f\br{t}}{t^c}
        =c.
    \end{equation*}
    To address the \(c=0\) case, we follow a similar argument to show that \(tf'\br{t}\) is definable, and that the limit \(\lim_{t\to\infty}tf'\br{t}\) exists in the extended real. Since \(\lim_{t\to\infty}\log t=\infty\), we may apply the L'H\^{o}pital's rule to deduce
    \begin{equation*}
        0=\lim_{t\to\infty} \frac{f\br{t}}{\log t}
        =\lim_{t\to\infty} \frac{f'\br{t}}{1/t}
        =\lim_{t\to\infty}tf'\br{t}.
    \end{equation*}
    Finally, we conclude:
    \begin{equation*}
        \lim_{t\to\infty}\frac{tf'\br{t}}{f\br{t}}
        =\lim_{t\to\infty}tf'\br{t}
        /\lim_{t\to\infty}f\br{t}
        =0/C=c.
    \end{equation*}
\end{proof}

Given \textbf{Proposition \ref{prop_homo_deg_ext}}, we may rephrase \textbf{Theorem \ref{thm_THT_asym_homo_curves}} for the polynomially bounded setting:

\begin{corollary}[Polynomially bounded settings]\label{cor_THT_poly_bdd_omini}
    Let \(\cR\) be a polynomially bounded o-minimal expansion of \(\ang{\R,+,\cdot,<}\). Let \(\gamma_j:\R\setminus\BR{0}\to\R\) be definable near \(\ast\) and satisfy \(\lim_{\abs{t}\to\ast}\frac{t\gamma'_j\br{t}}{\gamma_j\br{t}}\in \R\setminus\BR{0}\) for both \(\ast \in\BR{0,\infty}\).
    % and either even or odd. 
    If \eqref{eq_THT_ass_curvature_neq_0} and \eqref{eq_THT_ass_growth_match} hold, the operator \(T_\vgamma\) satisfies the estimates \eqref{eq_THT_Lp_bdds}.
\end{corollary}

We shall also introduce one specific application of \textbf{Corollary \ref{cor_THT_poly_bdd_omini}} for the polynomially bounded o-minimal structure \(\cR=\R_\an\).

\begin{corollary}[Real analytic settings]\label{cor_THT_real_ana}
    Let \(\gamma_j:\R\setminus\BR{0}\to\R\) be real meromorphic\footnote{That is, the function coincides with its Laurent series expansion, and all of its singularities near \(0\) and \(\infty\) are of finite order.} on both \(U_0\) and \(U_\infty\). If \eqref{eq_THT_ass_curvature_neq_0} and \eqref{eq_THT_ass_growth_match} hold, the operator \(T_\vgamma\) satisfies the estimates \eqref{eq_THT_Lp_bdds}.
\end{corollary}

Indeed, take the Laurent series expansion of \(\gamma_j\) on \(U_0\br{r_0}\) centered at \(0\). Either \(\dif{\gamma_j}_{U_0}\equiv C\) for some constant \(C\in\R\), which will be ruled out by \eqref{eq_THT_ass_curvature_neq_0}, or there is an integer \(n\in\Z\) such that the function \(f\br{t}:=t^n\gamma_j\br{t}\) is real analytic on \(U_0\br{r_0}\) with removable singularity at \(0\). Via a harmless translation that preserves the operator norm, we may assume without loss of generality that \(n\neq 0\) and impose that \(f\br{0}\neq 0\).
Taking \(I_0:=\mr{-r_0/2,r_0/2}\), we have that \(\dif{f}_{I_0}\) is definable in \(\R_\an\). Consequently, it follows that \(\gamma_j\br{t}=t^{-n}f\br{t}\) as a function on \(U_0\br{r_0/2}\) is definable in \(\R_\an\). Directly compute the following limit:
\begin{equation*}
    \lim_{\abs{t}\to 0}
        \frac{t\gamma'_j\br{t}}{\gamma_j\br{t}}
    =\lim_{t\to 0}
        \frac{
            -nt^{-n}f\br{t}
            +t^{1-n}f'\br{t}
        }{
            t^{-n}f\br{t}
        }
    =
    -n
        +
    \lim_{t\to 0}
        \frac{
            tf'\br{t}
        }{
            f\br{t}
        }
    =-n\in\R\setminus\BR{0}.
\end{equation*}
Similarly, take the Laurent series expansion of \(\gamma_j\) on \(U_\infty\br{r_\infty}\) centered at \(\infty\). There is another integer \(m\in\Z\) such that the function \(g\br{t}:=t^m\gamma_j\br{1/t}\) is real analytic on \(U_0\br{1/r_\infty}\) with removable singularity at \(0\). Again, via translation,\footnote{The translation here can be different from the previous one. This discrepancy does not cause any issue since one may decompose the operator to separate the contribution near \(\infty\) from that near \(0\).} we may assume \(m\neq 0\) and extend with \(g\br{0}\neq 0\). Taking \(I_\infty:=\mr{-\frac{1}{2r_\infty},\frac{1}{2r_\infty}}\), we have that \(\dif{g}_{I_\infty}\) is definable in \(\R_\an\). Thus, we argue that \(\gamma_j\br{t}=t^m g\br{1/t}\) as a function on \(U_\infty\br{2r_\infty}\) is definable in \(\R_\an\). Finally, by direct computation, we show that:
\begin{equation*}
    \lim_{\abs{t}\to\infty}
    \frac{t\gamma'_j\br{t}}{\gamma_j\br{t}}
    =
    \lim_{\abs{t}\to\infty}
    \frac{
        mt^mg\br{1/t}-t^{m-1}g'\br{1/t}
    }{
        t^mg\br{1/t}
    }
    =m
    -
    \lim_{t\to 0}
    \frac{tg'\br{t}}{g\br{t}}
    =m\in\R\setminus\BR{0}.
\end{equation*}
Applying now \textbf{Corollary \ref{cor_THT_poly_bdd_omini}} yields the conclusion of \textbf{Corollary \ref{cor_THT_real_ana}}.

\subsection{Maximal operator}\label{subsec_Max}
In this section, we explain how we derive  the maximal operator result \textbf{Theorem \ref{thm_max}} for \eqref{eq_M_vgamma} analogous to \textbf{Theorem \ref{thm_THT_asym_homo_curves}}. In view of \textbf{Proposition \ref{prop_HF_tame}, \ref{prop_HF_dec}, and \ref{prop_LF_M}} and the fact that the requirement for applying \textbf{Proposition \ref{prop_LF_M}} is less stringent than \textbf{Proposition \ref{prop_LF_SIO}}, the estimates stated in \textbf{Theorem \ref{thm_max}} for the smaller range
\begin{equation*}
    1<p_1,p_2<\infty,\quad \frac{1}{p_3}=\frac{1}{p_1}+\frac{1}{p_2}\leq 1
\end{equation*}
thus follow from a treatment identical to those presented in \textsc{Section \ref{subsubsec_tht_asym_homo}} with the application of \textbf{Proposition \ref{prop_LF_SIO}} replaced by that of \textbf{Proposition \ref{prop_LF_M}}. To recover the full range \eqref{eq_M_vgamma_range}. By symmetry, it remains to address the \(p_1=\infty\) case.

We begin with the trivial bound
\begin{equation*}
    M_\vgamma\br{f_1,f_2}\br{\vx}\leq
    \nrm{f_1}_{L^\infty}
    \sup_{R>0}
    \fint^R_0
        \abs{f_2\br{x_1,x_2+\gamma_2\br{t}}}
    dt
\end{equation*}
Via standard argument, we dominate the second factor with its lacunary analog
\begin{equation*}
    \sup_{R>0}
    \fint^R_0
        \abs{f_2\br{x_1,x_2+\gamma_2\br{t}}}
    dt
    \lesssim
    \sup_{k\in\Z}
    \int
        \abs{f_2\br{x_1,x_2+\gamma_2\br{2^k t}}}
    \rho\br{t}
    dt,
\end{equation*}
where \(\rho\) is a suitably chosen bump function satisfying \(\rho\br{t}\neq 0\implies t\eqsim 1\). Recall now \eqref{eq_Max_Gamma_rho}. We may equate the right-hand side of the above with \(M^{\br{2}}_{\BR{\Dil^\infty_{2^{-k}}\gamma_2}_k,\rho}f_2\br{\vx}\). Invoking notations \eqref{eq_vM_k_4_THT} and \eqref{eq_vgamma_s_4_THT}, we identify the two operators
\begin{equation*}
    M_{\BR{\Dil^\infty_{2^{-k}}\gamma_2}_k,\rho}f\br{x}=
    M_{\BR{M_{k,2}\gamma_{k,2}}_k,\rho}f\br{x},
\end{equation*}
whose full range \(L^p\) boundedness has been demonstrated in \textsc{Section \ref{subsubsec_LF_M_lin}}. Consequently, it follows that
\begin{equation*}
    \nrm{
        M_\vgamma\br{f_1,f_2}
    }_{L^p}
    \lesssim
    \nrm{f_1}_{L^\infty}
    \nrm{
        M^{\br{2}}_{\BR{M_{k,2}\gamma_{k,2}}_k,\rho}f_2
    }_{L^p}
    \lesssim
    \nrm{f_1}_{L^\infty}
    \nrm{f_2}_{L^p},\quad 
    \forall p\in\bR{1,\infty}.
\end{equation*}
This addresses the \(p_1=\infty\) endpoint case and completes the proof of \textbf{Theorem \ref{thm_max}}.

% \MH{add proof of \textbf{Theorem \ref{thm_max}}}

% In view of \textbf{Proposition \ref{prop_HF_tame}, \ref{prop_HF_dec}, and \ref{prop_LF_M}}, and the fact that the requirement for applying \textbf{Proposition \ref{prop_LF_M}} is more stringent than \textbf{Proposition \ref{prop_LF_SIO}}, we provide the following result for the maximal operators \eqref{eq_M_vgamma} analogous to \textbf{Theorem \ref{thm_THT_asym_homo_curves}}.

On the same note, due to \textbf{Theorem \ref{mainthm_global}}, we remark that \textbf{Proposition \ref{prop_HF_tame}, \ref{prop_HF_dec}, and \ref{prop_LF_M}} in combination actually provide stronger but more technical statements regarding the general form \eqref{eq_SIO_M_models}. 
\begin{theorem}
    Let \(\BR{\vM_k}_k\) be diagonal matrices satisfying \eqref{eq_lacu_seqs_condi}, and let \(\BR{\vgamma_k:\R\supset I\simto\Gamma_k\subseteq \R^2}_k\) be a sequence of \(A\)-admissible curves. If \(N_\ast:=\sup_k \Xi_{\Gamma_k}<\infty\), we have for \(\rho\in C^2_c\br{I}\) the estimates
    \begin{equation*}
        \nrm{
            \sup_k
            \abs{
                \cA_{\vM_k\vgamma_k,\rho}\br{f_1,f_2}
            }\br{\vx}
        }_{L^{p_3}\br{d\vx}}
        \underset{\substack{M_\ast,N_\ast,A,\\
        \rho,p_j}}{\lesssim}
        \nrm{f_1}_{L^{p_1}}
        \nrm{f_2}_{L^{p_2}}
    \end{equation*}
    for the same range as in \eqref{eq_M_vgamma_range}.
\end{theorem}
This follows from a direct modification of the proof presented above. Notably, the treatment of the endpoint \(p_1=\infty\) case follows exactly from the result in \textsc{Section \ref{subsubsec_LF_M_lin}}. We thus omit the standard details.

\subsection{Spherical maximal operator}

Continuing on the maximal operator theme from the previous section, we explain how we derive the maximal operator result \textbf{Theorem \ref{thm_spherMax}} for \eqref{eq_SM_vgamma}. Combining \textbf{Proposition \ref{prop_HF_tame}, \ref{prop_HF_dec}, and \ref{prop_LF_SM}} with \textbf{Theorem \ref{thm_real_ana_smoothing}}, coupled with an endpoint discussion similar to those presented in \textsc{Section \ref{subsec_Max}} but with the result in \textsc{Section \ref{subsubsec_LF_SM_lin}} in place of \textsc{Section \ref{subsubsec_LF_M_lin}}, we derive
\begin{theorem}
    Let \(\vgamma:I\to\R^2\) be a real analytic function with \(\img \vgamma\) not contained in the four classes \eqref{eq_line}, \eqref{eq_exp}, \eqref{eq_log}, or \eqref{eq_psi_line}. For any diagonal matrices \(\BR{\vM_k}_k\) satisfying \eqref{eq_lacu_seqs_condi} and any \(\rho\in C^2_c\br{I}\), we have
    \begin{equation*}
        \nrm{\sup_k\abs{\cA_{\vM_k\vgamma,\rho}\br{f_1,f_2}}\br{\vx}}_{L^{p_3}\br{d\vx}}
        \underset{M_\ast,\vgamma,p_j,\rho}{\lesssim}
        \nrm{f_1}_{L^{p_1}}
        \nrm{f_2}_{L^{p_2}}
    \end{equation*}
    for the same range as in \eqref{eq_M_vgamma_range}.
\end{theorem}
 In particular, the above holds with \(\sup_k\abs{\cA_{\vM_k\vgamma,\rho}\br{f_1,f_2}}\br{\vx}\) replaced by \eqref{eq_SM_vgamma}. That is, \textbf{Theorem \ref{thm_spherMax}} follows as a direct consequence.

\subsection{Corner-type pattern}

% In this subsection, we study the existence of patterns
% \begin{equation*}
%     \br{x_1,x_2},\quad
%     \br{x_1+\gamma_1\br{t},x_2},\quad
%     \br{x_1,x_2+\gamma_2\br{t}}
% \end{equation*}
% in a given set \(E\subseteq \R^2\) for some \(\vgamma:=\br{\gamma_1,\gamma_2}\) real analytic in a neighborhood of \(\Br{0,1}\) with \(\vgamma\br{t}=\vnull\). 

% Typically, it suffices to consider
% Suggested by the discussion at the beginning of the section, it may
% will inherit languages established at the beginning of the section to 

% present an adaptation of Bourgain's energy pigeonholing argument 

% \begin{equation*}
%     \br{x_1,x_2},\quad
%     \br{x_1+M_{k,1}\gamma_{k,1}\br{t}, x_2},\quad
%     \br{x_1, x_2+M_{k,2}\gamma_{k,2}\br{t}}
% \end{equation*}
% sequences of functions \(\BR{\vgamma_\sigma:I\to\R^2}_k\) and scalars \(\BR{M_{\sigma,j}}_k\) and some \(t\in I\).

In this subsection, we aim to prove \textbf{Theorem} \ref{thm_pattern}. Without loss of generality, we may assume via changes of variable that \(\vgamma \) is analytic in a neighborhood of \(\Br{0,2}\). Recall \eqref{eq_psi_delta} and take \(\delta:=1/2\). Consider now the following inner product expression
\begin{equation}\label{eq_conting_inner_prod_sigma}
    \ang{
        \1_E,\cA_{\vgamma\br{\sigma\cdot },\psi_{1/2}}\br{\1_E,\1_E}
    }
    =
    \int
        \1_E\br{\vx}
        \1_E\br{x_1+\gamma_1\br{\sigma t},x_2}
        \1_E\br{x_1,x_2+\gamma_2\br{\sigma t}}
        \psi_{1/2}\br{t}
    dt
    d\vx.
\end{equation}
Note that whenever the above expression is non-zero, we deduce the existence of some \(\vx\in\R^2\) and \(\sigma/2< t \leq 2\sigma\) such that the conclusion \eqref{eq_corner_pattern_in_E} holds.\footnote{From \eqref{eq_conting_inner_prod_sigma} to the existence of the above mentioned \(\sigma/2< t \leq 2\sigma\) a change of variable \(\widetilde{t}=\sigma t\) is performed.} As a result, the main task of proving \textbf{Theorem \ref{thm_pattern}} becomes demonstrating such a number \(\sigma\in\br{0,\infty}\) can be found in the following range
\begin{equation}\label{eq_sigma_range_4_corner_pattern}
    \exp\br{-\exp\br{C_\vgamma\abs{E}^{-c_\vgamma}}}\lesssim \sigma<1.
\end{equation}

In fact, in view of the general discussion at the beginning of \textsc{Section \ref{sec_application}}, the formula \eqref{eq_conting_inner_prod_sigma} suggests that we study the following general formulation:
\begin{equation}\label{eq_counting_inner_prod_sigma_gen}
    \ang{
        \1_E,\cA_{\vM_\sigma\vgamma_\sigma,\rho}\br{\1_E,\1_E}
    }
    =
    \int_{\R^2}
        \int
            \1_E\br{\vx}
            \1_E\br{x_1+M_{\sigma,1}\gamma_{\sigma,1}\br{t},x_2}
            \1_E\br{x_1,x_2+M_{\sigma,2}\gamma_{\sigma,2}\br{t}}
            \rho\br{t}
        dt
    d\vx,
\end{equation}
where we are given a family of positive diagonal matrices \(\BR{\vM_\sigma}_{\sigma\in\Sigma}\), a family of functions \(\BR{\vgamma_\sigma:I\to\R^2}_{\sigma\in\Sigma}\), and a fixed non-negative \(\rho\in C^2_c\br{I}\). We aim to identify properties of the pairings \(\BR{\br{\vM_\sigma,\vgamma_\sigma}}_{\sigma\in\Sigma}\) that guarantee the existence of \(\sigma\in\Sigma\) such that \eqref{eq_counting_inner_prod_sigma_gen} is non-trivial.

To proceed, we begin with a closer inspection of \eqref{eq_counting_inner_prod_sigma_gen} for a fixed \(\sigma\in\Sigma\). For a moment, we omit the \(\sigma\) dependence and simply write
\begin{equation}\label{eq_counting_inner_prod_gen}
    \ang{
        \1_E,\cA_{\vM\vgamma,\rho}\br{\1_E,\1_E}
    }=
    \int_{\R^3}
            \1_E\br{x_1,x_2}
            \1_E\br{x_1+M_1\gamma_1\br{t},x_2}
            \1_E\br{x_1,x_2+M_2\gamma_2\br{t}}
            \rho\br{t}
        dt
    dx_1dx_2.
\end{equation}
Standard consideration suggests that we decompose \eqref{eq_counting_inner_prod_gen} into the three parts
\begin{equation*}
    \ang{
        \1_E,\cA_{\vM\vgamma,\rho}\br{\1_E,\1_E}
    }
    =I_L+I_M+I_H,
\end{equation*}
where the exact definition of the low-frequency part \(I_L\), the intermediate frequency part \(I_M\), and the high frequency part \(I_H\) will be specified in a moment. A standard strategy to show the non-triviality of \eqref{eq_counting_inner_prod_gen} is to demonstrate the following inequality
\begin{equation}\label{eq_counting_inner_prod_gen_nontriv}
    \abs{\ang{
        \1_E,\cA_{\vM\vgamma,\rho}\br{\1_E,\1_E}
    }}\geq \abs{I_L}-\abs{I_M}-\abs{I_H}>0.
\end{equation}
This strategy hinges on a non-triviality statement for \(\abs{I_L}\) and the smallness of \(\abs{I_M},\abs{I_H}\). 

\subsubsection{Non-triviality of \(\abs{I_L}\)}
We begin with a paraphrase of \cite[\textbf{Lemma 5.1}]{MR4295087}.
\begin{lemma}\label{lem_Bourgain_pos}
Recall \eqref{eq_phi_uni_even} and denote \(\Phi_s g:=g\ast\Dil^1_s\phi/\nrm{\phi}_{L^1}\) for \(g\in L^1_{\operatorname{loc}}\br{\R}\).
There is a universal constant \(c_0>0\) such that for any non-negative \(f \in L^1\br{\R^2}\) with \(\supp f\subseteq \Br{0,1}^2\), the uniform lower bounds hold
\begin{equation}\label{eq_Bourgain_pos}
    \int_{\R^2}
        f\br{\vx}
        \Phi^{\br{1}}_{s_1}f\br{\vx}
        \Phi^{\br{2}}_{s_2}f\br{\vx}
    d\vx
    \geq
    c_0
    \nrm{f}^3_{L^1},\quad \forall s_j\in\bR{0,1}.
\end{equation}
\end{lemma}
Given \textbf{Lemma \ref{lem_Bourgain_pos}} and parameters \(R_j\in\bR{0,1}\) to be decided later, we set
\begin{equation}\label{eq_def_I_L_R_j}
    I_L:=\ang{
        \1_E,\cA_{\vM\vgamma,\rho}
        \br{
            \Phi^{\br{1}}_{R_1}\1_E,
            \Phi^{\br{2}}_{R_2}\1_E
        }
    }
\end{equation}
and compare \eqref{eq_def_I_L_R_j} with the quantity
\begin{equation}\label{eq_Bourgain_pos_E}
    \int_{\R^2}
        \1_E \br{\vx}
        \Phi^{\br{1}}_{R_1}\1_E\br{\vx}
        \Phi^{\br{2}}_{R_2}\1_E\br{\vx}
    d\vx.
\end{equation}
On the one hand, invoking the fact for \(f\in L^1_{\operatorname{loc}}\br{\R}\) that
\begin{equation*}
    \Phi_R f\br{x+z}-\Phi_R f\br{x}
    =
    \frac{z}{R\nrm{\phi}_{L^1}}
    \int^1_0
        f\ast \Dil^1_R\phi'\br{x+\theta z}
    d\theta,
\end{equation*}
we deduce via telescoping the following error estimate
\begin{equation*}
    \abs{\eqref{eq_def_I_L_R_j}- \nrm{\rho}_{L^1}
    \eqref{eq_Bourgain_pos_E}}
    \leq
    \int
        \br{
        \frac{M_1\abs{\gamma_1\br{t}}}{R_1\nrm{\phi}_{L^1}}+
        \frac{M_2\abs{\gamma_2\br{t}}}{R_2\nrm{\phi}_{L^1}}
        }\rho\br{t}
    dt
    \cdot
    \nrm{\phi'}_{L^1}\abs{E}
\end{equation*}
\begin{equation*}
    \leq
    \br{\frac{M_1}{R_1}+\frac{M_2}{R_2}}
    \nrm{\vgamma}_{L^\infty\br{I}}
    \br{\nrm{\phi'}_{L^1}/\nrm{\phi}_{L^1}}
    \nrm{\rho}_{L^1}\abs{E}.
\end{equation*}
On the other hand, \textbf{Lemma \ref{lem_Bourgain_pos}} yields the lower bound \(\eqref{eq_Bourgain_pos_E}\geq c_0\abs{E}^3\). Consequently, the lower bound
\begin{equation}\label{eq_I_L_lower_bd}
    \abs{I_L}\geq c_0\nrm{\rho}_{L^1}\abs{E}^3/2>0
\end{equation}
holds whenever \eqref{eq_def_I_L_R_j} can be taken with parameters satisfying
\begin{equation}\label{eq_I_L_pos_conds}
    R_j:=M_j/\lambda_L,\quad
    \max\br{M_1,M_2} \leq \lambda_L \leq \frac{c_0\nrm{\phi}_{L^1}\abs{E}^2}{4\nrm{\vgamma}_{L^\infty\br{I}}\nrm{\phi'}
    _{L^1}}.
\end{equation}
As a brief remark, the above condition also constrains the sizes of \(M_j\) that enables the aforementioned strategy.

\subsubsection{Smallness of \(\abs{I_H}\)}
Let again \(r_j\in\bR{0,1}\) to be decided later. Suggested by the heuristic
\begin{equation*}
    \widehat{\cQ_\lambda f}\br{\xi}=
    \varphi\br{\xi/\lambda}\widehat{f}\br{\xi}
    \approx
    \br{\widehat{\phi}\br{\xi/\lambda}/\widehat{\phi}\br{0}}\widehat{f}\br{\xi}
    =
    \widehat{\Phi_{1/\lambda}f}\br{\xi},
\end{equation*}
we shall formulate \(I_H\) in the following manner
\begin{equation}\label{eq_I_H_r_j}
    I_H:=\ang{
        \1_E,
        \cA_{\vM_\vgamma,\rho}\br{
            \1_E,
            \1_E
        }
        -
        \cA_{\vM_\vgamma,\rho}\br{
            \Phi^{\br{1}}_{r_1}\1_E,
            \Phi^{\br{2}}_{r_2}\1_E
        }
    }
\end{equation}
to cancel most of the low-frequency contribution. With the Sobolev smoothing inequality \eqref{eq_soblev_smoothing_std} associated with \(\cA_{\vgamma,\rho}\), we expect to derive power saving via making suitable choice of \(r_j\).

To justify the above heuristic, let \(\Psi_s:= \operatorname{Id}-\Phi_s\) and consider the following telescoping identity
\begin{equation}\label{eq_I_H_Avg_tele}
    \cA_{\vM_\vgamma,\rho}\br{
        \1_E,
        \1_E
    }
    -
    \cA_{\vM_\vgamma,\rho}\br{
        \Phi^{\br{1}}_{r_1}\1_E,
        \Phi^{\br{2}}_{r_2}\1_E
    }
    =
    \cA_{\vM\vgamma,\rho}\br{\Psi^{\br{1}}_{r_1}\1_E,\1_E}
    +
    \cA_{\vM\vgamma,\rho}\br{
        \Phi^{\br{1}}_{r_1}\1_E,
        \Psi^{\br{2}}_{r_2}\1_E
    }.
\end{equation}
To utilize \eqref{eq_soblev_smoothing_std} associated with \(\cA_{\vgamma,\rho}\), we invoke \eqref{eq_A_Mk_2_Dil_A} and an identities for \(\Phi_s\) and \(\Psi_s\) analogous to \eqref{eq_HF_dec_LP_Dil_com}
\begin{equation*}
    \Dil^2_{\vM^{-1}}\Phi^{\br{j}}_{r_j}\1_E=\Phi^{\br{j}}_{r_j/M_j}\Dil^2_{\vM^{-1}}\1_E,\quad
    \Dil^2_{\vM^{-1}}\Psi^{\br{j}}_{r_j}\1_E=\Psi^{\br{j}}_{r_j/M_j}\Dil^2_{\vM^{-1}}\1_E
\end{equation*}
to deduce the following chain of estimates
\begin{equation*}
    \abs{I_H}\leq
    \nrm{
        \eqref{eq_I_H_Avg_tele}
    }_{L^1}
    =
    \nrm{
        \cA_{\vgamma,\rho}
        \br{
            \Psi^{\br{1}}_{r_1/M_1}\Dil^2_{\vM^{-1}}\1_E,
            \Dil^2_{\vM^{-1}}\1_E
        }
    +
    \cA_{\vgamma,\rho}
        \br{
            \Phi^{\br{1}}_{r_1/M_1}\Dil^2_{\vM^{-1}}\1_E,
            \Psi^{\br{2}}_{r_2/M_2}\Dil^2_{\vM^{-1}}\1_E
        }
    }_{L^1}
\end{equation*}
\begin{equation}\label{eq_I_H_soblev_bd}
    \leq
    C_{\vgamma,\rho,\varepsilon}
    \abs{E}^{1/2}
    \br{
    \nrm{
        \Psi^{\br{1}}_{r_1/M_1}\Dil^2_{\vM^{-1}}\1_E
    }_{H^{\br{-\varepsilon,0}}}
    +
    \nrm{
        \Psi^{\br{2}}_{r_2/M_2}\Dil^2_{\vM^{-1}}\1_E
    }_{H^{\br{0,-\varepsilon}}}
    }
\end{equation}
To obtain the key power saving, it remains to study the Sobolev norm expression \(\nrm{\Psi_s f}_{H^{-\varepsilon}}\) for \(f\in L^2\br{\R}\). Indeed, Plancherel theorem yields the following estimate 
\begin{equation}\label{eq_Psi_sobolev_Plancherel}
    \nrm{\Psi_s f}_{H^{-\varepsilon}}
    =\nrm{
        \ang{\xi}^{-\varepsilon}
        \br{1-\widehat{\phi}\br{s\xi}/\widehat{\phi}\br{0}}
        \widehat{f}\br{\xi}
    }_{L^2\br{d\xi}}
    \leq
    \nrm{
        \frac{\widehat{\phi}\br{0}-\widehat{\phi}\br{s\xi}}{
        \ang{\xi}^\varepsilon\widehat{\phi}\br{0}
        }
    }_{L^\infty\br{d\xi}}
    \nrm{f}_{L^2}.
\end{equation}
Furthermore, we combine the trivial estimate with the estimate provided by mean value theorem
\begin{equation*}
    \abs{\widehat{\phi}\br{0}-\widehat{\phi}\br{\xi}}\leq \min\br{2\big\Vert\widehat{\phi}\big\Vert_{L^\infty}, \abs{\xi}\cdot\big\Vert\widehat{\phi}'\big\Vert_{L^\infty}}
    \leq
    \min\br{2\widehat{\phi}\br{0}, \abs{\xi}\cdot\big\Vert\widehat{\phi}'\big\Vert_{L^\infty}}
\end{equation*}
to further dominate the factor in \eqref{eq_Psi_sobolev_Plancherel} with
\begin{equation*}
    \nrm{
        \frac{\widehat{\phi}\br{0}-\widehat{\phi}\br{s\xi}}{
        \ang{\xi}^\varepsilon\widehat{\phi}\br{0}
        }
    }_{L^\infty\br{d\xi}}
    \leq
    \max
    \br{
        2
        \sup_{\abs{\xi}>\frac{2\widehat{\phi}\br{0}}{s\Vert\widehat{\phi}'\Vert_{L^\infty}}}
        \ang{\xi}^{-\varepsilon}
        ,
        \frac{s
        \big\Vert\widehat{\phi}'\big\Vert_{L^\infty}}{\widehat{\phi}\br{0}}
        \sup_{\abs{\xi}\leq \frac{2\widehat{\phi}\br{0}}{s\Vert\widehat{\phi}'\Vert_{L^\infty}}
        }
        \abs{\xi}/\ang{\xi}^\varepsilon
    }
    =2\ang{
        \frac{2\widehat{\phi}\br{0}}{s\Vert\widehat{\phi}'\Vert_{L^\infty}}
    }^{-\varepsilon}
\end{equation*}
In combination, we conclude
\begin{equation*}
    \nrm{\Psi_s f}_{H^{-\varepsilon}}
    \leq
    2\ang{
        \frac{2\widehat{\phi}\br{0}}{s\Vert\widehat{\phi}'\Vert_{L^\infty}}
    }^{-\varepsilon}
    \nrm{f}_{L^2}.
\end{equation*}
Applying the above estimate to \eqref{eq_I_H_soblev_bd} yields
\begin{equation*}
    \abs{I_H}\leq
    4C_{\vgamma,\rho,\varepsilon}
    \abs{E}
    \max_{j=1,2}
        \ang{\frac{2M_j\nrm{\phi}_{L^1}}{r_j\big\Vert\widehat{\phi}'\big\Vert_{L^\infty}}}^{-\varepsilon}
    .
\end{equation*}
In view of the above and the formerly established lower bound \eqref{eq_I_L_lower_bd}, we shall take
\begin{equation}\label{eq_I_H_small_cond}
    r_j:= M_j/\lambda_H,\quad \lambda_H \geq 
    \frac{1}{2}
    \br{\frac{
        16 C_{\vgamma,\rho,\varepsilon}
    }{
        c_0\nrm{\rho}_{L^1}\abs{E}^2
    }}^{1/\varepsilon}
    \big\Vert\widehat{\phi}'\big\Vert_{L^\infty}/\nrm{\phi}_{L^1}
\end{equation}
such that \(\abs{I_H}\leq\abs{I_L}/2\) and thus, the following lower bound holds
\begin{equation}\label{eq_I_L_I_H_diff_low_bd}
    \abs{I_L}-\abs{I_H}\geq c_0\nrm{\rho}_{L^1}\abs{E}^3/4.
\end{equation}

\subsubsection{Energy control on \(\abs{I_M}\)}
Let \(\lambda_L,\lambda_H\) be chosen such that \(R_j,r_j\) satisfy \eqref{eq_I_L_pos_conds} and \eqref{eq_I_H_small_cond} respectively. It remains to examine the intermediate frequency contribution
\begin{equation*}
    I_M:=\ang{\1_E,\cA_{\vM\vgamma,\rho}\br{\1_E,\1_E}}-I_L-I_H
\end{equation*}
\begin{equation}\label{eq_I_M_def}
    =
    \ang{
        \1_E,
        \cA_{\vM\vgamma,\rho}\br{
            \Phi^{\br{1}}_{M_1/\lambda_H}\1_E,
            \Phi^{\br{2}}_{M_2/\lambda_H}\1_E
        }
        -
        \cA_{\vM\vgamma,\rho}\br{
            \Phi^{\br{1}}_{M_1/\lambda_L}\1_E,
            \Phi^{\br{2}}_{M_2/\lambda_L}\1_E
        }
    }.
\end{equation}
In the unlikely scenarios, we could have 
\begin{equation}\label{eq_I_M_bd_goal}
    \abs{I_M}<c_0\nrm{\rho}_{L^1}\abs{E}^3/4
\end{equation} 
and deduce via \eqref{eq_counting_inner_prod_gen_nontriv} and \eqref{eq_I_L_I_H_diff_low_bd} the non-triviality of \eqref{eq_counting_inner_prod_gen}.
% This would demonstrate the presence of pattern:
% \begin{equation*}
%     \br{x_1,x_2},\quad\br{x_1+M_1\gamma_1\br{t},x_2},\quad\br{x_1,x_2+M_2\gamma_2\br{t}}\in E,\quad t\in \supp \rho.
% \end{equation*}
In most cases however, the upper bound \eqref{eq_I_M_bd_goal} is too demanding. Nonetheless, a trivial energy control on \(I_M\) is available. Observe the telescoping identity
\begin{equation*}
    \cA_{\vM\vgamma,\rho}\br{
        \Phi^{\br{1}}_{M_1/\lambda_H}\1_E,
        \Phi^{\br{2}}_{M_2/\lambda_H}\1_E
    }
    -
    \cA_{\vM\vgamma,\rho}\br{
        \Phi^{\br{1}}_{M_1/\lambda_L}\1_E,
        \Phi^{\br{2}}_{M_2/\lambda_L}\1_E
    }
\end{equation*}
\begin{equation*}
    =
    \cA_{\vM\vgamma,\rho}\br{
        \br{\Phi^{\br{1}}_{M_1/\lambda_H}-\Phi^{\br{1}}_{M_1/\lambda_L}}\1_E,
        \Phi^{\br{2}}_{M_2/\lambda_H}\1_E
    }
    +
    \cA_{\vM\vgamma,\rho}\br{
        \Phi^{\br{1}}_{M_1/\lambda_L}\1_E,
        \br{\Phi^{\br{2}}_{M_2/\lambda_H}-\Phi^{\br{2}}_{M_2/\lambda_L}}\1_E
    }
\end{equation*}
Applying the above identity to \eqref{eq_I_M_def} yields the following trivial bound
\begin{equation*}
    \abs{I_M}
    \leq \abs{E}^{1/2}\nrm{\rho}_{L^1}\br{
        \nrm{
            \br{\Phi^{\br{1}}_{M_1/\lambda_H}-\Phi^{\br{1}}_{M_1/\lambda_L}}\1_E
        }_{L^2}
        +
        \nrm{
            \br{\Phi^{\br{2}}_{M_2/\lambda_H}-\Phi^{\br{2}}_{M_2/\lambda_L}}\1_E
        }_{L^2}
    }.
\end{equation*}
Notice that the operator \(\Phi_A-\Phi_B\) possesses a Fourier multiplier representation with symbol \(m_{A,B}\) given by
\begin{equation*}
    m_{A,B}\br{\xi}:=
    \br{\widehat{\phi}\br{A\xi}-\widehat{\phi}\br{B\xi}}/\widehat{\phi}\br{0}
    =
    \int^A_B \frac{s\xi\widehat{\phi}'\br{s\xi}}{\widehat{\phi}\br{0}}\frac{ds}{s}.
\end{equation*}
With little effort, we deduce from the trivial estimate
\begin{equation*}
    \abs{m_{A,B}\br{\xi}}
    \leq 
    \min\br{
        2,
        \int^{\max\br{A,B}}_{\min\br{A,B}} \abs{s\xi\widehat{\phi}'\br{s\xi}/\widehat{\phi}\br{0}}\frac{ds}{s}
    }
\end{equation*}
the following more refined energy control on \(I_M\)
\begin{equation}\label{eq_I_M_square_bd}
    \abs{I_M}^2
    \leq 8 \abs{E}\nrm{\rho}^2_{L^1}
    \sum_{j=1,2}
    \int_{\R^2}
        \abs{\widehat{\1_E}\br{\vxi}}^2
            \int^{M_j/\lambda_L}_{M_j/\lambda_H} 
            \abs{
                s\xi_j\widehat{\phi}'\br{s\xi_j}/
                \widehat{\phi}\br{0}
            }
            \frac{ds}{s}
    d\vxi,
\end{equation}
which prepares us for the next step.

\subsubsection{Bourgain's energy pigeonholing argument}
Instead of aiming to show \eqref{eq_I_M_bd_goal} for a single pairing \(\br{\vM_\sigma,\vgamma_\sigma}=\br{\vM,\vgamma}\), we shall consider the whole families of \(\BR{\br{\vM_\sigma,\vgamma_\sigma}}_{\sigma\in\Sigma}\) and provide weaker conditions such that \eqref{eq_I_M_bd_goal} holds on average for certain sub-family \(\BR{\br{\vM_\omega,\vgamma_\omega}}_{\omega\in \Omega}\). This would demonstrate the existence of \(\sigma=\omega\in \Omega\subseteq \Sigma\) such that the inner product \eqref{eq_counting_inner_prod_sigma_gen} is non-trivial.

Following previous discussion, we shall decompose \eqref{eq_counting_inner_prod_sigma_gen} into \(I_{\sigma,L}\), \(I_{\sigma,M}\), and \(I_{\sigma,H}\), by taking \(R_{\sigma,j}:=M_{\sigma,j}/\lambda_{\sigma,L}\) and \(r_{\sigma,j}:=M_{\sigma,j}/\lambda_{\sigma,H}\) with \(\lambda_{\sigma,L}\) and \(\lambda_{\sigma,H}\) satisfying \eqref{eq_I_L_pos_conds} and \eqref{eq_I_H_small_cond} respectively
\begin{equation}\label{eq_lambda_L_H_conds}
    \max_{j=1,2} M_{\sigma,j} \leq \lambda_{\sigma,L} \leq \frac{c_0\nrm{\phi}_{L^1}\abs{E}^2}{4\nrm{\vgamma_\sigma}_{L^\infty\br{I}}\nrm{\phi'}
    _{L^1}},\quad
    \lambda_{\sigma,H} \geq 
    \frac{1}{2}
    \br{\frac{
        16 C_{\vgamma_\sigma,\rho,\varepsilon_\sigma}
    }{
        c_0\nrm{\rho}_{L^1}\abs{E}^2
    }}^{1/\varepsilon_\sigma}
    \big\Vert\widehat{\phi}'\big\Vert_{L^\infty}/\nrm{\phi}_{L^1}
\end{equation}
such that we may leverage the energy estimates \eqref{eq_I_M_square_bd} to obtain for any probability measure \(d\mathbb{P}_\Omega\br{\omega}\) on \(\Omega\subseteq \Sigma\), the following expectation bound
\begin{equation*}
    \E_\Omega
    \abs{I_{\omega,M}}^2
    \leq 8 \abs{E}\nrm{\rho}^2_{L^1}
    \sum_{j=1,2}
    \int_{\R^2}
        \abs{\widehat{\1_E}\br{\vxi}}^2
            \int^\infty_0
            \E_\Omega
            \1_{\bR{M_{\omega,j}/\lambda_{\omega,H}, M_{\omega,j}/\lambda_{\omega,L}}}\br{s}
            \abs{
                s\xi_j\widehat{\phi}'\br{s\xi_j}/
                \widehat{\phi}\br{0}
            }
            \frac{ds}{s}
    d\vxi
\end{equation*}
\begin{equation*}
    \leq 16 \abs{E}^2\nrm{\rho}^2_{L^1}
    \br{\big\Vert\widehat{\phi}'\big\Vert_{L^1}/
    \nrm{\phi}_{L^1}}
    \sup_{\substack{0<s<\infty\\j=1,2}}
        \mathbb{P}_\Omega\br{M_{\omega,j}/\lambda_{\omega,H}<s\leq M_{\omega,j}/\lambda_{\omega,L}}.
\end{equation*}
Compare now the above with the desired estimate \eqref{eq_I_M_bd_goal}. Once we have the following
\begin{equation}\label{eq_counting_cond}
    \sup_{\substack{0<s<\infty\\j=1,2}}
    \mathbb{P}_\Omega\br{M_{\omega,j}/\lambda_{\omega,H}<s\leq M_{\omega,j}/\lambda_{\omega,L}}
    <
    \frac{c^2_0\abs{E}^4\nrm{\phi}_{L^1}}{256\big\Vert\widehat{\phi}'\big\Vert_{L^1}},
\end{equation}
we may infer the existence of \(\sigma=\omega\in \Omega\subseteq \Sigma\) such that \(\abs{I_{\sigma,M}}\) satisfies \eqref{eq_I_M_bd_goal}, and thus, the inner product \eqref{eq_counting_inner_prod_sigma_gen} is non-trivial.
% In fact, we may refine the above statement by conditioning the probability to a subset \(\Omega\subseteq \Sigma\).
In other words, with \eqref{eq_lambda_L_H_conds}, the condition \eqref{eq_counting_cond}
% , a condition of the following form
% \begin{equation}\label{eq_counting_cond_restricted}
%     \sup_{\substack{0<s<\infty\\j=1,2}}
%     \mathbb{P}\br{
%         M_{\sigma,j}/\lambda_{\sigma,H}<s\leq M_{\sigma,j}/\lambda_{\sigma,L}
%         \mid
%         \sigma\in \Omega
%     }
%     <
%     \frac{c^2_0\abs{E}^4\nrm{\phi}_{L^1}}{256\big\Vert\widehat{\phi}'\big\Vert_{L^1}}
% \end{equation}
guarantees the existence of \(\omega\in\Omega\subseteq \Sigma\), \(t\in \supp \rho\), and \(\br{x_1,x_2}\in \R^2\) such that the following corner pattern exists
\begin{equation}\label{eq_corner_pattern_omega}
    \br{x_1,x_2},\quad
    \br{x_1+M_{\omega,1}\gamma_{\omega,1}\br{t},x_2},\quad
    \br{x_1,x_2+M_{\omega,2}\gamma_{\omega,2}\br{t}}\in E.
\end{equation}

\subsubsection{Proof of \textbf{Theorem \ref{thm_pattern}}}
To adopt the general strategy displayed above, the first task is to choose the representing pairings \(\BR{\br{\vM_\sigma,\vgamma_\sigma}}_{\sigma\in\br{0,1}}\) for the family of functions \(\BR{\vgamma\br{\sigma\cdot}}_{\sigma\in\br{0,1}}\). On the one hand, due to the analyticity of \(\vgamma\) and the fact that \(\vgamma\br{0}=\vnull\), we may assume without loss of generality that\footnote{Here, a change of variable is performed to rescale the domain of \(\vgamma\) for convenience.}
\begin{equation}\label{eq_gamma_j_taylor}
    \gamma_j\br{t}=a_j t^{n_j}+g_j\br{t},\quad
    \abs{g^{\br{k}}_j\br{t}}\leq C t^{n_j+\delta-k},\quad \forall k\in\BR{0,1,\dots,5},\,\forall t\in\Br{0,2}
\end{equation}
for some fixed \(a_j\neq 0\), \(1\leq n_1\leq n_2 \in\N\), and \(C,\delta>0\). As a direct consequence, we have
\begin{equation*}
    \gamma_j\br{\sigma t}= \sigma^{n_j}\br{a_jt^{n_j}+\sigma^{-n_j}g_j\br{\sigma t}},
\end{equation*}
with the error term satisfying
\begin{equation*}
    \abs{\partial^k_t \sigma^{-n_j}g_j\br{\sigma t}}\leq 
    \sigma^\delta C t^{n_j+\delta-k}.
\end{equation*}
This suggests the natural setting to achieve \(\vM_\sigma\vgamma_\sigma\br{t}=\vgamma\br{\sigma t}\) would be to consider
\begin{equation*}
    \vM_\sigma:=
    \begin{pmatrix}
        M_{\sigma,1} & 0\\
        0 & M_{\sigma,2}
    \end{pmatrix}
    :=
    \begin{pmatrix}
        \sigma^{n_1} & 0\\
        0 & \sigma^{n_2}
    \end{pmatrix}
    ,\quad
    \vgamma_\sigma: \br{1/2,2}\to\R^2:
    t\mapsto\vM^{-1}_\sigma\vgamma\br{\sigma t}.
\end{equation*}
It is an easy computation to show that, for some fixed constant \(C_{\vgamma,L}\in\br{0,1}\), the following holds
\begin{equation*}
    C_{\vgamma,L}\abs{E}^2
    \leq
    \frac{
        c_0\nrm{\phi}_{L^1}\abs{E}^2
    }{
        4\nrm{\vgamma_\sigma}_{L^\infty}\nrm{\phi'}_{L^1}
    },\quad \forall\sigma\in\br{0,1}.
\end{equation*}
We may thus make the following choice to fulfill the left-hand side of \eqref{eq_lambda_L_H_conds}
\begin{equation}\label{eq_lambda_L_conds_ana}
    \lambda_{\sigma,L}:=C_{\vgamma,L}\abs{E}^2,\quad \forall \sigma \leq 
    C^{\frac{1}{n_1}}_{\vgamma,L}\abs{E}^{\frac{2}{n_1}}.
\end{equation}
On the other hand, by the assumption that none of the four classes of forbidden sets \eqref{eq_line}, \eqref{eq_exp}, \eqref{eq_log}, and \eqref{eq_psi_line} contains \(\img \vgamma\), we may proceed as discussed in \textsc{Section \ref{subsec_pf_thm_real_ana_smoothing}} to obtain a lower bound similar to \eqref{eq_ana_bad_at_0}
\begin{equation}\label{eq_ana_bad_at_0_4_pattern}
    \abs{\gamma'_j}\br{t},\,\abs{\vgamma'\wedge \vgamma''}\br{t},\,\abs{\vchi'_\vgamma\wedge\vchi''_\vgamma}\br{t}\gtrsim t^n,\quad \forall t\in\Br{0,2},\,j=1,2
\end{equation}
with some \(n=n_\vgamma\in\N\sqcup\BR{0}\). To simplify the discussion, we may take \(n>0\). By direct computation, we infer that
\(\vgamma_\sigma\) is \(O\br{\sigma^{-n}}\)-admissible.
Using also the fact that \(\BR{\vgamma_\sigma:\br{1/2,2}\to\R^2}_{\sigma\in\br{0,1}}\) forms a definable family in \(\R_{\operatorname{an}}\), we deduce via \textbf{Theorem \ref{thm_smoothing_o_mini_uni}} the estimate on the implicit constant for \eqref{eq_soblev_smoothing_std} associated to \(\cA_{\vgamma_\sigma,\psi_{1/2}}\)
\begin{equation*}
    C_{\vgamma_\sigma,\psi_{1/2},\varepsilon}\underset{\vgamma,\psi_{1/2}}{\lesssim} \sigma^{-cn},\quad \forall \sigma\in\br{0,1}.
\end{equation*}
With further computation, we may show for a fixed constant \(C_{\vgamma,H}\in\br{1,\infty}\) that
\begin{equation*}
    \frac{1}{2}\br{
        \frac{
            16C_{\vgamma_\sigma,\psi_{1/2},\varepsilon}
        }{
            c_0\nrm{\psi_{1/2}}_{L^1}\abs{E}^2
        }
    }^{1/\varepsilon}
    \big\Vert\widehat{\phi}'\big\Vert_{L^\infty}/\nrm{\phi}_{L^1}
    \leq 
    C_{\vgamma,H}\sigma^{-cn/\varepsilon}\abs{E}^{-2/\varepsilon}.
\end{equation*}
Finally, by taking
\begin{equation*}
    \lambda_{\sigma,H}:=C_{\vgamma,H}\sigma^{-cn/\varepsilon}\abs{E}^{-2/\varepsilon},\quad \forall \sigma\in\br{0,1},
\end{equation*}
we fulfill the right-hand side of \eqref{eq_lambda_L_H_conds}.
It remains to find a probability measure \(\mathbb{P}_\Omega\) on some subset \(\Omega\subseteq\br{0,1}\) to achieve \eqref{eq_counting_cond}. One way to simplify the matter would be to choose \(\Omega:=\BR{\omega_k}^N_{k=1}\subseteq\br{0,1}\) such that
\begin{equation*}
    \BR{
    I_{k,j}:=
        \bR{
            M_{\omega_k,j}/\lambda_{\omega_k,H},
            M_{\omega_k,j}/\lambda_{\omega_k,L}
        }
    }^N_{k=1}
\end{equation*}
form disjoint collections of intervals for both \(j=1,2\). We may further require \(\mathbb{P}_\Omega\) to assign equal probability on all \(\omega\in\Omega\). The condition \eqref{eq_counting_cond} thus suggests that we take
\begin{equation}\label{eq_counting_ana_bd_N}
    \left\lfloor
    \frac{
        256\big\Vert\widehat{\phi}'\big\Vert_{L^1}
    }{
        c^2_0\abs{E}^4\nrm{\phi}_{L^1}
    }
    \right\rfloor+1
    = \#\Omega=N.
\end{equation}
As such, we show the existence of the corner pattern \eqref{eq_corner_pattern_omega} for some \(\omega=\omega_k\) with \(k\leq N\). By unwrapping the identity, it justifies the existence of the corner pattern \eqref{eq_corner_pattern_in_E} for some \(t\eqsim \omega_k\). As for the gap estimate \eqref{eq_corner_pattern_gap_est}, we must show the constructed \(\Omega:=\BR{\omega_k}^N_{k=1}\) satisfies a lower bound of the form
\begin{equation}\label{eq_corner_omegas_gap_est}
    \min_{k\in\BR{1,\dots,N}}\omega_k\geq \exp\br{-\exp\br{C_\vgamma\abs{E}^{-c_\vgamma}}}.
\end{equation}

With the aforementioned plan, we shall design \(\Omega:=\BR{\omega_k}^N_{k=1}\subseteq\br{0,1}\) such that
\begin{equation}\label{eq_omegas_cond_4_pattern}
    \omega_{k+1}^{n_j}\abs{E}^{-2}/C_{\vgamma,L}=M_{\omega_{k+1},j}/\lambda_{\omega_{k+1},L}
    \leq
    M_{\omega_k,j}/\lambda_{\omega_k,H}=
    \omega_k^{n_j+cn/\varepsilon}\abs{E}^{2/\varepsilon}/C_{\vgamma,H}.
\end{equation}
To further simplify the discussion, we introduce the parameter \(\alpha\underset{\vgamma}{\eqsim}\abs{E}\) given by
\begin{equation}\label{eq_E_2_alpha_w_rel}
    \alpha:=\frac{C_{\vgamma,L}}{C_{\vgamma,H}}\cdot\abs{E}\leq \min\br{
        C_{\vgamma,L},C^{-1}_{\vgamma,H},\abs{E}
    }
    < 1
\end{equation}
and restrict \(\omega_k\) to the form \(\omega_k:=\alpha^{C_k}\). Suggests by \eqref{eq_lambda_L_conds_ana}, we shall take 
\begin{equation}\label{eq_omega_1_4_pattern}
    \omega_1:=\alpha^{C_1}\leq C^{\frac{1}{n_1}}_{\vgamma,L}\abs{E}^{\frac{2}{n_1}},\quad
    C_1:=\frac{3}{n_1}
\end{equation}
as the starting point of our inductive process. Next, substituting \(\omega_k:=\alpha^{C_k}\) turns the requirement \eqref{eq_omegas_cond_4_pattern} into
\begin{equation*}
    \alpha^{C_{k+1}n_j}\leq \alpha^{C_k\br{n_j+\frac{cn}{\varepsilon}}}\cdot \frac{C_{\vgamma,L}}{C_{\vgamma,H}}\cdot \abs{E}^{2+\frac{2}{\varepsilon}}.
\end{equation*}
Using \eqref{eq_E_2_alpha_w_rel}, we may insert an intermediate term to simplify the computation
\begin{equation*}
    \alpha^{C_{k+1}n_j}\leq 
    \alpha^{C_k\br{n_j + \frac{cn}{\varepsilon}}+ 2+ \frac{2}{\varepsilon}}\leq
    \alpha^{C_k\br{n_j + \frac{cn}{\varepsilon}}}\cdot \frac{C_{\vgamma,L}}{C_{\vgamma,H}}\cdot \abs{E}^{2+\frac{2}{\varepsilon}}.
\end{equation*}
To make sure the left-hand side holds for both \(j=1,2\), we require 
\begin{equation}\label{eq_C_k_induct_4_pattern}
    C_{k+1} := C_k\br{1+\frac{cn}{\varepsilon n_1}}+\frac{2}{n_1}+\frac{2}{\varepsilon n_1}.
\end{equation}
Solving the recurrence relation given by \eqref{eq_omega_1_4_pattern} and \eqref{eq_C_k_induct_4_pattern}, we conclude
\begin{equation}\label{eq_C_k_close_form_4_pattern}
    C_k = \br{1+\frac{cn}{\varepsilon n_1}}^{k-1}\br{\frac{3}{n_1}+\frac{2+2\varepsilon}{cn}}-\frac{2+2\varepsilon}{cn}
    =O\br{\br{1+\frac{cn}{\varepsilon n_1}}^k}.
\end{equation}
Finally, with \eqref{eq_C_k_close_form_4_pattern} and \eqref{eq_counting_ana_bd_N}, we conclude \eqref{eq_corner_omegas_gap_est}
\begin{equation*}
    \min_{k\in\BR{1,\dots,N}}\omega_k
    =\alpha^{C_N}
    =\br{\frac{C_{\vgamma,L}}{C_{\vgamma,H}}\cdot \abs{E}}^{O\br{\br{1+\frac{cn}{\varepsilon n_1}}^{O\br{\abs{E}^{-4}}}}}
    \geq \exp\br{-\exp\br{C_\vgamma\abs{E}^{-c_\vgamma}}}
\end{equation*}
for some \(C_\vgamma,\,c_\vgamma>0\). This completes the proof of \textbf{Theorem \ref{thm_pattern}}.

% \phantomsection\section[phantom]{}

% \appendixautorefname
\appendix

\section{Van der Corput lemma}\label{sec_van_der_corput}

\begin{proposition}[Van der Corput lemma]\label{prop_oscillatory}
Let $\varphi : [a,b] \rightarrow \mathbb{R}$ be a $C^{k}$ function and suppose that $\varphi^{(k)}(x) \neq 0$ for some $k\geq 1$ and for all $x\in [a,b]$. If $k=1$, we further assume that $\varphi'$ is monotonic. Let $\psi:[a,b]\rightarrow \mathbb{R}$ be a $C^{1}$ function. Then we have the estimate
    \begin{equation}\label{oscillatory_est2}
        \abs{
            \int^b_a
                e\br{\varphi\br{t}}
                \psi\br{t}
            dt
        }
        \underset{k}{\lesssim}
        \frac{1}{\inf_{t\in\mr{a,b}} \abs{\varphi^{\br{k}}}^{\frac{1}{k}}}\cdot 
        \br{
            \nrm{\psi}_{L^\infty\br{\mr{a,b}}}
            +\nrm{\psi'}_{L^1\br{\mr{a,b}}}
        }.
    \end{equation}
\end{proposition}

Recall that we denote \(C\br{B}\) the number of connected components of \(B\subseteq \R\) and set \(C\br{\varnothing}=0\).
\begin{remark}[Basic facts about connected components in \(\R\)]\label{rmk_connect_rel}
    Given sets \(A,B\subseteq \R\), it's easy to show that:
    \begin{equation*}
        C\br{A\cup B},C\br{A\cap B}\leq C\br{A}+C\br{B}
        ,\quad
        C\br{A^c}\leq 1+C\br{A}.
    \end{equation*}
    % If \(A,B\) are disconnected from each other, we have the equality:
    % \begin{equation*}
    %     C\br{A\sqcup B}=C\br{A}+C\br{B}.
    % \end{equation*}
\end{remark}

\begin{definition}\label{def_M_phi}
    Given a function \(\varphi:\R\supset \dom \varphi \to \R\), we define the following two quantities:
    \begin{equation}\label{eq_M_phi}
        M_w\br{\varphi}:=
        C\br{
            \BR{
                t\in \dom \varphi
            \::\:
                \varphi\br{t}=w
            }
        },
        \quad
        M\br{\varphi}:=\sup_w
        M_w\br{\varphi}.
    \end{equation}
    Note that we always have \(1,M_w\br{\varphi}\leq M\br{\varphi}\) whenever \(\dom \varphi\neq \varnothing\).
\end{definition}

\begin{remark}[Relation to continuous functions]\label{rmk_connect_rel_to_func}
    Given \(\varphi:I \to \R\) continuous and \(w\in \R\), we have:
    \begin{equation*}
        C\br{
            \BR{
                t\in I
            \::\:
                \varphi\br{t}>w
            }
        }
        +
        C\br{
            \BR{
                t\in I
            \::\:
                \varphi\br{t}<w
            }
        }
        =
        C\br{
            \BR{
                t\in I
            \::\:
                \varphi\br{t}\neq w
            }
        }
        \leq 1+M_w\br{\varphi}.
    \end{equation*}
    % If \(A,B\) are disconnected from each other, we have the equality:
    % \begin{equation*}
    %     C\br{A\sqcup B}=C\br{A}+C\br{B}.
    % \end{equation*}
\end{remark}

\begin{proposition}\label{prop_M_0_to_M}
    Let \(\varphi: I\to \R\) be a \(C^1\) function. We have the following inequality:
    \begin{equation*}
        M\br{\varphi}\leq 2+ 2M_0\br{\varphi'}.
    \end{equation*}
\end{proposition}
\begin{proof}
    Without loss of generality, we assume \(M_0\br{\varphi'}<\infty\). 
    % Observe that we always have:
    % \begin{equation*}
    %     \max\br{1,M_0\br{\varphi'}}\leq M\br{\varphi'}.
    % \end{equation*}
    % The right-hand side inequality thus follows. It remains to prove the left-hand side inequality. 
    By \textbf{Remark \ref{rmk_connect_rel}}, we have:
    \begin{equation*}
        C\br{
            \BR{
                t\in I
            \::\:
                \varphi'\br{t}\neq 0
            }
        }
        =
        C\br{
            I\setminus
            \BR{
                t\in I
            \::\:
                \varphi'\br{t}=0
            }
        }
        \leq
        2+M_0\br{\varphi'}.
    \end{equation*}
    As a result, there is a partition of \(I\) consisting of \(N\leq 2+2M_0\br{\varphi'}\) many disjoint intervals \(\BR{I_j}^N_{j=1}\) such that on each \(I_j\), \(\varphi\) is either monotone or constant. As a result, for any \(w\in\R\), we deduce that:
    \begin{equation*}
        C\br{
            \BR{
                t\in I_j
            \::\:
                \varphi\br{t}=w
            }
        }
        \leq 1.
    \end{equation*}
    Apply \textbf{Remark \ref{rmk_connect_rel}} again. We deduce:
    \begin{equation*}
        M_w\br{\varphi}
        \leq 
        \sum^N_{j=1}
        C\br{
            \BR{
                t\in I_j
            \::\:
                \varphi\br{t}=w
            }
        }
        \leq N\leq 2+2M_0\br{\varphi'}.
    \end{equation*}
    Taking supremum over \(w\in\R\) gives the desired conclusion.
\end{proof}

\begin{proposition}[Reformulation of Van der Corput]\label{prop_VdC_reform}
    Let \(k\geq 2\). For any \(C^k\) function \(\varphi:\mr{a,b}\to\R\) and \(C^1\) function \(\psi:\mr{a,b}\to\R\), we have the following estimate:
    \begin{equation*}
        \abs{
            \int^b_a
                e\br{\varphi\br{t}}
                \psi\br{t}
            dt
        }
        \underset{k}{\lesssim}
        \frac{
            M\br{\varphi^{\br{k}}}
        }{
            \displaystyle{\inf_{t\in\mr{a,b}}
            \max_{1\leq j\leq k}}
            \abs{\varphi^{\br{j}}\br{t}}^{\frac{1}{j}}
        }
        \cdot 
        \br{
            \nrm{\psi}_{L^\infty\br{\mr{a,b}}}
            +\nrm{\psi'}_{L^1\br{\mr{a,b}}}
        }
        .
    \end{equation*}
\end{proposition}
\begin{proof}
    Let \(\lambda >0\) be such that:
    \begin{equation*}
        \lambda <
        \inf_{t\in\mr{a,b}}
            \max_{1\leq j\leq k}
                \abs{\varphi^{\br{j}}\br{t}}^{\frac{1}{j}}.
    \end{equation*}
    Consider the sets:
    \begin{equation*}
        A_\ast:=
        \BR{
            t\in \mr{a,b}
        \::\:
            \varphi''\br{t}=0
        },\quad
        A_j:=
        \BR{
            t\in \mr{a,b}
        \::\:
            \abs{\varphi^{\br{j}}\br{t}}>\lambda^j
        },\quad
        1\leq j \leq k.
    \end{equation*}
    By the choice of \(\lambda\), we have \(\mr{a,b}=\bigcup^k_{j=1}A_j\). We now construct the following sets:
    \begin{equation*}
        B_0:=A_1\cap A_\ast,\quad
        B_1:=A_1\setminus A_\ast,\quad 
        B_j:=A_j\setminus\bigcup_{i<j}A_i=
        A_j,\quad 
        1<j\leq k.
    \end{equation*}
    By design, we have \(\mr{a,b}=\bigsqcup^k_{j=0}B_j\).
    By \textbf{Remark \ref{rmk_connect_rel} and \ref{rmk_connect_rel_to_func}}, we deduce the following
    \begin{equation*}
        C\br{B_0},\, C\br{B_1}\leq 1+ C\br{A_1} + C\br{A_\ast} \leq 3+M_\lambda\br{\varphi'}+M_{-\lambda}\br{\varphi'}+M_0\br{\varphi''}
    \end{equation*}
    and for \(1<j\leq k\), the following
    \begin{equation*}
        C\br{B_j}=
        C\Big(A_j\setminus\bigcup_{i<j} A_i\Big)
        \leq C\br{A_j} +C\bigg(\Big(\bigcup_{i<j} A_i\Big)^c\bigg)
        \leq 1+ \sum^j_{i=1} C\br{A_i}
    \end{equation*}
    \begin{equation*}
        \leq 2j+ 1+ \sum^j_{i=1}
        M_{\lambda^i}\br{\varphi^{\br{i}}}+
        M_{-\lambda^i}\br{\varphi^{\br{i}}}.
    \end{equation*}
    Apply now \textbf{Proposition \ref{prop_M_0_to_M}}. We deduce the chain of relations:
    \begin{equation*}
        M_w\br{\varphi^{\br{j}}}
        \leq 
        M\br{\varphi^{\br{j}}}
        \leq
        2+2
        M_0\br{\varphi^{\br{j+1}}}
        \lesssim
        M\br{\varphi^{\br{j+1}}},\quad
        w\in\R, \; 1\leq j <k.
    \end{equation*}
    In combination with the trivial fact \(M\br{\varphi^{\br{k}}}\geq 1\), we establish the following estimates:
    \begin{equation*}
        C\br{B_j}\underset{k}{\lesssim} M\br{\varphi^{\br{k}}},\quad
        \forall j\in\BR{0,1,2,\dots, k}.
    \end{equation*}
    As a result, for each \(j\), we can decompose \(B_j\) into disjoint unions of \(N_j\lesssim M\br{\varphi^{\br{k}}}\) many intervals \(\BR{I_{j,l}}^{N_j}_{l=1}\).
    This allows us to apply the standard Van der Corput lemma \textbf{\ref{prop_oscillatory}} on each \(I_{j,l}\) and deduce:
    \begin{equation*}
        \abs{
        \int^b_a
            e\br{\varphi\br{t}}
            \psi\br{t}
        dt
        }
        \leq
        \sum^k_{j=0}
        \sum^{N_j}_{l=1}
        \abs{
        \int_{I_{j,l}}
            e\br{\varphi\br{t}}
            \psi\br{t}
        dt
        }
        \underset{k}{\lesssim}
        \frac{M\br{\varphi^{\br{k}}}}{\lambda}
        \br{\nrm{\psi}_{L^\infty\br{\mr{a,b}}}+\nrm{\psi'}_{L^1\br{\mr{a,b}}}},
    \end{equation*}
    where the issue with the monotonicity of \(\varphi'\) is addressed by the splitting of \(B_0\) and \(B_1\).
    Sending now \(\lambda\) to \(\displaystyle{\inf_{t\in\mr{a,b}}
            \max_{1\leq j\leq k}}
            \abs{\varphi^{\br{j}}\br{t}}^{\frac{1}{j}}\) completes the proof.
\end{proof}

\section{O-minimal theory}\label{sec_o_min}

\begin{definition}[Structure on \(M\)]\label{def_struct_over_M}
\(\cR:=\br{\cR_n}_{n\in\N}\) is a structure on a set \(M\neq\varnothing\) if:
\begin{enumerate}
    \item \(\cR_n\) is a boolean algebra of subsets of \(M^n\);
    \item \(A\in\cR_n\) implies \(A\times M,M\times A\in\cR_{n+1}\);
    \item \(\BR{\br{x_1,\dots,x_n}\in M^n\::\: x_1=x_n}\in \cR_n\);
    \item \(A\in \cR_{n+1}\) implies \(\pi\br{A}\in \cR_n\), where \(\pi\br{x_1,\dots,x_{n+1}}:=\br{x_1,\dots,x_n}\).
\end{enumerate}
\end{definition}

\begin{definition}[Definable sets, relations, and functions]
Given a structure \(\cR\) on \(M\), a set \(A\subseteq M^n\), an \(n\)-ary relation \(R\) on \(M\), and a function \(f:A\to M^m\), we say:
\begin{itemize}
    \item \(A\) is definable if \(A\in\cR_n\);
    \item \(R\) is definable if the set \(\BR{\br{x_1,\dots,x_n}\in M^n\::\: R\br{x_1,\dots,x_n}}\) is definable;
    \item \(f\) is definable if the \(n+m\)-ary relation given by \(f\br{x}=y\) is definable.
\end{itemize}
\end{definition}

\begin{definition}[Definable families]\label{def_definable_family}
    Given a structure \(\cR\) on \(M\) and a definable set \(\Sigma\subseteq M^k\), we say a family of sets \(\BR{S_\sigma\subseteq M^n}_{\sigma\in \Sigma}\) is a definable family if there is a definable set \(S\subseteq \Sigma\times M^n\) such that for each \(\sigma\in \Sigma\), the set \(S_\sigma\subseteq M^n\) coincides with the following fiber of \(S\):
    \begin{equation*}
        S_\sigma=\BR{
            x\in M^n
        \::\:
            \br{\sigma,x}\in S
        }.
    \end{equation*}
    We say a family of functions \(
    \BR{
        f_\sigma:M^m\supset A_\sigma\to M^n
    }_{\sigma\in\Sigma}
    \) is a definable family if the following family of sets:
    \begin{equation*}
        \BR{
            S_\sigma:=\BR{
                \br{a,b}\in A_\sigma\times M^n
        \::\:
                f_\sigma\br{a}=b
            }
            \subseteq M^{m+n}
        }_{\sigma\in\Sigma}
    \end{equation*}
    is a definable family of sets. In other words, there is a definable function \(f:M^{k+m}\supset A\to M^n\) such that for all \(\sigma\in \Sigma\), we have:
    \begin{equation*}
        \dom f_\sigma=
        \BR{
            a\in M^n
        \::\:
            \br{\sigma,a}\in A
        }
        =\br{\dom f}_\sigma;
    \end{equation*}
    \begin{equation*}
        f_\sigma\br{a}=f\br{\sigma,a},\quad
        \forall a\in \dom f_\sigma=\br{\dom f}_\sigma.
    \end{equation*}
\end{definition}

% \begin{definition}[Definable families of functions]
%     The family of functions \(\BR{f_a:I\to M^n}_{a\in A}\) is a definable family if there is a definable function \(F:A\times I\to M^n\) such that \(f_a\br{x}=F\br{a,x}\) for all \(\br{a,x}\in A\times I\).
% \end{definition}

The concept of structure and definability originates from the model theory. We choose to follow the presentation in \cite{MR1633348} and only mention the bare minimal concepts required from model theory. The key fact we will be using is the following:

\begin{remark}[Relation to first order formula]\label{rmk_def_set_construct}
Fix a structure \(\cR\) on \(M\). Let \(\phi\br{x}\) be a first-order formula with free variables \(x\in M^n\) built from predicates \(P\br{a}\) of the form \(a\in A\) for some definable set \(A\) (and thus, including formulas \(R\br{f_1\br{x},\dots,f_n\br{x}}\) composed of definable functions \(f_i\) and definable relations \(R\)). We have that
\begin{equation*}
    \BR{x\in M^n\::\: \phi\br{x}}\in \cR_n.
\end{equation*}
Technically speaking, the term ``definable" means ``definable with parameters" throughout the current paper.
\end{remark}
As a direct consequence, we have the following closure properties of definable functions:
\begin{proposition}\label{prop_func_def_closure}
    Fix a structure \(\cR\) on \(M\). Let functions \(f,g\) be definable.
    \begin{itemize}
        \item If \(A\subseteq \operatorname{dom}\br{f}\) is definable, \(f\vert_A\) is definable.
        % \item If \(f\br{x,y}=0\) defines an implicit function \(y=h\br{x}\), \(h\) is definable.
        \item If \(f\) is invertible, \(f^{-1}\) is definable.
        \item If \(\operatorname{dom}\br{f}\supset \operatorname{img}\br{g}\), \(f\circ g\) is definable.
    \end{itemize}
\end{proposition}

\begin{definition}[O-minimal expansion of \(\ang{\R,+,\cdot,<}\)]
    \(\cR\) is an o-minimal expansion of \(\ang{\R,+,\cdot,<}\) if it is a structure on \(\R\) that satisfies the following additional properties:
    \begin{enumerate}
    \setcounter{enumi}{4} 
        \item \(+\) and \(\cdot\) as bivariate functions are definable;
        \item \(<\) is a definable relation;
        \item \(A\in\cR_1\) if and only if \(A\) is a finite union of open intervals (bounded or unbounded) and isolated points in \(\R\).
    \end{enumerate}
\end{definition}

Henceforth, we fix an o-minimal expansion of \(\ang{\R,+,\cdot,<}\).

\begin{remark}\label{rmk_def_func}
    One may easily verify that the following are definable functions:
    \begin{itemize}
        \item the additive inverse \(\R\to\R: x\mapsto -x\),
        \item the multiplicative inverse \(\R\setminus\BR{0}\to\R: x\mapsto 1/x\),
        \item the binary maximum \(\R\times\R\to\R: \br{x,y}\mapsto \max\br{x,y}\),
        \item the absolute value \(\R\to\R:x\mapsto\abs{x}\),
        \item power functions \(\br{0,\infty}\to\R:x\mapsto x^r\) for \(r\in \mathbb{Q}\).
    \end{itemize}
\end{remark}

\begin{proposition}[Derivatives are definable]\label{prop_diff_def}
    Given a definable function \(\vf: \R^m\supset A\to\R^n\), the set of points \(A_j\subseteq A\) where \(\partial_j \vf\) exists is definable.\footnote{Recall the convention stated in \textbf{Remark \ref{rmk_dom_func}}.} Additionally, the function \(\partial_j \vf:A_j\to \R^n\) is also definable.
\end{proposition}
\begin{proof}
Without loss of generality, we may assume \(j=1\) and write \(\ve_1:=\br{1,0,\dots,0}\).
Observe that the relation \(\partial_1 \vf\br{\vx}=\vy\) can be expressed as a first-order formula \(\cD_1\br{\vx,\vy}\) composed of only definable relations and functions:
 \begin{equation*}
    \forall \epsilon>0,\exists \delta>0,\forall 
    h\in\R ,\;  
    \Br{
        0<\abs{h}<\delta \implies 
        \br{
           \vx, \vx+h \ve_1\in A
        \;\land\;
            \abs{\frac{\vf\br{\vx+h \ve_1}-\vf\br{\vx}}{h}-\vy}<\epsilon
        }
    }.
 \end{equation*}
 This implies that the following set is definable:
\begin{equation*}
    D_1:=\BR{
        \br{\vx,\vy}\in A\times \R^n
    \::\:
        \cD_1\br{\vx,\vy}
    }\in \cR_{m+n},
\end{equation*}
and thus, both the projection to the first \(m\) variables
\begin{equation*}
    A_1:=\pi_m D_1
    =\BR{
        \vx\in A\::\:
        \exists \vy\in \R^n,\,
        \cD_1\br{\vx,\vy}
    }
\end{equation*} 
and the function \(\partial_1 \vf: A_1\to\R^n\) are definable.
\end{proof}

\begin{theorem}[Monotonicity (\textsc{Chapter 3} \textbf{(1.2)} in \cite{MR1633348})]\label{thm_mono}
    Let \(-\infty\leq a <b\leq \infty\).\footnote{\(a,b\) are numbers in the extended real \(\R\sqcup\BR{-\infty,\infty}\).} Given a definable function \(f:\br{a,b}\to\R\), there are finite points:
    \begin{equation*}
        -\infty\leq a=a_0<a_1<\cdots<a_N=b\leq \infty
    \end{equation*}
    such that for all \(1\leq n\leq N\), the restriction \(f:\br{a_{n-1},a_n}\to\R\) is continuous and is either constant or strictly monotone.
\end{theorem}

\begin{corollary}[Existence of one-sided limits]\label{cor_lim}
    Let \(-\infty\leq a <b \leq \infty\). Given a definable function \(f:\br{a,b}\to\R\) and \(x\in\br{a,b}\), the following four limits exist in \(\R\sqcup\BR{-\infty,\infty}\):
    \begin{equation*}
        \lim_{t\to a^+} f\br{t},\quad
        \lim_{t\to b^-} f\br{t},\quad
        \lim_{t\to x^+} f\br{t},\quad
        \lim_{t\to x^-} f\br{t}.
    \end{equation*}
\end{corollary}

\begin{proposition}[Uniform finiteness (\textsc{Chapter 3} \textbf{(3.6)} in \cite{MR1633348})]\label{prop_uni_fin}
    Given a definable set \(A\subseteq \R^n\times\R\), there is a number \(M_A\in\N\) depending on the set \(A\) such that each fiber of the form:
    \begin{equation*}
        A_x:=\BR{y\in\R\::\: \br{x;y}\in A},\quad
        x\in\R^n
    \end{equation*}
    has at most \(M_A\) many connected components.
\end{proposition}

\begin{proposition}[Differentiability (\textsc{Chapter 7} \textbf{(2.5)} in \cite{MR1633348})]\label{prop_diff_ae}
Let \(-\infty\leq a <b \leq \infty\). A definable function \(f:\br{a,b}\to\R\) is differentiable outside a finite set \(Z\subseteq\br{a,b}\). 
\end{proposition}

\begin{remark}
The above collection of statements can be viewed under the natural framework of the cell decomposition theorem. We only state the bare minimum required for our application. Interested readers may refer to \cite{MR1633348} and \cite{coste2000introduction}.
\end{remark}

\begin{definition}[Polynomially bounded o-minimal expansion of \(\ang{\R,+,\cdot,<}\)]\label{def_omin_poly_bdd}
\(\cR\) is a polynomially bounded o-minimal expansion of \(\ang{\R,+,\cdot,<}\) if it is an o-minimal expansion of \(\ang{\R,+,\cdot,<}\) that satisfies the following additional properties:
\begin{enumerate}
\setcounter{enumi}{7} 
    \item For all definable \(f:\br{a,\infty}\to\R\), there is \(t_0\in\br{a,\infty}\) and \(N\in\N\) such that for all \(t>t_0\), the inequality \(\abs{f\br{t}}\leq \abs{t}^N\) holds.
\end{enumerate}
\end{definition}

One such structure is \(\R^\R_{\mathrm{an}}\) introduced in \cite{MR1278550}. In \(\R^\R_{\mathrm{an}}\), a large class of functions is definable. This includes all functions constructed from combinations of power functions \(x^\alpha\) and truncated real analytic functions through arithmetic, compositions, and inverses on their natural definable domains. Two typical examples are the class of generalized polynomials \(p\br{t}:=\sum_{j}c_jt^{\alpha_j}\) and the class of rational functions defined by generalized polynomials.

\section{Variants of shifted maximal operators and Fefferman-Stein inequalities}\label{sec_var_shift_n_feff_stein}
\begin{definition}[Variants of maximal operators with lacunary scales and bounded shifts]
    Given a lacunary sequence \(\Lambda:=\BR{\lambda_k}_k\subseteq \br{0,\infty}\) with \(1<\lambda_\ast:=\inf_k \lambda_{k+1}/\lambda_k\), a dyadic system \(\I:=\bigsqcup_k\I_k\) of intervals with \(I\in\I_k\implies \abs{I}=2^k\), and bounded sequence \(\Sigma:=\BR{\sigma_k}_k\), we introduce the following variants of shifted maximal operators:
    \begin{equation}\label{eq_MN_lac_bd}
        \cM_{\Lambda,\Sigma} f\br{x}:=
        \sup_k
            \int
                \abs{f}\br{x+\lambda_k\br{y+\sigma_k}}
                \ang{y}^{-N}
                % \varphi\br{y}
            dy,\quad N\gg 1\text{ fixed},
    \end{equation}
    \begin{equation}\label{eq_M_lac_bd}
        M_{\Lambda,\Sigma} f\br{x}:=
        \sup_k
            \int^1_0
                \abs{f}\br{x+\lambda_k\br{y+\sigma_k}}
                % \varphi\br{y}
            dy
        =\sup_k
            \fint_{x+\lambda_k\br{\sigma_k+\left[0,1\right)}}
                \abs{f}\br{y}
            dy,
    \end{equation}
    \begin{equation}\label{eq_Md_bd}
        M_{\I,\Sigma}f\br{x}:=\sup_k \sup_{\substack{
            I\in \I_k\\
            x+2^k\sigma_k \in I
        }}
        \fint_I \abs{f}\br{y}dy.
    \end{equation}
\end{definition}

The following results are slight modifications of \cite[\textbf{Theorem 3.1}]{MR3669936}.
\begin{proposition}\label{prop_FS_MN_lac_bd}
    \begin{equation*}
        \nrm{
            \nrm{
                \cM_{\Lambda,\Sigma}
                f_j
            }_{\ell^2\br{j}}
        }_{L^p}
        \underset{p}{\lesssim}
        \log^2_{\min\br{2,\lambda_\ast}}
        \br{e+\nrm{\sigma_k}_{\ell^\infty\br{k}}}
        \nrm{
            \nrm{f_j}_{\ell^2\br{j}}
        }_{L^p},\quad
        \forall p\in\br{1,\infty}.
    \end{equation*}
\end{proposition}

\begin{proposition}\label{prop_FS_M_lac_bd}
    \begin{equation*}
        \nrm{
            \nrm{
                M_{\Lambda,\Sigma}
                f_j
            }_{\ell^2\br{j}}
        }_{L^p}
        \underset{p}{\lesssim}
        \log^2_{\min\br{2,\lambda_\ast}}
        \br{e+\nrm{\sigma_k}_{\ell^\infty\br{k}}}
        \nrm{
            \nrm{f_j}_{\ell^2\br{j}}
        }_{L^p},\quad
        \forall p\in\br{1,\infty}.
    \end{equation*}
\end{proposition}

\begin{proposition}\label{prop_FS_Md_bd}
    \begin{equation*}
        \nrm{
            \nrm{
                M_{\I,\Sigma}
                f_j
            }_{\ell^2\br{j}}
        }_{L^p}
        \underset{p}{\lesssim}
        \log^2
        \br{e+\nrm{\sigma_k}_{\ell^\infty\br{k}}}
        \nrm{
            \nrm{f_j}_{\ell^2\br{j}}
        }_{L^p},\quad
        \forall p\in\br{1,\infty}.
    \end{equation*}
\end{proposition}
We observe that \textbf{Proposition \ref{prop_FS_MN_lac_bd}} follows from \textbf{Proposition \ref{prop_FS_M_lac_bd}} via the relation between \eqref{eq_MN_lac_bd} and \eqref{eq_M_lac_bd}:
\begin{equation*}
    \cM_{\Lambda,\Sigma}f\br{x}\lesssim \sum_{z\in\Z}\ang{z}^{-N}M_{\Lambda,z+\Sigma}f\br{x},\quad
    \nrm{z+\sigma_k}_{\ell^\infty\br{k}}\leq \nrm{\sigma_k}_{\ell^\infty}+\abs{z}.
\end{equation*}
Next, we claim that \textbf{Proposition \ref{prop_FS_M_lac_bd}} follows from \textbf{Proposition \ref{prop_FS_Md_bd}} via a standard trick on dyadic system:
\begin{definition}[One-third shifted dyadic system \cite{MR3065022,MR4007575}]
    Let \(u\in\BR{0,1/3,2/3}\). Define
    \begin{equation*}
        \I^{\br{u}}:=
        \bigsqcup_k \I^{\br{u}}_k,\quad
        \I^{\br{u}}_k:=
        \BR{
            2^k\br{m+\br{-1}^k u +\left[0,1\right)}\subseteq \R
            \::\:
            m\in\Z
        }.
    \end{equation*}
\end{definition}
\begin{remark}\label{rmk_3latice_trick}
    For all finite interval \(J\subseteq \R\), there is \(u\in\BR{0,1/3,2/3}\) and \(I\in \I^{\br{u}}\) such that:
    \begin{equation*}
        J\subseteq I,\quad
        3\abs{J}< \abs{I}\leq 6\abs{J}.
    \end{equation*}
    For details, see
    \cite[\textbf{Lemma 2.5}]{MR3065022} and \cite[\textbf{Remark 3.2}]{MR4007575}.
\end{remark}
We now relate \eqref{eq_M_lac_bd} and \eqref{eq_Md_bd} via \textbf{Remark \ref{rmk_3latice_trick}}.
Fix \(\Lambda\) lacunary and \(\Sigma\) bounded. By splitting \(\Lambda\) into \(O\br{1+\lceil \log_{\lambda_\ast} 2\rceil}\) many subsequences, we may assume \(\lambda_\ast > 2\). 
For \(k\in\Z\), we let \(n_k\in \Z\) and \(s_{n_k}\in\R\) be such that:
\begin{equation}\label{eq_lac_2_dya}
    2^{n_k-1}\leq 3\lambda_k<2^{n_k},\quad
    \lambda_k\sigma_k=2^{n_k}s_{n_k}.
\end{equation}
As a direct consequence of \textbf{Remark \ref{rmk_3latice_trick}} and \eqref{eq_lac_2_dya}, we may find for \(x\in\R\) and \(k\in\Z\) a choice of \(u\in\BR{0,1/3,2/3}\) such that there is \(I\in\I^{\br{u}}_{n_k}\) satisfying
\begin{equation*}
    x+2^{n_k}s_{n_k}=x+\lambda_k\sigma_k\in
    x+\lambda_k\br{\sigma_k+\left[0,1\right)}\subseteq I.
\end{equation*}
By setting \(S:=\BR{s_n}_n\) with
\begin{equation*}
    s_n:=
    \begin{cases}
        s_{n_k} & n=n_k \text{ for some }k\in\Z;\\
        0 & \text{else},
    \end{cases}
\end{equation*}
we have the following relations:
\begin{equation*}
    M_{\Lambda,\Sigma}f\br{x}\lesssim \sum_{u\in\BR{0,1/3,2/3}}
    M_{\I^{\br{u}},S}f\br{x},\quad
    \nrm{s_n}_{\ell^\infty\br{n}}\eqsim \nrm{\sigma_k}_{\ell^\infty\br{k}}.
\end{equation*}
Consequently, \textbf{Proposition \ref{prop_FS_M_lac_bd}} follows from \textbf{Proposition \ref{prop_FS_Md_bd}} via triangle inequality.
To address \textbf{Proposition \ref{prop_FS_Md_bd}}, we begin with an observation:
\begin{remark}
    Given bounded sequence \(S:=\BR{s_k}_k\) in \(\R\), we associate \(S\) to a pair of bounded integer sequences \(\lfloor S\rfloor:=\BR{\lfloor s_k\rfloor}_k\) and \(\lceil S\rceil:=\BR{\lceil s_k\rceil}_k\). By design \eqref{eq_Md_bd}, we have the following trivial relations:
    \begin{equation*}
        M_{\I,S}f\br{x}\leq M_{\I,\lfloor S\rfloor}f\br{x}+
        M_{\I,\lceil S\rceil}f\br{x},\quad
        \nrm{\lfloor s_k\rfloor}_{\ell^\infty\br{k}},\nrm{\lceil s_k\rceil}_{\ell^\infty\br{k}}\leq\nrm{s_k}_{\ell^\infty\br{k}}+1.
    \end{equation*}
\end{remark}

Hence, for the purpose of proving \textbf{Proposition \ref{prop_FS_Md_bd}}, we may assume \(S\subseteq \Z\).
Following the treatment of \cite[\textbf{Theorem 3.1}]{MR3669936}, \textbf{Proposition \ref{prop_FS_Md_bd}} reduces to a weighted estimate analogous to \cite[\textbf{Lemma 3.2}]{MR3669936}:
\begin{lemma}\label{lem_Md_weight}
    Let \(S:=\BR{s_k}_k\subseteq \Z\). For any locally integrable \(\omega \geq 0\), we have:
    \begin{equation*}
        \nrm{M_{\I,S} f}_{L^{1,\infty}\br{\omega}}
        \lesssim 
        \log
        \br{e+\nrm{s_k}_{\ell^\infty\br{k}}}
        \nrm{f}_{L^1\br{M_{\I,-S}\omega}}
        +\nrm{f}_{L^1\br{M\omega}}.
    \end{equation*}
    % with the implicit constant independent of the choice of \(\I\) and \(S\).
\end{lemma}
The proof of \textbf{Lemma \ref{lem_Md_weight}} follows from a direct modification of the proof of \cite[\textbf{Lemma 3.2}]{MR3669936}.

\section{On twisted paraproducts}\label{sec_on_twist_para_prod}
Given a symbol \(m=\widehat{K}\in L^\infty\br{\R^2}\), we associate to \(m\) the following operator:
\begin{equation}\label{eq_T_m_twisted_para}
    T_m\br{f_1,f_2}\br{\vx}:=
    \int
        f_1\br{x_1-y_1,x_2}
        f_2\br{x_1,x_2-y_2}
        K\br{\vy}
    d\vy
\end{equation}
\begin{equation}\label{eq_T_m_twisted_multiplier}
    =
    \int
        \4{1}f_1\br{\xi_1,x_2}
        \4{2}f_2\br{x_1,\xi_2}
        m\br{\vxi}
        e\br{\vx\cdot\vxi}
    d\vxi.
\end{equation}
Observe for positive definite diagonal matrix \(\vLambda\in GL_2\br{\R}\) the following rescaling identity:
\begin{equation}\label{eq_T_m_twisted_rescaling_id}
    \Dil^{p_3}_\vLambda T_{\Dil^\infty_\vLambda m}\br{f_1,f_2}=T_m\br{\Dil^{p_1}_\vLambda f_1, \Dil^{p_2}_\vLambda f_2},\quad
    \frac{1}{p_3}=\frac{1}{p_1}+\frac{1}{p_2}.
\end{equation}
As a direct consequence, any H\"{o}lder type bound for \(T_m\) will be preserved under dilation \(m\mapsto \Dil^\infty_\vLambda m\).
One sufficient condition for such H\"{o}lder type bounds that respect in part the above rescaling properties is the following symbol condition:
\begin{definition}[Anisotropic Mikhlin-type symbol condition]\label{def_ansio_mikh}
    Given \(\valpha\in\br{0,\infty}^2\) and \(m\in L^\infty\br{\R^2}\), we say \(m\) satisfies the \(\valpha\)-Mikhlin-type symbol condition up to order \(N\in\N\sqcup\BR{0}\) if:
    \begin{equation}\label{eq_ansio_mikh}
        \abs{
            \partial^\vbeta_\vxi
            m\br{\vxi} 
        }
        \underset{N}{\lesssim}1/
        \nrm{\vxi}^\vbeta_\valpha
        ,\quad
        \vbeta\in\BR{0,1,\dots,N}^2
        ,\quad
        \nrm{\vxi}^\vbeta_\valpha:=
        \br{\abs{\xi_1}^{1/\alpha_1}+\abs{\xi_2}^{1/\alpha_2}}^{\vbeta\cdot\valpha}.
    \end{equation}
\end{definition}
\begin{theorem}[Twisted paraproduct estimate with anisotropic scaling \cite{MR4295087}]\label{thm_aniso_twist}
For any symbol \(m=\widehat{K}\) satisfying the \(\valpha\)-Mikhlin-type symbol condition up to some sufficiently large order \(N\), the following estimate:
\begin{equation}\label{eq_T_m_twisted_multipli_est}
    \nrm{T_m\br{f_1,f_2}}_{L^{p_3}}\underset{m,p_j}{\lesssim}
    \nrm{f_1}_{L^{p_1}}
    \nrm{f_2}_{L^{p_2}}
\end{equation}
holds whenever the exponents satisfy the following condition:
\begin{equation}\label{eq_T_m_twisted_multiplier_cond}
    p_1,p_2\in\br{1,\infty}
    ,\quad
    \frac{1}{2}<
    \frac{1}{p_3}=\frac{1}{p_1}+\frac{1}{p_2}.
\end{equation}
\end{theorem}

For our purposes, we consider symbols that possess additional structure. To elaborate, let \(\vPhi:=\BR{\Phi_k}_k\subseteq C^\infty_c\br{\R^2}\) be a sequence of smooth bump functions and \(\vM:=\BR{\vM_k}_k\) be a sequence of positive definite diagonal matrices.
For any sequence of complex numbers 
\(\vepsilon:=\BR{\epsilon_k}_k\in \ell^\infty\), we define:
\begin{equation}\label{eq_symb_w_seq_struct}
    m_{\vepsilon,\vPhi,\vM}\br{\vxi}:=\sum_k\epsilon_k \Phi_k\br{\vM_k \vxi}.
\end{equation}
We remark that the proof of \textbf{Theorem \ref{thm_aniso_twist}} involves a reduction step that morally reduces \(m\) to symbol of the form \eqref{eq_symb_w_seq_struct} with \(1\)-bounded sequence \(\vepsilon\), \(\vPhi\) satisfying the regularity condition
\begin{equation}\label{eq_vPhi_seq_reg_cond}
    \abs{\partial^\vbeta\Phi_k\br{\vxi}}\lesssim 
    \1_{\br{1/2,2}}\br{\abs{\vxi}},
    \quad
    \forall k,\, \forall \vbeta\in\BR{0,1,\dots,N}^2,
\end{equation}
and \(\vM\) of the following particular form:
\begin{equation}\label{eq_vM_valpha_scaling}
    \vM_k:=
    \begin{pmatrix}
        2^{k \alpha_1 } & 0\\
        0 & 2^{k \alpha_2}
    \end{pmatrix}.
\end{equation}
In fact, direct computation shows that the symbol \(m_{\vepsilon,\vPhi,\vM}\) satisfies the \(\valpha\)-Mikhlin-type symbol condition up to order \(N\).
As one natural extension, we propose the following conjecture for general \(\vM\):
\begin{conjecture}[Twisted paraproduct associated to \(\vM\)]\label{conj_twist_vM}
Given a sequence of positive definite diagonal matrices \(\vM:=\BR{\vM_k}_k\) that satisfies \eqref{eq_lacu_seqs_condi},
there is a constant \(N_\vM\in\N\sqcup\BR{0}\) such that for any \(\vepsilon:=\BR{\epsilon_k}_k\in\ell^\infty\) and \(\vPhi:=\BR{\Phi_k}_k\) satisfying \eqref{eq_vPhi_seq_reg_cond} with \(N=N_\vM\), the following estimate
\begin{equation}\label{eq_T_m_seq_est}
    \nrm{
    T_{m_{\vepsilon,\vPhi,\vM}}\br{f_1,f_2}
    }_{L^{p_3}}\underset{p_j,M_\ast}{\lesssim}
    \nrm{\epsilon_k}_{\ell^\infty\br{k}}
    \nrm{f_1}_{L^{p_1}}
    \nrm{f_2}_{L^{p_2}}
\end{equation}
holds for the range mentioned in \eqref{eq_T_m_twisted_multiplier_cond}.
\end{conjecture}
Via previous discussions and the rescaling identity \eqref{eq_T_m_twisted_rescaling_id}, we have that \textbf{Conjecture \ref{conj_twist_vM}} holds and is morally equivalent to \textbf{Theorem \ref{thm_aniso_twist}} for sequences \(\vM\) of the following form:
\begin{equation}\label{eq_vM_hold_for_twist}
    \vM_k:= 
    \begin{pmatrix}
        2^{k \alpha_1 } /\lambda_1 & 0\\
        0 & 2^{k \alpha_2} /\lambda_2
    \end{pmatrix},\quad
    \text{ for some fixed }
    \valpha:=\br{\alpha_1,\alpha_2},\vlambda:=\br{\lambda_1,\lambda_2} \in\br{0,\infty}^2.
\end{equation}
% In fact, we believe the following shifted variant statement is true.
% \begin{conjecture}[Shifted twisted paraproduct associated to \(\vM\)]\label{conj_twist_vM_shifted}
% Given a sequence of positive definite diagonal matrices \(\vM:=\BR{\vM_k}_k\) that satisfies \eqref{eq_lacu_seqs_condi},
% there is a constant \(N_\vM\in\N\sqcup\BR{0}\) and \(C>0\) such that for any sequence of numbers \(\vepsilon:=\BR{\epsilon_k}_k \subseteq \C\), vectors \(\vU:=\BR{\vu_k}_k  \subseteq \R^2\), and functions \(\vPhi:=\BR{\Phi_k}_k\) satisfying \eqref{eq_vPhi_seq_reg_cond} with \(N=N_\vM\), the operator \(T_{m_{\vepsilon,\vU,\vPhi,\vM}}\) given by the symbol
% \begin{equation}\label{eq_symb_w_seq_struct_shifted}
%     m_{\vepsilon,\vU,\vPhi,\vM}\br{\vxi}:=
%     \sum_k \epsilon_k e\br{ \vu_k^\top \vM_k \vxi}\Phi_k\br{\vM_k \vxi}
% \end{equation}
% satisfies the following estimate
% \begin{equation}\label{eq_T_m_seq_est_shifted}
%     \nrm{
%     T_{m_{\vepsilon,\vU,\vPhi,\vM}}\br{f_1,f_2}
%     }_{L^{p_3}}\underset{p_j,M_\ast}{\lesssim}
%     \log^C\br{e+\nrm{\vu_k}_{\ell^\infty}}
%     \nrm{\epsilon_k}_{\ell^\infty\br{k}}
%     \nrm{f_1}_{L^{p_1}}
%     \nrm{f_2}_{L^{p_2}}
% \end{equation}
% holds for the range mentioned in \eqref{eq_T_m_twisted_multiplier_cond}.
% \end{conjecture}

% Such a statement may extend the boundedness range of the high frequency contribution
% \begin{equation*}
%     \nrm{\sum_k\cA^H_{\vM_k,\vgamma_k,\varphi}\br{f_1,f_2}}_{L^{p_3}}\lesssim\nrm{f_1}_{L^{p_1}}\nrm{f_2}_{L^{p_2}},\quad
%     p_3\in\br{1/2,\infty}.
% \end{equation*}

\newpage

%\printbibliography[title={References}]
\bibliographystyle{plain}
\bibliography{ref}

\begin{align*}
&\textsc{ Martin Hsu, Rutgers University, USA}\\
 &\textit{Email address:}\: \:\textbf{martin.hsu@rutgers.edu}\\
 \\
&\textsc{ Fred Yu-Hsiang Lin, purdue university, USA}\\
 &\textit{Email address:}\: \:\textbf{ lin2311@purdue.edu}\\
\end{align*}

\end{document}